\documentclass[12pt]{amsart}
\usepackage{amssymb}
\usepackage{hyperref}
\usepackage{thmtools}
\usepackage{thm-restate}
\usepackage{multirow,makecell} 
\usepackage[all]{xy}
\usepackage{mathrsfs}
\usepackage{enumerate}
\usepackage{algorithm}
\usepackage{algpseudocode}
\usepackage{geometry}
\usepackage{todonotes}
\usepackage{cleveref}
\usepackage{comment}
\usepackage[normalem]{ulem} 
\usepackage{bm} 

\usepackage{rotating}
\usepackage{tikz,graphicx}
\usetikzlibrary{matrix,shapes,arrows,fit,calc, math}
\usetikzlibrary{arrows.meta}

\definecolor{ForestGreen}{cmyk}{0.91,0,0.88,0.12}
\definecolor{myYellow}{rgb}{249,166,3}
\definecolor{Cerulean}{RGB}{20,0,150}
\definecolor{NavyBlue}{cmyk}{0.94,0.54,0,0}
\definecolor{cof}{RGB}{219,144,71}
\definecolor{pur}{RGB}{186,146,162}
\definecolor{greeo}{RGB}{91,173,69}
\definecolor{greet}{RGB}{52,111,72}
\definecolor{LimeGreen}{cmyk}{0.50,0,1,0}
\definecolor{RawSienna}{cmyk}{0,0.72,1,0.45}
\definecolor{RedViolet}{cmyk}{0.07,0.90,0,0.34}

\definecolor{blue}{rgb}{0, 0.445, 0.695}
\definecolor{bluishgreen}{rgb}{0, 0.626, 0.456}
\definecolor{red}{rgb}{0.896, 0.395, 0}
\definecolor{purple}{rgb}{0.783, 0.464, 0.640}
\definecolor{skyblue}{rgb}{0.359, 0.752, 0.973}
\definecolor{orange}{rgb}{0.999, 0.706, 0.0}
\definecolor{yellow}{rgb}{0.937, 0.890, 0.258}

\definecolor{olive}{RGB}{116,141,19}
\definecolor{green}{RGB}{108,208,48}
\definecolor{teal}{RGB}{47,77,62}
\definecolor{turquoise}{RGB}{86,235,211}
\definecolor{lightblue}{RGB}{150,178,153}
\definecolor{blue2}{RGB}{25,50,191}
\definecolor{indigo}{RGB}{142,128,251}
\definecolor{indigo2}{RGB}{114,32,246}
\definecolor{lightpurple}{RGB}{243,197,250}
\definecolor{purple2}{RGB}{105,66,131}
\definecolor{magenta}{RGB}{206,43,188}
\definecolor{brown}{RGB}{110,57,13} 
\definecolor{darkblue}{rgb}{0.0, 0.0, 0.7}

\newtheorem{theorem}{Theorem}[section]

\renewcommand{\thethmx}{\Alph{thmx}} 
\newtheorem{lemma}[theorem]{Lemma}
\newtheorem{proposition}[theorem]{Proposition}
\newtheorem{corollary}[theorem]{Corollary}
\theoremstyle{definition}
\newtheorem{definition}[theorem]{Definition}
\newtheorem{example}[theorem]{Example}
\newtheorem{question}[theorem]{Question}
\newtheorem{remark}[theorem]{Remark}
\newtheorem{conjecture}[theorem]{Conjecture}

\newcommand\ba{\mathbf{a}}
\newcommand\bb{\mathbf{b}}
\newcommand\bd{\mathbf{d}}
\newcommand\be{\mathbf{e}}
\newcommand\bp{\mathbf{p}}
\newcommand\bq{\mathbf{q}}
\newcommand\br{\mathbf{r}}
\newcommand\bs{\mathbf{s}}
\newcommand\bt{\mathbf{t}}
\newcommand\bu{\mathbf{u}}
\newcommand\bv{\mathbf{v}}
\newcommand\bw{\mathbf{w}}
\newcommand\bx{\mathbf{x}}
\newcommand\by{\mathbf{y}}
\newcommand\bz{\mathbf{z}}

\newcommand\balpha{\bm{\alpha}}

\newcommand\bbB{\mathbb{B}}
\newcommand\bbC{\mathbb{C}}
\newcommand\bbF{\mathbb{F}}
\newcommand\bbG{\mathbb{G}}
\newcommand\bbI{\mathbb{I}}
\newcommand\bbK{\mathbb{K}}
\newcommand\bbN{\mathbb{N}}
\newcommand\bbP{\mathbb{P}}
\newcommand\bbQ{\mathbb{Q}}
\newcommand\bbR{\mathbb{R}}
\newcommand\bbS{\mathbb{S}}
\newcommand\bbZ{\mathbb{Z}}

\newcommand\calA{\mathcal{A}}
\newcommand\calB{\mathcal{B}}
\newcommand\calC{\mathcal{C}}
\newcommand\calD{\mathcal{D}}
\newcommand\calE{\mathcal{E}}
\newcommand\calF{\mathcal{F}}
\newcommand\calG{\mathcal{G}}
\newcommand\calH{\mathcal{H}}
\newcommand\calI{\mathcal{I}}
\newcommand\calJ{\mathcal{J}}
\newcommand\calK{\mathcal{K}}
\newcommand\calL{\mathcal{L}}
\newcommand\calM{\mathcal{M}}
\newcommand\calN{\mathcal{N}}
\newcommand\calO{\mathcal{O}}
\newcommand\calP{\mathcal{P}}
\newcommand\calQ{\mathcal{Q}}
\newcommand\calR{\mathcal{R}}
\newcommand\calS{\mathcal{S}}
\newcommand\calT{\mathcal{T}}
\newcommand\calU{\mathcal{U}}
\newcommand\calV{\mathcal{V}}
\newcommand\calW{\mathcal{W}}
\newcommand\calX{\mathcal{X}}
\newcommand\calY{\mathcal{Y}}
\newcommand\calZ{\mathcal{Z}}

\newcommand\scrB{\mathscr{B}}
\newcommand\scrC{\mathscr{C}}
\newcommand\scrE{\mathscr{E}}
\newcommand\scrF{\mathscr{F}}
\newcommand\scrI{\mathscr{I}}
\newcommand\scrL{\mathscr{L}}
\newcommand\scrM{\mathscr{M}}
\newcommand\scrN{\mathscr{N}}
\newcommand\scrO{\mathscr{O}}
\newcommand\scrP{\mathscr{P}}
\newcommand\scrR{\mathscr{R}}
\newcommand\scrS{\mathscr{S}}

\newcommand\frakj{\mathfrak{j}} 
\newcommand\frakm{\mathfrak{m}} 

\newcommand\horiz{\mathrm{horiz}}
\newcommand\inedge{\mathrm{in}}
\newcommand\outedge{\mathrm{out}}
\newcommand\In{\mathrm{In}}
\newcommand\Out{\mathrm{Out}}
\newcommand\inv{^{-1}}
\newcommand\precccwrot{\prec_{\mathrm{rot}}^{\mathrm{ccw}}}
\newcommand\leqccwrot{\leq_{\mathrm{rot}}^{\mathrm{ccw}}}
\newcommand\cw{\mathrm{cw}}
\newcommand\ccw{\mathrm{ccw}}
\newcommand\down{\mathrm{down}}
\newcommand\up{\mathrm{up}}
\newcommand\mapcw{\varphi^\mathrm{cw}}
\newcommand\mapccw{\varphi^\mathrm{ccw}}

\newcommand{\oru}[1]{\mathrm{oru}(#1)}
\newcommand{\multioruga}[1]{G_{#1}}
\newcommand{\cambrianCaracolGraph}[1]{G_{#1}}
\newcommand{\crossGridGraph}[1]{G_{#1}}

\newcommand{\defn}[1]{{\color{green!50!black}\textbf{\emph{#1}}}}

\newcommand{\cesar}[1]{\todo[size=\footnotesize,color=orange!30,inline]{#1  \hfill --- C.}}
\newcommand{\cesarm}[1]{\todo[size=\footnotesize,color=orange!30]{#1  \hfill --- C.}}
\newcommand{\cesarformatias}[1]{\todo[size=\footnotesize,color=red!30,inline]{For Matias: #1  \hfill --- C.}}

\newcommand{\cesaredit}[1]{{\color{red} #1}}

\newcommand{\matias}[1]{\todo[size=\footnotesize,color=blue!20,inline]{#1  \hfill --- M.}}
\newcommand{\matiasm}[1]{\todo[size=\footnotesize,color=blue!20]{#1  \hfill --- M.}}
\newcommand{\matiasforcesar}[1]{\todo[size=\footnotesize,color=green!20,inline]{For Cesar: #1  \hfill --- M.}}

\DeclareMathOperator{\hor}{hor}
\DeclareMathOperator{\ver}{ver}
\DeclareMathOperator{\car}{car}
\DeclareMathOperator{\Cat}{Cat}
\DeclareMathOperator{\Nar}{Nar}
\DeclareMathOperator{\vol}{vol}
\DeclareMathOperator{\col}{col}
\DeclareMathOperator{\row}{row}
\DeclareMathOperator{\grid}{grid}
\DeclareMathOperator{\Tam}{Tam}
\DeclareMathOperator{\Dyck}{Dyck}
\DeclareMathOperator{\Rem}{Rem}
\DeclareMathOperator{\conv}{conv}
\DeclareMathOperator{\rev}{rev}
\DeclareMathOperator{\Top}{Top}
\DeclareMathOperator{\Bot}{Bot}

\usepackage{etoolbox}
\DeclareRobustCommand{\verylongrightarrow}{\joinrel\relbar\joinrel\relbar\joinrel\relbar\joinrel\relbar\joinrel\relbar\joinrel\rightarrow}
\newcommand{\includeSymbol}[1]{\ensuremath{%
	\mathchoice
		{\raisebox{-.7mm}{\includegraphics[height=2.2ex]{#1}}}	
		{\raisebox{-.7mm}{\includegraphics[height=2.2ex]{#1}}}
		{\raisebox{-.6mm}{\includegraphics[height=1.6ex]{#1}}}
		{\raisebox{-.5mm}{\includegraphics[height=1ex]{#1}}}
}}
\robustify{\includeSymbol}
\newcommand{\noneCirc}{\includeSymbol{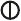}}
\newcommand{\upCirc}{\includeSymbol{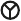}}
\newcommand{\downCirc}{\includeSymbol{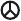}}
\newcommand{\upDownCirc}{\includeSymbol{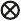}}

\makeatletter
\renewcommand\part[1]{%
  \refstepcounter{part}
  \addcontentsline{toc}{part}{Part~\thepart:\ #1}
  \Hy@raisedlink{\hypertarget{part\thepart}{}}
  \markboth{}{}%
  \vspace*{.7\linespacing}%
  \par\noindent
  \begin{center}
    \LARGE\sffamily
    Part~\thepart:\ #1
  \end{center}%
  \vspace*{.8\linespacing}%
  \@addtoreset{section}{part}
}
\makeatother

\makeatletter
\@addtoreset{section}{part}
\makeatother
\renewcommand{\thesection}{\arabic{part}.\arabic{section}}

\makeatletter
\def\l@section{\@tocline{1}{5pt}{0pc}{}{}}
\makeatother

\let\oldtocpart=\tocpart
\renewcommand{\tocpart}[2]{\bf\large\oldtocpart{#1}{#2}}

\let\oldtocsection=\tocsection
\renewcommand{\tocsection}[2]{\bf\oldtocsection{#1}{#2}}

\title[Framing lattices]{Framing lattices and flow polytopes}

\author[von Bell]{Matias von Bell}
\address[von Bell]{Institute of Geometry\\
         Graz University of Technology\\
}
\email{matias.vonbell@gmail.com}

\author[Ceballos]{Cesar Ceballos}
\address[Ceballos]{Institute of Geometry\\
         Graz University of Technology\\
}
\email{cesar.ceballos@tugraz.at}

\usepackage[backend=bibtex]{biblatex}
\begin{document}
\parskip=5pt

\begin{abstract}
Flow polytopes of acyclic oriented graphs arise naturally in combinatorial
optimization, and the study of their volumes and triangulations has revealed
intriguing connections across combinatorics, geometry, algebra, and representation
theory.

In this work, we introduce the framing lattice associated with a framed graph, whose Hasse diagram is dual to a framed triangulation of the corresponding flow polytope.
Framing lattices are remarkable in that they provide a unifying framework encompassing many classical and well-studied lattice structures, including the Boolean lattice, the Tamari lattice, and the weak order on permutations. 
They further subsume a broad array of examples such as all type-$A$ Cambrian lattices, the Grassmann and grid-Tamari lattices, the alt-$\nu$-Tamari and cross-Tamari lattices, the permutree lattices, and the $\tau$-tilting posets of certain gentle algebras.

We show, among several foundational structural properties,  that the framing lattice is a semidistributive, congruence uniform, and polygonal lattice, with its polygons consisting of squares, pentagons, and hexagons. 
We study its connections to noncrossing partitions via Reading’s core label orders, simple representations of its join and meet irreducible elements, and several of its lattice congruences and quotients induced by a graph operation called an M-move.

\vspace{.5cm}
\noindent{\bf \keywordsname}: Framing lattice, flow polytope,  DKK triangulation, weak order, Tamari lattice, Cambrian lattice, cross-Tamari lattice.
\end{abstract}

\maketitle 
\tableofcontents

\section*{Introduction}\label{sec.intro}

Flow polytopes encode spaces of flows in directed acyclic graphs and form a central class of objects in algebraic and geometric combinatorics, with connections to representation theory~\cite{BV08}, toric geometry~\cite{Hil03}, diagonal harmonics~\cite{MMR17}, and gentle algebras~\cite{BBBHPSY22, BS24}. A particularly rich combinatorial structure arises from their triangulations. Among these, the framing (Danilov--Karzanov--Koshevoy) triangulations~\cite{DKK12} are defined by equipping a directed acyclic graph $G$ with a linear ordering of the incoming and outgoing edges at each vertex.

Recent work has revealed that the dual graphs of such triangulations often carry the structure of well-known lattices, including Tamari-like lattices~\cite{BGMY23}, generalizations of the weak order~\cite{GMPTY23}, and lattices arising in the representation theory of gentle algebras~\cite{BBBHPSY22}. Despite this diversity, these lattice structures have so far been understood only through case-by-case identifications with previously known models. In particular, fundamental structural properties such as semidistributivity and congruence uniformity have required intricate and specialized arguments in each setting.

In this paper, we show that these phenomena admit a uniform explanation. To each framed directed acyclic graph $(G,F)$, we associate a canonical partial order on the facets of the corresponding triangulation and prove that it always forms a lattice, which we call a \emph{framing lattice}. These lattices lie in a highly constrained class: they are semidistributive, congruence-uniform, and polygonal, with local structure governed by polygons of size four, five, or six. Thus, a broad family of lattices arising across combinatorics is obtained from a single combinatorial construction.

A key feature of this framework is that it replaces a range of case-by-case arguments in the literature with a single conceptual mechanism. Properties that were previously established individually for specific families now follow uniformly from the general theory. At the same time, the construction admits precise combinatorial control, including descriptions of lattice operations, join- and meet-irreducible elements, and associated auxiliary orders such as the core label order, which connects these lattices to the theory of noncrossing partitions.

The theory is compatible with natural operations on framed graphs. In particular, an edge-cutting operation ($M$-move) induces lattice congruences and yields canonical lattice quotients; iterating this operation produces a distributive lattice depending only on the underlying graph $G$, and not on the choice of framing $F$.

Taken together, these results initiate a general theory of framing lattices and identify them as a common combinatorial structure underlying several previously disparate phenomena. In particular, they arise naturally at the intersection of the geometry of flow polytopes, lattice-theoretic combinatorics, and the representation theory of gentle algebras. This perspective not only unifies existing examples, but also establishes a conceptual bridge between these areas.

The framework developed here also points toward further directions. Framing lattices appear to admit natural geometric realizations as edge graphs of polytopal complexes~\cite{framingtopes}, and their connection to flow polytopes suggests potential applications in enumerative questions, including the Ehrhart theory of these polytopes. More broadly, this provides a systematic method for constructing and studying large families of lattices, and suggests that many further examples remain to be explored. This theory also opens new directions for applying lattice-theoretic methods to geometric, combinatorial, and representation-theoretic problems arising from flow polytopes.

\section*{Acknowledgements}
The authors were partially supported by the Austrian Science Fund FWF, grants P 33278 and I 5788.
We are especially grateful for discussions with Jonah Berggren, Cl\'ement Chenevi\`ere, Sergio Fernandez de soto, Rafael Goz\'alez D'Le\'on, Eva Philippe, Germain Poullot, Daniel Tamayo Jim\'enez, and Martha Yip. Additionally, we wish to thank Germain Poullot for helpful feedback on the first draft and suggesting Questions \ref{question.doubling} and \ref{question.simplecriterion}.

\part{Framing lattices}

In this first part, we present the seminal theory of framing lattices. We start with some preliminary background on flow polytopes and framed triangulations (\Cref{sec.background}), followed by the introduction of the framing poset (\Cref{sec_framing_poset}) and the proof of several lattice properties, including polygonality, semidistrivitubity and congruence uniformity (\Cref{sec_lattice_properties}). We also propose a framing generalization of the poset of noncrossing partitions via the core label order (\Cref{sec_core_label_order}), and study certain lattice quotients of the framing lattice (\Cref{sec.quotients}).

\section{Background and terminology}\label{sec.background}

In this section we recall some preliminaries about flow polytopes and their triangulations, and present some useful lemmas that will be used throughout the paper. 
We refer to~\cite{B22,DKK12,MM19} for further details and background on flow polytopes. 

\subsection{Flow polytopes}
Let $G$ be a directed acyclic graph on vertex set $V(G) = [n]$ and edge multiset $E(G)$ such that 
all edges are directed from smaller vertices to larger vertices and $G$ has a unique source $s=1$ and sink $t=n$.
We call such a graph $G$ a \defn{flow graph}.
A path from the source to the sink is said to be a \defn{route}.

Given a flow graph $G$ with vertex set $[n]$, a \defn{unit flow} on $G$ is a tuple $(x_e)_{e\in E(G)} \in \bbR^{|E(G)|}_{\geq0}$ satisfying
$$\sum_{e \in \Out(j)} x_e - \sum_{e \in \In(j)}  x_e = u_j,$$
where $u_1 = 1$, $u_n = -1$, and $u_j = 0$ for $1 < 
j < n$.
The \defn{flow polytope} of $G$ is the set $\calF_G$ of unit flows on $G$.
The dimension of a flow polytope $\calF_{G}$ is given by the formula $|E(G)| - |V(G)| +1$. 
The vertices of $\calF_G$ can be characterized as the unit flows on $G$ which have value one on the edges of a route and value zero on the remaining edges.
Thus $\calF_G$ can be described as the convex hull of the indicator vectors of the routes of $G$.

Flow polytopes for general net flows $(u_1,\dots,u_n)\in \mathbb{R}^n$  arise naturally in combinatorial optimization, and include many interesting renown examples, such as Tesler polytopes \cite{MMR17}, Gelfand--Tsetlin polytopes \cite{LMS19}, certain order polytopes and faces of the alternating sign matrix polytope \cite{MM19}, the Chan--Robbins--Yuen polytope \cite{CRY00}, and the Pitman--Stanley polytope~\cite{PS02}.

\begin{remark}
    We remark that the flow polytopes we consider are limited to those with unit flows. We also assume that $G$ has a unique source and sink. While this assumption can be omitted without affecting our lattice theoretic results (by connecting all sources to a ``global source'' and all sinks to a ``global sink''), we keep it to retain the connection to triangulated unit flow polytopes explained further below.
\end{remark}

\begin{example}[The oruga graph and the cube]\label{example_running_oruga}
    Let $G_n=\oru{n}$ be the \defn{oruga graph} on the vertex set~$[n+1]$ containing two edges between $i$ and $i+1$ for $i \in [n]$. These two edges are oriented from smallest vertex to largest and we label them by $e_{2i-1}$ and $e_{2i}$. Some examples of the oruga graph are illustrated on the top of~\Cref{fig_oruga_graph}. 
    The name ``oruga" was given in \cite{GMPTY23}, meaning caterpillar in Spanish.

\begin{figure}[htb]
        \centering
    \begin{tikzpicture}[scale=0.8,
    fnode/.style={circle, draw, inner sep=1.4pt, fill},
    fnodesmall/.style={circle, draw, inner sep=0.7pt, fill},
    fnodeextrasmall/.style={circle, draw, inner sep=0.5pt, fill},
    ]
    %
    %
    \begin{scope}[shift={(0,0)}]
        \draw[thick, color=black] (0,0) .. controls (0.4, 0.6) and (1.6, 0.6) .. (2,0);
        \draw[thick, color=black] (0,0) .. controls (0.4, -0.6) and (1.6, -0.6) .. (2,0);
        \node[fnode] (1) at (0,0) {};
        \node[fnode] (2) at (2,0) {};
    
        \node[] (a) at (1,0.8) {\scriptsize $e_1$};
        \node[] (a) at (1,-0.8) {\scriptsize $e_2$};

        \node[] at (1,-1.5) {\footnotesize $G_1=\oru{1}$};    
    \end{scope}
    %
    %
    \begin{scope}[shift={(5,0)}]
        \draw[thick, color=black] (0,0) .. controls (0.4, 0.6) and (1.6, 0.6) .. (2,0);
        \draw[thick, color=black] (0,0) .. controls (0.4, -0.6) and (1.6, -0.6) .. (2,0);    
        \draw[thick, color=black] (2,0) .. controls (2.4, 0.6) and (3.6, 0.6) .. (4,0);
        \draw[thick, color=black] (2,0) .. controls (2.4, -0.6) and (3.6, -0.6) .. (4,0);
        
        \node[fnode] (1) at (0,0) {};
        \node[fnode] (2) at (2,0) {};
        \node[fnode] (3) at (4,0) {};
    
        \node[] (a) at (1,0.8) {\scriptsize $e_1$};
        \node[] (a) at (1,-0.8) {\scriptsize $e_2$};
        \node[] (a) at (3,0.8) {\scriptsize $e_3$};    
        \node[] (a) at (3,-0.8) {\scriptsize $e_4$};    

        \node[] at (2,-1.5) {\footnotesize $G_2=\oru{2}$};    
    \end{scope}
    %
    %
    \begin{scope}[shift={(12,0)}]
        \draw[thick, color=black] (0,0) .. controls (0.4, 0.6) and (1.6, 0.6) .. (2,0);
        \draw[thick, color=black] (0,0) .. controls (0.4, -0.6) and (1.6, -0.6) .. (2,0);    
        \draw[thick, color=black] (2,0) .. controls (2.4, 0.6) and (3.6, 0.6) .. (4,0);
        \draw[thick, color=black] (2,0) .. controls (2.4, -0.6) and (3.6, -0.6) .. (4,0);
        \draw[thick, color=black] (4,0) .. controls (4.4, 0.6) and (5.6, 0.6) .. (6,0);
        \draw[thick, color=black] (4,0) .. controls (4.4, -0.6) and (5.6, -0.6) .. (6,0);          
        \node[fnode] (1) at (0,0) {};
        \node[fnode] (2) at (2,0) {};
        \node[fnode] (3) at (4,0) {};
        \node[fnode] (4) at (6,0) {};    
    
        \node[] (a) at (1,-0.8) {\scriptsize $e_2$};
        \node[] (a) at (3,-0.8) {\scriptsize $e_4$};    
        \node[] (a) at (5,-0.8) {\scriptsize $e_6$};
        \node[] (a) at (1,0.8) {\scriptsize $e_1$};
        \node[] (a) at (3,0.8) {\scriptsize $e_3$};    
        \node[] (a) at (5,0.8) {\scriptsize $e_5$};

        \node[] at (3,-1.5) {\footnotesize $G_3=\oru{3}$};
    \end{scope}
    %
    %
    \begin{scope}[shift={(1,-6)},scale=2]
        \draw[thick] (0,0) -- (0,1);
        \node[fnode] at (0,0) {};
        \node[fnode] at (0,1) {};
        \node[] at (0,-0.5) {\footnotesize $\calF_{G_1}$};           
        %
        %
        \begin{scope}[shift={(0.3,0)}, scale=0.15]
            \draw[thick, color=black] (0,0) .. controls (0.4, 0.6) and (1.6, 0.6) .. (2,0);
            \node[] at (1,1.3) {\scriptsize $10$};

            \node[fnodeextrasmall] (1) at (0,0) {};
            \node[fnodeextrasmall] (2) at (2,0) {};
        
        \end{scope}
        \begin{scope}[shift={(0.3,1)}, scale=0.15]
            \draw[thick, color=black] (0,0) .. controls (0.4, -0.6) and (1.6, -0.6) .. (2,0);    
            \node[] at (1,1.3) {\scriptsize $01$};

            \node[fnodeextrasmall] (1) at (0,0) {};
            \node[fnodeextrasmall] (2) at (2,0) {};
        
        \end{scope}
    \end{scope}
    %
    %
    \begin{scope}[shift={(6,-6)},scale=2]
        \draw[thick] (0,0) -- (0,1) -- (1,1) -- (1,0) -- cycle ;
        \node[fnode] at (0,0) {};
        \node[fnode] at (0,1) {};
        \node[fnode] at (1,1) {};
        \node[fnode] at (1,0) {};
        \node[] at (0.5,-0.5) {\footnotesize $\calF_{G_2}$};      
        %
        %
        \begin{scope}[shift={(-0.8,0)}, scale=0.15]
            \draw[thick, color=black] (0,0) .. controls (0.4, 0.6) and (1.6, 0.6) .. (2,0);
            \draw[thick, color=black] (2,0) .. controls (2.4, 0.6) and (3.6, 0.6) .. (4,0);
            \node[] at (2,1.3) {\scriptsize $1010$};
            
            \node[fnodeextrasmall] (1) at (0,0) {};
            \node[fnodeextrasmall] (2) at (2,0) {};
            \node[fnodeextrasmall] (3) at (4,0) {};
        
        \end{scope}        
        \begin{scope}[shift={(1.2,0)}, scale=0.15]
            \draw[thick, color=black] (0,0) .. controls (0.4, 0.6) and (1.6, 0.6) .. (2,0);
            \draw[thick, color=black] (2,0) .. controls (2.4, -0.6) and (3.6, -0.6) .. (4,0);
            \node[] at (2,1.3) {\scriptsize $1001$};
            
            \node[fnodeextrasmall] (1) at (0,0) {};
            \node[fnodeextrasmall] (2) at (2,0) {};
            \node[fnodeextrasmall] (3) at (4,0) {};
        
        \end{scope}        
        \begin{scope}[shift={(-0.8,1)}, scale=0.15]
            \draw[thick, color=black] (0,0) .. controls (0.4, -0.6) and (1.6, -0.6) .. (2,0);    
            \draw[thick, color=black] (2,0) .. controls (2.4, 0.6) and (3.6, 0.6) .. (4,0);
            \node[] at (2,1.3) {\scriptsize $0110$};
            
            \node[fnodeextrasmall] (1) at (0,0) {};
            \node[fnodeextrasmall] (2) at (2,0) {};
            \node[fnodeextrasmall] (3) at (4,0) {};
        
        \end{scope}        
        \begin{scope}[shift={(1.2,1)}, scale=0.15]
            \draw[thick, color=black] (0,0) .. controls (0.4, -0.6) and (1.6, -0.6) .. (2,0);    
            \draw[thick, color=black] (2,0) .. controls (2.4, -0.6) and (3.6, -0.6) .. (4,0);
            \node[] at (2,1.3) {\scriptsize $0101$};
            
            \node[fnodeextrasmall] (1) at (0,0) {};
            \node[fnodeextrasmall] (2) at (2,0) {};
            \node[fnodeextrasmall] (3) at (4,0) {};
        
        \end{scope}        
    \end{scope}
    %
    %
    \begin{scope}[shift={(13.5,-6)},
        x={(2cm, 0cm)},
        y={(0cm, 2cm)},
        z={(1cm, 1cm)}]
        \draw[thick] (0,0,0) -- (1,0,0) -- (1,1,0) -- (0,1,0) -- cycle;
        \draw[thick] (1,0,0) -- (1,0,1) -- (1,1,1) -- (1,1,0);
        \draw[thick] (0,1,0) -- (0,1,1) -- (1,1,1);
        \draw[dashed] (0,0,0) -- (0,0,1);
        \draw[dashed] (0,1,1) -- (0,0,1) -- (1,0,1);
        \node[fnode] at (0,0,0) {};
        \node[circle, draw, inner sep=1.4pt,fill=white] at (0,0,1) {};
        \node[fnode] at (0,1,0) {};
        \node[fnode] at (0,1,1) {};
        \node[fnode] at (1,0,0) {};
        \node[fnode] at (1,0,1) {};
        \node[fnode] at (1,1,0) {};
        \node[fnode] at (1,1,1) {};
        \node[] at (0.75,-0.5) {\footnotesize $\calF_{G_3}$};    
        %
        %
        \begin{scope}[shift={(-1.1,0,0)}, scale=0.15]
            \draw[thick, color=black] (0,0) .. controls (0.4, 0.6) and (1.6, 0.6) .. (2,0);
            \draw[thick, color=black] (2,0) .. controls (2.4, 0.6) and (3.6, 0.6) .. (4,0);
            \draw[thick, color=black] (4,0) .. controls (4.4, 0.6) and (5.6, 0.6) .. (6,0);
            \node[] at (3,1.3) {\scriptsize $101010$};

            \node[fnodeextrasmall] (1) at (0,0) {};
            \node[fnodeextrasmall] (2) at (2,0) {};
            \node[fnodeextrasmall] (3) at (4,0) {};
            \node[fnodeextrasmall] (4) at (6,0) {};    
        
        \end{scope}
        \begin{scope}[shift={(1.2,0,0)}, scale=0.15]
            \draw[thick, color=black] (0,0) .. controls (0.4, 0.6) and (1.6, 0.6) .. (2,0);
            \draw[thick, color=black] (2,0) .. controls (2.4, -0.6) and (3.6, -0.6) .. (4,0);
            \draw[thick, color=black] (4,0) .. controls (4.4, 0.6) and (5.6, 0.6) .. (6,0);
            \node[] at (3,1.3) {\scriptsize $100110$};

            \node[fnodeextrasmall] (1) at (0,0) {};
            \node[fnodeextrasmall] (2) at (2,0) {};
            \node[fnodeextrasmall] (3) at (4,0) {};
            \node[fnodeextrasmall] (4) at (6,0) {};    
        
        \end{scope}
        \begin{scope}[shift={(-1.1,1,0)}, scale=0.15]
            \draw[thick, color=black] (0,0) .. controls (0.4, -0.6) and (1.6, -0.6) .. (2,0);    
            \draw[thick, color=black] (2,0) .. controls (2.4, 0.6) and (3.6, 0.6) .. (4,0);
            \draw[thick, color=black] (4,0) .. controls (4.4, 0.6) and (5.6, 0.6) .. (6,0);
            \node[] at (3,1.3) {\scriptsize $011010$};

            \node[fnodeextrasmall] (1) at (0,0) {};
            \node[fnodeextrasmall] (2) at (2,0) {};
            \node[fnodeextrasmall] (3) at (4,0) {};
            \node[fnodeextrasmall] (4) at (6,0) {};    
        
        \end{scope}
        \begin{scope}[shift={(1.2,1,0)}, scale=0.15]
            \draw[thick, color=black] (0,0) .. controls (0.4, -0.6) and (1.6, -0.6) .. (2,0);    
            \draw[thick, color=black] (2,0) .. controls (2.4, -0.6) and (3.6, -0.6) .. (4,0);
            \draw[thick, color=black] (4,0) .. controls (4.4, 0.6) and (5.6, 0.6) .. (6,0);
            \node[] at (3.15,1.3) {\scriptsize $010110$};

            \node[fnodeextrasmall] (1) at (0,0) {};
            \node[fnodeextrasmall] (2) at (2,0) {};
            \node[fnodeextrasmall] (3) at (4,0) {};
            \node[fnodeextrasmall] (4) at (6,0) {};    
        
        \end{scope}
        \begin{scope}[shift={(-1.1,0,1)}, scale=0.15]
            \draw[thick, color=black] (0,0) .. controls (0.4, 0.6) and (1.6, 0.6) .. (2,0);
            \draw[thick, color=black] (2,0) .. controls (2.4, 0.6) and (3.6, 0.6) .. (4,0);
            \draw[thick, color=black] (4,0) .. controls (4.4, -0.6) and (5.6, -0.6) .. (6,0);       
            \node[] at (3,1.3) {\scriptsize $101001$};

            \node[fnodeextrasmall] (1) at (0,0) {};
            \node[fnodeextrasmall] (2) at (2,0) {};
            \node[fnodeextrasmall] (3) at (4,0) {};
            \node[fnodeextrasmall] (4) at (6,0) {};    
        
        \end{scope}
        \begin{scope}[shift={(1.2,0,1)}, scale=0.15]
            \draw[thick, color=black] (0,0) .. controls (0.4, 0.6) and (1.6, 0.6) .. (2,0);
            \draw[thick, color=black] (2,0) .. controls (2.4, -0.6) and (3.6, -0.6) .. (4,0);
            \draw[thick, color=black] (4,0) .. controls (4.4, -0.6) and (5.6, -0.6) .. (6,0);          
            \node[] at (3,1.3) {\scriptsize $100101$};

            \node[fnodeextrasmall] (1) at (0,0) {};
            \node[fnodeextrasmall] (2) at (2,0) {};
            \node[fnodeextrasmall] (3) at (4,0) {};
            \node[fnodeextrasmall] (4) at (6,0) {};    
        
        \end{scope}
        \begin{scope}[shift={(-1.1,1,1)}, scale=0.15]
            \draw[thick, color=black] (0,0) .. controls (0.4, -0.6) and (1.6, -0.6) .. (2,0);    
            \draw[thick, color=black] (2,0) .. controls (2.4, 0.6) and (3.6, 0.6) .. (4,0);
            \draw[thick, color=black] (4,0) .. controls (4.4, -0.6) and (5.6, -0.6) .. (6,0);          
            \node[] at (3,1.3) {\scriptsize $011001$};

            \node[fnodeextrasmall] (1) at (0,0) {};
            \node[fnodeextrasmall] (2) at (2,0) {};
            \node[fnodeextrasmall] (3) at (4,0) {};
            \node[fnodeextrasmall] (4) at (6,0) {};    
        
        \end{scope}
        \begin{scope}[shift={(1.2,1,1)}, scale=0.15]
            \draw[thick, color=black] (0,0) .. controls (0.4, -0.6) and (1.6, -0.6) .. (2,0);    
            \draw[thick, color=black] (2,0) .. controls (2.4, -0.6) and (3.6, -0.6) .. (4,0);
            \draw[thick, color=black] (4,0) .. controls (4.4, -0.6) and (5.6, -0.6) .. (6,0);          
            \node[] at (3,1.3) {\scriptsize $010101$};

            \node[fnodeextrasmall] (1) at (0,0) {};
            \node[fnodeextrasmall] (2) at (2,0) {};
            \node[fnodeextrasmall] (3) at (4,0) {};
            \node[fnodeextrasmall] (4) at (6,0) {};    
        
        \end{scope}        
    \end{scope}
\end{tikzpicture}
        \caption{Some examples of the oruga graph and their flow polytopes.}
        \label{fig_oruga_graph}
\end{figure}

    The flow polytope $\calF_{G_n}$ is the set of points $(x_1,\dots,x_{2n})\in \bbR^{2n}_{\geq 0}$ (that is $x_i\geq 0$ for all $i$) such that 
    $
    x_{2i-1}+x_{2i} = 1
    $
    for every $i\in [n]$.
    Combinatorially, this flow polytope is a cube of dimension $n$ in $\bbR^{2n}$. 
    Its vertices are of the form 
    \[
    e_{i_1}+\dots+e_{i_n},
    \]
    where $e_i\in \bbR^{2n}$ denote the standard basis vectors and $i_k$ has two possibilities, $i_k=2k-1$ or~$i_k=2k$, for each value $k\in [n]$. 
    These are precisely the indicator vectors of the routes of~$G_n$ consisting of the edges $e_{i_1},\dots,e_{i_n}$.

    For instance, the flow polytope $\calF_{G_1}$ of the oruga graph for $n=1$ is the convex hull 
    \[
    \calF_{G_1} = \conv\{
        (1,0),(0,1)
    \}.
    \]
    It is a one dimensional segment in $\bbR^{2}$, and is illustrated on the bottom left of~\Cref{fig_oruga_graph}. 
    In general, $\calF_{G_n}$ is the product of $n$ segments, an $n$-dimensional cube. For $n=2$ we get a square 
    \[
    \calF_{G_2} = \conv\{
        (1,0,1,0),(1,0,0,1),(0,1,1,0),(0,1,0,1)
    \}.
    \]
    For $n=3$ we get a cube
    \begin{align*}
    \calF_{G_3} = \conv\{ &
        (1,0,1,0,1,0),(1,0,1,0,0,1), \\
       &(1,0,0,1,1,0),(1,0,0,1,0,1), \\
       &(0,1,1,0,1,0),(0,1,1,0,0,1), \\
       &(0,1,0,1,1,0),(0,1,0,1,0,1)
    \}.        
    \end{align*}
    These flow polytopes are also illustrated on the bottom of~\Cref{fig_oruga_graph}, where we omit parentheses and commas for the coordinates of the vertices for simplicity.     
\end{example}

\subsection{framed triangulations}
A substantial amount of research on flow polytopes in the literature has been focused on their volumes and triangulations. A remarkable family of known triangulations are the framed triangulations (also called DKK triangulations), which are induced by a framing on the underlying graph. We briefly recall these triangulations in this section and refer to~\cite{DKK12} for more details.

Let $G$ be a flow graph as above. 
For each vertex $v$, let $\In(v)$ and $\Out(v)$ respectively denote the (possibly empty) sets of incoming and outgoing edges at $v$.
A \defn{framing} at the vertex $v$ is a pair of linear orders $(\leq_{\In(v)}, \leq_{\Out(v)})$ on the incoming and outgoing edges at $v$.
A \defn{framed graph}, denoted $(G,F)$, is a flow graph with a framing $F$ at every vertex.
Two different framings of the $\oru{2}$ graph are shown in~\Cref{fig_oruga2_frame_triangulations}, where the labels indicate the order of the incoming and outgoing edges at every vertex.
 
\begin{figure}[htb]
    \centering
    \begin{tikzpicture}[scale=0.8,
    fnode/.style={circle, draw, inner sep=1.4pt, fill},
    fnodesmall/.style={circle, draw, inner sep=0.7pt, fill},
    fnodeextrasmall/.style={circle, draw, inner sep=0.5pt, fill},
    Gedge/.style={thick, color=black},
    facet/.style={fill=red!95!black,fill opacity=0.800000}]

    %
    %
    \begin{scope}[shift={(0,0)}]
        %
        %
        \begin{scope}[shift={(0,0)}]
            \draw[Gedge] (0,0) .. controls (0.4, 0.6) and (1.6, 0.6) .. (2,0);
            \draw[Gedge] (0,0) .. controls (0.4, -0.6) and (1.6, -0.6) .. (2,0);    
            \draw[Gedge] (2,0) .. controls (2.4, 0.6) and (3.6, 0.6) .. (4,0);
            \draw[Gedge] (2,0) .. controls (2.4, -0.6) and (3.6, -0.6) .. (4,0);
            
            \node[fnode] (1) at (0,0) {};
            \node[fnode] (2) at (2,0) {};
            \node[fnode] (3) at (4,0) {};
        
            \node[] (a) at (1.5,0.6) {\scriptsize $1$};
            \node[] (a) at (1.5,-0.6) {\scriptsize $2$};
            \node[] (a) at (3.5,0.6) {\scriptsize $1$};    
            \node[] (a) at (3.5,-0.6) {\scriptsize $2$};    
    
            \node[] (a) at (0.5,0.6) {\scriptsize $1$};
            \node[] (a) at (0.5,-0.6) {\scriptsize $2$};
            \node[] (a) at (2.5,0.6) {\scriptsize $1$};    
            \node[] (a) at (2.5,-0.6) {\scriptsize $2$};    
        \end{scope}
        %
        %
        \begin{scope}[shift={(0,-6)},scale=4]
            \draw[thick, fill=orange] (0,0) -- (1,1) -- (1,0) -- cycle;
            \draw[thick, fill=purple] (0,0) -- (1,1) -- (0,1) -- cycle;
            \node[fnode] at (0,0) {};
            \node[fnode] at (0,1) {};
            \node[fnode] at (1,1) {};
            \node[fnode] at (1,0) {};
            %
            %
            \begin{scope}[shift={(0.55,0.4)}]
                \begin{scope}[shift={(0,0.0)}, scale=0.09]
                    \draw[Gedge] (0,0) .. controls (0.4, -0.6) and (1.6, -0.6) .. (2,0);    
                    \draw[Gedge] (2,0) .. controls (2.4, -0.6) and (3.6, -0.6) .. (4,0);

                    \node[fnodesmall] (1) at (0,0) {};
                    \node[fnodesmall] (2) at (2,0) {};
                    \node[fnodesmall] (3) at (4,0) {};
                
                \end{scope}      
                \begin{scope}[shift={(0,-0.1)}, scale=0.09]
                    \draw[Gedge] (0,0) .. controls (0.4, 0.6) and (1.6, 0.6) .. (2,0);
                    \draw[Gedge] (2,0) .. controls (2.4, -0.6) and (3.6, -0.6) .. (4,0);

                    \node[fnodesmall] (1) at (0,0) {};
                    \node[fnodesmall] (2) at (2,0) {};
                    \node[fnodesmall] (3) at (4,0) {};
                
                \end{scope}      
                \begin{scope}[shift={(0,-0.2)}, scale=0.09]
                    \draw[Gedge] (0,0) .. controls (0.4, 0.6) and (1.6, 0.6) .. (2,0);
                    \draw[Gedge] (2,0) .. controls (2.4, 0.6) and (3.6, 0.6) .. (4,0);

                    \node[fnodesmall] (1) at (0,0) {};
                    \node[fnodesmall] (2) at (2,0) {};
                    \node[fnodesmall] (3) at (4,0) {};
                
                \end{scope}      
            \end{scope}   
            %
            %
            \begin{scope}[shift={(0.1,0.8)}]
                \begin{scope}[shift={(0,0)}, scale=0.09]
                    \draw[Gedge] (0,0) .. controls (0.4, -0.6) and (1.6, -0.6) .. (2,0);    
                    \draw[Gedge] (2,0) .. controls (2.4, -0.6) and (3.6, -0.6) .. (4,0);

                    \node[fnodesmall] (1) at (0,0) {};
                    \node[fnodesmall] (2) at (2,0) {};
                    \node[fnodesmall] (3) at (4,0) {};
                
                \end{scope}      
                \begin{scope}[shift={(0,-0.1)}, scale=0.09]
                    \draw[Gedge] (0,0) .. controls (0.4, -0.6) and (1.6, -0.6) .. (2,0);    
                    \draw[Gedge] (2,0) .. controls (2.4, 0.6) and (3.6, 0.6) .. (4,0);

                    \node[fnodesmall] (1) at (0,0) {};
                    \node[fnodesmall] (2) at (2,0) {};
                    \node[fnodesmall] (3) at (4,0) {};
                
                \end{scope}      
                \begin{scope}[shift={(0,-0.2)}, scale=0.09]
                    \draw[Gedge] (0,0) .. controls (0.4, 0.6) and (1.6, 0.6) .. (2,0);
                    \draw[Gedge] (2,0) .. controls (2.4, 0.6) and (3.6, 0.6) .. (4,0);

                    \node[fnodesmall] (1) at (0,0) {};
                    \node[fnodesmall] (2) at (2,0) {};
                    \node[fnodesmall] (3) at (4,0) {};
                
                \end{scope}      
            \end{scope}   
        \end{scope}
    \end{scope}
    
    %
    %
    \begin{scope}[shift={(8,0)}]
        %
        %
        \begin{scope}[shift={(0,0)}]
            \draw[Gedge] (0,0) .. controls (0.4, 0.6) and (1.6, 0.6) .. (2,0);
            \draw[Gedge] (0,0) .. controls (0.4, -0.6) and (1.6, -0.6) .. (2,0);    
            \draw[Gedge] (2,0) .. controls (2.4, 0.6) and (3.6, 0.6) .. (4,0);
            \draw[Gedge] (2,0) .. controls (2.4, -0.6) and (3.6, -0.6) .. (4,0);
            
            \node[fnode] (1) at (0,0) {};
            \node[fnode] (2) at (2,0) {};
            \node[fnode] (3) at (4,0) {};
        
            \node[] (a) at (1.5,0.6) {\scriptsize $1$};
            \node[] (a) at (1.5,-0.6) {\scriptsize $2$};
            \node[] (a) at (3.5,0.6) {\scriptsize $1$};    
            \node[] (a) at (3.5,-0.6) {\scriptsize $2$};    
    
            \node[] (a) at (0.5,0.6) {\scriptsize $1$};
            \node[] (a) at (0.5,-0.6) {\scriptsize $2$};
            \node[] (a) at (2.5,0.6) {\scriptsize $2$};    
            \node[] (a) at (2.5,-0.6) {\scriptsize $1$};    
        \end{scope}
        %
        %
        \begin{scope}[shift={(0,-6)},scale=4]
            \draw[thick, fill=orange] (0,0) -- (0,1) -- (1,0) -- cycle;
            \draw[thick, fill=purple] (1,1) -- (0,1) -- (1,0) -- cycle;
            \node[fnode] at (0,0) {};
            \node[fnode] at (0,1) {};
            \node[fnode] at (1,1) {};
            \node[fnode] at (1,0) {};
            %
            %
            \begin{scope}[shift={(0.1,0.4)}]
                \begin{scope}[shift={(0,0)}, scale=0.09]
                    \draw[Gedge] (0,0) .. controls (0.4, -0.6) and (1.6, -0.6) .. (2,0);    
                    \draw[Gedge] (2,0) .. controls (2.4, 0.6) and (3.6, 0.6) .. (4,0);

                    \node[fnodesmall] (1) at (0,0) {};
                    \node[fnodesmall] (2) at (2,0) {};
                    \node[fnodesmall] (3) at (4,0) {};
                
                \end{scope}      
                \begin{scope}[shift={(0,-0.1)}, scale=0.09]
                    \draw[Gedge] (0,0) .. controls (0.4, 0.6) and (1.6, 0.6) .. (2,0);
                    \draw[Gedge] (2,0) .. controls (2.4, -0.6) and (3.6, -0.6) .. (4,0);

                    \node[fnodesmall] (1) at (0,0) {};
                    \node[fnodesmall] (2) at (2,0) {};
                    \node[fnodesmall] (3) at (4,0) {};
                
                \end{scope}      
                \begin{scope}[shift={(0,-0.2)}, scale=0.09]
                    \draw[Gedge] (0,0) .. controls (0.4, 0.6) and (1.6, 0.6) .. (2,0);
                    \draw[Gedge] (2,0) .. controls (2.4, 0.6) and (3.6, 0.6) .. (4,0);

                    \node[fnodesmall] (1) at (0,0) {};
                    \node[fnodesmall] (2) at (2,0) {};
                    \node[fnodesmall] (3) at (4,0) {};
                
                \end{scope}      
            \end{scope}   
            %
            %
            \begin{scope}[shift={(0.55,0.8)}]
                \begin{scope}[shift={(0,0)}, scale=0.09]
                    \draw[Gedge] (0,0) .. controls (0.4, -0.6) and (1.6, -0.6) .. (2,0);    
                    \draw[Gedge] (2,0) .. controls (2.4, -0.6) and (3.6, -0.6) .. (4,0);

                    \node[fnodesmall] (1) at (0,0) {};
                    \node[fnodesmall] (2) at (2,0) {};
                    \node[fnodesmall] (3) at (4,0) {};
                
                \end{scope}      
                \begin{scope}[shift={(0,-0.1)}, scale=0.09]
                    \draw[Gedge] (0,0) .. controls (0.4, -0.6) and (1.6, -0.6) .. (2,0);    
                    \draw[Gedge] (2,0) .. controls (2.4, 0.6) and (3.6, 0.6) .. (4,0);

                    \node[fnodesmall] (1) at (0,0) {};
                    \node[fnodesmall] (2) at (2,0) {};
                    \node[fnodesmall] (3) at (4,0) {};
                
                \end{scope}      
                \begin{scope}[shift={(0,-0.2)}, scale=0.09]
                    \draw[Gedge] (0,0) .. controls (0.4, 0.6) and (1.6, 0.6) .. (2,0);
                    \draw[Gedge] (2,0) .. controls (2.4, -0.6) and (3.6, -0.6) .. (4,0);

                    \node[fnodesmall] (1) at (0,0) {};
                    \node[fnodesmall] (2) at (2,0) {};
                    \node[fnodesmall] (3) at (4,0) {};
                
                \end{scope}      
            \end{scope}   
        \end{scope}
    \end{scope}
    
\end{tikzpicture}
    \caption{Two framings of the $G_2=\oru{2}$ graph and the framing triangulations of the corresponding flow polytope $\calF_{G_2}$.}
    \label{fig_oruga2_frame_triangulations}
\end{figure}

For a path $P$ containing a vertex $v$, let $Pv$ (resp. $vP$) denote the maximal subpath of~$P$ ending (resp. beginning) at $v$.
Furthermore, let $\scrI(v)$ (resp. $\scrO(v)$) denote the set of paths in~$G$ ending (resp. beginning) at $v$.
Our notation $\scrI$ stands for Incoming and $\scrO$ for Outgoing.
We consider $\scrI(s)$ as containing only the path of length $0$ at vertex $s$, and $\scrO(t)$ as containing only the path of length $0$ at vertex $t$. 
We define the \defn{relations $\leq_{\scrI(v)}$} and \defn{$\leq_{\scrO(v)}$} on $\scrI(v)$ and~$\scrO(v)$ as follows.

Given paths $Pv, Qv\in\scrI(v)$, let $w\leq v$ be the first vertex after which $Pv$ and $Qv$ coincide. 
If $w$ is the first vertex of $Pv$ or $Qv$, we say that $Pv=_{\scrI(v)} Qv$.
Otherwise let $e_P$ be the edge of $P$ entering $w$ and let $e_Q$ be the edge of $Q$ entering $w$.
Then $Pv <_{\scrI(v)} Qv$ if and only if~$e_P <_{\In(w)} e_Q$.
Similarly for $vP, vQ\in\scrO(v)$, let $w'\geq v$ be the last vertex before which $vP$ and $vQ$ coincide.
If $w'$ is the largest vertex of $vP$ or $vQ$, then $vP =_{\scrO(v)} vQ$.
Otherwise let $e_P'$ be the edge of $P$ leaving $w'$ and let~$e_Q'$ be the edge of $Q$ leaving $w'$.
Then $vP <_{\scrO(v)} vQ$ if and only if $e_P' <_{\Out(w')} e_Q'$.

\begin{figure}[htb]
    \centering
\begin{tikzpicture}
\begin{scope}[xshift=0]
    \draw[ultra thick,color=blue] (2.5,0.8) .. controls (2.6, 0.8) and (3.6, 0.8) .. (4,0);
    \draw[ultra thick,color=blue] (4,0.07)--(6,0.07);
    \draw[ultra thick,color=blue, dashed] (2,0.8)--(2.5,0.8);
    
    \draw[ultra thick,color=red] (2.5,-0.8) .. controls (2.6, -0.8) and (3.6, -0.8) .. (4,0);
    \draw[ultra thick,color=red] (4,-0.07)--(6,-0.07);
    \draw[ultra thick,color=red, dashed] (2,-0.8)--(2.5,-0.8);    

    \draw[ultra thick,color=black!20] (4,0)--(6,0);
    
    \node[circle, draw, inner sep=2pt, fill,label=above:\small{$w$}] (s) at (4,0) {};    
    \node[circle, draw, inner sep=2pt, fill,label=above:\small{$v$}] (c) at (6,0) {};
    
    \node[] (P) at (2.8,1) {\color{blue} $P$};
    \node[] (Q) at (2.8,-0.5) {\color{red} $Q$};
    \node[] (R) at (3.6,0) {\color{black!30} $R$};

\end{scope}
\end{tikzpicture}
    \caption{The relation $\leq_{\scrI(v)}$ is a partial order on incoming paths to $v$ only if they all begin at the source.}
    \label{fig.PQR}
\end{figure}

Note that if $Rv$ is a subpath of $Pv$, then $Rv =_{\scrI(v)} Pv$. But, if they do not start at the same vertex, then they are different paths. Therefore, the relation $\leq_{\scrI(v)}$ is not even a partial order. For example, in Figure~\ref{fig.PQR} we have $Pv =_{\scrI(v)} Rv =_{\scrI(v)} Qv$, but $Pv <_{\scrI(v)} Qv$. However, if we restrict $\leq_{\scrI(v)}$ (resp. $\leq_{\scrO(v)}$) to the set of paths starting at the source~$s$ (resp. $v$) and ending at $v$ (resp. the sink $t$), then we get a linear order.

\begin{lemma}\label{lem_linear_order}
    The following hold:
    \begin{enumerate}
        \item The restriction of $\leq_{\scrI(v)}$ to the set of paths starting at the source $s$ and ending at $v$ is a linear order. 
        \item The restriction of $\leq_{\scrO(v)}$ to the set of paths starting at $v$ and ending at the sink $t$ is a linear order. 
    \end{enumerate}
\end{lemma}
\begin{proof}
    For the proof of (1), it is straight forward to check that the restriction of $\leq_{\scrI(v)}$ to the set of paths starting at the source $s$ and ending at $v$ is reflexive, transitive and antisymmetric, which implies that the relation is a partial order.
    The transitivity is illustrated in Figure~\ref{fig.transitivity}.
    Furthermore, since every pair of paths $Pv$ and $Qv$ starting at the source~$s$ are comparable in the order $\leq_{\scrI(v)}$, we deduce that $\leq_{\scrI(v)}$ is a linear order. 
    A similar argument shows part~(2).
\end{proof}

\begin{figure}[htb]
    \centering
    \begin{tikzpicture}
\begin{scope}[xshift=0]
    \draw[ultra thick,color=blue] (2.5,0.8) .. controls (2.6, 0.8) and (3.6, 0.8) .. (4,0);
    \draw[ultra thick,color=blue] (4,0.04)--(5,0.04);
    \draw[ultra thick,color=blue] (5,0.04)--(6,0.04);
    \draw[ultra thick,color=blue, dashed] (2,0.8)--(2.5,0.8);
    
    \draw[ultra thick,color=red] (2,-0.6) .. controls (2.4, -0.6) and (3.6, -0.5) .. (4,0);
    \draw[ultra thick,color=red] (4,-0.02)--(5,-0.02);
    \draw[ultra thick,color=red] (5,-0.02)--(6,-0.02);
    \draw[ultra thick,color=red, dashed] (1.5,-0.6)--(2,-0.6);    

    \draw[ultra thick,color=black!30] (4.6,-0.08)--(6,-0.08);
    \draw[ultra thick,color=black!30] (3,-1) .. controls (3.3, -1) and (4.2, -0.7) .. (4.6,-0.08);
    \draw[ultra thick,color=black!30, dashed] (2.5,-1)--(3,-1);    
    
    \node[circle, draw, inner sep=2pt, fill,label=above:\small{$w$}] (s) at (4,0) {};
    \node[circle, draw, inner sep=2pt, fill,label=above:\small{$w'$}] (s) at (4.6,0) {};    
    \node[circle, draw, inner sep=2pt, fill,label=above:\small{$v$}] (c) at (6,0) {};
    
    \node[] (P) at (2.8,1) {\color{blue} $P$};
    \node[] (Q) at (2.8,-0.2) {\color{red} $Q$};
    \node[] (R) at (4,-1) {\color{black!30} $R$};

\end{scope}

\begin{scope}[xshift=240]
    \draw[ultra thick,color=blue] (2.5,0.8) .. controls (2.6, 0.8) and (3.6, 0.8) .. (4,0);
    \draw[ultra thick,color=blue] (4,0.04)--(5,0.04);
    \draw[ultra thick,color=blue] (5,0.04)--(6,0.04);
    \draw[ultra thick,color=blue, dashed] (2,0.8)--(2.5,0.8);
    
    \draw[ultra thick,color=red] (2,-0.6) .. controls (2.4, -0.6) and (3.6, -0.5) .. (4,0);
    \draw[ultra thick,color=red] (4,-0.02)--(5,-0.02);
    \draw[ultra thick,color=red] (5,-0.02)--(6,-0.02);
    \draw[ultra thick,color=red, dashed] (1.5,-0.6)--(2,-0.6);    

    \draw[ultra thick,color=black!30] (4,-0.08)--(6,-0.08);
    \draw[ultra thick,color=black!30] (3.2,-0.5) .. controls (3.3, -0.48) and (3.8, -0.35) .. (4,-0.08);
    \draw[ultra thick,color=black!30] (2,-1) .. controls (2.3, -1) and (2.8, -1) .. (3.2,-0.49);
    
    \draw[ultra thick,color=black!30, dashed] (1.5,-1)--(2,-1);    
    
    \node[circle, draw, inner sep=2pt, fill,label=above:\small{$w$}] (s) at (4,0) {};
    \node[circle, draw, inner sep=2pt, fill,label=above:\small{$w'$}] (s) at (3.2,-0.45) {};    
    \node[circle, draw, inner sep=2pt, fill,label=above:\small{$v$}] (c) at (6,0) {};
    
    \node[] (P) at (2.8,1) {\color{blue} $P$};
    \node[] (Q) at (2.5,-0.2) {\color{red} $Q$};
    \node[] (R) at (3,-1) {\color{black!30} $R$};

\end{scope}

\end{tikzpicture}
    \caption{Two cases in the proof of the transitivity of $\leq_{\scrI(v)}$ in~\Cref{lem_linear_order}.}
    \label{fig.transitivity}
\end{figure}

We say that a vertex $v$ of a path $P$ is an \defn{inner vertex} if $v$ is not the first or last vertex of the path.
If $v$ is an inner vertex of paths $P$ and $Q$, we say that $P$ and $Q$ are \defn{incoherent at $v$} if $Pv <_{\scrI(v)} Qv$ and $vQ <_{\scrO(v)} vP$, or if $Qv <_{\scrI(v)} Pv$ and $vP <_{\scrO(v)} vQ$, and we say that they are \defn{coherent at $v$} otherwise.
Paths $P$ and $Q$ are then said to be \defn{coherent} if they are coherent at each common inner vertex and they are \defn{incoherent} otherwise.
A set of pairwise coherent routes is called a \defn{clique}.
We denote by \defn{$\calC$} the \defn{collection of maximal cliques}. 
Examples of these concepts are illustrated in~\Cref{fig_oruga2_coherent_example}.

\begin{figure}[h]
    \centering
    \begin{tikzpicture}[scale=0.7,
    fnode/.style={circle, draw, inner sep=1.4pt, fill},
    fnodesmall/.style={circle, draw, inner sep=0.7pt, fill},
    fnodeextrasmall/.style={circle, draw, inner sep=0.5pt, fill},
    Gedge/.style={thick, color=black},
    facet/.style={fill=red!95!black,fill opacity=0.800000}]

    %
    %
    \begin{scope}[shift={(0,0)}]
        \draw[Gedge] (0,0) .. controls (0.4, 0.6) and (1.6, 0.6) .. (2,0);
        \draw[Gedge] (0,0) .. controls (0.4, -0.6) and (1.6, -0.6) .. (2,0);    
        \draw[Gedge] (2,0) .. controls (2.4, 0.6) and (3.6, 0.6) .. (4,0);
        \draw[Gedge] (2,0) .. controls (2.4, -0.6) and (3.6, -0.6) .. (4,0);
        
        \node[fnode] (1) at (0,0) {};
        \node[fnode] (2) at (2,0) {};
        \node[fnode] (3) at (4,0) {};
    
        \node[] (a) at (1.5,0.6) {\scriptsize $1$};
        \node[] (a) at (1.5,-0.6) {\scriptsize $2$};
        \node[] (a) at (3.5,0.6) {\scriptsize $1$};    
        \node[] (a) at (3.5,-0.6) {\scriptsize $2$};    

        \node[] (a) at (0.5,0.6) {\scriptsize $1$};
        \node[] (a) at (0.5,-0.6) {\scriptsize $2$};
        \node[] (a) at (2.5,0.6) {\scriptsize $1$};    
        \node[] (a) at (2.5,-0.6) {\scriptsize $2$};    
    \end{scope}

    %
    %
    \begin{scope}[shift={(6,0)}]
        \begin{scope}[shift={(0,0.6)}]
            \draw[Gedge] (0,0) .. controls (0.4, 0.6) and (1.6, 0.6) .. (2,0);
            \draw[Gedge] (2,0) .. controls (2.4, 0.6) and (3.6, 0.6) .. (4,0);
            \node[fnode] (1) at (0,0) {};
            \node[fnode] (2) at (2,0) {};
            \node[fnode] (3) at (4,0) {};
        \end{scope}
        \begin{scope}[shift={(0,-0.6)}]
            \draw[Gedge] (0,0) .. controls (0.4, -0.6) and (1.6, -0.6) .. (2,0);    
            \draw[Gedge] (2,0) .. controls (2.4, -0.6) and (3.6, -0.6) .. (4,0);       %
            \node[fnode] (1) at (0,0) {};
            \node[fnode] (2) at (2,0) {};
            \node[fnode] (3) at (4,0) {};
            \node[] at (2,-1) {Coherent};
        \end{scope}
    \end{scope}

    %
    %
    \begin{scope}[shift={(12,0)}]
        \begin{scope}[shift={(0,0.6)}]
            \draw[Gedge] (0,0) .. controls (0.4, 0.6) and (1.6, 0.6) .. (2,0);
            \draw[Gedge] (2,0) .. controls (2.4, -0.6) and (3.6, -0.6) .. (4,0);       %
            \node[fnode] (1) at (0,0) {};
            \node[fnode] (2) at (2,0) {};
            \node[fnode] (3) at (4,0) {};
        \end{scope}
        \begin{scope}[shift={(0,-0.6)}]
            \draw[Gedge] (0,0) .. controls (0.4, -0.6) and (1.6, -0.6) .. (2,0);    
            \draw[Gedge] (2,0) .. controls (2.4, 0.6) and (3.6, 0.6) .. (4,0);
            \node[fnode] (1) at (0,0) {};
            \node[fnode] (2) at (2,0) {};
            \node[fnode] (3) at (4,0) {};
            \node[] at (2,-1) {Incoherent};
        \end{scope}
    \end{scope}

    %
    %
    \begin{scope}[shift={(18,0)}]
        \begin{scope}[shift={(0,1.2)}]
            \draw[Gedge] (0,0) .. controls (0.4, 0.6) and (1.6, 0.6) .. (2,0);
            \draw[Gedge] (2,0) .. controls (2.4, 0.6) and (3.6, 0.6) .. (4,0);
            \node[fnode] (1) at (0,0) {};
            \node[fnode] (2) at (2,0) {};
            \node[fnode] (3) at (4,0) {};
        \end{scope}
        \begin{scope}[shift={(0,0)}]
            \draw[Gedge] (0,0) .. controls (0.4, 0.6) and (1.6, 0.6) .. (2,0);
            \draw[Gedge] (2,0) .. controls (2.4, -0.6) and (3.6, -0.6) .. (4,0);       %
            \node[fnode] (1) at (0,0) {};
            \node[fnode] (2) at (2,0) {};
            \node[fnode] (3) at (4,0) {};
            
        \end{scope}        
        \begin{scope}[shift={(0,-1.2)}]
            \draw[Gedge] (0,0) .. controls (0.4, -0.6) and (1.6, -0.6) .. (2,0);    
            \draw[Gedge] (2,0) .. controls (2.4, -0.6) and (3.6, -0.6) .. (4,0);       %
            \node[fnode] (1) at (0,0) {};
            \node[fnode] (2) at (2,0) {};
            \node[fnode] (3) at (4,0) {};
            \node[] at (2,-1) {A maximal clique};
        \end{scope}        
    \end{scope}
\end{tikzpicture}
    \caption{Examples of coherent routes, incoherent routes, and a maximal clique for the given framing of the $\oru{2}$ graph.}
    \label{fig_oruga2_coherent_example}
\end{figure}

Given a framed graph $(G,F)$, the convex hull of the indicator vectors of the routes in a maximal clique induced by the framing form a simplex inside the flow polytope $\calF_{G}$. 
In the example in~\Cref{fig_oruga2_coherent_example}, there are exactly two maximal cliques;  they correspond to the two triangles of the triangulation of $\calF_{G_2}$ on the left of~\Cref{fig_oruga2_frame_triangulations}. 
Changing the framing changes the coherence relation and therefore the triangles.
The right of ~\Cref{fig_oruga2_frame_triangulations} shows the two triangles corresponding to the maximal cliques induced by the other shown framing of the graph.

The motivation for the definition of a framed graph $(G,F)$ and the induced coherent relation in the set of routes of $G$, is that the collection of simplices associated to cliques give a nice triangulation of the flow polytope $\calF_G$. 
This is stated in the following proposition, where $\Delta_C$ denotes the convex hull of the indicator vectors of the routes in a maximal clique~$C$.

\begin{proposition}[{Danilov et al.~\cite{DKK12}}]
Let $(G,F)$ be a framed graph. 
The set $\{\Delta_C \mid C \in \calC\}$ is the set of the top-dimensional simplices in a regular unimodular triangulation of $\calF_G$.
\qed
\end{proposition}



\begin{corollary}\label{cor.route_rotatable_only_in_one_way}
    Every codimension 1 clique is contained in at most two maximal cliques. In particular,
    if $C$, $C' = C\setminus R\cup R'$, and $C'' = C\setminus R \cup R''$ are maximal cliques, then $C' = C''$ and $R' = R''$. \qed
\end{corollary}

Recall that the \defn{dual graph} of a triangulation is the graph whose vertices are the facets of the triangulation, with edges between facets sharing a codimension $1$ face. 
The following corollary then follows from the previous corollary.

\begin{corollary}\label{cor.triangleFree}
    The dual graph of a framed triangulation of $\calF_G$ is triangle-free. \qed 
\end{corollary}

A triangulation of $\calF_G$ whose facets are the maximal cliques of $(G,F)$ for some framing~$F$ is called a \defn{framed triangulation of $\calF_G$}. 
The routes appearing in every maximal clique of $(G,F)$ are called the \defn{exceptional routes}, and we use $\calE$ to denote the set of exceptional routes. 
The exceptional routes of $(G,F)$ correspond to cone points in the framed triangulation, i.e. points contained in every facet of the triangulation.

From now on, unless otherwise specified, we draw the framed graphs $(G,F)$ in such a way that the order of the framing of the incoming and outgoing edges at every vertex is \defn{increasing from top to bottom}; for example, as in~\Cref{fig_oruga2_frame_triangulations} (left) and not as in~\Cref{fig_oruga2_frame_triangulations} (right). 

This has two advantages. First, we do not need to include the labels of a framing for the incoming and outgoing edges to the figure because they are just ordered from top to bottom. Second, the coherence relation becomes very intuitive because two paths are coherent at a vertex $v$ if they ``\defn{do not cross}'' at $v$, as illustrated in~\Cref{fig_coherent_crossing}. 

\begin{figure}[htb]
    \centering
    \begin{tikzpicture}
\begin{scope}[xshift=0]
    \draw[ultra thick,color=blue] (2.5,0.6) .. controls (2.6, 0.6) and (3.6, 0.6) .. (4,0);
    \draw[ultra thick,color=blue] (4,+0.03)--(5,+0.03);
    \draw[ultra thick,color=blue, dashed] (2,0.6)--(2.5,0.6);
    \draw[ultra thick,color=blue] (5,0) .. controls (5.4, 0.6) and (6.4, 0.6) .. (6.5,0.6);
    \draw[ultra thick,color=blue, dashed] (6.5,0.6)--(7,0.6);
    
    \draw[ultra thick,color=red] (2.5,-0.6) .. controls (2.6, -0.6) and (3.6, -0.5) .. (4,0);
    \draw[ultra thick,color=red] (4,-0.03)--(5,-0.03);
    \draw[ultra thick,color=red, dashed] (2,-0.6)--(2.5,-0.6);    
    \draw[ultra thick,color=red] (5,0) .. controls (5.4, -0.6) and (6.4, -0.6) .. (6.5,-0.6);
    \draw[ultra thick,color=red, dashed] (6.5,-0.6)--(7,-0.6);

    \node[circle, draw, inner sep=1.6pt, fill,label=below:\small{$v$}] (s) at (4.5,0) {};    
    
    \node[] (P) at (1.6,0.6) {\color{blue} $R$};
    \node[] (Q) at (1.6,-0.6) {\color{red} $R'$};
    
    \node[] (P) at (7.4,0.6) {\color{blue} $R$};
    \node[] (Q) at (7.4,-0.6) {\color{red} $R'$};

    \node[] at (4.5,-1.3) {\small coherent at $v$};
\end{scope}

\begin{scope}[xshift=240]
    \draw[ultra thick,color=blue] (2.5,0.6) .. controls (2.6, 0.6) and (3.6, 0.6) .. (4,0);
    \draw[ultra thick,color=blue] (4,-0.03)--(5,-0.03);
    \draw[ultra thick,color=blue, dashed] (2,0.6)--(2.5,0.6);
    \draw[ultra thick,color=blue] (5,0) .. controls (5.4, -0.6) and (6.4, -0.6) .. (6.5,-0.6);
    \draw[ultra thick,color=blue, dashed] (6.5,-0.6)--(7,-0.6);
    
    \draw[ultra thick,color=red] (2.5,-0.6) .. controls (2.6, -0.6) and (3.6, -0.5) .. (4,0);
    \draw[ultra thick,color=red] (4,0.03)--(5,0.03);
    \draw[ultra thick,color=red, dashed] (2,-0.6)--(2.5,-0.6);    
    \draw[ultra thick,color=red] (5,0) .. controls (5.4, 0.6) and (6.4, 0.6) .. (6.5,0.6);
    \draw[ultra thick,color=red, dashed] (6.5,0.6)--(7,0.6);
    
    \node[circle, draw, inner sep=1.6pt, fill,label=below:\small{$v$}] (s) at (4.5,0) {};    
    
    \node[] (P) at (1.6,0.6) {\color{blue} $R$};
    \node[] (Q) at (1.6,-0.6) {\color{red} $R'$};
    
    \node[] (P) at (7.4,0.6) {\color{red} $R'$};
    \node[] (Q) at (7.4,-0.6) {\color{blue} $R$};

    \node[] at (4.5,-1.3) {\small incoherent at $v$ = crossing at $v$};
    \node[] at (4.5,-1.9) {\small $R$ is cw from $R'$ at $v$};
\end{scope}

\end{tikzpicture}    
    \caption{The coherence relation coincides with the non-crossing relation when the framing of incoming and outgoing edges increases from top to bottom at every vertex. If two routes $R,R'$ are incoherent at $v$ as shown on the right, we say that $R$ is clockwise (cw) from $R'$ at $v$.}
    \label{fig_coherent_crossing}
\end{figure}

This convention motivates the following definition.
We say that a route $R$ is \defn{clockwise (cw)} from $R'$ at $v$ 
if $Rv <_{\scrI (v)} R'v$ and $vR' <_{\scrO (v)} vR$. 
We use the \defn{notation~$R <_v^{\cw}R'$} when $R$ is cw from $R'$ at~$v$, see an illustration in~\Cref{fig_coherent_crossing} (right). 
In particular, $R$ and $R'$ are incoherent at $v$ if and only if $R <_v^{\cw} R'$ or $R' <_v^{\cw} R$.

\begin{lemma}
\label{lem.transitive}
The relation $<_v^{\cw}$ is transitive. That is, if $R <_v^{\cw} R'$ and $R' <_v^{\cw} R''$, then $R <_v^{\cw} R''$. 
\end{lemma}

\begin{example}[A framed triangulation of the oruga graph]
\label{ex.oruga_with_planar_framing}
Let $G_n=\oru{n}$ be the oruga graph from~\Cref{example_running_oruga}, and let $F$ be the framing that orders the incoming and outgoing edges of $G_n$ from top to bottom. 
The maximal cliques of $(G,F)$ are in correspondence with permutations of $[n]$ as follows.

Given a permutation $[i_1,\dots,i_n]$ of $[n]$, construct a maximal clique consisting of $n+1$ routes $R_0,\dots,R_n$, where $R_k$ is the route containing the top edges $e_{2i_j-1}$ for $1\leq j \leq k$, and the bottom edges $e_{2i_j}$ for $k< j \leq n$. That is, $R_k$ is the route with top edges at positions $i_1,\dots, i_k$ and bottom edges at the positions $i_{k+1},\dots, i_n$. An example is illustrated in~\Cref{fig_oruga2_maxclique_permutation}; see~\Cref{fig_weakorder_framinglattice} for a bigger example.  

\begin{figure}[htb]
    \centering
    \begin{tikzpicture}[scale=0.8,
    fnode/.style={circle, draw, inner sep=1.4pt, fill},
    fnodesmall/.style={circle, draw, inner sep=0.7pt, fill},
    fnodeextrasmall/.style={circle, draw, inner sep=0.5pt, fill},
    Gedge/.style={thick, color=black},
    facet/.style={fill=red!95!black,fill opacity=0.800000}]

    %
    %
    \begin{scope}[shift={(0,0)}]
        \begin{scope}[shift={(0,0.0)}]
            \draw[Gedge] (0,0) .. controls (0.4, -0.6) and (1.6, -0.6) .. (2,0);    
            \draw[Gedge] (2,0) .. controls (2.4, -0.6) and (3.6, -0.6) .. (4,0);
            \node[] at (4.5,0) {\small $R_0$};
            
            \node[fnode] (1) at (0,0) {};
            \node[fnode] (2) at (2,0) {};
            \node[fnode] (3) at (4,0) {};
        
        \end{scope}      
        \begin{scope}[shift={(0,-1.3)}]
            \draw[Gedge] (0,0) .. controls (0.4, 0.6) and (1.6, 0.6) .. (2,0);
            \draw[Gedge] (2,0) .. controls (2.4, -0.6) and (3.6, -0.6) .. (4,0);
            \node[] at (1,0.2) {\footnotesize $1$};
            %
            \node[] at (4.5,0) {\small $R_1$};
            
            \node[fnode] (1) at (0,0) {};
            \node[fnode] (2) at (2,0) {};
            \node[fnode] (3) at (4,0) {};
        
        \end{scope}      
        \begin{scope}[shift={(0,-2.6)}]
            \draw[Gedge] (0,0) .. controls (0.4, 0.6) and (1.6, 0.6) .. (2,0);
            \draw[Gedge] (2,0) .. controls (2.4, 0.6) and (3.6, 0.6) .. (4,0);
            %
            \node[] at (1,0.2) {\footnotesize $1$};
            \node[] at (3,0.2) {\footnotesize $2$};
            \node[] at (4.5,0) {\small $R_2$};
            
            \node[fnode] (1) at (0,0) {};
            \node[fnode] (2) at (2,0) {};
            \node[fnode] (3) at (4,0) {};
        
        \end{scope}    
        %
        %
        \begin{scope}[shift={(-1,-1.3)}]
            \node[] at (0,0) {$12$};
        \end{scope}
    \end{scope}   
    %
    %
    \begin{scope}[shift={(8,0)}]
        \begin{scope}[shift={(0,0)}]
            \draw[Gedge] (0,0) .. controls (0.4, -0.6) and (1.6, -0.6) .. (2,0);    
            \draw[Gedge] (2,0) .. controls (2.4, -0.6) and (3.6, -0.6) .. (4,0);
            \node[] at (4.5,0) {\small $R_0$};
            
            \node[fnode] (1) at (0,0) {};
            \node[fnode] (2) at (2,0) {};
            \node[fnode] (3) at (4,0) {};
        
        \end{scope}      
        \begin{scope}[shift={(0,-1.3)}]
            \draw[Gedge] (0,0) .. controls (0.4, -0.6) and (1.6, -0.6) .. (2,0);    
            \draw[Gedge] (2,0) .. controls (2.4, 0.6) and (3.6, 0.6) .. (4,0);
            %
            \node[] at (3,0.2) {\footnotesize $2$};
            \node[] at (4.5,0) {\small $R_1$};
            
            \node[fnode] (1) at (0,0) {};
            \node[fnode] (2) at (2,0) {};
            \node[fnode] (3) at (4,0) {};
        
        \end{scope}      
        \begin{scope}[shift={(0,-2.6)}]
            \draw[Gedge] (0,0) .. controls (0.4, 0.6) and (1.6, 0.6) .. (2,0);
            \draw[Gedge] (2,0) .. controls (2.4, 0.6) and (3.6, 0.6) .. (4,0);
            %
            \node[] at (1,0.2) {\footnotesize $1$};
            \node[] at (3,0.2) {\footnotesize $2$};
            \node[] at (4.5,0) {\small $R_2$};
            
            \node[fnode] (1) at (0,0) {};
            \node[fnode] (2) at (2,0) {};
            \node[fnode] (3) at (4,0) {};
        
        \end{scope}      
        %
        %
        \begin{scope}[shift={(-1,-1.3)}]
            \node[] at (0,0) {$21$};
        \end{scope}
    \end{scope}   

\end{tikzpicture}
    \caption{Maximal cliques of the $\oru{n}$ graph are in bijection with permutations of $[n]$.}
    \label{fig_oruga2_maxclique_permutation}
\end{figure}

It is not hard to see that the resulting set of routes is a maximal clique, and that all the maximal cliques are of this form. Therefore, the maximal simplices of the framed triangulation of $\calF_{G_n}$ induced by the framing $F$ are naturally labeled by permutations of $[n]$. Moreover, two facets are adjacent if and only if the corresponding permutations can be obtained from each other by swapping two consecutive numbers. 
Thus, the dual graph of this framed triangulation of $\calF_{G_n}$ is the Hasse diagram of the classical weak order of permutations of $[n]$.
\end{example}

\subsection{Some useful lemmas}
The purpose of this paper is to introduce and study a partial order structure whose Hasse diagram is the dual graph of a framed triangulation of a flow polytope, for any flow graph and any framing on it. 
For this, we present a toolkit of useful lemmas in this section that will be used throughout the paper. 

\begin{lemma}[Extending coherent paths]
\label{lem.path_is_extendable}
    If a path $P$ is coherent with all the routes in a clique $C$, then $P$ can be extended to a route that is coherent with all the routes in $C$. 
    In particular, if $C$ is a maximal clique, then every edge is contained in some route of $C$. 
\end{lemma}

\begin{proof}
    If $P$ is a route, then we are done, so we assume that either $P$ does not begin at~$s$ or $P$ does not end at $t$.
    We consider first the latter case.
    Let $v$ be the largest vertex of $P$ and suppose that $v \neq t$.
    Let $e$ be an edge $(v,w) \in E(G)$ where $v< w$.
    Consider $P\cup \{e\}$. 
    If $P\cup \{e\}$ is coherent with all routes in $C$, then we extend $P$ by append $e$ to it.
    If $P\cup \{e\}$ is incoherent with a route $R$ in $C$, they must be incoherent at $v$.
    If $Pv <_{\scrI (v)} Rv$ and $vR <_{\scrO (v)} vP$ we can assume that $Rv$ is maximal with respect to $\leq_{\scrI(v)}$ and then assume that $vR$ is minimal with respect to $\leq_{\scrO(v)}$. 
    Now if $w'$ is the first vertex of $R$ after $v$, then by construction, appending the edge $e'=(v,w')$ to $P$ results in a path coherent with all routes in $C$.
    Similarly, if $Rv <_{\scrI (v)} Pv$ and $vP <_{\scrO (v)} vR$, we take $R$ to be minimal with respect to $\leq_{\scrI(v)}$ and then assume that $vR$ is maximal 
    with respect to $\leq_{\scrO(v)}$.
    In this case, we again extend $P$ by the edge $e'=(v,w')$ where $w'$ is the first vertex of $R$ after $v$.
    
    In the same way, if the first vertex of $P$ is not $s$, then we can extend it by an edge towards the source.
    The argument can be repeated until $P$ has been extended to a route.
\end{proof}

Two facets in a framed triangulation of $\calF_G$ are adjacent if and only if their corresponding maximal cliques differ by a single pair of routes.
The following lemma gives a necessary condition on such a pair of routes.

\begin{lemma}\label{lem.incoherent_only_on_one_path}
    Let $C_1\neq C_2$ be two adjacent maximal cliques satisfying $C_2=(C_1\setminus R_1)\cup R_2$, then the following hold.    
    \begin{itemize}
        \item[(i)] The routes $R_1$ and $R_2$ incoherent at some vertex $v$. Furthermore, they are incoherent at every vertex in the maximal path $P_v$ in $R_1\cap R_2$ that contains $v$. 
        \item[(ii)] The routes $R_1vR_2$ and $R_2vR_1$ are contained in $C_1\cap C_2$. 
        \item[(iii)] The routes $R_1$ and $R_2$ are incoherent only on the vertices of $P_v$.
    \end{itemize}
\end{lemma}

\begin{figure}[htb]
    \centering
    \begin{tikzpicture}
\begin{scope}[xshift=0, scale=0.8]

    \draw[ultra thick,color=blue, dashed] (2,0.6)--(2.5,0.6);
    \draw[ultra thick,color=blue] (2.5,0.6) .. controls (2.6, 0.6) and (3.6, 0.6) .. (4,0);
    \draw[ultra thick,color=blue] (4,0.035)--(5,0.035);    
    \draw[ultra thick,color=red] (5,0.035)--(6,0.035);
    \draw[ultra thick,color=red] (6,0) .. controls (6.4,0.6) and (7.4, 0.6) .. (7.5,0.6);
    \draw[ultra thick,color=red, dashed] (7.5,0.6)--(8.2,0.6);
    
    \draw[ultra thick,color=red, dashed] (2,-0.6)--(2.5,-0.6);
    \draw[ultra thick,color=red] (2.5,-0.6) .. controls (2.6, -0.6) and (3.6, -0.6) .. (4,0);
    \draw[ultra thick,color=red] (4,-0.035)--(5,-0.035);
    \draw[ultra thick,color=blue] (5,-0.035)--(6,-0.035);
    \draw[ultra thick,color=blue] (6,0) .. controls (6.4, -0.6) and (7.4, -0.6) .. (7.5,-0.6);
    \draw[ultra thick,color=blue, dashed] (7.5,-0.6)--(8.2,-0.6);        
    \node[circle, draw, inner sep=1.7pt, fill,label=below:\small{$v$}] (c) at (5,0) {};

    \node[] (r1) at (1.5,0.6) {\color{blue} $R_1$};
    \node[] (r1) at (8.6,-0.6) {\color{blue} $R_1$};
    \node[] (r1) at (1.5,-0.6) {\color{red} $R_2$};
    \node[] (r1) at (8.6,0.6) {\color{red} $R_2$};

    \node[] (r1) at (5,0.6) {$\overbrace{\text{\hspace{1.5cm}}}^{P_v}$};
\end{scope}

\begin{scope}[xshift=237, scale=0.8]

    \draw[ultra thick,color=bluishgreen, dashed] (2,0.6)--(2.5,0.6);
    \draw[ultra thick,color=bluishgreen] (2.5,0.6) .. controls (2.6, 0.6) and (3.6, 0.6) .. (4,0);
    \draw[ultra thick,color=bluishgreen] (4,0.035)--(5,0.035);    
    \draw[ultra thick,color=bluishgreen] (5,0.035)--(6,0.035);
    \draw[ultra thick,color=bluishgreen] (6,0) .. controls (6.4,0.6) and (7.4, 0.6) .. (7.5,0.6);
    \draw[ultra thick,color=bluishgreen, dashed] (7.5,0.6)--(8.2,0.6);
    
    \draw[ultra thick,color=orange, dashed] (2,-0.6)--(2.5,-0.6);
    \draw[ultra thick,color=orange] (2.5,-0.6) .. controls (2.6, -0.6) and (3.6, -0.6) .. (4,0);
    \draw[ultra thick,color=orange] (4,-0.035)--(5,-0.035);
    \draw[ultra thick,color=orange] (5,-0.035)--(6,-0.035);
    \draw[ultra thick,color=orange] (6,0) .. controls (6.4, -0.6) and (7.4, -0.6) .. (7.5,-0.6);
    \draw[ultra thick,color=orange, dashed] (7.5,-0.6)--(8.2,-0.6);        
    \node[circle, draw, inner sep=1.7pt, fill,label=below:\small{$v$}] (c) at (5,0) {};

    \node[] (r1) at (1,0.6) {\color{bluishgreen} $R_1vR_2$};
    \node[] (r1) at (9,-0.6) {\color{orange} $R_2vR_1$};
    \node[] (r1) at (1,-0.6) {\color{orange} $R_2vR_1$};
    \node[] (r1) at (9,0.6) {\color{bluishgreen} $R_1vR_2$};
    \node[] (r1) at (5,0.6) {$\overbrace{\text{\hspace{1.5cm}}}^{P_v}$};
\end{scope}
\end{tikzpicture}    
    \caption{The routes $R_1$, $R_2$, $R_1vR_2$, and $R_2vR_1$ of Lemma~\ref{lem.incoherent_only_on_one_path}.}
    \label{fig.R1R2}
\end{figure}

\begin{proof}
    We first prove (i).
    If the routes $R_1$ and $R_2$ are coherent, then since $R_2$ is coherent with all the routes in $C_1\cap C_2$, we have that $C_1 \cup R_2$ is a clique larger than $C_1$. 
    However, this contradicts the maximality of $C_1$, and so $R_1$ and $R_2$ must be incoherent at some $v$. 
    It follows from the definition of incoherent routes that $R_1$ and $R_2$ are also incoherent at every vertex of $P_v$.
    We assume without loss of generality that $R_1 <_v^{\cw}R_2$ as shown in~\Cref{fig.R1R2}.

    
    Next, we prove (ii). 
    First we prove that $R_1vR_2\in C_1\cap C_2$. Since $R_1vR_2$ is different to $R_1$ and $R_2$, by Corollary~\ref{cor.route_rotatable_only_in_one_way}, it is sufficient to show that $R_1vR_2$ is coherent with every route in $C_1\cap C_2$.
    Suppose toward a contradiction that there is a route $R \in C_1\cap C_2$ that is incoherent with $R_1vR_2$ at some $x$. 
    Let $P_x$ be the maximal path in $R_1vR_2 \cap R$ containing $x$.
    We consider three cases (see~\Cref{fig_lem_incoherent_only_on_one_path}): (1) $\min(P_x)\leq \max(P_x) < v$, (2) $\min(P_x) \leq v \leq \max(P_x)$, and (3) $v < \min(P_x) \leq \max(P_x)$. 
    In case (1), $R$ is incoherent with $R_1$.
    In case (3), $R$ is incoherent with $R_2$.
    In case (2), we consider two sub-cases: (2a) if $R_1vR_2 <_x^{\cw} R$, then $R$ incoherent with $R_1$, and (2b) if $R <_x^{\cw} R_1vR_2$, then $R$ is incoherent with $R_2$. 
    However, $R$ cannot be incoherent with $R_1$ or $R_2$ since $R \in C_1$ and $R \in C_2$. 
    The proof for $R_2vR_1$ is similar, and so~(ii) follows.

\begin{figure}[htb]
        \centering
        \begin{tikzpicture}
\begin{scope}[xshift=0, scale=0.8]
    \draw[ultra thick,color=black!30,dashed] (2.4,0)--(3,0);
    \draw[ultra thick,color=black!30] (3,0) .. controls (3.1,0) and (3.4,0.05) .. (3.6,0.35);
    \draw[ultra thick,color=black!30] (4.3,0) .. controls (4.6, 0.6) and (5,0.6) .. (5.1,0.6);
    \draw[ultra thick,color=black!30,dashed] (5.2,0.6)--(5.6,0.6);

    \draw[ultra thick,color=blue,dashed] (1.8,0.6)--(2.5,0.6);
    \draw[ultra thick,color=blue] (2.5,0.6) .. controls (2.6, 0.6) and (3.6, 0.6) .. (4,0);
    \draw[ultra thick,color=blue] (4,0.035)--(5,0.035);
    \draw[ultra thick,color=blue] (5,-0.035)--(6,-0.035);
    \draw[ultra thick,color=blue] (6,0) .. controls (6.4, -0.6) and (7.4, -0.6) .. (7.5,-0.6);
    \draw[ultra thick,color=blue,dashed] (7.5,-0.6)--(8.2,-0.6);

    \draw[ultra thick,color=red,dashed] (1.8,-0.6)--(2.5,-0.6);
    \draw[ultra thick,color=red] (2.5,-0.6) .. controls (2.6, -0.6) and (3.6, -0.6) .. (4,0);
    \draw[ultra thick,color=red] (4,-0.04)--(5,-0.04);
    \draw[ultra thick,color=red] (5,0.04)--(6,0.04);
    \draw[ultra thick,color=red] (6,0) .. controls (6.4, 0.6) and (7.4,0.6) .. (7.5,0.6);
    \draw[ultra thick,color=red,dashed] (7.5,0.6)--(8.2,0.6);

    \node[circle, draw, inner sep=1.4pt, fill,label=below:\small{$v$}](c) at (5,0) {};

    \node[] (r1) at (1.3,0.6) {\small $R_1$};
    \node[] (r1) at (8.7,-0.6) {\small $R_1$};
    \node[] (r1) at (1.3,-0.6) {\small $R_2$};
    \node[] (r1) at (8.7,0.6) {\small $R_2$};

    \node[] (r1) at (2,0) {\small $R$};
    \node[] (r1) at (5.9,0.6) {\small $R$};    


    \node[] (r1) at (5,-1.6) {\small Case 1};

\end{scope}

\begin{scope}[xshift=220, scale=0.8]
    \draw[ultra thick,color=black!30] (5.7,0) .. controls (5.4, -0.6) and (5,-0.6) .. (4.9,-0.6);
    \draw[ultra thick,color=black!30,dashed] (4.8,-0.6)--(4.4,-0.6);
    \draw[ultra thick,color=black!30] (6.4,0.35) .. controls (6.6, 1.1) and (7.2,1.1) .. (7.3,1.1);
    \draw[ultra thick,color=black!30,dashed] (7.4,1.1)--(7.85,1.1);

    \draw[ultra thick,color=blue,dashed] (1.8,0.6)--(2.5,0.6);
    \draw[ultra thick,color=blue] (2.5,0.6) .. controls (2.6, 0.6) and (3.6, 0.6) .. (4,0);
    \draw[ultra thick,color=blue] (4,0.035)--(5,0.035);
    \draw[ultra thick,color=blue] (5,-0.035)--(6,-0.035);
    \draw[ultra thick,color=blue] (6,0) .. controls (6.4, -0.6) and (7.4, -0.6) .. (7.5,-0.6);
    \draw[ultra thick,color=blue,dashed] (7.5,-0.6)--(8.2,-0.6);

    \draw[ultra thick,color=red,dashed] (1.8,-0.6)--(2.5,-0.6);
    \draw[ultra thick,color=red] (2.5,-0.6) .. controls (2.6, -0.6) and (3.6, -0.6) .. (4,0);
    \draw[ultra thick,color=red] (4,-0.04)--(5,-0.04);
    \draw[ultra thick,color=red] (5,0.04)--(6,0.04);
    \draw[ultra thick,color=red] (6,0) .. controls (6.4, 0.6) and (7.4,0.6) .. (7.5,0.6);
    \draw[ultra thick,color=red,dashed] (7.5,0.6)--(8.2,0.6);

    \node[circle, draw, inner sep=1.4pt, fill,label=below:\small{$v$}](c) at (5,0) {};

    \node[] (r1) at (1.3,0.6) {\small $R_1$};
    \node[] (r1) at (8.7,-0.6) {\small $R_1$};
    \node[] (r1) at (1.3,-0.6) {\small $R_2$};
    \node[] (r1) at (8.7,0.6) {\small $R_2$};

    \node[] (r1) at (4,-0.6) {\small $R$};
    \node[] (r1) at (8.2,1.1) {\small $R$};    


    \node[] (r1) at (5,-1.6) {\small Case 3};
\end{scope}
\begin{scope}[xshift=0, yshift=-110, scale=0.8]
    \draw[ultra thick,color=black!30,dashed] (2.4,0)--(3,0);
    \draw[ultra thick,color=black!30] (3,0) .. controls (3.1,0) and (3.4,0.05) .. (3.6,0.35);
    \draw[ultra thick,color=black!30] (6.4,0.35) .. controls (6.6, 1.1) and (7.2,1.1) .. (7.3,1.1);
    \draw[ultra thick,color=black!30,dashed] (7.4,1.1)--(7.85,1.1);

    \draw[ultra thick,color=blue,dashed] (1.8,0.6)--(2.5,0.6);
    \draw[ultra thick,color=blue] (2.5,0.6) .. controls (2.6, 0.6) and (3.6, 0.6) .. (4,0);
    \draw[ultra thick,color=blue] (4,0.035)--(5,0.035);
    \draw[ultra thick,color=blue] (5,-0.035)--(6,-0.035);
    \draw[ultra thick,color=blue] (6,0) .. controls (6.4, -0.6) and (7.4, -0.6) .. (7.5,-0.6);
    \draw[ultra thick,color=blue,dashed] (7.5,-0.6)--(8.2,-0.6);

    \draw[ultra thick,color=red,dashed] (1.8,-0.6)--(2.5,-0.6);
    \draw[ultra thick,color=red] (2.5,-0.6) .. controls (2.6, -0.6) and (3.6, -0.6) .. (4,0);
    \draw[ultra thick,color=red] (4,-0.04)--(5,-0.04);
    \draw[ultra thick,color=red] (5,0.04)--(6,0.04);
    \draw[ultra thick,color=red] (6,0) .. controls (6.4, 0.6) and (7.4,0.6) .. (7.5,0.6);
    \draw[ultra thick,color=red,dashed] (7.5,0.6)--(8.2,0.6);

    \node[circle, draw, inner sep=1.4pt, fill,label=below:\small{$v$}](c) at (5,0) {};

    \node[] (r1) at (1.3,0.6) {\small $R_1$};
    \node[] (r1) at (8.7,-0.6) {\small $R_1$};
    \node[] (r1) at (1.3,-0.6) {\small $R_2$};
    \node[] (r1) at (8.7,0.6) {\small $R_2$};

    \node[] (r1) at (2,0) {\small $R$};
    \node[] (r1) at (8.2,1.1) {\small $R$};    

    \node[] (r1) at (5,-2) {\small Case 2a};
\end{scope}
\begin{scope}[xshift=220, yshift=-110, scale=0.8]
    \draw[ultra thick,color=black!30,dashed] (7.6,0)--(7,0);
    \draw[ultra thick,color=black!30] (7,0) .. controls (6.9,0) and (6.6,0.05) .. (6.4,0.35);
    \draw[ultra thick,color=black!30] (3.6,0.35) .. controls (3.4, 1.1) and (2.8,1.1) .. (2.7,1.1);
    \draw[ultra thick,color=black!30,dashed] (2.6,1.1)--(2.15,1.1);

    \draw[ultra thick,color=blue,dashed] (1.8,0.6)--(2.5,0.6);
    \draw[ultra thick,color=blue] (2.5,0.6) .. controls (2.6, 0.6) and (3.6, 0.6) .. (4,0);
    \draw[ultra thick,color=blue] (4,0.035)--(5,0.035);
    \draw[ultra thick,color=blue] (5,-0.035)--(6,-0.035);
    \draw[ultra thick,color=blue] (6,0) .. controls (6.4, -0.6) and (7.4, -0.6) .. (7.5,-0.6);
    \draw[ultra thick,color=blue,dashed] (7.5,-0.6)--(8.2,-0.6);

    \draw[ultra thick,color=red,dashed] (1.8,-0.6)--(2.5,-0.6);
    \draw[ultra thick,color=red] (2.5,-0.6) .. controls (2.6, -0.6) and (3.6, -0.6) .. (4,0);
    \draw[ultra thick,color=red] (4,-0.04)--(5,-0.04);
    \draw[ultra thick,color=red] (5,0.04)--(6,0.04);
    \draw[ultra thick,color=red] (6,0) .. controls (6.4, 0.6) and (7.4,0.6) .. (7.5,0.6);
    \draw[ultra thick,color=red,dashed] (7.5,0.6)--(8.2,0.6);

    \node[circle, draw, inner sep=1.4pt, fill,label=below:\small{$v$}](c) at (5,0) {};

    \node[] (r1) at (1.3,0.6) {\small $R_1$};
    \node[] (r1) at (8.7,-0.6) {\small $R_1$};
    \node[] (r1) at (1.3,-0.6) {\small $R_2$};
    \node[] (r1) at (8.7,0.6) {\small $R_2$};

    \node[] (r1) at (1.8,1.1) {\small $R$};
    \node[] (r1) at (8,0) {\small $R$};    

    \node[] (r1) at (5,-2) {\small Case 2b};
\end{scope}
\end{tikzpicture}
        \caption{Examples of the three cases in the proof of~\Cref{lem.incoherent_only_on_one_path}.}
        \label{fig_lem_incoherent_only_on_one_path}
\end{figure}

    For (iii), we let $y \in R_1\cap R_2$ such that $y\notin P_v$.
    Note that $R_1$ and $R_2$ must be coherent at $y$, as otherwise $R_1vR_2$ and $R_2vR_1$ would be incoherent at $y$, which is not possible by (ii).
    Hence we have (iii). 
\end{proof}

The following proposition characterizes precisely the condition under which a route in a maximal clique can be replaced by another route to form an adjacent maximal clique in a framed triangulation.

\begin{definition}
    Let $R_1<_v^{\cw} R_2$ be incoherent routes at a vertex $v$.
    We say that a route~$R$ is \defn{weakly in between $R_1$ and $R_2$ at $v$} if 
    \begin{itemize}
        \item[(i)] $R_1v \leq_{\scrI (v)} Rv\leq_{\scrI (v)} R_2v$ and $vR_2 \leq_{\scrO (v)} vR \leq_{\scrO (v)} vR_1$.        
    \end{itemize}
    Furthermore, $R$ is \defn{in between $R_1$ and $R_2$ at $v$} if one of the following conditions holds.  
    \begin{itemize}
        \item[(i)] $R_1v <_{\scrI (v)} Rv <_{\scrI (v)} R_2v$ and $vR_2 \leq_{\scrO (v)} vR \leq_{\scrO (v)} vR_1$; or
        \item[(ii)] $R_1v \leq_{\scrI (v)} Rv \leq_{\scrI (v)} R_2v$ and $vR_2 <_{\scrO (v)} vR <_{\scrO (v)} vR_1$.        
    \end{itemize}
\end{definition}

\begin{figure}[htb]
    \centering
    \begin{tikzpicture}
\begin{scope}[xshift=0, yshift=0, scale=0.8]
    \draw[ultra thick,color=black!30,dashed] (2.4,0)--(3,0);
    \draw[ultra thick,color=black!30] (3,0) .. controls (3.1,0) and (3.4,0.05) .. (3.6,0.35);
    \draw[ultra thick,color=black!30,dashed] (7.6,0)--(7,0);
    \draw[ultra thick,color=black!30] (7,0) .. controls (6.9,0) and (6.6,0.05) .. (6.4,0.35);    

    \draw[ultra thick,color=blue,dashed] (1.8,0.6)--(2.5,0.6);
    \draw[ultra thick,color=blue] (2.5,0.6) .. controls (2.6, 0.6) and (3.6, 0.6) .. (4,0);
    \draw[ultra thick,color=blue] (4,0.035)--(5,0.035);
    \draw[ultra thick,color=blue] (5,-0.035)--(6,-0.035);
    \draw[ultra thick,color=blue] (6,0) .. controls (6.4, -0.6) and (7.4, -0.6) .. (7.5,-0.6);
    \draw[ultra thick,color=blue,dashed] (7.5,-0.6)--(8.2,-0.6);

    \draw[ultra thick,color=red,dashed] (1.8,-0.6)--(2.5,-0.6);
    \draw[ultra thick,color=red] (2.5,-0.6) .. controls (2.6, -0.6) and (3.6, -0.6) .. (4,0);
    \draw[ultra thick,color=red] (4,-0.04)--(5,-0.04);
    \draw[ultra thick,color=red] (5,0.04)--(6,0.04);
    \draw[ultra thick,color=red] (6,0) .. controls (6.4, 0.6) and (7.4,0.6) .. (7.5,0.6);
    \draw[ultra thick,color=red,dashed] (7.5,0.6)--(8.2,0.6);

    \node[circle, draw, inner sep=1.4pt, fill,label=below:\small{$v$}](c) at (5,0) {};

    \node[] (r1) at (1.3,0.6) {\small $R_1$};
    \node[] (r1) at (8.7,-0.6) {\small $R_1$};
    \node[] (r1) at (1.3,-0.6) {\small $R_2$};
    \node[] (r1) at (8.7,0.6) {\small $R_2$};

    \node[] (r1) at (2.0,0) {\small $R$};
    \node[] (r1) at (8.1,0) {\small $R$};    

\end{scope}
\end{tikzpicture}
    \caption{The route $R$ is in between $R_1$ and $R_2$ at $v$. The same is true for routes $RvR_1$, $RvR_2$, $R_1vR$, and $R_2vR$.}
    \label{fig_lem_inbetween}
\end{figure}

\begin{proposition}\label{lem.rotation}
    Let $R_1$ be a route in a maximal clique $C_1$, and let $R_2$ be a route such that $R_1 <_v^{\cw} R_2$ for some $v$.
    Then $C_2=(C_1\setminus R_1)\cup R_2$ is a maximal clique if and only if the following statements hold. 
    \begin{itemize}
        \item[(i)] (The ``Top-Bottom Property'') The routes $\Top(R_1,R_2) := R_1vR_2$ and $\Bot(R_1,R_2) := R_2vR_1$ are contained in $C_1\cap C_2$.
        \item[(ii)] (The ``In Between Property'') No route in $C_1$ is in between $R_1$ and $R_2$ at $v$.
    \end{itemize} 
\end{proposition}

\begin{figure}[htb]
    \centering
        \begin{tikzpicture}[scale=0.95]
\begin{scope}[xshift=0, scale=0.8]

    \draw[ultra thick,color=blue, dashed] (2,0.6)--(2.5,0.6);
    \draw[ultra thick,color=blue] (2.5,0.6) .. controls (2.6, 0.6) and (3.6, 0.6) .. (4,0);
    \draw[ultra thick,color=blue] (4,0.035)--(5,0.035);    
    \draw[ultra thick,color=red] (5,0.035)--(6,0.035);
    \draw[ultra thick,color=red] (6,0) .. controls (6.4,0.6) and (7.4, 0.6) .. (7.5,0.6);
    \draw[ultra thick,color=red, dashed] (7.5,0.6)--(8.2,0.6);
    
    \draw[ultra thick,color=red, dashed] (2,-0.6)--(2.5,-0.6);
    \draw[ultra thick,color=red] (2.5,-0.6) .. controls (2.6, -0.6) and (3.6, -0.6) .. (4,0);
    \draw[ultra thick,color=red] (4,-0.035)--(5,-0.035);
    \draw[ultra thick,color=blue] (5,-0.035)--(6,-0.035);
    \draw[ultra thick,color=blue] (6,0) .. controls (6.4, -0.6) and (7.4, -0.6) .. (7.5,-0.6);
    \draw[ultra thick,color=blue, dashed] (7.5,-0.6)--(8.2,-0.6);        
    \node[circle, draw, inner sep=1.7pt, fill,label=below:$v$] (c) at (5,0) {};

    \node[] (r1) at (1.5,0.6) {\color{blue} $R_1$};
    \node[] (r1) at (8.6,-0.6) {\color{blue} $R_1$};
    \node[] (r1) at (1.5,-0.6) {\color{red} $R_2$};
    \node[] (r1) at (8.6,0.6) {\color{red} $R_2$};

    \node[] (r1) at (5,0.6) {$\overbrace{\text{\hspace{1.5cm}}}^{\text{\normalsize $P_v$}}$};
\end{scope}

\begin{scope}[xshift=237, scale=0.8]

    \draw[ultra thick,color=bluishgreen, dashed] (2,0.6)--(2.5,0.6);
    \draw[ultra thick,color=bluishgreen] (2.5,0.6) .. controls (2.6, 0.6) and (3.6, 0.6) .. (4,0);
    \draw[ultra thick,color=bluishgreen] (4,0.035)--(5,0.035);    
    \draw[ultra thick,color=bluishgreen] (5,0.035)--(6,0.035);
    \draw[ultra thick,color=bluishgreen] (6,0) .. controls (6.4,0.6) and (7.4, 0.6) .. (7.5,0.6);
    \draw[ultra thick,color=bluishgreen, dashed] (7.5,0.6)--(8.2,0.6);
    
    \draw[ultra thick,color=orange, dashed] (2,-0.6)--(2.5,-0.6);
    \draw[ultra thick,color=orange] (2.5,-0.6) .. controls (2.6, -0.6) and (3.6, -0.6) .. (4,0);
    \draw[ultra thick,color=orange] (4,-0.035)--(5,-0.035);
    \draw[ultra thick,color=orange] (5,-0.035)--(6,-0.035);
    \draw[ultra thick,color=orange] (6,0) .. controls (6.4, -0.6) and (7.4, -0.6) .. (7.5,-0.6);
    \draw[ultra thick,color=orange, dashed] (7.5,-0.6)--(8.2,-0.6);        
    \node[circle, draw, inner sep=1.7pt, fill,label=below:$v$] (c) at (5,0) {};

    \node[] (r1) at (0.5,0.6) {\color{bluishgreen} $\Top(R_1,R_2)$};
    \node[] (r1) at (9.7,-0.6) {\color{orange} $\Bot(R_1,R_2)$};
    \node[] (r1) at (0.5,-0.6) {\color{orange} $\Bot(R_1,R_2)$};
    \node[] (r1) at (9.7,0.6) {\color{bluishgreen} $\Top(R_1,R_2)$};
    
    \node[] (r1) at (5,0.6) {$\overbrace{\text{\hspace{1.5cm}}}^{\text{\normalsize $P_v$}}$};
\end{scope}
\end{tikzpicture}
    \caption{Routes $R_1$, $R_2$, $\Top(R_1,R_2)$, and $\Bot(R_1,R_2)$ in~\Cref{lem.rotation}.}
    \label{fig.TopBot}
\end{figure}

\begin{proof}
    First, we show that properties (i) and (ii) imply that $C_2=(C_1\setminus R_1)\cup R_2$ is a maximal clique. 
    We proceed by contradiction. Suppose that there is a route $R\in C_1\cap C_2$ that is incoherent with $R_2$ at some $x$. 
    Let $P_x$ be the maximal path in $R \cap R_2$ containing $x$.
    Similarly, as in the proof of the previous lemma, we consider three cases: (1) $\min(P_x)\leq \max(P_x) < v$, (2) $\min(P_x) \leq v \leq \max(P_x)$, and (3) $v < \min(P_x) \leq \max(P_x)$. 
    In case (1), $R$ is incoherent with $\Bot(R_1,R_2)$.
    In case (3), $R$ is incoherent with $\Top(R_1,R_2)$. 
    In case (2), note that if $R_1 <_v^{\cw} R_2<_v^{\cw} R$ then $R$ would be incoherent with $R_1$, which is not possible. Therefore we can assume that $R<_v^{\cw} R_2$, that is 
    \begin{align*}
    Rv <_{\scrI (v)} R_2v      \,  \text{ and } \,
    vR_2 <_{\scrO (v)} vR.
    \end{align*}
    We consider three sub-cases (see~\Cref{fig_lem_rotation}): 
    
    (2a) $Rv <_{\scrI (v)} R_1v$

    (2b) $vR_1 <_{\scrO (v)} vR$

    (2c) $R_1v \leq_{\scrI (v)} Rv <_{\scrI (v)} R_2v$ and $vR_2 <_{\scrO (v)} vR \leq_{\scrO (v)} vR_1$

    \smallskip
    \noindent
    In case (2a), $R$ is incoherent with $\Top(R_1,R_2)$. 
    In case (2b), $R$ is incoherent with $\Bot(R_1,R_2)$.
    In case (2c), one of the two weak inequalities must be strict (otherwise $R=R_1$), and therefore $R$ would be a route in between $R_1$ and $R_2$, which is not possible by assumption of the In Between Property (ii).    
    In either case we have a contradiction, since $R$ can not be incoherent with $\Top(R_1,R_2)$ or $\Bot(R_1,R_2)$ because of the Top-Bottom Property (i).   
    This finishes the proof of the backward direction.

\begin{figure}[htb]
        \centering
        \begin{tikzpicture}
\begin{scope}[xshift=0, yshift=0, scale=0.8]
    \draw[ultra thick,color=black!30,dashed] (7.6,0)--(7,0);
    \draw[ultra thick,color=black!30] (7,0) .. controls (6.9,0) and (6.6,0.05) .. (6.4,0.35);
    \draw[ultra thick,color=black!30] (3.6,0.35) .. controls (3.4, 1.1) and (2.8,1.1) .. (2.7,1.1);
    \draw[ultra thick,color=black!30,dashed] (2.6,1.1)--(2.15,1.1);

    \node[] (r1) at (1.9,1.1) {\small $R$};
    \node[] (r1) at (7.9,0) {\small $R$};

    \draw[ultra thick,color=blue,dashed] (1.8,0.6)--(2.5,0.6);
    \draw[ultra thick,color=blue] (2.5,0.6) .. controls (2.6, 0.6) and (3.6, 0.6) .. (4,0);
    \draw[ultra thick,color=blue] (4,0.035)--(5,0.035);
    \draw[ultra thick,color=blue] (5,-0.035)--(6,-0.035);
    \draw[ultra thick,color=blue] (6,0) .. controls (6.4, -0.6) and (7.4, -0.6) .. (7.5,-0.6);
    \draw[ultra thick,color=blue,dashed] (7.5,-0.6)--(8.2,-0.6);

    \draw[ultra thick,color=red,dashed] (1.8,-0.6)--(2.5,-0.6);
    \draw[ultra thick,color=red] (2.5,-0.6) .. controls (2.6, -0.6) and (3.6, -0.6) .. (4,0);
    \draw[ultra thick,color=red] (4,-0.04)--(5,-0.04);
    \draw[ultra thick,color=red] (5,0.04)--(6,0.04);
    \draw[ultra thick,color=red] (6,0) .. controls (6.4, 0.6) and (7.4,0.6) .. (7.5,0.6);
    \draw[ultra thick,color=red,dashed] (7.5,0.6)--(8.2,0.6);

    \node[circle, draw, inner sep=1.4pt, fill,label=below:\small{$v$}](c) at (5,0) {};

    \node[] (r1) at (1.4,0.6) {\small $R_1$};
    \node[] (r1) at (8.5,-0.6) {\small $R_1$};
    \node[] (r1) at (1.4,-0.6) {\small $R_2$};
    \node[] (r1) at (8.5,0.6) {\small $R_2$};    

    \node[] (r1) at (5,-2) {Case (2a)};
\end{scope}

\begin{scope}[xshift=220, yshift=0, scale=0.8]
    \draw[ultra thick,color=black!30,dashed] (2.4,0)--(3,0);
    \draw[ultra thick,color=black!30] (3,0) .. controls (3.1,0) and (3.4,-0.05) .. (3.6,-0.35);
    \draw[ultra thick,color=black!30] (6.4,-0.35) .. controls (6.6, -1.1) and (7.2,-1.1) .. (7.3,-1.1);
    \draw[ultra thick,color=black!30,dashed] (7.4,-1.1)--(7.85,-1.1);

    \draw[ultra thick,color=blue,dashed] (1.8,0.6)--(2.5,0.6);
    \draw[ultra thick,color=blue] (2.5,0.6) .. controls (2.6, 0.6) and (3.6, 0.6) .. (4,0);
    \draw[ultra thick,color=blue] (4,0.035)--(5,0.035);
    \draw[ultra thick,color=blue] (5,-0.035)--(6,-0.035);
    \draw[ultra thick,color=blue] (6,0) .. controls (6.4, -0.6) and (7.4, -0.6) .. (7.5,-0.6);
    \draw[ultra thick,color=blue,dashed] (7.5,-0.6)--(8.2,-0.6);

    \draw[ultra thick,color=red,dashed] (1.8,-0.6)--(2.5,-0.6);
    \draw[ultra thick,color=red] (2.5,-0.6) .. controls (2.6, -0.6) and (3.6, -0.6) .. (4,0);
    \draw[ultra thick,color=red] (4,-0.04)--(5,-0.04);
    \draw[ultra thick,color=red] (5,0.04)--(6,0.04);
    \draw[ultra thick,color=red] (6,0) .. controls (6.4, 0.6) and (7.4,0.6) .. (7.5,0.6);
    \draw[ultra thick,color=red,dashed] (7.5,0.6)--(8.2,0.6);

    \node[circle, draw, inner sep=1.4pt, fill,label=below:\small{$v$}](c) at (5,0) {};

    \node[] (r1) at (1.4,0.6) {\small $R_1$};
    \node[] (r1) at (8.5,-0.6) {\small $R_1$};
    \node[] (r1) at (1.4,-0.6) {\small $R_2$};
    \node[] (r1) at (8.5,0.6) {\small $R_2$};

    \node[] (r1) at (2.1,0) {\small $R$};
    \node[] (r1) at (8.1,-1.1) {\small $R$};    

    \node[] (r1) at (5,-2) {Case (2b)};
\end{scope}

\begin{scope}[xshift=110, yshift=-110, scale=0.8]
    \draw[ultra thick,color=black!30,dashed] (2.4,0)--(3,0);
    \draw[ultra thick,color=black!30] (3,0) .. controls (3.1,0) and (3.4,-0.05) .. (3.6,-0.35);
    \draw[ultra thick,color=black!30,dashed] (7.6,0)--(7,0);
    \draw[ultra thick,color=black!30] (7,0) .. controls (6.9,0) and (6.6,0.05) .. (6.4,0.35);

    \draw[ultra thick,color=blue,dashed] (1.8,0.6)--(2.5,0.6);
    \draw[ultra thick,color=blue] (2.5,0.6) .. controls (2.6, 0.6) and (3.6, 0.6) .. (4,0);
    \draw[ultra thick,color=blue] (4,0.035)--(5,0.035);
    \draw[ultra thick,color=blue] (5,-0.035)--(6,-0.035);
    \draw[ultra thick,color=blue] (6,0) .. controls (6.4, -0.6) and (7.4, -0.6) .. (7.5,-0.6);
    \draw[ultra thick,color=blue,dashed] (7.5,-0.6)--(8.2,-0.6);

    \draw[ultra thick,color=red,dashed] (1.8,-0.6)--(2.5,-0.6);
    \draw[ultra thick,color=red] (2.5,-0.6) .. controls (2.6, -0.6) and (3.6, -0.6) .. (4,0);
    \draw[ultra thick,color=red] (4,-0.04)--(5,-0.04);
    \draw[ultra thick,color=red] (5,0.04)--(6,0.04);
    \draw[ultra thick,color=red] (6,0) .. controls (6.4, 0.6) and (7.4,0.6) .. (7.5,0.6);
    \draw[ultra thick,color=red,dashed] (7.5,0.6)--(8.2,0.6);

    \node[circle, draw, inner sep=1.4pt, fill,label=below:\small{$v$}](c) at (5,0) {};

    \node[] (r1) at (1.4,0.6) {\small $R_1$};
    \node[] (r1) at (8.5,-0.6) {\small $R_1$};
    \node[] (r1) at (1.4,-0.6) {\small $R_2$};
    \node[] (r1) at (8.5,0.6) {\small $R_2$};

    \node[] (r1) at (2.1,0) {\small $R$};
    \node[] (r1) at (7.9,0) {\small $R$};    


    \node[] (r1) at (5,-1.6) {Case (2c)};
\end{scope}
\end{tikzpicture}
        \caption{Examples of the cases (2a), (2b) and (2c) in the proof of~\Cref{lem.rotation}.}
        \label{fig_lem_rotation}
\end{figure}

    For the proof of the other direction assume that $C_2=(C_1\setminus R_1)\cup R_2$ is a maximal clique. By~\Cref{lem.incoherent_only_on_one_path} (ii), the Top-Bottom Property (i) in this lemma holds. 
    For part (ii), consider a route $R$ in between $R_1$ and $R_2$ at $v$. Then, $R$ is necessarily incoherent with either $R_1$ or~$R_2$ at $v$. But every route in $C_1\cap C_2$ is coherent with $R_1$ and $R_2$, and so $R \notin C_1\cap C_2$. Since~$R\neq R_1$, we have that $R\notin C_1$. This proves the In Between Property (ii).
\end{proof}

\begin{remark}\label{rem_notInBetween}
    Under the conditions of~\Cref{lem.rotation}, no route in $C_1$ is between~$R_1$ and~$R_2$. However, we remark that this does not necessarily holds for routes not in~$C_1$. 
    Figure~\ref{fig_notInBetween} shows two maximal cliques $C_1$ and $C_2$ satisfying 
    $C_2=(C_1\setminus R_1)\cup R_2$ with 
    $R_1 <_v^{\cw} R_2$, 
    and a route $R'\notin C_1$ in between $R_1$ and $R_2$.
    In this case, the route $R''\in C_1$ is incoherent with $R'$. 
\end{remark}

\begin{figure}[htb]
    \centering

    \caption{Example of a route $R'\notin C_1$ that is in between the routes $R_1$ and $R_2$ involved in a rotation.}
    \label{fig_notInBetween}
\end{figure}

\begin{figure}[h]
    \centering
    %
    \caption{The weak order as a framing lattice.}
    \label{fig_weakorder_framinglattice}
\end{figure}

\section{The Framing poset}\label{sec_framing_poset}

In this section, we introduce the main object of study of this paper, the framing poset. This poset is indexed by a framed graph $(G,F)$, and its Hasse diagram is dual to the corresponding framed triangulation. 

\subsection{The framing poset}
Let $C\neq C'$ be two maximal cliques in $(G,F)$ such that $C'=(C\setminus R)\cup R'$. 
Then $R$ and $R'$ are incoherent at some vertex $v$.
If $R <_v^{\cw} R'$, then we say that $R'$ is obtained from $R$ by a ccw rotation at $v$. 
In this case we also say that $C'$ is obtained from $C$ by a ccw rotation.
Define the cover relation $C \precccwrot C'$ if $C'$ is obtainable from $C$ by a ccw rotation. 
The \defn{framing poset} $\scrL_{G,F} = (\calC, \leqccwrot)$ is the poset of ccw rotations of maximal cliques induced by the transitive closure of the cover relation $\precccwrot$.
We simply write $C \leq C'$ when the partial order is clear from context.
Reflexivity and transitivity are immediate, and antisymmetry follows from the next lemma.

\begin{lemma}\label{lem.ccwExists}
    Let $C$ and $C'$ be maximal cliques such that $C\leq C'$ and let $R$ be a route in $C$.
    If $R^*$ is a route in $(G,F)$ satisfying $R^* <_v^{\cw} R$ for some $v$, then there exists a route $R' \in C'$ such that  $R^* <_w^{\cw} R'$ at some node $w$ 
    in the maximal path in $R\cap R^*$ containing $v$.
\end{lemma}
\begin{proof}
    Let $P$ be the maximal path in $R \cap R^*$ containing $v$.
    Let $e$ be the incoming edge of $R$ to $P$ and let $f$ be the outgoing edge of $P$ to $R$.
    Let $R_2$ be a ccw rotation of a route~$R_1$ containing~$ePf$.
    If $R_1 <_w^{\cw} R_2$ for some $w \notin P$, then $\Top(R_1,R_2)$ or $\Bot(R_1,R_2)$ contains $ePf$. 
    If $R_1 <_w^{\cw} R_2$ for some $w \in P$, then let $P_w$ be the maximal path in $R_1\cap R_2$ containing~$w$, and let $w_1=\min P_w$ and $w_2=\max P_w$. 
    We consider three cases: 
    (1) If $w_1\in P$, then~$R^*<_{w_1}^{\cw} R_1 <_{w_1}^{\cw} R_2$. 
    (2) If $w_2\in P$, then~$R^*<_{w_2}^{\cw} R_1 <_{w_2}^{\cw} R_2$.
    (3) If $w_1\notin P$ and $w_2\notin P$, then $ePf$ is contained in $P_w$, and so both  $\Top(R_1,R_2)$ and $\Bot(R_1,R_2)$ contain $ePf$.
    It follows that any maximal clique of routes obtained from $C$ by ccw rotations must contain a route $R'$ for which $R^* <_w^{\cw} R'$ for some node $w\in P$.     
\end{proof}

\begin{corollary}
    If $C \leq C'$ and $C' \leq C$, then $C = C'$. 
\end{corollary}

\begin{proof}
    Suppose toward a contradiction that $C\neq C'$. 
    Then, there are routes $R \in C$ and $R'\in C'$ such that $R' <_v^{\cw} R$ or $R <_v^{\cw} R'$ at some $v$.
    If $R' <_v^{\cw} R$, then since $C \leq C'$ it follows by Lemma~\ref{lem.ccwExists} that there is a route $R'' \in C'$ such that $R' <_w^{\cw} R''$ at some node $w$.
    However, then $R'$ and $R''$ are incohrent routes in $C'$, which is a contradiction.
    If $R <_v^{\cw} R'$, then we similarly reach a contradiction.
\end{proof}

\begin{corollary}
    $\scrL_{G,F}$ is a poset. \qed 
\end{corollary}

\begin{example}[The weak order]\label{ex_weak_via_framing}
    Let $G_n = \oru{n}$ from Example~\ref{example_running_oruga} and let $F$ be the framing of $G_n$ in Example~\ref{ex.oruga_with_planar_framing}. 
    The poset $\scrL_{G_n,F}$ is the weak order on permutations of length $n$. 
    The case when $n=3$ is shown in Figure~\ref{fig_weakorder_framinglattice}.
\end{example}

\begin{example}
    Let $G$ be the graph $\oru{3}$ but with the added edge $(1,3)$.
    Let the framing~$F$ be induced by the drawing in Figure~\ref{fig.framing_poset_example}. 
    The poset $\scrL_{G,F}$ is shown on the right of Figure~\ref{fig.framing_poset_example}
\end{example}

\begin{figure}[h]
    \centering

    \caption{A framed graph $(G,F)$ and its framing poset.}
    \label{fig.framing_poset_example}
\end{figure}

The length of an edge $(i,j)$ is $|j-i|$. 
For a graph $G$ define the \defn{length framing} as the framing induced by ordering the incoming and outgoing edges at each vertex in increasing order of length.

\begin{example}[The Tamari lattice] \label{ex_caracol3_Tamari}
    The \defn{caracol graph} $\car(n)$ is the path graph on $n$ vertices with the added edges $(1,i)$ for $2 < i \leq n-1$ and $(j,n)$ for $2 \leq j < n-1$. 
    Let~$F$ be the length framing of $\car(n)$, as drawn on the left of Figure~\ref{fig.caracol3}.
    The framing poset $\scrL_{\car(n),F}$ is the \defn{Tamari lattice} $\Tam(n-3)$, whose elements are triangulations of a convex $(n-1)$-gon and cover relations are increasing slope diagonal flips. 

    This connection can be argued as follows. Each route of the $(\car(n),F)$ is uniquely determined by its first and last edge. Label the outgoing edges of the source from top to bottom by the numbers $1,2,3,\dots, n-2$, and the incoming edges of the sink from top to bottom by $n-1,\dots, 2$. Label each route by the pair $(i,j)$ of its first and last edges. This gives a correspondence between routes of the graph and edges and diagonals of the $(n-1)$-gon, with the exceptional routes corresponding to the edges of the $(n-1)$-gon. Under this correspondence, two routes are incoherent if and only if the corresponding segments in the $(n-1)$-gon cross, so maximal cliques of $(\car(n),F)$ correspond to triangulations of the $(n-1)$-gon. Moreover, $\ccw$ rotations correspond to increasing slope diagonal flips. Thus, the framing poset $\scrL_{\car(n),F}$ is the Tamari lattice $\Tam(n-3)$.
\end{example}

\begin{example}[The Dyck lattice] \label{ex_caracol3_Dyck}
    Consider the caracol graph $\car(n)$ as in the previous example, but with the framing $\widetilde{F}$ obtained from length framing $F$ by reversing the order of the incoming edges at each inner vertex. 
    The case $n=6$ is depicted on the right in Figure~\ref{fig.caracol3}. 
    The framing poset~$\scrL_{\car(n),\widetilde{F}}$ is the \defn{Dyck lattice} $\Dyck(n-3)$, whose elements are lattice paths in an $(n-3)\times (n-3)$ grid using north and east steps from the left-bottom corner to the top-right corner that stay weakly above the main diagonal, i.e. \defn{Dyck paths} in an $(n-3) \times (n-3)$ square. A path $\pi_2$ covers $\pi_1$ in the Dyck lattice if $\pi_2$ can be obtained from~$\pi_1$ by adding one box.  

    The bijection between maximal cliques of $(\car(n),\widetilde{F})$ and Dyck paths in an $(n-3)\times (n-3)$ grid is given as follows. 
    As before, the routes of $(\car(n),\widetilde{F})$ are determined by their first and last edges.
    Label the outgoing edges of the source from bottom to top by the numbers $0,1,\dots, n-3$, and the incoming edges of the sink from top to bottom by $n-3,\dots, 0$. Label a route by the pair $(i,j)$ of its first and last edges. Such a pair satisfies $0\leq i\leq j \leq n-3$.
    Therefore, the pair $(i,j)$ can also be regarded as a lattice point weakly above the main diagonal of the $(n-3)\times (n-3)$ square, where $i$ denotes the column number (from left to right) and $j$ the row number (from bottom to top).  
    Thus, we have a bijection between routes of $(\car(n),\widetilde{F})$ and lattice points weakly above the main diagonal of an $(n-3)\times (n-3)$ grid.
    Taking the image of the routes of a maximal clique gives the set of lattice points of a Dyck path. Moreover, two maximal cliques are related by a $\ccw$ rotation if and only if the corresponding two paths are related by adding one box. 
    As a consequence, the framing poset~$\scrL_{\car(n),\widetilde{F}}$ is the \defn{Dyck lattice} $\Dyck(n-3)$.
    
\end{example}

\begin{remark}
    The caracol graph $\car(n)$ has been considered several times in the literature in connection with flow polytopes and Catalan structures, see for example \cite{BGMY23, BGHHKMY19, MM19}.
    Its name comes from the resemblance of the drawing of $(\car(n),\widetilde{F})$ to a snail, which is referred to as ``caracol'' in Spanish (see the right of Figure~\ref{fig.caracol3}).
\end{remark}

\begin{remark}
    The Tamari lattice and Dyck lattice are particular cases of a more general family of lattices that we call cross-Tamari lattices, which are studied in depth in~\Cref{sec_crossTamari}. 
    The graphs we use in~\Cref{sec_crossTamari} are slightly different to the caracol graph, and have the advantage of making the description a bit clearer. 
\end{remark}

\begin{figure}[h]
    \centering

    \caption{The Tamari lattice and the Dyck lattice as framing lattices of the caracol graph $\car(6)$ (with two different framings).}
    \label{fig.caracol3}
\end{figure}

\begin{example}
    Consider the directed complete graph $K_6$ with edges directed toward the larger vertex. 
    Its corresponding framing poset under the length framing is shown in~\Cref{fig.K6}.
\end{example}

\begin{figure}
    \centering
    %
    \caption{The graph $K_6$ drawn with the length framing, its six exceptional routes, nine non-exceptional routes, and framing poset.}
    \label{fig.K6}
\end{figure}

\subsection{Some graph operations}

There are two natural operations that we can apply to a directed graph $G$ and a framing~$F$ of it. We denote by $G^{\rev}$ the directed acyclic graph~$G$ but with the direction of its edges reversed. The framing $F$ of~$G$ naturally induces a framing on $G^{\rev}$, which by abuse of notation we also denote by $F$, producing a framed graph $(G^{\rev},F)$. 
On the other hand, we denote by $F^{\rev}$ the framing of $G$ obtained from $F$ by reversing the orders of the incoming and outgoing edges at each vertex. 
These operations and the corresponding framing posets are illustrated for an example in~\Cref{fig.graph_operations_ex}.

\begin{lemma}
    \label{lem.dualizing_operations}
    The following hold:
    \begin{itemize}
        \item[(1)] $\scrL_{G^{\rev},F} \cong \mathrm{dual}(\scrL_{G,F})$;
        \item[(2)] $\scrL_{G,F^{\rev}} \cong \mathrm{dual}(\scrL_{G,F})$; and
        \item[(3)] $\scrL_{G^{\rev},F^{\rev}} \cong \scrL_{G,F}$.
    \end{itemize}
\end{lemma}

\begin{proof}
    Observe that reversing the flow in $G$ preserves the pairwise coherence of routes.
    Thus it does not change the maximal cliques beyond reversing the direction of the routes in the maximal cliques.
    However, counterclockwise rotations of maximal cliques in $(G,F)$ become clockwise rotations in $(G^{\rev},F)$ (and vice versa), and hence (1) follows.

    Reversing the orders of the incoming and outgoing edges at each vertex preserves the pairwise coherence of routes, and hence the maximal cliques in $(G,F)$ and $(G,F^{\rev})$ are the same.
    The counterclockwise rotations of maximal cliques in $(G,F)$ then become clockwise rotations in $(G,F^{\rev})$ (and vice versa), and hence (2) follows.

    Finally, (3) follows from (1) and (2).
\end{proof}

\begin{figure}
    \centering
    \begin{tikzpicture}[scale=1]

\begin{scope}[scale=0.4, xshift=0, yshift=80]
    \draw[thick, color=black] (0,0) .. controls (0.4, 0.6) and (1.6, 0.6) .. (2,0);
    \draw[thick, color=black] (0,0) .. controls (0.4, -0.6) and (1.6, -0.6) .. (2,0);    
    \draw[thick, color=black] (0,0) .. controls (0.4, 1.3) and (3.6, 1.3) .. (4,0);       
    \draw[thick, color=black] (2,0) .. controls (2.4, 0.6) and (3.6, 0.6) .. (4,0);
    \draw[thick, color=black] (2,0) .. controls (2.4, -0.6) and (3.6, -0.6) .. (4,0);
    \draw[thick, color=black] (4,0) .. controls (4.4, 0.6) and (5.6, 0.6) .. (6,0);
    \draw[thick, color=black] (4,0) .. controls (4.4, -0.6) and (5.6, -0.6) .. (6,0);        
    
    \node[circle, draw, inner sep=1.4pt, fill] (1) at (0,0) {};
    \node[circle, draw, inner sep=1.4pt, fill] (2) at (2,0) {};
    \node[circle, draw, inner sep=1.4pt, fill] (3) at (4,0) {};
    \node[circle, draw, inner sep=1.4pt, fill] (4) at (6,0) {};    

    \node[] at (3,-1.8) {\footnotesize $(G,F)$};   
\end{scope}	

\begin{scope}[scale=0.25, xshift=150, yshift=-320]
	\draw[thick, color=NavyBlue]    (0,0) -- (4,3) -- (4,7) -- (0,10) -- (-4,7) -- (-4,3) -- (0,0);		
	\node[circle,fill,inner sep=1.5pt,color=NavyBlue] at (0,0)  {};
	\node[circle,fill,inner sep=1.5pt,color=NavyBlue]  at (4,3) {};
	\node[circle,fill,inner sep=1.5pt,color=NavyBlue]  at (4,7) {}; 
	\node[circle,fill,inner sep=1.5pt,color=NavyBlue]  at (0,10) {};
	\node[circle,fill,inner sep=1.5pt,color=NavyBlue]  at (-4,7) {};
	\node[circle,fill,inner sep=1.5pt,color=NavyBlue]  at (-4,3) {};
	\node[circle,fill,inner sep=1.5pt,color=NavyBlue]  at (-7,-1) {};
 	\node[circle,fill,inner sep=1.5pt,color=NavyBlue]  at (-3,-4) {};

    \draw[thick, color=NavyBlue]    (-4,3) -- (-7,-1) -- (-3,-4) -- (0,0);  
    \node[] at (0,-6.5) {$\mathscr{L}_{G,F}$}; 
\end{scope}

\begin{scope}[xshift=100]
\begin{scope}[scale=0.4, xshift=0, yshift=80]
    \draw[thick, color=black] (0,0) .. controls (0.4, 0.6) and (1.6, 0.6) .. (2,0);
    \draw[thick, color=black] (0,0) .. controls (0.4, -0.6) and (1.6, -0.6) .. (2,0);    
    \draw[thick, color=black] (2,0) .. controls (2.4, 1.3) and (5.6, 1.3) .. (6,0);       
    \draw[thick, color=black] (2,0) .. controls (2.4, 0.6) and (3.6, 0.6) .. (4,0);
    \draw[thick, color=black] (2,0) .. controls (2.4, -0.6) and (3.6, -0.6) .. (4,0);
    \draw[thick, color=black] (4,0) .. controls (4.4, 0.6) and (5.6, 0.6) .. (6,0);
    \draw[thick, color=black] (4,0) .. controls (4.4, -0.6) and (5.6, -0.6) .. (6,0);        
    
    \node[circle, draw, inner sep=1.4pt, fill] (1) at (0,0) {};
    \node[circle, draw, inner sep=1.4pt, fill] (2) at (2,0) {};
    \node[circle, draw, inner sep=1.4pt, fill] (3) at (4,0) {};
    \node[circle, draw, inner sep=1.4pt, fill] (4) at (6,0) {};    

    \node[] at (3,-1.8) {\footnotesize $(G^{\rev},F)$};   
\end{scope}
\begin{scope}[scale=0.25, xshift=180, yshift=-420]
	\draw[thick, color=NavyBlue]    (0,0) -- (4,3) -- (4,7) -- (0,10) -- (-4,7) -- (-4,3) -- (0,0);		
	\node[circle,fill,inner sep=1.5pt,color=NavyBlue] at (0,0)  {};
	\node[circle,fill,inner sep=1.5pt,color=NavyBlue]  at (4,3) {};
	\node[circle,fill,inner sep=1.5pt,color=NavyBlue]  at (4,7) {}; 
	\node[circle,fill,inner sep=1.5pt,color=NavyBlue]  at (0,10) {};
	\node[circle,fill,inner sep=1.5pt,color=NavyBlue]  at (-4,7) {};
	\node[circle,fill,inner sep=1.5pt,color=NavyBlue]  at (-4,3) {};
	\node[circle,fill,inner sep=1.5pt,color=NavyBlue]  at (-7,11) {};
 	\node[circle,fill,inner sep=1.5pt,color=NavyBlue]  at (-3,14) {};

    \draw[thick, color=NavyBlue]    (-4,7) -- (-7,11) -- (-3,14) -- (0,10);  
    
    \node[] at (0,-3) {$\scrL_{G^{\rev},F}$}; 
\end{scope}
\end{scope}

\begin{scope}[xshift=200]
\begin{scope}[scale=0.4, xshift=0, yshift=80]
    \draw[thick, color=black] (0,0) .. controls (0.4, 0.6) and (1.6, 0.6) .. (2,0);
    \draw[thick, color=black] (0,0) .. controls (0.4, -0.6) and (1.6, -0.6) .. (2,0);    
    \draw[thick, color=black] (0,0) .. controls (0.4, -1.3) and (3.6, -1.3) .. (4,0);       
    \draw[thick, color=black] (2,0) .. controls (2.4, 0.6) and (3.6, 0.6) .. (4,0);
    \draw[thick, color=black] (2,0) .. controls (2.4, -0.6) and (3.6, -0.6) .. (4,0);
    \draw[thick, color=black] (4,0) .. controls (4.4, 0.6) and (5.6, 0.6) .. (6,0);
    \draw[thick, color=black] (4,0) .. controls (4.4, -0.6) and (5.6, -0.6) .. (6,0);        
    
    \node[circle, draw, inner sep=1.4pt, fill] (1) at (0,0) {};
    \node[circle, draw, inner sep=1.4pt, fill] (2) at (2,0) {};
    \node[circle, draw, inner sep=1.4pt, fill] (3) at (4,0) {};
    \node[circle, draw, inner sep=1.4pt, fill] (4) at (6,0) {};    

    \node[] at (3,-1.8) {\footnotesize $(G,F^{\rev})$};   
\end{scope}
\begin{scope}[scale=0.25, xshift=180, yshift=-420]
	\draw[thick, color=NavyBlue]    (0,0) -- (4,3) -- (4,7) -- (0,10) -- (-4,7) -- (-4,3) -- (0,0);		
	\node[circle,fill,inner sep=1.5pt,color=NavyBlue] at (0,0)  {};
	\node[circle,fill,inner sep=1.5pt,color=NavyBlue]  at (4,3) {};
	\node[circle,fill,inner sep=1.5pt,color=NavyBlue]  at (4,7) {}; 
	\node[circle,fill,inner sep=1.5pt,color=NavyBlue]  at (0,10) {};
	\node[circle,fill,inner sep=1.5pt,color=NavyBlue]  at (-4,7) {};
	\node[circle,fill,inner sep=1.5pt,color=NavyBlue]  at (-4,3) {};
	\node[circle,fill,inner sep=1.5pt,color=NavyBlue]  at (-7,11) {};
 	\node[circle,fill,inner sep=1.5pt,color=NavyBlue]  at (-3,14) {};

    \draw[thick, color=NavyBlue]    (-4,7) -- (-7,11) -- (-3,14) -- (0,10);  

    \node[] at (0,-3) {$\scrL_{G,F^{\rev}}$}; 
\end{scope}
\end{scope}

\begin{scope}[xshift=300]
\begin{scope}[scale=0.4, xshift=0, yshift=80]
    \draw[thick, color=black] (0,0) .. controls (0.4, 0.6) and (1.6, 0.6) .. (2,0);
    \draw[thick, color=black] (0,0) .. controls (0.4, -0.6) and (1.6, -0.6) .. (2,0);    
    \draw[thick, color=black] (2,0) .. controls (2.4,-1.3) and (5.6,-1.3) .. (6,0);       
    \draw[thick, color=black] (2,0) .. controls (2.4, 0.6) and (3.6, 0.6) .. (4,0);
    \draw[thick, color=black] (2,0) .. controls (2.4, -0.6) and (3.6, -0.6) .. (4,0);
    \draw[thick, color=black] (4,0) .. controls (4.4, 0.6) and (5.6, 0.6) .. (6,0);
    \draw[thick, color=black] (4,0) .. controls (4.4, -0.6) and (5.6, -0.6) .. (6,0);        
    
    \node[circle, draw, inner sep=1.4pt, fill] (1) at (0,0) {};
    \node[circle, draw, inner sep=1.4pt, fill] (2) at (2,0) {};
    \node[circle, draw, inner sep=1.4pt, fill] (3) at (4,0) {};
    \node[circle, draw, inner sep=1.4pt, fill] (4) at (6,0) {};    

    \node[] at (3,-1.8) {\footnotesize $(G^{\rev},F^{\rev})$};   
\end{scope}	
\begin{scope}[scale=0.25, xshift=180, yshift=-320]
	\draw[thick, color=NavyBlue]    (0,0) -- (4,3) -- (4,7) -- (0,10) -- (-4,7) -- (-4,3) -- (0,0);		
	\node[circle,fill,inner sep=1.5pt,color=NavyBlue] at (0,0)  {};
	\node[circle,fill,inner sep=1.5pt,color=NavyBlue]  at (4,3) {};
	\node[circle,fill,inner sep=1.5pt,color=NavyBlue]  at (4,7) {}; 
	\node[circle,fill,inner sep=1.5pt,color=NavyBlue]  at (0,10) {};
	\node[circle,fill,inner sep=1.5pt,color=NavyBlue]  at (-4,7) {};
	\node[circle,fill,inner sep=1.5pt,color=NavyBlue]  at (-4,3) {};
	\node[circle,fill,inner sep=1.5pt,color=NavyBlue]  at (-7,-1) {};
 	\node[circle,fill,inner sep=1.5pt,color=NavyBlue]  at (-3,-4) {};

    \draw[thick, color=NavyBlue]    (-4,3) -- (-7,-1) -- (-3,-4) -- (0,0);  
    
    \node[] at (0,-6.5) {$\mathscr{L}_{G^{\rev},F^{\rev}}$};     
\end{scope}
\end{scope}
\end{tikzpicture}
    \caption{An example of the graph operations in Lemma~\ref{lem.dualizing_operations} and the corresponding framing lattices.}
    \label{fig.graph_operations_ex}
\end{figure}

An edge $(v,w)$ in $G$ is said to be \defn{idle} if $v$ has out-degree one or $w$ has in-degree one.
If a graph $G'$ is obtained from $G$ by contracting an idle edge, then the flow polytopes $\calF_{G}$ and~$\calF_{G'}$ are integrally equivalent, meaning they are affinely equivalent and have the same Ehrhart polynomial \cite[Lemma 2.2]{MS20}. 
Contracting idle edges in $G$ also preserves the framing poset, along with the operations mentioned in the following lemma.

\begin{lemma}
\label{lem.iso_operations}
    The following operations of $(G,F)$ are poset isomorphisms between framing posets.
    \begin{itemize}
        \item[(1)] Contracting an idle edge in $G$.
        \item[(2)] Changing the framing at the source $s$ or sink $t$ of $G$.
        \item[(3)] Reversing the incoming order of two edges $e = (s,v)$ and $e' = (s,v)$ at a vertex $v$ that are consecutive in $\leq_{\In(v)}$.
        \item[(4)] Reversing the outgoing order of two edges $e = (w,t)$ and $e' = (w,t)$ at a vertex $w$ that are consecutive in $\leq_{\Out(v)}$.
    \end{itemize}
\end{lemma}

\begin{proof}
    Operations (1) and (2) are immediate.
    We consider (3), as (4) is shown analogously. 
    Let $F_1$ denote the original framing with order $\leq_{\In(v)}$ and let $F_2$ be the framing otherwise the same as $F_1$ but with $\leq_{\In(v)}^{\mathrm{rev}}$ instead of $\leq_{\In(v)}$.
    The $\varphi$ be the map sending a route $R$ to the route $R'$ obtained from $R$ by replaceing $e$ with $e'$.
    Observe that two routes are coherent at $v$ under $F_1$ if and only if their images under $\varphi$ are coherent at $v$ under $F_2$.
    Thus $\varphi$ extends to a bijection $\Phi$ between maximal cliques.
    Furthermore, it is easy to verify that $\Phi$ is order-preserving.
\end{proof}

\subsection{The existence of a ccw rotation}

The following lemma gives a sufficient condition for the existence of a ccw rotation in a maximal clique.
It is a helpful tool that we will use later in order to give a global characterization of the partial order of the framing poset.

\begin{lemma}\label{lem.rotationExists}
    Let $R$ be a route in a maximal clique $C_1$. 
    If there is a route $R^*$ such that $R <_v^{\cw} R^*$ at some $v$, then there exists routes $R_1 \in C_1$ and $R_2 \notin C_1$ with $R_1 <_w^{\cw} R^*$ for some~$w$ $\in R^*$ such that $C_2 = (C_1\setminus R_1) \cup R_2$ is a maximal clique and $C_1\precccwrot C_2$. 
\end{lemma}


\begin{proof}
    Of the routes in $C_1$ containing $vR$ consider the route $R'$ minimal with respect to~$\leq_{\scrI (v)}$. 
    Then from the routes in $C_1$ containing $R'v$ take $Q$ to be the route maximal with respect to~$\leq_{\scrO (v)}$. 
    Now $Q \in C_1$ and by construction $Q <_v^{\cw} R^*$.    
    Moreover:
    \begin{itemize}
        \item[(A)] There is no $Q'\in C_1$ such that $Qv=Q'v$ and $vQ<_{\scrO (v)} vQ'$. 
        \item[(B)] There is no $Q'\in C_1$ such that $vQ=vQ'$ and $Q'v<_{\scrI (v)} Qv$.
    \end{itemize}
    Item (A) follows by construction. The proof of (B) requires further arguments. We proceed by contradiction, assuming that there is $Q'\in C_1$ such that $vQ=vQ'$ and $Q'v<_{\scrI (v)} Qv$.
    First, note that under this assumption we have that $Q$ is the only route in $C_1$ that contains $R'v=Qv$, because by the $\leq_{\scrO (v)}$-maximality of $Q$ any other such a route would be incoherent with $Q'$. 
    This means that $Q=R'$ and $vQ=vR'=vR$ by construction. 
    But then, $Q'$ contradicts the $\leq_{\scrI (v)}$-minimality of $R'=Q$.

    Now, We consider the following three possible cases,

    \begin{enumerate}
        \item There are routes $V$ and $V'$ in $C_1$ through $v$ such that $Qv <_{\scrI(v)} Vv$ and $vV' <_{\scrO(v)} vQ$. 
        \item There does not exist a route $V$ in $C_1$ through $v$ such that $Qv <_{\scrI(v)} Vv$.
        \item There does not exist a route $V'$ in $C_1$ through $v$ such that $vV' <_{\scrO(v)} vQ$.
    \end{enumerate}

    First consider Case (1).
    We can assume that $V$ is minimal with respect to $\leq_{\scrI(v)}$, and that~$V'$ is maximal with respect to $\leq_{\scrO(v)}$.
    Then take $w = v$, and let $R_1 = Q$ and $R_2 = VwV'$, as depicted in Figure~\ref{fig.rotationExists}.
    The following properties hold:
    \begin{enumerate}[(i)]
        \item $R_1\in C_1$. \newline
        This is clear because $R_1=Q\in C_1$. \smallskip
        \item $\Top(R_1,R_2)=R_1vR_2\in C_1$. \newline
        We will show that $R_1vR_2=QvV'$ is coherent with every route in $C_1$; since $C_1$ is a maximal clique, this implies that $QvV'\in C_1$ as desired. 
        We proceed by contradiction. 
        Assume that there exist $Q'\in C_1$ that is incoherent with $QvV'$.
        Then, they must be incoherent at $v$, otherwise $Q'$ would be incoherent with either $Q\in C_1$ or $V'\in C_1$.
        
        If $QvV'<_v^\cw Q'$ then $Q<_v^\cw Q'$ which is a contradiction. Therefore, $Q'<_v^\cw QvV'$. 
        But then $vQ=vQ'$, otherwise by the $\leq_{\scrO(v)}$-maximality of $V'$ we would have $vQ<_{\scrO(v)}vQ'$ which implies that $Q'$ is incoherent with $Q$.
        Thus we have $vQ=vQ'$, but also $Q'v<_{\scrI(v)}Qv$, contradicting Item (B) above.
        \item $\Bot(R_1,R_2)=R_2vR_1\in C_1$. \newline
        Again, we need to show that $R_2vR_1=VvQ$ is coherent with every route in $C_1$. We proceed by contradiction. 
        If there is $Q'\in C_1$ incoherent with $VvQ$ then they must be incoherent at $v$, otherwise it would be incoherent with either $V\in C_1$ or $Q\in C_1$.
        If $VvQ<_v^\cw Q'$ then $Q<_v^\cw Q'$ which is a contradiction. 
        Therefore, $Q'<_v^\cw VvQ$.
        But then $Qv=Q'v$, otherwise by the $\leq_{\scrI(v)}$-minimality of $V$ we would have $Q'v<_{\scrI(v)} Qv$ which implies that $Q'$ is incoherent with $Q$. 
        Thus we have $Qv=Q'v$, but also $vQ<_{\scrO(v)} vQ'$, contradicting Item (A) above.
        \item There is not route in $C_1$ in between $R_1=Q$ and $R_2=VvV'$. \newline
        We proceed by contradiction. Assume that there is $Q'\in C_1$ in between $Q$ and $VvV'$. 
        Then either $Qv=Q'v$ or $vQ=vQ'$, otherwise $Q'$ would be incoherent with $Q$.  
        If $Qv=Q'v$, then $Q'$ contradicts the $\leq_{\scrO(v)}$-maximality of $V'$. 
        If $vQ=vQ'$, then $Q'$ contradicts the $\leq_{\scrI(v)}$-minimality of $V$. 
    \end{enumerate}

    Properties (i)-(iv) show that 
    the Top-Bottom Property and the In Between Property of~\Cref{lem.rotation} hold for $R_1$ and $R_2$. 
    Thus, $C_2 = (C_1\setminus R_1) \cup R_2$ is a maximal clique.
    Since $R_2 <_w^{\cw} R_1$, we have that $C_1 \precccwrot C_2$.   
    We also have $R_1 <_w^{\cw} R^*$ because $Q=R_1$. 

    Next we consider Case (2). 
    Let $P_v$ be the maximal path containing $v$ in $Q\cap R^*$, and let~$w$ be the minimal vertex in $P_v$.
    By Lemma~\ref{lem.path_is_extendable}, there is a route in $C_1$ through the edge of~$R^*$ incoming to $w$. 
    In particular, there is a route $W$ in $C_1$ such that $Qw <_{\scrI(w)} Ww$. 
    We assume that $W$ is a minimal such route with respect to $\leq_{\scrI(w)}$, and that of such routes $W$ is the minimal route with respect to $\leq_{\scrO(w)}$.
    Note that $W$ cannot pass through~$v$, as otherwise we would have $Qv <_{\scrI(v)} Wv$, which contradicts our assumption for Case (2).
    Thus $w<v$ and $wQ <_{\scrO(w)} wW$.

    Let $w'$ be the maximal vertex in the path of $W\cap Q$ containing $w$. 
    Let $e$ be the edge of $Q$ incoming to $w$ and let $e'$ be the edge of $W$ outgoing from $w'$.
    Let $\widehat{P}$ be the path of $W\cap Q$ containing $w$, together with the edges $e$ and $e'$. 
    
    Suppose towards a contradiction that there is a route $R''\in C_1$ that is incoherent with the path $\widehat{P}$. 
    Since $R''$ must be coherent with $Q$ and $W$, it must be in between $QwW$ and $WwQ$ at $w$. However, if $R''w <_{\scrI(w)} Ww$, then that contradicts the minimality of the choice of $W$ with respect to $\leq_{\scrI(w)}$. Similarly, if $wR'' <_{\scrO(w)} wW$, then that contradicts the minimality of the choice of $W$ with respect to $\leq_{\scrO(w)}$. 
    
    Thus the path $\widehat{P}$ is coherent with all routes in $C_1$ and by Lemma~\ref{lem.path_is_extendable}, it can be extended to a route \mbox{$\widehat{Q} \in C_1$}. 
    Furthermore, since $\widehat{Q}$ cannot be in between $QwW$ and $WwQ$, it follows that $\widehat{Q}w=Qw$ and $w\widehat{Q}=wW$. 
    
    Now, since $Q,W\in C_1$ must be coherent with $\widehat{Q}$, we have that $\widehat{Q}$ satisfies conditions (A) and (B) above at $w$. 
    Since $\widehat{Q}$ is clockwise from $R^*$ at $w$, We can now repeat the argument in Case (1) using $\widehat{Q}$ instead of $Q$, $w$ instead of $v$, and $V:=W \in C_1$ and $V':=Q \in C_1$.
    
    Case (3) can be argued symmetrically to Case (2). 
\end{proof}

\begin{figure}[htb]
    \centering
    \begin{tikzpicture}
\begin{scope}[xshift=0, scale=0.8]
    \draw[ultra thick,color=black!30,dashed] (2.4,-1.2)--(3,-1.2);
    \draw[ultra thick,color=black!30] (3,-1.2) .. controls (3.3, -1.2) and (4.0, -1) .. (4.5,-0.035);
    \draw[ultra thick,color=black!30] (6,0) .. controls (6.4, 1.2) and (7.4,1.2) .. (7.5,1.2);
    \draw[ultra thick,color=black!30,dashed] (7.5,1.2)--(8.2,1.2);
    
    \draw[ultra thick,color=blue,dashed] (1.8,0.6)--(2.5,0.6);
    \draw[ultra thick,color=blue] (2.5,0.6) .. controls (2.6, 0.6) and (3.6, 0.6) .. (4,0);
    \draw[ultra thick,color=blue] (4,0.035)--(5.6,0.035);
    \draw[ultra thick,color=blue] (5.6,-0.035)--(6,-0.035);
    \draw[ultra thick,color=blue] (6,0) .. controls (6.4, -0.6) and (7.4, -0.6) .. (7.5,-0.6);
    \draw[ultra thick,color=blue,dashed] (7.5,-0.6)--(8.2,-0.6);
    
    \draw[ultra thick,color=red,dashed] (1.8,-0.6)--(2.5,-0.6);
    \draw[ultra thick,color=red] (2.5,-0.6) .. controls (2.6, -0.6) and (3.6, -0.6) .. (4,0);
    \draw[ultra thick,color=red] (4,-0.04)--(5.6,-0.04);
    \draw[ultra thick,color=red] (5.6,0.04)--(6,0.04);
    \draw[ultra thick,color=red] (6,0) .. controls (6.4, 0.6) and (7.4,0.6) .. (7.5,0.6);
    \draw[ultra thick,color=red,dashed] (7.5,0.6)--(8.2,0.6);

    \node[circle, draw, inner sep=1.4pt, fill,label=above:\scriptsize{$v$}](c) at (5.6,0) {};

    \node[] (r1) at (1.5,0.6) {\scriptsize $Q$};
    \node[] (r1) at (8.5,-0.6) {\scriptsize $Q$};
    \node[] (r1) at (1.5,-0.6) {\scriptsize $V$};
    \node[] (r1) at (8.5,0.6) {\scriptsize $V'$};

    \node[] (r1) at (2,-1.2) {\scriptsize $R^*$};
    \node[] (r1) at (8.5,1.2) {\scriptsize $R^*$};    


    \node[] (r1) at (2.5,-2.2) {\scriptsize $R_1 = Q$};
    \node[] (r1) at (7.5,-2.2) {\scriptsize $R_2 = VvV'$};

\end{scope}
\begin{scope}[xshift=220, scale=0.8]
    \draw[ultra thick,color=black!30,dashed] (2.4,-1.2)--(3,-1.2);
    \draw[ultra thick,color=black!30] (3,-1.2) .. controls (3.3, -1.2) and (4.0, -0.8) .. (4.3,-0.035);
    \draw[ultra thick,color=black!30] (6,0) .. controls (6.4, 1.2) and (7.4,1.2) .. (7.5,1.2);
    \draw[ultra thick,color=black!30,dashed] (7.5,1.2)--(8.2,1.2);

    \draw[ultra thick,color=blue,dashed] (1.8,0.6)--(2.5,0.6);
    \draw[ultra thick,color=blue] (2.5,0.6) .. controls (2.6, 0.6) and (3.6, 0.6) .. (4,0.035);
    \draw[ultra thick,color=blue] (3.97,0.035)--(4.3,0.035);
    \draw[ultra thick,color=blue] (5.6,0.0)--(6.02,0.0);
    \draw[ultra thick,color=blue] (6,0.0) .. controls (6.4, -0.6) and (7.4, -0.6) .. (7.5,-0.6);
    \draw[ultra thick,color=blue,dashed] (7.5,-0.6)--(8.2,-0.6);
    
    \draw[ultra thick,color=red,dashed] (1.8,-0.6)--(2.5,-0.6);
    \draw[ultra thick,color=red] (2.5,-0.6) .. controls (2.6,-0.6) and (3.6,-0.6) .. (4,-0.04);
    \draw[ultra thick,color=red] (3.97,-0.035)--(4.3,-0.035);
    \draw[ultra thick,color=blue] (4.3,0.035)--(5.1,0.035);   
    \draw[ultra thick,color=blue] (5.1,0)--(5.6,0);   
    \draw[ultra thick,color=red] (4.3,-0.035)--(5.1,-0.035);    

    \draw[ultra thick,color=red] (5.1,-0.04) .. controls (5.5, -1) and (7.3, -1.2) .. (7.5,-1.2);
    \draw[ultra thick,color=red,dashed] (7.5,-1.2)--(8.2,-1.2); 
    
    \node[circle, draw, inner sep=1.4pt, fill,label=above:\scriptsize{$v$}](c) at (5.6,0.01) {};
    \node[circle, draw, inner sep=1.4pt, fill,label=above:\scriptsize{$w$}](c) at (4.3,0) {};
    \node[circle, draw, inner sep=1.4pt, fill,label=above:\scriptsize{$w'$}](c) at (5.1,0) {};

    \node[] (r1) at (1.5,0.6) {\scriptsize $Q$};
    \node[] (r1) at (8.5,-0.6) {\scriptsize $Q$};
    \node[] (r1) at (1.5,-0.6) {\scriptsize $W$};
    \node[] (r1) at (8.6,-1.2) {\scriptsize $W$};  

    \node[] (r1) at (2,-1.2) {\scriptsize $R^*$};
    \node[] (r1) at (8.5,1.2) {\scriptsize $R^*$};    

\end{scope}
\end{tikzpicture}
    \caption{Cases (1) and (2) in the proof of~\Cref{lem.rotationExists}.}
    \label{fig.rotationExists}
\end{figure}

For completeness, we also include the analog of~\Cref{lem.rotationExists} for the existence of $\cw$ rotation, which follows by symmetry. 

\begin{lemma}\label{lem.rotationExists_cw_version}
    Let $R$ be a route in a maximal clique $C_1$. 
    If there is a route $R^*$ such that $R^* <_v^{\cw} R$ at some $v$, then there exists routes $R_1 \in C_1$ and $R_2 \notin C_1$ with $R^* <_w^{\cw} R_1$ for some~$w$ $\in R^*$ such that $C_2 = (C_1\setminus R_1) \cup R_2$ is a maximal clique and $C_2\precccwrot C_1$. 
\end{lemma}

\subsection{Characterization of the partial order}

The purpose of this section is to prove the following characterization of the partial order relation of the framing poset. 

Given two sets of pairwise coherent paths (for example two maximal cliques) $C$ and $C'$, we say that \defn{$C$ is cw from $C'$} if for all paths $P \in C$, $P'\in C'$, and $v \in P\cap P'$, we have that $P$ and $P'$ are coherent at $v$ or $P <_v^{\cw} P'$.

\begin{theorem}[Characterization of the partial order]\label{thm.ccw}
    Let $C$ and $C'$ be maximal cliques. 
    Then $C \leq C'$ if and only if $C$ is cw from $C'$.    
\end{theorem}

In order to show this, we will introduce two algorithms that we call the $C_{\max}$ algorithm and the $C_{\min}$ algorithm. These algorithms will also be very useful later to prove that the framing poset is a lattice.  

\subsubsection{The $C_{\max}$ and $C_{\min}$ algorithms}

Let $S$ be a set of pairwise coherent paths.
We construct a maximal clique $C_{\max}(S)$ containing the ccw-most routes that are coherent with~$S$. 
Informally, this is done by adding the ccw-most routes that are coherent with $S$ recursively at each vertex until a maximal clique is formed.
If $S$ is a set of pairwise coherent \emph{routes}, then $S$ will be contained in $C_{\max}(S)$; this case will be very important for us. 
The formal construction is described in Algorithm~\ref{alg.CmaxS}, where $\leq_{\scrI(v)}^{\rev}$ denotes the reverse order to $\leq_{\scrI(v)}$.

\begin{algorithm}
  \caption{The construction of $C_{\max}(S)$}\label{alg.CmaxS}
  \begin{algorithmic}[1]
    \State $C_{\max}(S) := \emptyset$
    \For{$v \in V(G)$ (in increasing order)}
        \For{$Pv \in \scrI(v)$ (in the order $\leq_{\scrI(v)}^{\rev}$)} \Comment{$Pv$ possibly empty}
            \For{$vQ \in \scrO(v)$ (in the order $\leq_{\scrO(v)}$)} \Comment{$vQ$ possibly empty}
                \If{$PvQ$ is coherent with all paths in $C_{\max}(S) \cup S$}
                    \State $C_{\max}(S) := C_{\max}(S) \cup \{PvQ\}$
                    \State \textbf{break} \Comment{This terminates the innermost loop}
                \EndIf
            \EndFor
        \EndFor
    \EndFor  
  \end{algorithmic}
\end{algorithm}

It is worth emphasizing that the same route may be added to $C_{\max}(S)$ several times during its construction, making the algorithm inefficient. 
It should also be noted that whenever a route $PvQ$ is coherent with $S$ and all routes already added to $C_{\max}(S)$, the innermost loop terminates regardless of whether $PvQ$ is already in the set $C_{\max}(S)$ or not. 
The following example illustrates these points explicitly.

\begin{example} \label{ex.alg1}
    We apply Algorithm~\ref{alg.CmaxS} to the graph $\oru{3}$ in Example~\ref{example_running_oruga} with \mbox{$S = \emptyset$}.
    At vertex $1$, $P1$ is the empty path, and $1Q$ is the path $e_1e_3e_5$.
    Thus the first path added to $C_{\max}(S)$ is $e_1e_3e_5$.
    We have considered all paths incoming to $1$, so we proceed to vertex $2$.
    The path $e_2e_3e_5$ is considered first and added to $C_{\max}(S)$.
    The next path to be considered is the path $e_1e_3e_5$, which is already contained in $C_{\max}(S)$.
    Nevertheless, it is added again and the innermost loop in the algorithm terminates, and we proceed to vertex $3$. 
    The first route considered at vertex $3$ is $e_2e_4e_5$, which is added.
    The routes $e_2e_3e_5$ and $e_1e_3e_5$ are added again, and no more new routes are added at vertex $3$.
    Finally, at vertex $4$ every route is considered, but the only new route added is $e_2e_4e_6$. 
    We have that $C_{\max}(\emptyset) = \{e_1e_3e_5, e_2e_3e_5, e_2e_4e_5, e_2e_4e_6\}$.
    See the top row in Figure~\ref{fig.Cmax_and_Cmin_ex} for the order in which each route is added.
\end{example}

\begin{figure}
    \centering

    \caption{Computing $C_{\max}(\emptyset)$ and $C_{\min}(\emptyset)$ for the graph $\oru{3}$.}
    \label{fig.Cmax_and_Cmin_ex}
\end{figure}

In a similar way, we construct a maximal clique $C_{\min}(S)$ whose routes are as clockwise as possible while coherent with $S$.
The algorithm is otherwise the same, but the orders in which $\scrI(v)$ and~$\scrO(v)$ are read in Algorithm~\ref{alg.CmaxS} are reversed. That is, we replace $\leq_{\scrI(v)}^{\rev}$ with $\leq_{\scrI(v)}$ and replace $\leq_{\scrO(v)}$ with $\leq_{\scrO(v)}^{\rev}$, where $\leq_{\scrO(v)}^{\rev}$ denotes the reverse order to $\leq_{\scrO(v)}$.
The precise construction is described in Algorithm~\ref{alg.CminS}.

\begin{algorithm}
  \caption{The construction of $C_{\min}(S)$}\label{alg.CminS}
  \begin{algorithmic}[1]
    \State $C_{\min}(S) := \emptyset$
    \For{$v \in V(G)$ (in increasing order)}
        \For{$Pv \in \scrI(v)$ (in the order $\leq_{\scrI(v)}$)} \Comment{$Pv$ possibly empty}
            \For{$vQ \in \scrO(v)$ (in the order $\leq_{\scrO(v)}^{\rev}$)} \Comment{$vQ$ possibly empty}
                \If{$PvQ$ is coherent with all paths in $C_{\max}(S) \cup S$}
                    \State $C_{\min}(S) := C_{\min}(S) \cup \{PvQ\}$
                    \State \textbf{break} \Comment{This terminates the innermost loop}
                \EndIf
            \EndFor
        \EndFor
    \EndFor  
  \end{algorithmic}
\end{algorithm}

Similar to Example~\ref{ex.alg1}, we demonstrate Algorithm~\ref{alg.CminS} on the graph $\oru{3}$ below.

\begin{example} \label{ex.alg2}
    We apply Algorithm~\ref{alg.CminS} to the graph $\oru{3}$ in Example~\ref{example_running_oruga} with $S = \emptyset$.
    At vertex $1$, $P1$ is the empty path, and $1Q$ is the path $e_2e_4e_6$.
    Thus the first path added to $C_{\min}(S)$ is $e_2e_4e_6$.
    We have considered all paths incoming to $1$, so we proceed to vertex $2$.
    The path $e_1e_4e_6$ is considered first and added to $C_{\min}(S)$.
    The next path to be considered is the path $e_2e_4e_6$, which is already contained in $C_{\min}(S)$.
    Nevertheless, the innermost loop in the algorithm terminates, and we proceed to vertex $3$. 
    The first route considered at vertex $3$ is $e_1e_3e_6$, which is added.
    The routes $e_1e_4e_6$ and $e_2e_4e_6$ are added again, and no more new routes are added at vertex $3$.
    Finally, at vertex $4$ every route is considered, but only the new route $e_1e_3e_5$ is added. 
    We have that $C_{\min}(\emptyset) = \{e_2e_4e_6, e_1e_4e_6, e_1e_3e_6, e_1e_3e_5\}$.
    See the bottom row of Figure~\ref{fig.Cmax_and_Cmin_ex} for the order in which each route is added.
\end{example}

\begin{lemma}\label{lem.CmaxS_is_ccw_maximal}
    The cliques $C_{\max}(S)$ and $C_{\min}(S)$ are the unique maximal cliques with the following properties. If a route $R$ is coherent with all paths in $S$, then 
    \begin{enumerate}
        \item for any $R' \in C_{\max}(S)$ and $v \in R\cap R'$ either $R$ and $R'$ are coherent at $v$ or $R <_v^{\cw} R'$.
        \item for any $R' \in C_{\min}(S)$ and $v \in R\cap R'$ either $R$ and $R'$ are coherent at $v$ or $R' <_v^{\cw} R$.
    \end{enumerate} 
\end{lemma}
\begin{proof}
    We prove property (1) for $C_{\max}(S)$. The proof of property (2) for $C_{\min}(S)$ is similar. 
    
    First note that $C_{\max}(S)$ is a maximal clique since all routes are considered in the for-loop in Line~4 of Algorithm~\ref{alg.CmaxS} when $vQ$ is the empty path with $v=t$.
    Suppose toward a contradiction that there is a route $R$ such that $R$ is coherent with all paths of $S$ but incoherent with a route $R' \in C_{\max}(S)$ such that $R' <_v^{\cw} R$ at some $v \in R\cap R'$.
    We assume that $R'$ is the first route added in the algorithm that is not coherent with $R$ and take $v$ to be the minimal point at which $R$ and $R'$ are incoherent.
    Let $v'$ be the vertex in Algorithm~\ref{alg.CmaxS} at which the route $R'$ is added.
    We consider two cases, namely when $v < v'$ and when $v \geq v'$.
    These cases are depicted in Figure~\ref{fig.CmaxS_proof_cases}. 

    First, if $v < v'$, then consider the step in the algorithm where the extensions of the path $Rv$ are considered.
    At this point $R'$ has not been added since $v < v'$.
    Since the extension $R = RvR$ is considered before $RvR'$ and $R$ is coherent with all the routes added before this step and the paths in $S$, it follows that in this step we add a route $\widetilde{R} = Rv\widetilde{R}$ to $C_{\max}(S)$ where $v\widetilde{R} \leq_{\scrO(v)} vR$. 
    However, this route $\widetilde{R}$ is incoherent with $R'$, which contradicts the fact that $R'$ is added to~$C_{\max}(S)$ at a later step. 


    Next, if $v \geq v'$, then consider the step in the algorithm when $R'$ is added.
    Since $R$ and $R'$ are coherent with all the previously added routes to $C_{\max}(S)$ and the paths in $S$, the route $R'vR$ must also be coherent with all the routes previously added to $C_{\max}(S)$ and the paths in $S$.
    However, $v'R$ is read before $v'R'$ in the order $\leq_{\scrO(v)}$, which contradicts the fact that $R'$ is added to $C_{\max}(S)$ at this step.

    It remains to check uniqueness.
    Let $C\neq C_{\max}(S)$ be another maximal clique satisfying property (1).
    Then, there is a route $R_1 \in C$ and $R_2 \in C_{\max}(S)$ such that $R_1$ and $R_2$ are coherent with all paths in $S$, but they are incoherent at some $v$.
    Note that $R_1$ and $R_2$ cannot be in $S$ and thus we must have $R_1 <_v^{\cw} R_2$ and $R_2 <_v^{\cw} R_1$, which is a contradiction.
    It follows that $C_{\max}(S)$ is unique.
\end{proof}

\begin{figure}
    \centering    
    \begin{tikzpicture}
\begin{scope}[xshift=0]
    \draw[ultra thick,color=blue] (2.5,0.6) .. controls (2.6, 0.6) and (3.6, 0.6) .. (4,0);
    \draw[ultra thick,color=blue] (4,-0.03)--(5,-0.03);
    \draw[ultra thick,color=blue, dashed] (2,0.6)--(2.5,0.6);
    \draw[ultra thick,color=blue] (5,0) .. controls (5.4, -0.6) and (6.4, -0.6) .. (6.5,-0.6);
    \draw[ultra thick,color=blue, dashed] (6.5,-0.6)--(7,-0.6);
    
    \draw[ultra thick,color=red] (2.5,-0.6) .. controls (2.6, -0.6) and (3.6, -0.5) .. (4,0);
    \draw[ultra thick,color=red] (4,0.03)--(5,0.03);
    \draw[ultra thick,color=red, dashed] (2,-0.6)--(2.5,-0.6);    
    \draw[ultra thick,color=red] (5,0) .. controls (5.4, 0.6) and (6.4, 0.6) .. (6.5,0.6);
    \draw[ultra thick,color=red, dashed] (6.5,0.6)--(7,0.6);

    \node[circle, draw, inner sep=1.6pt, fill,label=below:\small{$v$}] (s) at (4,0) {};    
    \node[circle, draw, inner sep=1.6pt, fill] (c) at (5,0) {};
    \node[circle, draw, inner sep=1.6pt, fill,label=below:\small{$v'$}] (s) at (5.75,-0.5) {};    
    
    \node[] (P) at (1.6,0.6) {\color{blue} $R'$};
    \node[] (Q) at (1.6,-0.6) {\color{red} $R$};
    
    \node[] (P) at (7.4,0.6) {\color{red} $R$};
    \node[] (Q) at (7.4,-0.6) {\color{blue} $R'$};

\end{scope}

\begin{scope}[xshift=240]
    \draw[ultra thick,color=blue] (2.5,0.6) .. controls (2.6, 0.6) and (3.6, 0.6) .. (4,0);
    \draw[ultra thick,color=blue] (4,-0.03)--(5,-0.03);
    \draw[ultra thick,color=blue, dashed] (2,0.6)--(2.5,0.6);
    \draw[ultra thick,color=blue] (5,0) .. controls (5.4, -0.6) and (6.4, -0.6) .. (6.5,-0.6);
    \draw[ultra thick,color=blue, dashed] (6.5,-0.6)--(7,-0.6);
    
    \draw[ultra thick,color=red] (2.5,-0.6) .. controls (2.6, -0.6) and (3.6, -0.5) .. (4,0);
    \draw[ultra thick,color=red] (4,0.03)--(5,0.03);
    \draw[ultra thick,color=red, dashed] (2,-0.6)--(2.5,-0.6);    
    \draw[ultra thick,color=red] (5,0) .. controls (5.4, 0.6) and (6.4, 0.6) .. (6.5,0.6);
    \draw[ultra thick,color=red, dashed] (6.5,0.6)--(7,0.6);

    
    \node[circle, draw, inner sep=1.6pt, fill,label=below:\small{$v$}] (s) at (4,0) {};    
    \node[circle, draw, inner sep=1.6pt, fill] (c) at (5,0) {};
    \node[circle, draw, inner sep=1.6pt, fill,label=below:\small{$v'$}] (s) at (3.2,0.5) {};    
    
    \node[] (P) at (1.6,0.6) {\color{blue} $R'$};
    \node[] (Q) at (1.6,-0.6) {\color{red} $R$};
    
    \node[] (P) at (7.4,0.6) {\color{red} $R$};
    \node[] (Q) at (7.4,-0.6) {\color{blue} $R'$};

\end{scope}

\end{tikzpicture}
    \caption{The two cases in the proof of Lemma~\ref{lem.CmaxS_is_ccw_maximal}.}
    \label{fig.CmaxS_proof_cases}
\end{figure}

When $S = \emptyset$, 
the next lemma implies that the framing poset $\scrL_{G,F}$ has a unique minimal element $\widehat{0}$ and a unique maximal element $\widehat{1}$.

\begin{lemma}\label{lem.rotationExists2}
    Let $S$ be a set of pairwise coherent paths and $C$ be a maximal clique whose routes are coherent with the paths in $S$. The following hold:
    \begin{enumerate}
        \item If $C\neq C_{\max}(S)$ then there is a maximal clique $C'$ whose routes are coherent with the paths in $S$ such that $C\precccwrot C'$. 
        \item If $C\neq C_{\min}(S)$ then there is a maximal clique $C'$ whose routes are coherent with the paths in $S$ such that $C'\precccwrot C$. 
    \end{enumerate}
    In particular, $C_{\min} := C_{\min}(\emptyset)$ and $C_{\max} := C_{\max}(\emptyset)$ are respectively the $\widehat{0}$ and $\widehat{1}$ of $\scrL_{G,F}$.
\end{lemma}

\begin{proof}
    We prove property (1) for $C_{\max}(S)$. The proof of property (2) for $C_{\min}(S)$ is similar. 

    Since $C\neq C_{\max}(S)$, it follows from Lemma~\ref{lem.CmaxS_is_ccw_maximal} that there is a route $R \in C\setminus S$, a route $R'\in C_{\max}(S)$, and a vertex $v \in R\cap R'$ such that $R<_v^{\cw} R'$.
    By~\Cref{lem.rotationExists}, there exists a route $R_1 \in C$ such that $R_1 <_w^{\cw} R'$ for some $w$ and a route $R_2\notin C$ such that $C' = (C \setminus R_1)\cup R_2$ is a maximal clique. 
    Since $R'$ is coherent with the paths in $S$ and $R_1 <_w^{\cw} R'$, we have that $R_1 \notin S$.
    We need to show that $R_2$ is coherent with all paths in $S$.
    Assume to the contrary that $R_2$ is incoherent with a path $P$ in $S$.
    Since $P$ is coherent with~$C$, it is coherent with $R_1$, $\Top(R_1,R_2)$, and $\Bot(R_1,R_2)$, and therefore $P$ must be weakly in between $R_1$ and $R_2$ at~$w$. 
    If $P\subseteq R_1$, then $P$ is incoherent with $R' \in C_{\max}(S)$, which is a contradiction.
    If $P\not \subseteq R_1$, then $P$ is in between $R_1$ and $R_2$.
    Then, by Lemma~\ref{lem.path_is_extendable} the path $P$ is extendable to a route $R_P\in C$ that is in between $R_1$ and $R_2$, which is again a contradiction.
\end{proof}

The following corollaries are straight forward consequences of~\Cref{lem.rotationExists2}.

\begin{corollary}
\label{cor_Cmin_properties}
    Let $S$ be a set of pairwise coherent paths. The following hold:
    \begin{enumerate}
        \item $C_{\max}(S)$ is the unique maximal clique that is bigger in the order $\leq_{\mathrm{rot}}^{\ccw}$ than all the maximal cliques whose routes are coherent with the paths in $S$. 
        \item $C_{\min}(S)$ is the unique maximal clique that is smaller in the order $\leq_{\mathrm{rot}}^{\ccw}$ than all the maximal cliques whose routes are coherent with the paths in $S$.
        \qed
    \end{enumerate}
\end{corollary}

\begin{corollary}
    \label{cor.interval}
    Let $S$ be a set of pairwise coherent paths.
    The set of maximal cliques whose routes are coherent with the paths in $S$ is the interval $[C_{\min}(S), C_{\max}(S)]$ of the framing poset $\scrL_{G,F}$. \qed
\end{corollary}

Specializing to the case where $S$ is a set of pairwise coherent routes we get:

\begin{corollary}
    \label{cor_interval_max_cliques}
    Let $S$ be a set of pairwise coherent routes.
    The set of maximal cliques containing $S$ is the interval $[C_{\min}(S), C_{\max}(S)]$ of the framing poset $\scrL_{G,F}$. \qed
\end{corollary}

\subsubsection{Proof of~\Cref{thm.ccw}}
We now have all the necessary ingredients to prove the Characterization Theorem~\ref{thm.ccw} of the partial order of the framing poset, which states that two maximal cliques satisfy $C \leq C'$ if and only if $C$ is cw from $C'$. 

    We first show the forward direction and assume $C\leq C'$. 
    Suppose $C$ is not cw from $C'$, i.e. that there is a route $R \in C$ and $R'\in C'$ such that $R' <_v^{\cw} R$ at some $v$.
    By Lemma~\ref{lem.ccwExists} there is then a route $R'' \in C'$ 
    such that $R'<_w^{\cw} R''$ at some $w$. But this would imply that~$R''$ is incoherent with $R'$, which is a contradiction.
        
    To show the backward direction, we assume that $C$ is cw from $C'$ and show that we can apply a sequence of rotations to $C$ until obtaining $C'$.
    If all routes of $C$ are coherent with all routes in $C'$, then $C = C'$ and we are done. 
    Thus we consider the case when there is a route $R \in C$ and $R' \in C'$ such that $R <_v^{\cw} R'$ for some $v$. 
    Take $S = C\cap C'$.
    By Lemma~\ref{lem.CmaxS_is_ccw_maximal}, we have $C \neq C_{\max}(S)$. 
    By~\Cref{lem.rotationExists2}, there are routes $R_1$ and $R_2$ through some $w$ with $R_1 <_w^{\cw} R_2$ and $R_1\in C\setminus S$ such that $(C\setminus R_1) \cup R_2$ is a maximal clique.
    
    Next, we show that $(C\setminus R_1)\cup R_2$ is also cw from $C'$.
    Suppose toward a contradiction that~$R_2$ satisfies $\widetilde{R}' <_x^{\cw} R_2$ for some $\widetilde{R}'$ in $C'$ and $x \in R_2\cap \widetilde{R}'$; see~\Cref{fig_case_example_proof_theorem} for two examples of such a route $\widetilde{R}'$.
    Let $P_x$ be the maximal path in $\widetilde{R}'\cap R_2$ containing $x$, and consider $\min(P_x)$ and $\max(P_x)$. 
    We cannot have $\max(P_x) < w$, as then $\widetilde{R}' <_w^{\cw} R_2wR_1$ and $R_2wR_1 = \Bot(R_1,R_2)\in C$, which contradicts the fact that $C$ is cw from $C'$.
    Similarly, we cannot have $w < \min(P_x)$, as then $\widetilde{R}' <_x^{\cw} R_1wR_2$ and $R_1wR_2 = \Top(R_1,R_2)\in C$, which contradicts the fact that $C$ is cw from $C'$.

\begin{figure}[htb]
        \centering
        \begin{tikzpicture}

\begin{scope}[xshift=0, yshift=0, scale=0.8]
    \draw[ultra thick,color=black!30,dashed] (2.4,0)--(3,0);
    \draw[ultra thick,color=black!30] (3,0) .. controls (3.1,0) and (3.4,-0.05) .. (3.6,-0.35);
    \draw[ultra thick,color=black!30,dashed] (7.6,0)--(7,0);
    \draw[ultra thick,color=black!30] (7,0) .. controls (6.9,0) and (6.6,0.05) .. (6.4,0.35);

    \draw[ultra thick,color=blue,dashed] (1.8,0.6)--(2.5,0.6);
    \draw[ultra thick,color=blue] (2.5,0.6) .. controls (2.6, 0.6) and (3.6, 0.6) .. (4,0);
    \draw[ultra thick,color=blue] (4,0.035)--(5,0.035);
    \draw[ultra thick,color=blue] (5,-0.035)--(6,-0.035);
    \draw[ultra thick,color=blue] (6,0) .. controls (6.4, -0.6) and (7.4, -0.6) .. (7.5,-0.6);
    \draw[ultra thick,color=blue,dashed] (7.5,-0.6)--(8.2,-0.6);

    \draw[ultra thick,color=red,dashed] (1.8,-0.6)--(2.5,-0.6);
    \draw[ultra thick,color=red] (2.5,-0.6) .. controls (2.6, -0.6) and (3.6, -0.6) .. (4,0);
    \draw[ultra thick,color=red] (4,-0.04)--(5,-0.04);
    \draw[ultra thick,color=red] (5,0.04)--(6,0.04);
    \draw[ultra thick,color=red] (6,0) .. controls (6.4, 0.6) and (7.4,0.6) .. (7.5,0.6);
    \draw[ultra thick,color=red,dashed] (7.5,0.6)--(8.2,0.6);

    \node[circle, draw, inner sep=1.4pt, fill,label=below:\scriptsize{$w$}](c) at (5,0) {};

    \node[circle, draw, inner sep=1.4pt, fill,label=below:\scriptsize{$\min(P_x)$}](c) at (3.6,-0.35) {};

    \node[circle, draw, inner sep=1.4pt, fill,label=above:\scriptsize{$\max(P_x)$}](c) at (6.4,0.35) {};

    \node[] (r1) at (1.5,0.6) {\scriptsize $R_1$};
    \node[] (r1) at (8.5,-0.6) {\scriptsize $R_1$};
    \node[] (r1) at (1.5,-0.6) {\scriptsize $R_2$};
    \node[] (r1) at (8.5,0.6) {\scriptsize $R_2$};

    \node[] (r1) at (2.1,0) {\scriptsize $\widetilde R'$};
    \node[] (r1) at (7.9,0) {\scriptsize $\widetilde R'$};    

\end{scope}

\begin{scope}[xshift=220, yshift=0, scale=0.8]
    \draw[ultra thick,color=black!30,dashed] (2.4,0)--(3,0);
    \draw[ultra thick,color=black!30] (3,0) .. controls (3.1,0) and (3.4,0.05) .. (3.6,0.35);
    \draw[ultra thick,color=black!30,dashed] (7.6,0)--(7,0);
    \draw[ultra thick,color=black!30] (7,0) .. controls (6.9,0) and (6.6,0.05) .. (6.4,0.35);

    \draw[ultra thick,color=blue,dashed] (1.8,0.6)--(2.5,0.6);
    \draw[ultra thick,color=blue] (2.5,0.6) .. controls (2.6, 0.6) and (3.6, 0.6) .. (4,0);
    \draw[ultra thick,color=blue] (4,0.035)--(5,0.035);
    \draw[ultra thick,color=blue] (5,-0.035)--(6,-0.035);
    \draw[ultra thick,color=blue] (6,0) .. controls (6.4, -0.6) and (7.4, -0.6) .. (7.5,-0.6);
    \draw[ultra thick,color=blue,dashed] (7.5,-0.6)--(8.2,-0.6);

    \draw[ultra thick,color=red,dashed] (1.8,-0.6)--(2.5,-0.6);
    \draw[ultra thick,color=red] (2.5,-0.6) .. controls (2.6, -0.6) and (3.6, -0.6) .. (4,0);
    \draw[ultra thick,color=red] (4,-0.04)--(5,-0.04);
    \draw[ultra thick,color=red] (5,0.04)--(6,0.04);
    \draw[ultra thick,color=red] (6,0) .. controls (6.4, 0.6) and (7.4,0.6) .. (7.5,0.6);
    \draw[ultra thick,color=red,dashed] (7.5,0.6)--(8.2,0.6);

    \node[circle, draw, inner sep=1.4pt, fill,label=below:\scriptsize{$w$}](c) at (5,0) {};

    \node[circle, draw, inner sep=1.4pt, fill,label=above:\scriptsize{$\min(P_x)$}](c) at (3.6,0.35) {};

    \node[circle, draw, inner sep=1.4pt, fill,label=above:\scriptsize{$\max(P_x)$}](c) at (6.4,0.35) {};

    \node[] (r1) at (1.5,0.6) {\scriptsize $R_1$};
    \node[] (r1) at (8.5,-0.6) {\scriptsize $R_1$};
    \node[] (r1) at (1.5,-0.6) {\scriptsize $R_2$};
    \node[] (r1) at (8.5,0.6) {\scriptsize $R_2$};

    \node[] (r1) at (2.1,0) {\scriptsize $\widetilde R'$};
    \node[] (r1) at (7.9,0) {\scriptsize $\widetilde R'$};    

\end{scope}
\end{tikzpicture}
        \caption{Two examples for the route $\widetilde{R}'$ in the proof of~\Cref{thm.ccw}.}
        \label{fig_case_example_proof_theorem}
    \end{figure}
    
    If $\min(P_x) \leq w \leq \max(P_x)$, observe that that since $C$ is cw from $C'$ we cannot have $\widetilde{R}'$ be cw from $R_1wR_2$, $R_2wR_1$, or $R_1$ as they are all routes in $C$.   
    Thus $R_1w \leq_{\scrI(w)} \widetilde{R}'w <_{\scrI(w)} R_2w$ and $wR_2 <_{\scrO(w)} w\widetilde{R}' \leq_{\scrO(w)} wR_1$.
    Consider the case when $R_1w <_{\scrI(w)} \widetilde{R}'w$.
    Let $q$ be the largest vertex after which $R_1w$ and $\widetilde{R}'w$ coincide. 
    Let $Q$ denote the path formed by the edge of $\widetilde{R}'w$ entering $q$ and the subpath of $\widetilde{R}'$ from $q$ to $w$.
    Now $Q$ must be coherent with all routes in $C$, as otherwise 
    $Q<_y^{\cw} \widetilde{R}$ or $\widetilde{R} <_y^{\cw} Q$ for some route $R_1 \neq \widetilde{R}\in C$ at some node $y$; 
    In the first case, $\widetilde{R}'<_y^{\cw} \widetilde{R}$ which contradicts that $C$ is cw from $C'$, and in the second case $\widetilde{R} <_y^{\cw} R_2$ which contradicts that the only route of $C$ incoherent with $R_2$ is $R_1$. 

    Since the path $Q$ is coherent with all the routes in $C$, by Lemma~\ref{lem.path_is_extendable}, $Q$ is extendable to a route $R_Q$ in $C$.
    Note that $R_1w <_{\scrI(w)} R_Q w <_{\scrI(w)} R_2w$, 
    and since $R_Q$ must be coherent with $\Top(R_1,R_2)$ and $\Bot(R_1,R_2)$, then $R_Q$ must be a route in between $R_1$ and~$R_2$. This contradicts the In Between Property of~\Cref{lem.rotation} for $R_1$ and $R_2$ at $w$. 
    Therefore $R_1w =_{\scrI(w)} \widetilde{R}'w$.
    By a similar argument, we cannot have $w\widetilde{R}' <_{\scrO(w)} wR_1$, and hence $w\widetilde{R}' =_{\scrO(w)} wR_1$. 
    It follows that $R_1 = \widetilde{R}' \in C\cap C'$ and therefore $R_1\in S$, which is a contradiction.
    
    We have now shown that if $C'\neq C$ and $C$ is cw from $C'$, then we can apply a rotation to~$C$ in order to obtain a maximal clique $C^*:= (C\setminus R_1) \cup R_2$ such that $C \precccwrot C^*$ and $C^*$ is cw from $C'$.
    If $C^* = C'$, then we are finished.
    Otherwise, we can repeat the process until we reach $C'$. 
    \qed

\section{Lattice properties}\label{sec_lattice_properties}
\subsection{The framing poset is a polygonal lattice}

To show that $\scrL_{G,F}$ is a lattice, we rely on the BEZ lemma, which is stated as follows.

\begin{lemma}(BEZ Lemma \cite[Lemma 9-2.2]{Rea16}) \label{lem.BEZ}
Suppose $P$ is a finite poset with $\widehat{0}$. 
Suppose also that for all $x$ and $y$ in $P$ such that $x$ and $y$ cover a common element $z$, the join $x\vee y$ exists. 
Then $P$ is a lattice. \qed
\end{lemma}

To apply the BEZ lemma on $\scrL_{G,F}$, we investigate the case when two maximal cliques cover a single maximal clique in the following lemmas.

\begin{lemma}
\label{lem.R2Q}
    Let $Q$ be a maximal clique covered by two distinct maximal cliques $C_1 = Q\setminus R_1^Q \cup R_1$ and $C_2 = Q\setminus R_2^Q\cup R_2$.
    For $i\in \{1,2\}$ let $w_i$ be a point at which $R_i^Q$ and $R_i$ are incoherent and let $P_{w_i}$ denote the maximal path on which $R_i^Q$ and $R_i$ are incoherent.
    If $R_1$ and $R_2$ are incoherent at $x$ and $P_x$ denotes the maximal path containing $x$ at which $R_1$ and $R_2$ are incoherent, then one of the following holds.
    \begin{itemize}
        \item[(1)] If $R_1 <_x^{\cw} R_2$, then $P_x \cap P_{w_1} = \emptyset$. 
        Furthermore,
        $R_2^Q = \Bot(R_1^Q,R_1)$ when $x < w_1$ and $R_2^Q = \Top(R_1^Q,R_1)$ when $w_1 < x$.
        \item[(2)] If $R_2 <_x^{\cw} R_1$, then $P_x \cap P_{w_2} = \emptyset$.
        Furthermore,
        $R_1^Q = \Bot(R_2^Q,R_2)$ when $x < w_2$ and $R_1^Q = \Top(R_2^Q,R_2)$ when $w_2 < x$.
    \end{itemize}
\end{lemma}

\begin{proof}
    By symmetry it suffices to prove (1). 
    First, suppose toward a contradiction that $P_x\cap P_{w_1} \neq \emptyset$.
    Then there is a vertex $v$ at which we have $R_1^Q <_v^{\cw} R_1 <_v^{\cw} R_2$.
    In particular, $R_1^Q$ and $R_2$ are incoherent.
    Since the only route in $Q$ incoherent with $R_2$ is $R_2^Q$, it follows that $R_1^Q=R_2^Q$, which contradicts the fact that $C_1\neq C_2$ (Corollary~\ref{cor.route_rotatable_only_in_one_way}).
    Hence $P_x\cap P_{w_1} = \emptyset$.

    Now if $x < w_1$, then we have that $R_2$ is incoherent with $\Bot(R_1^Q,R_1) = R_1w_1R_1^Q \in Q$.
    The only route in $Q$ that is incoherent with $R_2$ is $R_2^Q$, and so necessarily $R_2^Q = \Bot(R_1^Q,R_1)$. 
    On the other hand, if $w_1 < x$, then we have that $R_2$ is incoherent with $\Top(R_1^Q,R_1) = R_1^Qw_1R_1 \in Q$.
    Therefore $R_2^Q = \Top(R_1^Q,R_1)$.
\end{proof}

\begin{proposition}\label{prop.intervalJoin}
    Let $C_1$ and $C_2$ be distinct maximal cliques covering a maximal clique $Q$ in $\scrL_{G,F}$ and let $S = C_1\cap C_2$. 
    Then, the following statements hold.
    \begin{itemize}
        \item[(i)] For $i\in \{1,2\}$, let $R_i = C_i\setminus Q$ and $R_i^Q = Q\setminus C_i$. 
        Then $R_1^Q \in C_2\setminus C_1$ and $R_2^Q\in C_1\setminus C_2$.
        \item[(ii)] The maximal clique $Q$ is the minimal element in the interval $I_S = [C_{\min}(S),C_{\max}(S)]$.
        \item[(iii)] The interval $I_S$ is the union of two chains $K_1$ and $K_2$, with $C_1\in K_1$ and $C_2\in K_2$, such that $K_1\cap K_2 = \{C_{\min}(S), C_{\max}(S)\}$.
        \item[(iv)] Furthermore, the chains $K_1$ and $K_2$ forming $I_S$ are of length $2$ or $3$. 
        In other words, $I_S$ is a square, pentagon, or a hexagon.
        \item[(v)] $C_1 \vee C_2$ exists and is $C_{\max}(S)$.
    \end{itemize}
\end{proposition}

\begin{proof}
    For (i) it suffices to simply observe that $R_1^Q \neq R_2^Q$.
    By Corollary~\ref{cor.interval}, $I_S$ is an interval, and $Q\in I_S$ since $S\subseteq Q$. 
    In addition, the routes in $Q\setminus S = \{R_1^Q,R_2^Q\}$ allow only ccw rotations and hence (ii) follows.
    To see (iii), we note that the routes in $S$ form a codimension $2$ inner face of a triangulation of $\calF_G$.
    Thus projecting along $S$ all the faces containing $S$ gives a triangulated polygon.
    Alternatively, a purely combinatorial proof follows from our proof of~(iv) and~(v) below.


    Next we show (v) by investigating the possible structures of the interval $I_S$. 
    Statement~(iv) will be proven in the process. 
    There are two cases to consider, namely when 
    
    Case (1): $R_1$ and $R_2$ are coherent, and 
    
    Case (2): when $R_1$ and $R_2$ are incoherent.
    
    We will see that Case (1) gives rise to a square, while Case (2) is split in several sub-cases, some of which give pentagons and some of which are hexagons. 
    In each case, we show that~$C_1 \vee C_2$ exists and is $C_{\max}(S)$.
    We start by considering the first case.

    \textbf{Case (1)}: $R_1$ and $R_2$ are coherent. 

    In this case, we have necessarily that $C_{\max}(S) \setminus S = \{R_1, R_2\}$, with the interval $[Q, C_{\max}(S)]$ being a square (see Figure~\ref{fig.square}). The maximal cliques of this square are, $Q$, $C_1$, $C_2$ and~$C_{\max}(S)$.

    It remains to show that~$C_1 \vee C_2$ exists and is equal to $C_{\max}(S)$. For this, 
    let $M$ be a maximal clique satisfying $C_1\leq M$ and $C_2 \leq M$. We want to show that $C_{\max}(S)\leq M$. 
    Let~$R^M$ be a route in $M$, $R$ be a route in $C_{\max}(S)$, and $v \in R^M\cap R$. 
    By Theorem~\ref{thm.ccw}, it suffices to show that $R^M$ is coherent with~$R$ at $v$ or that $R <_v^{\cw} R^M$. Every route $R$ in $S\cup \{R_1,R_2\}$ satisfies this condition, because $C_1\leq M$ and $C_2\leq M$. Since $C_{\max}(S) = S\cup \{R_1, R_2\}$ then we are done. 


\begin{figure}[htb]
    \centering
    \begin{tikzpicture}[scale=0.9]
\begin{scope}[xshift=0, scale=0.5]
    \draw[ultra thick,color=blue] (2.5,0.6) .. controls (2.6, 0.6) and (3.6, 0.6) .. (4,0);
    \draw[ultra thick,color=blue] (4,0.045)--(5,0.045);
    \draw[ultra thick,color=blue] (5,-0.045)--(6,-0.045);
    \draw[ultra thick,color=blue, dashed] (1.5,0.6)--(2.5,0.6);
    \draw[ultra thick,color=blue] (6,0) .. controls (6.4, -0.6) and (7.4, -0.6) .. (7.5,-0.6);
    \draw[ultra thick,color=blue, dashed] (7.5,-0.6)--(8.8,-0.6);
    
    \draw[ultra thick,color=black!30] (2.5,-0.6) .. controls (2.6, -0.6) and (3.6, -0.5) .. (4,0);
    \draw[ultra thick,color=black!30] (4,-0.045)--(5,-0.045);
    \draw[ultra thick,color=black!30] (5,0.045)--(6,0.045);
    \draw[ultra thick,color=black!30, dashed] (1.5,-0.6)--(2.5,-0.6);    
    \draw[ultra thick,color=black!30] (6,0) .. controls (6.4, 0.6) and (7.4, 0.6) .. (7.5,0.6);
    \draw[ultra thick,color=black!30, dashed] (7.5,0.6)--(8.8,0.6);

    \node[circle, draw, inner sep=1.4pt, fill,label=below:\scriptsize{$v_1$}] (s) at (5,0) {};    
        
    \node[] (P) at (0.6,0.5) {\scriptsize \color{blue} $R_1^Q$};
    \node[] (Q) at (0.6,-0.6) {\scriptsize \color{black!30} $R_1$};
    
    \node[] (P) at (9.5,0.6) {\scriptsize \color{black!30} $R_1$};
    \node[] (Q) at (9.5,-0.5) {\scriptsize \color{blue} $R_1^Q$};

    \draw[ultra thick,color=red] (2.5,2.6) .. controls (2.6, 2.6) and (3.6, 2.6) .. (4,2);
    \draw[ultra thick,color=red] (4,2.045)--(5,2.045);
    \draw[ultra thick,color=red] (5,1.955)--(6,1.955);
    \draw[ultra thick,color=red, dashed] (1.5,2.6)--(2.5,2.6);
    \draw[ultra thick,color=red] (6,2) .. controls (6.4, 1.4) and (7.4, 1.4) .. (7.5,1.4);
    \draw[ultra thick,color=red, dashed] (7.5,1.4)--(8.8,1.4);
    
    \draw[ultra thick,color=black!30] (2.5,1.4) .. controls (2.6, 1.4) and (3.6, 1.5) .. (4,2);
    \draw[ultra thick,color=black!30] (4,1.955)--(5,1.955);
    \draw[ultra thick,color=black!30] (5,2.045)--(6,2.045);
    \draw[ultra thick,color=black!30, dashed] (1.5,1.4)--(2.5,1.4);    
    \draw[ultra thick,color=black!30] (6,2) .. controls (6.4, 2.6) and (7.4, 2.6) .. (7.5,2.6);
    \draw[ultra thick,color=black!30, dashed] (7.5,2.6)--(8.8,2.6);

    \node[circle, draw, inner sep=1.4pt, fill,label=below:\scriptsize{$v_2$}] (s) at (5,2) {};    
        
    \node[] (P) at (0.6,2.5) {\scriptsize \color{red} $R_2^Q$};
    \node[] (Q) at (0.6,1.4) {\scriptsize \color{black!30} $R_2$};
    
    \node[] (P) at (9.5,2.6) {\scriptsize \color{black!30} $R_2$};
    \node[] (Q) at (9.5,1.5) {\scriptsize \color{red} $R_2^Q$};
\end{scope}

\begin{scope}[xshift=-140, yshift = 80, scale=0.5]
    \draw[ultra thick,color=black!30] (2.5,0.6) .. controls (2.6, 0.6) and (3.6, 0.6) .. (4,0);
    \draw[ultra thick,color=black!30] (4,0.045)--(5,0.045);
    \draw[ultra thick,color=black!30] (5,-0.045)--(6,-0.045);
    \draw[ultra thick,color=black!30, dashed] (1.5,0.6)--(2.5,0.6);
    \draw[ultra thick,color=black!30] (6,0) .. controls (6.4, -0.6) and (7.4, -0.6) .. (7.5,-0.6);
    \draw[ultra thick,color=black!30, dashed] (7.5,-0.6)--(8.8,-0.6);
    
    \draw[ultra thick,color=blue] (2.5,-0.6) .. controls (2.6, -0.6) and (3.6, -0.5) .. (4,0);
    \draw[ultra thick,color=blue] (4,-0.045)--(5,-0.045);
    \draw[ultra thick,color=blue] (5,0.045)--(6,0.045);
    \draw[ultra thick,color=blue, dashed] (1.5,-0.6)--(2.5,-0.6);    
    \draw[ultra thick,color=blue] (6,0) .. controls (6.4, 0.6) and (7.4, 0.6) .. (7.5,0.6);
    \draw[ultra thick,color=blue, dashed] (7.5,0.6)--(8.8,0.6);

    \node[circle, draw, inner sep=1.4pt, fill,label=below:\scriptsize{$v_1$}] (s) at (5,0) {};    
        
    \node[] (P) at (0.6,0.5) {\scriptsize \color{black!30} $R_1^Q$};
    \node[] (Q) at (0.6,-0.6) {\scriptsize \color{blue} $R_1$};
    
    \node[] (P) at (9.5,0.6) {\scriptsize \color{blue} $R_1$};
    \node[] (Q) at (9.5,-0.5) {\scriptsize \color{black!30} $R_1^Q$};

    \draw[ultra thick,color=red] (2.5,2.6) .. controls (2.6, 2.6) and (3.6, 2.6) .. (4,2);
    \draw[ultra thick,color=red] (4,2.045)--(5,2.045);
    \draw[ultra thick,color=red] (5,1.955)--(6,1.955);
    \draw[ultra thick,color=red, dashed] (1.5,2.6)--(2.5,2.6);
    \draw[ultra thick,color=red] (6,2) .. controls (6.4, 1.4) and (7.4, 1.4) .. (7.5,1.4);
    \draw[ultra thick,color=red, dashed] (7.5,1.4)--(8.8,1.4);
    
    \draw[ultra thick,color=black!30] (2.5,1.4) .. controls (2.6, 1.4) and (3.6, 1.5) .. (4,2);
    \draw[ultra thick,color=black!30] (4,1.955)--(5,1.955);
    \draw[ultra thick,color=black!30] (5,2.045)--(6,2.045);
    \draw[ultra thick,color=black!30, dashed] (1.5,1.4)--(2.5,1.4);    
    \draw[ultra thick,color=black!30] (6,2) .. controls (6.4, 2.6) and (7.4, 2.6) .. (7.5,2.6);
    \draw[ultra thick,color=black!30, dashed] (7.5,2.6)--(8.8,2.6);

    \node[circle, draw, inner sep=1.4pt, fill,label=below:\scriptsize{$v_2$}] (s) at (5,2) {};    
        
    \node[] (P) at (0.6,2.5) {\scriptsize \color{red} $R_2^Q$};
    \node[] (Q) at (0.6,1.4) {\scriptsize \color{black!30} $R_2$};
    
    \node[] (P) at (9.5,2.6) {\scriptsize \color{black!30} $R_2$};
    \node[] (Q) at (9.5,1.5) {\scriptsize \color{red} $R_2^Q$};
\end{scope}

\begin{scope}[xshift=140, yshift = 80, scale=0.5]
    \draw[ultra thick,color=blue] (2.5,0.6) .. controls (2.6, 0.6) and (3.6, 0.6) .. (4,0);
    \draw[ultra thick,color=blue] (4,0.045)--(5,0.045);
    \draw[ultra thick,color=blue] (5,-0.045)--(6,-0.045);
    \draw[ultra thick,color=blue, dashed] (1.5,0.6)--(2.5,0.6);
    \draw[ultra thick,color=blue] (6,0) .. controls (6.4, -0.6) and (7.4, -0.6) .. (7.5,-0.6);
    \draw[ultra thick,color=blue, dashed] (7.5,-0.6)--(8.8,-0.6);
    
    \draw[ultra thick,color=black!30] (2.5,-0.6) .. controls (2.6, -0.6) and (3.6, -0.5) .. (4,0);
    \draw[ultra thick,color=black!30] (4,-0.045)--(5,-0.045);
    \draw[ultra thick,color=black!30] (5,0.045)--(6,0.045);
    \draw[ultra thick,color=black!30, dashed] (1.5,-0.6)--(2.5,-0.6);    
    \draw[ultra thick,color=black!30] (6,0) .. controls (6.4, 0.6) and (7.4, 0.6) .. (7.5,0.6);
    \draw[ultra thick,color=black!30, dashed] (7.5,0.6)--(8.8,0.6);

    \node[circle, draw, inner sep=1.4pt, fill,label=below:\scriptsize{$v_1$}] (s) at (5,0) {};    
        
    \node[] (P) at (0.6,0.5) {\scriptsize \color{blue} $R_1^Q$};
    \node[] (Q) at (0.6,-0.6) {\scriptsize \color{black!30} $R_1$};
    
    \node[] (P) at (9.5,0.6) {\scriptsize \color{black!30} $R_1$};
    \node[] (Q) at (9.5,-0.5) {\scriptsize \color{blue} $R_1^Q$};

    \draw[ultra thick,color=black!30] (2.5,2.6) .. controls (2.6, 2.6) and (3.6, 2.6) .. (4,2);
    \draw[ultra thick,color=black!30] (4,2.045)--(5,2.045);
    \draw[ultra thick,color=black!30] (5,1.955)--(6,1.955);
    \draw[ultra thick,color=black!30, dashed] (1.5,2.6)--(2.5,2.6);
    \draw[ultra thick,color=black!30] (6,2) .. controls (6.4, 1.4) and (7.4, 1.4) .. (7.5,1.4);
    \draw[ultra thick,color=black!30, dashed] (7.5,1.4)--(8.8,1.4);
    
    \draw[ultra thick,color=red] (2.5,1.4) .. controls (2.6, 1.4) and (3.6, 1.5) .. (4,2);
    \draw[ultra thick,color=red] (4,1.955)--(5,1.955);
    \draw[ultra thick,color=red] (5,2.045)--(6,2.045);
    \draw[ultra thick,color=red, dashed] (1.5,1.4)--(2.5,1.4);    
    \draw[ultra thick,color=red] (6,2) .. controls (6.4, 2.6) and (7.4, 2.6) .. (7.5,2.6);
    \draw[ultra thick,color=red, dashed] (7.5,2.6)--(8.8,2.6);

    \node[circle, draw, inner sep=1.4pt, fill,label=below:\scriptsize{$v_2$}] (s) at (5,2) {};    
        
    \node[] (P) at (0.6,2.5) {\scriptsize \color{black!30} $R_2^Q$};
    \node[] (Q) at (0.6,1.4) {\scriptsize \color{red} $R_2$};
    
    \node[] (P) at (9.5,2.6) {\scriptsize \color{red} $R_2$};
    \node[] (Q) at (9.5,1.5) {\scriptsize \color{black!30} $R_2^Q$};
\end{scope}

\begin{scope}[xshift=0, yshift = 160, scale=0.5]
    \draw[ultra thick,color=black!30] (2.5,0.6) .. controls (2.6, 0.6) and (3.6, 0.6) .. (4,0);
    \draw[ultra thick,color=black!30] (4,0.045)--(5,0.045);
    \draw[ultra thick,color=black!30] (5,-0.045)--(6,-0.045);
    \draw[ultra thick,color=black!30, dashed] (1.5,0.6)--(2.5,0.6);
    \draw[ultra thick,color=black!30] (6,0) .. controls (6.4, -0.6) and (7.4, -0.6) .. (7.5,-0.6);
    \draw[ultra thick,color=black!30, dashed] (7.5,-0.6)--(8.8,-0.6);
    
    \draw[ultra thick,color=blue] (2.5,-0.6) .. controls (2.6, -0.6) and (3.6, -0.5) .. (4,0);
    \draw[ultra thick,color=blue] (4,-0.045)--(5,-0.045);
    \draw[ultra thick,color=blue] (5,0.045)--(6,0.045);
    \draw[ultra thick,color=blue, dashed] (1.5,-0.6)--(2.5,-0.6);    
    \draw[ultra thick,color=blue] (6,0) .. controls (6.4, 0.6) and (7.4, 0.6) .. (7.5,0.6);
    \draw[ultra thick,color=blue, dashed] (7.5,0.6)--(8.8,0.6);

    \node[circle, draw, inner sep=1.4pt, fill,label=below:\scriptsize{$v_1$}] (s) at (5,0) {};    
        
    \node[] (P) at (0.6,0.5) {\scriptsize \color{black!30} $R_1^Q$};
    \node[] (Q) at (0.6,-0.6) {\scriptsize \color{blue} $R_1$};
    
    \node[] (P) at (9.5,0.6) {\scriptsize \color{blue} $R_1$};
    \node[] (Q) at (9.5,-0.5) {\scriptsize \color{black!30} $R_1^Q$};

    \draw[ultra thick,color=black!30] (2.5,2.6) .. controls (2.6, 2.6) and (3.6, 2.6) .. (4,2);
    \draw[ultra thick,color=black!30] (4,2.045)--(5,2.045);
    \draw[ultra thick,color=black!30] (5,1.955)--(6,1.955);
    \draw[ultra thick,color=black!30, dashed] (1.5,2.6)--(2.5,2.6);
    \draw[ultra thick,color=black!30] (6,2) .. controls (6.4, 1.4) and (7.4, 1.4) .. (7.5,1.4);
    \draw[ultra thick,color=black!30, dashed] (7.5,1.4)--(8.8,1.4);
    
    \draw[ultra thick,color=red] (2.5,1.4) .. controls (2.6, 1.4) and (3.6, 1.5) .. (4,2);
    \draw[ultra thick,color=red] (4,1.955)--(5,1.955);
    \draw[ultra thick,color=red] (5,2.045)--(6,2.045);
    \draw[ultra thick,color=red, dashed] (1.5,1.4)--(2.5,1.4);    
    \draw[ultra thick,color=red] (6,2) .. controls (6.4, 2.6) and (7.4, 2.6) .. (7.5,2.6);
    \draw[ultra thick,color=red, dashed] (7.5,2.6)--(8.8,2.6);

    \node[circle, draw, inner sep=1.4pt, fill,label=below:\scriptsize{$v_2$}] (s) at (5,2) {};    
        
    \node[] (P) at (0.6,2.5) {\scriptsize \color{black!30} $R_2^Q$};
    \node[] (Q) at (0.6,1.4) {\scriptsize \color{red} $R_2$};
    
    \node[] (P) at (9.5,2.6) {\scriptsize \color{red} $R_2$};
    \node[] (Q) at (9.5,1.5) {\scriptsize \color{black!30} $R_2^Q$};
\end{scope}

\begin{scope}[xshift=-10, yshift=30, scale = 0.5]
    \node[] (a) at (0,0) {};
    \node[] (b) at (-2,2) {};
    \draw[very thick] (a) -- (b);
\end{scope}
\begin{scope}[xshift=150, yshift=30, scale = 0.5]
    \node[] (a) at (0,0) {};
    \node[] (b) at (2,2) {};
    \draw[very thick] (a) -- (b);
\end{scope}
\begin{scope}[xshift=-35, yshift=125, scale = 0.5]
    \node[] (a) at (0,0) {};
    \node[] (b) at (2,2) {};
    \draw[very thick] (a) -- (b);
\end{scope}
\begin{scope}[xshift=175, yshift=125, scale = 0.5]
    \node[] (a) at (0,0) {};
    \node[] (b) at (-2,2) {};
    \draw[very thick] (a) -- (b);
\end{scope}

\begin{scope}[xshift=-155, yshift=95, scale = 0.5]
    \node[] (a) at (0,0) {$C_1$};
\end{scope}
\begin{scope}[xshift=125, yshift=95, scale = 0.5]
    \node[] (a) at (0,0) {$C_2$};
\end{scope}
\begin{scope}[xshift=-15, yshift=15, scale = 0.5]
    \node[] (a) at (0,0) {$Q$};
\end{scope}
\begin{scope}[xshift=-30, yshift=175, scale = 0.5]
    \node[] (a) at (0,0) {$C_1\vee C_2$};
\end{scope}

\end{tikzpicture}
    \caption{Case (1) in~\Cref{prop.intervalJoin}.}
\label{fig.square}
\end{figure}

    \textbf{Case (2)}: $R_1$ and $R_2$ are incoherent.

    Suppose that $R_1$ and $R_2$ are incoherent at a vertex $x$.
    Without lost of generality, we can assume that $R_1<_x^{\cw} R_2$.
    Let $P_x$ denote the maximal path containing $x$ in $R_1\cap R_2$ at which~$R_1$ and~$R_2$ are incoherent. 
    Let $w_1$ be a vertex at which $R_1^Q$ and $R_1$ are incoherent, and let $P_{w_1}$ be the path at which they are incoherent. 
    We assume that $x < w_1$ (the case $w_1 < x$ is obtained by a $180^\circ$ rotation symmetry).
    By Lemma~\ref{lem.R2Q}, we have that $P_x\cap P_{w_1} = \emptyset$, and $$R_2^Q = \Bot(R_1^Q,R_1) = R_1w_1R_1^Q.$$ 

    To simplify notation, let $P_{w_2}$ be the path at which $R_2^Q$ and $R_2$ are incoherent.
    For $i\in \{1,2\}$ let $b_i = \min(P_{w_i})$ and $c_i = \max(P_{w_i})$. 
    Let $d_2$ denote the smallest point after which~$R_2$ and~$R_2^Q$ coincide. 
    We consider the following sub-cases, which are depicted in Figure~\ref{fig.join_proof_cases}:

\begin{figure}
    \centering 
    \resizebox{\textwidth}{!}{%
    \begin{tikzpicture}
\begin{scope}[xshift=-13, yshift=-350, scale=1.7]
    \draw[thick, fill=blue!05] (0,0) -- (0,8) -- (5,8) -- (5,2) -- (10.3,2) -- (10.3,0) -- cycle;
\end{scope}

\begin{scope}[xshift=240, yshift=-145, scale=1.7]
    \draw[thick, fill=red!05] (0,0) -- (5,0) -- (5,1.7) -- (0,1.7) -- cycle;
\end{scope}

\begin{scope}[xshift=0]
    \draw[ultra thick,color=black!30] (2.5,0.6) .. controls (2.6, 0.6) and (3.6, 0.6) .. (4,0);
    \draw[ultra thick,color=blue] (4,0)--(5,0);
    \draw[ultra thick,color=black!30, dashed] (2,0.6)--(2.5,0.6);
    \draw[ultra thick,color=blue] (5,0) .. controls (5.4, -0.6) and (6.4, -0.6) .. (6.5,-0.6);
    \draw[ultra thick,color=blue, dashed] (6.5,-0.6)--(7,-0.6);
    
    \draw[ultra thick,color=blue] (2.5,-0.6) .. controls (2.6, -0.6) and (3.6, -0.5) .. (4,0);
    \draw[ultra thick,color=blue, dashed] (1,-0.6)--(1.5,-0.6);
    \draw[ultra thick,color=blue] (1.5,-0.6)--(2.5,-0.6);    
    \draw[ultra thick,color=black!30] (5,0) .. controls (5.4, 0.6) and (6.4, 0.6) .. (6.5,0.6);
    \draw[ultra thick,color=black!30, dashed] (6.5,0.6)--(7,0.6);

    \draw[ultra thick,color=red, dashed] (1,-1) --(1.5,-1);    
    \draw[ultra thick,color=red] (1.5,-1) .. controls (1.6, -1) and (1.9, -0.8) .. (2,-0.6);
    \draw[ultra thick,color=red] (2,-0.6)--(2.5,-0.6);  
    \draw[ultra thick,color=red] (2.5,-0.6) .. controls (2.7, -0.3) and (2.8, -0.3) .. (2.9,-0.3); 
    \draw[ultra thick,color=red, dashed] (2.9,-0.3) .. controls (3,-0.3) and (3.3, -0.7) .. (3.4,-0.6); 
    \draw[ultra thick,color=red] (3.4,-0.6) .. controls (3.5, -0.7) and (3.6, -0.4) .. (3.7,-0.25); 

    \node[circle, draw, inner sep=1.6pt, fill,label=below:\small{$b_2$}] (s) at (2,-0.6) {};
    \node[circle, draw, inner sep=1.6pt, fill,label=below:\small{$c_2$}] (s) at (2.5,-0.6) {};

    \node[circle, draw, inner sep=1.6pt, fill,label=below:\small{$d_2$}] (s) at (3.7,-0.25) {};
    
    \node[circle, draw, inner sep=1.6pt, fill,label=below:\small{$b_1$}] (s) at (4,0) {};    
    \node[circle, draw, inner sep=1.6pt, fill, label=below:\small{$c_1$}] (c) at (5,0) {};
       
    \node[] (P) at (1.6,0.6) {$R_1^Q$};
    \node[] (Q) at (0.6,-0.5) {$R_1$};
    \node[] (Q) at (0.6,-1.1) {$R_2$};
    
    \node[] (P) at (7.4,0.6) {$R_1$};
    \node[] (Q) at (7.4,-0.6) {$R_1^Q$};

    \node[] (P) at (0.4,0.8) {I.(i).};
    
\end{scope}

\begin{scope}[xshift=250]
    \draw[ultra thick,color=black!30] (2.5,0.6) .. controls (2.6, 0.6) and (3.6, 0.6) .. (4,0);
    \draw[ultra thick,color=blue] (4,0)--(5,0);
    \draw[ultra thick,color=black!30, dashed] (2,0.6)--(2.5,0.6);
    \draw[ultra thick,color=blue] (5,0) .. controls (5.4, -0.6) and (6.4, -0.6) .. (6.5,-0.6);
    \draw[ultra thick,color=blue, dashed] (6.5,-0.6)--(7,-0.6);
    
    \draw[ultra thick,color=blue] (2.5,-0.6) .. controls (2.6, -0.6) and (3.6, -0.5) .. (4,0);
    \draw[ultra thick,color=blue, dashed] (1,-0.6)--(1.5,-0.6);
    \draw[ultra thick,color=blue] (1.5,-0.6)--(2.5,-0.6);    
    \draw[ultra thick,color=black!30] (5,0) .. controls (5.4, 0.6) and (6.4, 0.6) .. (6.5,0.6);
    \draw[ultra thick,color=black!30, dashed] (6.5,0.6)--(7,0.6);

    \draw[ultra thick,color=red, dashed] (1,-1) --(1.5,-1);    
    \draw[ultra thick,color=red] (1.5,-1) .. controls (1.6, -1) and (1.9, -0.8) .. (2,-0.6);
    \draw[ultra thick,color=red] (2,-0.6)--(2.5,-0.6);  
    \draw[ultra thick,color=red] (2.5,-0.6) .. controls (2.6, -0.4) and (2.7, -0.3) .. (2.8,-0.2); 
    \draw[ultra thick,color=red, dashed] (2.8,-0.2) .. controls (3,0) and (3.2, 0) .. (3.3,0); 
    \draw[ultra thick,color=red] (3.3,0) .. controls (3.4, 0) and (3.6, -0.1) .. (3.7,-0.25); 

    \node[circle, draw, inner sep=1.6pt, fill,label=below:\small{$b_2$}] (s) at (2,-0.6) {};
    \node[circle, draw, inner sep=1.6pt, fill,label=below:\small{$c_2$}] (s) at (2.5,-0.6) {};

    \node[circle, draw, inner sep=1.6pt, fill,label=below:\small{$d_2$}] (s) at (3.7,-0.25) {};
    
    \node[circle, draw, inner sep=1.6pt, fill,label=below:\small{$b_1$}] (s) at (4,0) {};    
    \node[circle, draw, inner sep=1.6pt, fill, label=below:\small{$c_1$}] (c) at (5,0) {};
       
    \node[] (P) at (1.6,0.6) {$R_1^Q$};
    \node[] (Q) at (0.6,-0.5) {$R_1$};
    \node[] (Q) at (0.6,-1.1) {$R_2$};
    
    \node[] (P) at (7.4,0.6) {$R_1$};
    \node[] (Q) at (7.4,-0.6) {$R_1^Q$};

    \node[] (P) at (0.4,0.8) {I.(ii).};
\end{scope}

\begin{scope}[xshift=0,yshift=-200]
    \draw[ultra thick,color=black!30] (2.5,0.6) .. controls (2.6, 0.6) and (3.6, 0.6) .. (4,0);
    \draw[ultra thick,color=blue] (4,0)--(5,0);
    \draw[ultra thick,color=black!30, dashed] (2,0.6)--(2.5,0.6);
    \draw[ultra thick,color=blue] (5,0) .. controls (5.4, -0.6) and (6.4, -0.6) .. (6.5,-0.6);
    \draw[ultra thick,color=blue, dashed] (6.5,-0.6)--(7,-0.6);
    
    \draw[ultra thick,color=blue] (2.5,-0.6) .. controls (2.6, -0.6) and (3.6, -0.5) .. (4,0);
    \draw[ultra thick,color=blue, dashed] (1,-0.6)--(1.5,-0.6);
    \draw[ultra thick,color=blue] (1.5,-0.6)--(2.5,-0.6);    
    \draw[ultra thick,color=black!30] (5,0) .. controls (5.4, 0.6) and (6.4, 0.6) .. (6.5,0.6);
    \draw[ultra thick,color=black!30, dashed] (6.5,0.6)--(7,0.6);

    \draw[ultra thick,color=red, dashed] (1,-1) --(1.5,-1);    
    \draw[ultra thick,color=red] (1.5,-1) .. controls (1.6, -1) and (1.9, -0.8) .. (2,-0.6);
    \draw[ultra thick,color=red] (2,-0.6)--(2.5,-0.6);  
    \draw[ultra thick,color=red] (2.5,-0.6) .. controls (2.7, -0.3) and (2.8, -0.3) .. (2.9,-0.3); 
    \draw[ultra thick,color=red, dashed] (2.9,-0.3) .. controls (3,-0.3) and (3.3, -0.7) .. (3.8,-0.8); 
    \draw[ultra thick,color=red] (3.8,-0.8) .. controls (4, -0.8) and (4.3, -0.5) .. (4.5,0); 

    \node[circle, draw, inner sep=1.6pt, fill,label=below:\small{$b_2$}] (s) at (2,-0.6) {};
    \node[circle, draw, inner sep=1.6pt, fill,label=below:\small{$c_2$}] (s) at (2.5,-0.6) {};

    \node[circle, draw, inner sep=1.6pt, fill,label=below:\small{$d_2$}] (s) at (4.5,0) {};
    
    \node[circle, draw, inner sep=1.6pt, fill,label=below:\small{$b_1$}] (s) at (4,0) {};    
    \node[circle, draw, inner sep=1.6pt, fill, label=below:\small{$c_1$}] (c) at (5,0) {};
       
    \node[] (P) at (1.6,0.6) {$R_1^Q$};
    \node[] (Q) at (0.6,-0.5) {$R_1$};
    \node[] (Q) at (0.6,-1.1) {$R_2$};
    
    \node[] (P) at (7.4,0.6) {$R_1$};
    \node[] (Q) at (7.4,-0.6) {$R_1^Q$};

    \node[] (P) at (0.4,0.8) {III.(i).};
\end{scope}

\begin{scope}[xshift=250, yshift=-200]
    \draw[ultra thick,color=black!30] (2.5,0.6) .. controls (2.6, 0.6) and (3.6, 0.6) .. (4,0);
    \draw[ultra thick,color=blue] (4,0)--(5,0);
    \draw[ultra thick,color=black!30, dashed] (2,0.6)--(2.5,0.6);
    \draw[ultra thick,color=blue] (5,0) .. controls (5.4, -0.6) and (6.4, -0.6) .. (6.5,-0.6);
    \draw[ultra thick,color=blue, dashed] (6.5,-0.6)--(7,-0.6);
    
    \draw[ultra thick,color=blue] (2.5,-0.6) .. controls (2.6, -0.6) and (3.6, -0.5) .. (4,0);
    \draw[ultra thick,color=blue, dashed] (1,-0.6)--(1.5,-0.6);
    \draw[ultra thick,color=blue] (1.5,-0.6)--(2.5,-0.6);    
    \draw[ultra thick,color=black!30] (5,0) .. controls (5.4, 0.6) and (6.4, 0.6) .. (6.5,0.6);
    \draw[ultra thick,color=black!30, dashed] (6.5,0.6)--(7,0.6);

    \draw[ultra thick,color=red, dashed] (1,-1) --(1.5,-1);    
    \draw[ultra thick,color=red] (1.5,-1) .. controls (1.6, -1) and (1.9, -0.8) .. (2,-0.6);
    \draw[ultra thick,color=red] (2,-0.6)--(2.5,-0.6);  
    \draw[ultra thick,color=red] (2.5,-0.6) .. controls (2.6, -0.35) and (2.6, -0.35) .. (2.7,-0.2); 
    \draw[ultra thick,color=red, dashed] (2.7,-0.2) .. controls (2.9,0.1) and (3.5, 0.7) .. (4,0.5); 
    \draw[ultra thick,color=red] (4,0.5) .. controls (4.2, 0.4) and (4.3, 0.3) .. (4.5,0); 

    \node[circle, draw, inner sep=1.6pt, fill,label=below:\small{$b_2$}] (s) at (2,-0.6) {};
    \node[circle, draw, inner sep=1.6pt, fill,label=below:\small{$c_2$}] (s) at (2.5,-0.6) {};

    \node[circle, draw, inner sep=1.6pt, fill,label=below:\small{$d_2$}] (s) at (4.5,0) {};
    
    \node[circle, draw, inner sep=1.6pt, fill,label=below:\small{$b_1$}] (s) at (4,0) {};    
    \node[circle, draw, inner sep=1.6pt, fill, label=below:\small{$c_1$}] (c) at (5,0) {};
       
    \node[] (P) at (1.6,0.6) {$R_1^Q$};
    \node[] (Q) at (0.6,-0.5) {$R_1$};
    \node[] (Q) at (0.6,-1.1) {$R_2$};
    
    \node[] (P) at (7.4,0.6) {$R_1$};
    \node[] (Q) at (7.4,-0.6) {$R_1^Q$};

    \node[] (P) at (0.4,0.8) {III.(ii).};
\end{scope}

\begin{scope}[xshift=0, yshift=-300]
    \draw[ultra thick,color=black!30] (2.5,0.6) .. controls (2.6, 0.6) and (3.6, 0.6) .. (4,0);
    \draw[ultra thick,color=blue] (4,0)--(5,0);
    \draw[ultra thick,color=black!30, dashed] (2,0.6)--(2.5,0.6);
    \draw[ultra thick,color=blue] (5,0) .. controls (5.4, -0.6) and (6.4, -0.6) .. (6.5,-0.6);
    \draw[ultra thick,color=blue, dashed] (6.5,-0.6)--(7,-0.6);
    
    \draw[ultra thick,color=blue] (2.5,-0.6) .. controls (2.6, -0.6) and (3.6, -0.5) .. (4,0);
    \draw[ultra thick,color=blue, dashed] (1,-0.6)--(1.5,-0.6);
    \draw[ultra thick,color=blue] (1.5,-0.6)--(2.5,-0.6);    
    \draw[ultra thick,color=black!30] (5,0) .. controls (5.4, 0.6) and (6.4, 0.6) .. (6.5,0.6);
    \draw[ultra thick,color=black!30, dashed] (6.5,0.6)--(7,0.6);

    \draw[ultra thick,color=red, dashed] (1,-1) --(1.5,-1);    
    \draw[ultra thick,color=red] (1.5,-1) .. controls (1.6, -1) and (1.9, -0.8) .. (2,-0.6);
    \draw[ultra thick,color=red] (2,-0.6)--(2.5,-0.6);  
    \draw[ultra thick,color=red] (2.5,-0.6) .. controls (2.7, -0.3) and (2.8, -0.3) .. (2.9,-0.2); 
    \draw[ultra thick,color=red, dashed] (2.9,-0.2) .. controls (3.4,-0.2) and (3.8, -0.9) .. (5.3,-0.8); 
    \draw[ultra thick,color=red] (5.3,-0.8) .. controls (5.2, -0.8) and (5.8, -0.7) .. (6,-0.55); 

    \node[circle, draw, inner sep=1.6pt, fill,label=below:\small{$b_2$}] (s) at (2,-0.6) {};
    \node[circle, draw, inner sep=1.6pt, fill,label=below:\small{$c_2$}] (s) at (2.5,-0.6) {};

    \node[circle, draw, inner sep=1.6pt, fill,label=below:\small{$d_2$}] (s) at (6,-0.55) {};
    
    \node[circle, draw, inner sep=1.6pt, fill,label=below:\small{$b_1$}] (s) at (4,0) {};    
    \node[circle, draw, inner sep=1.6pt, fill, label=below:\small{$c_1$}] (c) at (5,0) {};
       
    \node[] (P) at (1.6,0.6) {$R_1^Q$};
    \node[] (Q) at (0.6,-0.5) {$R_1$};
    \node[] (Q) at (0.6,-1.1) {$R_2$};
    
    \node[] (P) at (7.4,0.6) {$R_1$};
    \node[] (Q) at (7.4,-0.6) {$R_1^Q$};

    \node[] (P) at (0.4,0.8) {IV.(i).};

\end{scope}

\begin{scope}[xshift=250, yshift=-300]
    \draw[ultra thick,color=black!30] (2.5,0.6) .. controls (2.6, 0.6) and (3.6, 0.6) .. (4,0);
    \draw[ultra thick,color=blue] (4,0)--(5,0);
    \draw[ultra thick,color=black!30, dashed] (2,0.6)--(2.5,0.6);
    \draw[ultra thick,color=blue] (5,0) .. controls (5.4, -0.6) and (6.4, -0.6) .. (6.5,-0.6);
    \draw[ultra thick,color=blue, dashed] (6.5,-0.6)--(7,-0.6);
    
    \draw[ultra thick,color=blue] (2.5,-0.6) .. controls (2.6, -0.6) and (3.6, -0.5) .. (4,0);
    \draw[ultra thick,color=blue, dashed] (1,-0.6)--(1.5,-0.6);
    \draw[ultra thick,color=blue] (1.5,-0.6)--(2.5,-0.6);    
    \draw[ultra thick,color=black!30] (5,0) .. controls (5.4, 0.6) and (6.4, 0.6) .. (6.5,0.6);
    \draw[ultra thick,color=black!30, dashed] (6.5,0.6)--(7,0.6);

    \draw[ultra thick,color=red, dashed] (1,-1) --(1.5,-1);    
    \draw[ultra thick,color=red] (1.5,-1) .. controls (1.6, -1) and (1.9, -0.8) .. (2,-0.6);
    \draw[ultra thick,color=red] (2,-0.6)--(2.5,-0.6);  
    \draw[ultra thick,color=red] (2.5,-0.6) .. controls (2.6, -0.3) and (2.7, -0.2) .. (2.8,-0.1); 
    \draw[ultra thick,color=red, dashed] (2.8,-0.1) .. controls (3.4,0.7) and (5.2, 0.7) .. (5.7,0); 
    \draw[ultra thick,color=red] (5.7,0) .. controls (5.8, -0.1) and (5.9, -0.3) .. (6,-0.55); 

    \node[circle, draw, inner sep=1.6pt, fill,label=below:\small{$b_2$}] (s) at (2,-0.6) {};
    \node[circle, draw, inner sep=1.6pt, fill,label=below:\small{$c_2$}] (s) at (2.5,-0.6) {};

    \node[circle, draw, inner sep=1.6pt, fill,label=below:\small{$d_2$}] (s) at (6,-0.55) {};
    
    \node[circle, draw, inner sep=1.6pt, fill,label=below:\small{$b_1$}] (s) at (4,0) {};    
    \node[circle, draw, inner sep=1.6pt, fill, label=below:\small{$c_1$}] (c) at (5,0) {};
       
    \node[] (P) at (1.6,0.6) {$R_1^Q$};
    \node[] (Q) at (0.6,-0.5) {$R_1$};
    \node[] (Q) at (0.6,-1.1) {$R_2$};
    
    \node[] (P) at (7.4,0.6) {$R_1$};
    \node[] (Q) at (7.4,-0.6) {$R_1^Q$};
    
    \node[] (P) at (0.4,0.8) {IV.(ii).};
\end{scope}

\begin{scope}[xshift=0, yshift=-100]
    \draw[ultra thick,color=black!30] (2.5,0.6) .. controls (2.6, 0.6) and (3.6, 0.6) .. (4,0);
    \draw[ultra thick,color=blue] (4,0)--(5,0);
    \draw[ultra thick,color=black!30, dashed] (2,0.6)--(2.5,0.6);
    \draw[ultra thick,color=blue] (5,0) .. controls (5.4, -0.6) and (6.4, -0.6) .. (6.5,-0.6);
    \draw[ultra thick,color=blue, dashed] (6.5,-0.6)--(7,-0.6);
    
    \draw[ultra thick,color=blue] (2.5,-0.6) .. controls (2.6, -0.6) and (3.6, -0.5) .. (4,0);
    \draw[ultra thick,color=blue, dashed] (1,-0.6)--(1.5,-0.6);
    \draw[ultra thick,color=blue] (1.5,-0.6)--(2.5,-0.6);    
    \draw[ultra thick,color=black!30] (5,0) .. controls (5.4, 0.6) and (6.4, 0.6) .. (6.5,0.6);
    \draw[ultra thick,color=black!30, dashed] (6.5,0.6)--(7,0.6);

    \draw[ultra thick,color=red, dashed] (1,-1) --(1.5,-1);    
    \draw[ultra thick,color=red] (1.5,-1) .. controls (1.6, -1) and (1.9, -0.8) .. (2,-0.6);
    \draw[ultra thick,color=red] (2,-0.6)--(2.5,-0.6);  
    \draw[ultra thick,color=red] (2.5,-0.6) .. controls (2.7, -0.3) and (2.8, -0.3) .. (2.9,-0.3); 
    \draw[ultra thick,color=red, dashed] (2.9,-0.3) .. controls (3,-0.3) and (3.3, -0.7) .. (3.5,-0.7); 
    \draw[ultra thick,color=red] (3.5,-0.7) .. controls (3.6, -0.7) and (3.8, -0.5) .. (4,0); 

    \node[circle, draw, inner sep=1.6pt,  fill,label=below:\small{$b_2$}] (s) at (2,-0.6) {};
    \node[circle, draw, inner sep=1.6pt, fill,label=below:\small{$c_2$}] (s) at (2.5,-0.6) {};

    \node[label=below:\small{$b_1$}] (s) at (4,0) {};
    \node[label=below:\small{$=$}] (s) at (4.3,-0.1) {};    
    \node[label=below:\small{$d_2$}] (s) at (4.7,0) {};
    \node[circle, draw, inner sep=1.6pt, fill] (s) at (4,0) {};    
    \node[circle, draw, inner sep=1.6pt, fill, label=above:\small{$c_1$}] (c) at (5,0) {};
       
    \node[] (P) at (1.6,0.6) {$R_1^Q$};
    \node[] (Q) at (0.6,-0.5) {$R_1$};
    \node[] (Q) at (0.6,-1.1) {$R_2$};
    
    \node[] (P) at (7.4,0.6) {$R_1$};
    \node[] (Q) at (7.4,-0.6) {$R_1^Q$};

    \node[] (P) at (0.4,0.8) {II.(i).};
\end{scope}

\begin{scope}[xshift=250, yshift=-100]
    \draw[ultra thick,color=black!30] (2.5,0.6) .. controls (2.6, 0.6) and (3.6, 0.6) .. (4,0);
    \draw[ultra thick,color=blue] (4,0)--(5,0);
    \draw[ultra thick,color=black!30, dashed] (2,0.6)--(2.5,0.6);
    \draw[ultra thick,color=blue] (5,0) .. controls (5.4, -0.6) and (6.4, -0.6) .. (6.5,-0.6);
    \draw[ultra thick,color=blue, dashed] (6.5,-0.6)--(7,-0.6);
    
    \draw[ultra thick,color=blue] (2.5,-0.6) .. controls (2.6, -0.6) and (3.6, -0.5) .. (4,0);
    \draw[ultra thick,color=blue, dashed] (1,-0.6)--(1.5,-0.6);
    \draw[ultra thick,color=blue] (1.5,-0.6)--(2.5,-0.6);    
    \draw[ultra thick,color=black!30] (5,0) .. controls (5.4, 0.6) and (6.4, 0.6) .. (6.5,0.6);
    \draw[ultra thick,color=black!30, dashed] (6.5,0.6)--(7,0.6);

    \draw[ultra thick,color=red, dashed] (1,-1) --(1.5,-1);    
    \draw[ultra thick,color=red] (1.5,-1) .. controls (1.6, -1) and (1.9, -0.8) .. (2,-0.6);
    \draw[ultra thick,color=red] (2,-0.6)--(2.5,-0.6);  
    \draw[ultra thick,color=red] (2.5,-0.6) .. controls (2.6, -0.35) and (2.6, -0.35) .. (2.7,-0.2); 
    \draw[ultra thick,color=red, dashed] (2.7,-0.2) .. controls (2.9,0.1) and (3.4, 1) .. (3.7,0.5); 
    \draw[ultra thick,color=red] (3.7,0.5) .. controls (3.8, 0.4) and (3.9, 0.3) .. (4,0); 

    \node[circle, draw, inner sep=1.6pt, fill,label=below:\small{$b_2$}] (s) at (2,-0.6) {};
    \node[circle, draw, inner sep=1.6pt, fill,label=below:\small{$c_2$}] (s) at (2.5,-0.6) {};

    \node[label=below:\small{$b_1$}] (s) at (4,0) {};
    \node[label=below:\small{$=$}] (s) at (4.3,-0.1) {};    
    \node[label=below:\small{$d_2$}] (s) at (4.7,0) {};
    \node[circle, draw, inner sep=1.6pt, fill] (s) at (4,0) {};    
    \node[circle, draw, inner sep=1.6pt, fill, label=above:\small{$c_1$}] (c) at (5,0) {};
       
    \node[] (P) at (1.6,0.6) {$R_1^Q$};
    \node[] (Q) at (0.6,-0.5) {$R_1$};
    \node[] (Q) at (0.6,-1.1) {$R_2$};
    
    \node[] (P) at (7.4,0.6) {$R_1$};
    \node[] (Q) at (7.4,-0.6) {$R_1^Q$};

    \node[] (P) at (0.4,0.8) {II.(ii).};
\end{scope}
\end{tikzpicture}   
    }
    \caption{The eight sub-cases of Case (2) in the proof of~\Cref{prop.intervalJoin}(v). The five cases in the light-blue region give rise to pentagons, the case in the light-orange region gives rise to the hexagon, and the remaining two cases are impossible.}
\label{fig.join_proof_cases}
\end{figure}
   
    \begin{itemize}
        \item[(I)] $d_2 < b_1$
        \begin{itemize}
            \item[(i)] $R_1d_2 <_{\scrI(d_2)} R_2d_2$ 
            \item[(ii)] $R_2d_2 <_{\scrI(d_2)} R_1d_2$
        \end{itemize}
        \item[(II)] $d_2 = b_1$
        \begin{itemize}
            \item[(i)] $R_1d_2 <_{\scrI(d_2)} R_2d_2$ 
            \item[(ii)] $R_2d_2 <_{\scrI(d_2)} R_1d_2$
        \end{itemize}            
        \item[(III)] $b_1 < d_2 \leq c_1$
        \begin{itemize}
            \item[(i)] $R_1d_2 <_{\scrI(d_2)} R_2d_2$ 
            \item[(ii)] $R_2d_2 <_{\scrI(d_2)} R_1d_2$
        \end{itemize}
        \item[(IV)] $c_1 < d_2$
        \begin{itemize}
            \item[(i)] $R_1^Qd_2 <_{\scrI(d_2)} R_2d_2$
            \item[(ii)] $R_2d_2 <_{\scrI(d_2)} R_1^Qd_2$ 
        \end{itemize}
    \end{itemize}

    In order to analyze each of these cases, we use the following strategy. 

        \textbf{Strategy.}
        In order find the polygon consisting of all the maximal cliques containing $S$, and proving that $C_1 \vee C_2$ exists and is $C_{\max}(S)$,  carry out the next steps. Let $\widetilde{G}$ be the subgraph of $G$ induced by $R_1^Q$, $R_2^Q$, $R_1$, and $R_2$.

        \begin{itemize}
            \item[\bf Step 1.] Find the set $\widetilde{S}$ of routes in $S$ inside $\widetilde{G}$. 

            \noindent
            These are the routes $\widetilde{S}=\{\Top_1,\Bot_1,\Top_2,\Bot_2\}\smallsetminus \{R_1^Q,R_2^Q,R_1,R_2\}$, where
            $\Top_i:=\Top(R_i^Q,R_i)$ and $\Bot_i:=\Bot(R_i^Q,R_i)$.

            \item[\bf Step 2.] Find the set $\widetilde{Z}$ of routes in $\widetilde{G}$ that are coherent with $\widetilde{S}$, but are itself not in $\widetilde{S}$. 

            \noindent
            These are $R_1^Q$, $R_2^Q$, $R_1$, $R_2$ and possibly others. (we will see that there are zero, one or two more in each case).

            \item[Step 2a.] Argue that the routes in $\widetilde{Z}$, from Step 2, are coherent with $S$.

            \noindent
            This step needs to be verified in each case (by a simple argument). If there is a route $R\in S$ incoherent with a route in $\widetilde{Z}$, then it would have to be incoherent with either $R_1^Q$, $R_2^Q$, $R_1$, or $R_2$, which is a contradiction.

            \item[\bf Step 3.] Find the coherence graph among the routes in $\widetilde{Z}$ (or equivalently $\widetilde{Z}\cup S$). 

            \noindent
            This step is straight forward. We will see that for each case this is a triangulated polygon: either a square, pentagon, or hexagon. 

            \item[Step 3a.] Argue that the maximal cliques determined by the polygon in Step 3 are exactly the maximal cliques containing $S$. 

            \noindent This step is automatic for the following reason. In each case, the maximal cliques containing $S$ determine a 2-dimensional triangulated polygon with a point representing $S$ in the middle. But we know that intersecting the maximal simplices of the framed triangulation of the flow polytope with the 2-dimensional space orthogonal to $S$ gives a triangulated polygon. So, the two triangulated polygons have to be the same.  

            \item[\bf Step 4.] The last step is to find $C_{\max}(S)$ among the maximal cliques in the polygon, and argue that $C_1 \vee C_2=C_{\max}(S)$.

            \noindent
            This step is relatively straightforward. 
            
        \end{itemize}

    Now, we proceed to analyze each of the cases using the previous strategy. We do this in an order that is convenient for us. First, we start by discarding an impossible case. 
    

    \textbf{Case I.(ii).} Impossible.

    Since $d_2 < b_1$, note that we cannot have $R_2d_2 <_{\scrI(d_2)} R_1d_2$ as otherwise $\Top(R_2^Q,R_2) \in Q$ violates the in between property of $R_1^Q$ and $R_1$ at $v_1$.
    Therefore case I.(i) is not possible.

    \textbf{Case I.(i).} Pentagon.

\begin{figure}[htb]
    \centering
    \scalebox{0.9}{

    }
    \caption{Steps 1-3 in the strategy applied to Case (2.I.i).}
    \label{fig_photo_caseI_i}
\end{figure}

    Next, we consider Case I.(i). 
    The outcome of the steps throughout our strategy are illustrated in~\Cref{fig_photo_caseI_i}. 
    In Step 1 of our strategy, we find the set $\widetilde{S}=\{\Top_1,\Top_2,\Bot_2\}$; note that $\Bot_1=\Bot(R_1^Q,R_1)=R_2^Q$ does not appear in this set because $R_2^Q\notin S$.
    In Step~2, we find the set $\widetilde{Z}=\{R_1^Q, R_2^Q, R_1, R_2, R_2b_2R_1\}$ of routes in $\widetilde{G}$ that are not in $\widetilde S$ and are coherent with $\widetilde{S}$.
    
    In Step~2a, we need to argue that the routes in $\widetilde{Z}$ are coherent with $S$. Since this holds for $R_1^Q, R_2^Q, R_1, R_2$, we only need to check that $R_2b_2R_1$ is coherent with $S$. Suppose there exists a route $R\in S$ that is incoherent with $R_2b_2R_1$. Then, they have to be incoherent at $b_2$, otherwise $R\in S$ would be incoherent with either $R_2$ or $R_1$ which is a contradiction. There are two cases to consider: 
    $R <_{b_2}^\cw R_2b_2R_1$ or 
    $R_2b_2R_1 <_{b_2}^\cw R$.
    In the first case, $R <_{b_2}^\cw R_2$ because 
    $Rb_2<_{\scrI(b_2)} R_2b_2$ and
    $b_2R_2<_{\scrO(b_2)} b_2R_1 <_{\scrO(b_2)} b_2R$.
    This implies that $R$ and $R_2$ are incoherent, which is a contradiction. 
    In the second case, $R_1 <_{b_2}^\cw R$ because 
    $R_1b_2 <_{\scrI(b_2)} R_2b_2 <_{\scrI(b_2)} Rb_2$ and 
    $b_2R <_{\scrO(b_2)} b_2R_1$.
    This implies that $R$ and $R_1$ are incoherent, which is a contradiction.

    In Step~3, we find the coherence graph between the routes in $\widetilde{Z}\cup S$, which is a triangulated pentagon, with a point in the middle representing the set $S$. The dashed arrows represent clockwise rotations from the smaller to the larger maximal clique. The pentagon is also illustrated in more detail in~\Cref{fig.pentagon}.

    In Step~4, we deduce that $C_{\max(S)}= S\cup \{R_2,R_2b_2R_1\}$ just by looking at the direction of the arrows in the triangulated pentagon to determine the largest maximal clique. 
    
    It remains to show that~$C_1 \vee C_2$ exists and is equal to $C_{\max}(S)$. For this, let $M$ be a maximal clique satisfying $C_1\leq M$ and $C_2 \leq M$. We want to show that $C_{\max}(S)\leq M$. 
    Let~$R^M$ be a route in $M$, $R$ be a route in $C_{\max}(S)$, and $v \in R^M\cap R$. 
    By Theorem~\ref{thm.ccw}, it suffices to show that $R^M$ is coherent with~$R$ at $v$ or that $R <_v^{\cw} R^M$. Every route $R$ in $S\cup \{R_1,R_2\}$ satisfies this condition, because $C_1\leq M$ and $C_2\leq M$. Since $C_{\max}(S) = S\cup \{R_2,R_2b_2R_1\}$ then we only need to check $R_2b_2R_1$. 
    We proceed by contradiction. 
    Suppose that $R^M <_v^\cw R_2b_2R_1$. 
    We can assume that $v=b_2$, because otherwise $R^M <_v^\cw R_1$ or $R^M <_v^\cw R_2$, which is a contradiction. 
    Now, since $R^M <_{b_2}^\cw R_2b_2R_1$ and $b_2R_2<_{\scrO(b_2)}b_2R_1$ then $R^M <_{b_2}^\cw R_2$, which is also a contradiction. 
    This finishes the argument for Step~4.


\begin{figure}
    \centering    
    \begin{tikzpicture}

\begin{scope}[xshift=55, yshift=0]
    \draw[ultra thick,color=black!30] (2.5,0.6) .. controls (2.6, 0.6) and (3.6, 0.6) .. (4,0);
    \draw[ultra thick,color=blue] (4,0)--(5,0);
    \draw[ultra thick,color=black!30, dashed] (2,0.6)--(2.5,0.6);
    \draw[ultra thick,color=blue] (5,0) .. controls (5.4, -0.6) and (6.4, -0.6) .. (6.5,-0.6);
    \draw[ultra thick,color=blue, dashed] (6.5,-0.6)--(7,-0.6);
    
    \draw[ultra thick,color=blue] (2.5,-0.6) .. controls (2.6, -0.6) and (3.6, -0.5) .. (4,0);
    \draw[ultra thick,color=blue, dashed] (1,-0.6)--(1.5,-0.6);
    \draw[ultra thick,color=blue] (1.5,-0.6)--(2.5,-0.6);    
    \draw[ultra thick,color=black!30] (5,0) .. controls (5.4, 0.6) and (6.4, 0.6) .. (6.5,0.6);
    \draw[ultra thick,color=black!30, dashed] (6.5,0.6)--(7,0.6);

    \draw[ultra thick,color=red, dashed] (1,-1) --(1.5,-1);    
    \draw[ultra thick,color=red] (1.5,-1) .. controls (1.6, -1) and (1.9, -0.8) .. (2,-0.6);
    \draw[ultra thick,color=red] (2,-0.6)--(2.5,-0.6);  
    \draw[ultra thick,color=red] (2.5,-0.6) .. controls (2.6, -0.35) and (2.6, -0.35) .. (2.7,-0.2); 
    \draw[ultra thick,color=red, dashed] (2.7,-0.2) .. controls (2.85,0) and (3.3,0.1) .. (3.4,0.1);

    \node[circle, draw, inner sep=1.6pt, fill,label=below:\small{$b_2$}] (s) at (2,-0.6) {};
    \node[circle, draw, inner sep=1.6pt, fill,label=below:\small{$c_2$}] (s) at (2.5,-0.6) {};

    \node[label=below:\small{$b_1$}] (s) at (4,0) {};
    
    \node[circle, draw, inner sep=1.6pt, fill] (s) at (4,0) {};    
    \node[circle, draw, inner sep=1.6pt, fill, label=below:\small{$c_1$}] (c) at (5,0) {};
       
    \node[] (P) at (1.6,0.6) {$R_1^Q$};
    \node[] (Q) at (0.6,-0.5) {$R_1$};
    \node[] (Q) at (0.6,-1.1) {$R_2$};
    
    \node[] (P) at (7.4,0.6) {$R_1$};
    \node[] (Q) at (7.4,-0.6) {$R_1^Q$};

\end{scope}

\begin{scope}[scale=0.25, xshift=1600, yshift=-150]
    \draw[thick, color=black]  (4,0) -- (4,10) -- (-2,10) -- (-5,5) -- (-2,0) -- (4,0);		
    \draw[thick, color=black]  (4,0) -- (0,5) -- (4,10);
    \draw[thick, color=black]  (-2,0) -- (0,5) -- (-5,5);
    \draw[thick, color=black]  (0,5) -- (4,10); 
    \draw[thick, color=black]  (0,5) -- (-2,10); 
    
	\node[circle,fill,inner sep=2pt,color=black] at (4,0)  {};
	\node[circle,fill,inner sep=2pt,color=black]  at (4,10) {};
	\node[circle,fill,inner sep=2pt,color=black]  at (-2,10) {};
	\node[circle,fill,inner sep=2pt,color=black]  at (-5,5) {};
	\node[circle,fill,inner sep=2pt,color=black]  at (-2,0) {};
 	\node[circle,fill,inner sep=2pt,color=black]  at (0,5) {};  

	\node[circle,fill,inner sep=1.3pt,color=blue] at (0,1.7)  {};
 	\node[circle,fill,inner sep=1.3pt,color=blue] at (-2.5,3.3)  {};
 	\node[circle,fill,inner sep=1.3pt,color=blue] at (3,5)  {};
 	\node[circle,fill,inner sep=1.3pt,color=blue] at (-2.5,6.5)  {};  
 	\node[circle,fill,inner sep=1.3pt,color=blue] at (0,8.3)  {};

    \draw[thick,dashed,color=blue,-{Stealth[scale=1]}]  (0,1.7) -- (-2.5,3.3);
    \draw[thick,dashed,color=blue,-{Stealth[scale=1]}]  (0,1.7) -- (3,5);    
    \draw[thick,dashed,color=blue,-{Stealth[scale=1]}]  (3,5) -- (0,8.3);
    \draw[thick,dashed,color=blue,-{Stealth[scale=1]}]  (-2.5,3.3) -- (-2.5,6.5);    
    \draw[thick,dashed,color=blue,-{Stealth[scale=1]}]  (-2.5,6.5) -- (0,8.3);

    \node[] (a) at (4,-1.5) {$R_1^Q$};
    \node[] (a) at (-2.4,-1.5) {$R_2^Q$};
    \node[] (a) at (4,11.4) {$R_2$};
    \node[] (a) at (-6.4,5) {$R_1$};  
    \node[] (a) at (1.07,5) {$S$};
    \node[] (a) at (-1.5,11.4) {$R_2b_2R_1$}; 
\end{scope}

\begin{scope}[xshift=0, yshift=-130]
    \draw[ultra thick,color=blue] (2.5,0.6) .. controls (2.6, 0.6) and (3.6, 0.6) .. (4,0);
    \draw[ultra thick,color=red] (4,0)--(5,0);
    \draw[ultra thick,color=blue, dashed] (2,0.6)--(2.5,0.6);
    \draw[ultra thick,color=red] (5,0) .. controls (5.4, 0.6) and (6.6, 0.6) .. (7,0);
    \draw[ultra thick,color=black!30] (8,0) .. controls (8.4, 0.6) and (9.4, 0.6) .. (9.5,0.6);
    \draw[ultra thick,color=black!30, dashed] (9.5,0.6)--(10,0.6);

    \draw[ultra thick,color=red] (2.5,-0.6) .. controls (2.6, -0.6) and (3.6, -0.5) .. (4,0);
    \draw[ultra thick,color=red, dashed] (2,-0.6)--(2.5,-0.6); 
    \draw[ultra thick,color=blue] (5,0) .. controls (5.4, -0.6) and (6.6, -0.6) .. (7,0);
    \draw[ultra thick,color=red] (7,0)--(8,0);
    \draw[ultra thick,color=red] (8,0) .. controls (8.4, -0.6) and (9.4, -0.6) .. (9.5,-0.6);
    \draw[ultra thick,color=red, dashed] (9.5,-0.6)--(10,-0.6);

    \node[circle, draw, inner sep=1.6pt, fill,label=below:\small{$b_2$}] (s) at (4,0) {};
    \node[circle, draw, inner sep=1.6pt, fill,label=below:\small{$c_2$}] (s) at (5,0) {};

    \node[circle, draw, inner sep=1.6pt, fill,label=below:\small{$b_1$}] (s) at (7,0) {};
    \node[label=below:\small{$=$}] (s) at (7.3,-0.1) {};    
    \node[label=below:\small{$d_2$}] (s) at (7.7,0.05) {};
    \node[circle, draw, inner sep=1.6pt, fill, label=above:\small{$c_1$}] (c) at (8,0) {};
       
    \node[] (P) at (2.6,1) {$R_1^Q, R_1, R_2^Q$};
    \node[] (Q) at (2.6,-1) {$R_2$};

    \node[] (P) at (6,1) {$R_1^Q, R_2$};
    \node[] (Q) at (6,-1) {$R_2^Q, R_1$};
    
    \node[] (P) at (9.4,1) {$R_1$};
    \node[] (Q) at (9.4,-1) {$R_1^Q, R_2^Q, R_2$};

\end{scope}

\begin{scope}[scale=0.25, xshift=1600, yshift=-660]
    \draw[thick, color=black]  (2,0) -- (5,5) -- (2,10) -- (-2,10) -- (-5,5) -- (-2,0) -- (2,0);		
    \draw[thick, color=black]  (2,0) -- (0,5) -- (5,5);
    \draw[thick, color=black]  (-2,0) -- (0,5) -- (-5,5);
    \draw[thick, color=black]  (-2,10) -- (0,5) -- (2,10);    
    
	\node[circle,fill,inner sep=2pt,color=black] at (2,0)  {};
	\node[circle,fill,inner sep=2pt,color=black]  at (5,5) {};
	\node[circle,fill,inner sep=2pt,color=black]  at (2,10) {};
        \node[circle,fill,inner sep=2pt,color=black]  at (-2,10) {}; 
	\node[circle,fill,inner sep=2pt,color=black]  at (-5,5) {};
	\node[circle,fill,inner sep=2pt,color=black]  at (-2,0) {};
 	\node[circle,fill,inner sep=2pt,color=black]  at (0,5) {};  

	\node[circle,fill,inner sep=1.3pt,color=blue] at (0,1.7)  {};
 	\node[circle,fill,inner sep=1.3pt,color=blue] at (-2.5,3.3)  {};
 	\node[circle,fill,inner sep=1.3pt,color=blue] at (2.5,3.3)  {};
 	\node[circle,fill,inner sep=1.3pt,color=blue] at (-2.5,6.6)  {};  
 	\node[circle,fill,inner sep=1.3pt,color=blue] at (2.5,6.6)  {};
 	\node[circle,fill,inner sep=1.3pt,color=blue] at (0,8.3)  {};  

    \draw[thick,dashed,color=blue,-{Stealth[scale=1]}]  (0,1.7) -- (-2.5,3.3);
    \draw[thick,dashed,color=blue,-{Stealth[scale=1]}]  (0,1.7) -- (2.5,3.3);    
    \draw[thick,dashed,color=blue,-{Stealth[scale=1]}]  (2.5,3.3) -- (2.5,6.6);
    \draw[thick,dashed,color=blue,-{Stealth[scale=1]}]  (-2.5,3.3) -- (-2.5,6.6);    
    \draw[thick,dashed,color=blue,-{Stealth[scale=1]}]  (-2.5,6.6) -- (0,8.3);
    \draw[thick,dashed,color=blue,-{Stealth[scale=1]}]  (2.5,6.6) -- (0,8.3);

    \node[] (a) at (2.4,-1.5) {$R_1^Q$};
    \node[] (a) at (-2.4,-1.5) {$R_2^Q$};
    \node[] (a) at (6.4,5) {$R_2$};
    \node[] (a) at (-6.4,5) {$R_1$};  
    \node[] (a) at (1,5.8) {$S$};
    \node[] (a) at (-3,11.4) {$R_2b_2R_1$}; 
    \node[] (a) at (3,11.4) {$R_2b_1R_1$};  
\end{scope}
\end{tikzpicture}   
    \caption{The pentagonal and hexagonal sub-cases and their respective coherence graphs.
}
\label{fig.polygons_and_coherence_graphs}
\end{figure} 

    \textbf{Cases II.(i), III.(i), IV.(i) and IV.(ii).} Pentagons.

    The cases II.(i), III.(i), IV.(i), and IV.(ii) are all argued identically to I.(i) above.

    \textbf{Case III.(ii).} Impossible.

    We consider the Case III.(ii) next.
    Observe that in this case the routes $\Top(R_2^Q,R_2)$ and~$R_1$ are incoherent at $d_2$.
    Since $\Top(R_2^Q,R_2)\in Q$, and the only route in $Q$ incoherent with~$R_1$ in~$Q$ is $R_1^Q$, we have that $\Top(R_2^Q,R_2) = R_1^Q$. 
    However, this would imply that $d_2 = b_1$, and so Case III.(ii) is impossible.

    \textbf{Case II.(ii).} Hexagon.
    
    \begin{figure}[htb]
    \centering
    \scalebox{0.8}{
    \begin{tikzpicture}
\begin{scope}[xshift=0, yshift=-10]
\begin{scope}[xshift=-10, yshift=-45, scale=1.6]
    \draw[thick, fill=ForestGreen!5] (-1,0) -- (5,0) -- (5,2.6) -- (-1,2.6) -- cycle;
\end{scope}
\begin{scope}[xshift=-85]
    \draw[ultra thick,color=blue] (2.5,0.6) .. controls (2.6, 0.6) and (3.6, 0.6) .. (4,0);
    \draw[ultra thick,color=red] (4,0)--(5,0);
    \draw[ultra thick,color=blue, dashed] (2,0.6)--(2.5,0.6);
    \draw[ultra thick,color=red] (5,0) .. controls (5.4, 0.6) and (6.6, 0.6) .. (7,0);
    \draw[ultra thick,color=black!30] (8,0) .. controls (8.4, 0.6) and (9.4, 0.6) .. (9.5,0.6);
    \draw[ultra thick,color=black!30, dashed] (9.5,0.6)--(10,0.6);

    \draw[ultra thick,color=red] (2.5,-0.6) .. controls (2.6, -0.6) and (3.6, -0.5) .. (4,0);
    \draw[ultra thick,color=red, dashed] (2,-0.6)--(2.5,-0.6); 
    \draw[ultra thick,color=blue] (5,0) .. controls (5.4, -0.6) and (6.6, -0.6) .. (7,0);
    \draw[ultra thick,color=red] (7,0)--(8,0);
    \draw[ultra thick,color=red] (8,0) .. controls (8.4, -0.6) and (9.4, -0.6) .. (9.5,-0.6);
    \draw[ultra thick,color=red, dashed] (9.5,-0.6)--(10,-0.6);

    \node[circle, draw, inner sep=1.6pt, fill,label=below:\small{$b_2$}] (s) at (4,0) {};
    \node[circle, draw, inner sep=1.6pt, fill,label=below:\small{$c_2$}] (s) at (5,0) {};

    \node[circle, draw, inner sep=1.6pt, fill,label=below:\small{$b_1$}] (s) at (7,0) {};
    \node[label=below:\small{$=$}] (s) at (7.3,-0.1) {};    
    \node[label=below:\small{$d_2$}] (s) at (7.7,0.05) {};
    \node[circle, draw, inner sep=1.6pt, fill, label=above:\small{$c_1$}] (c) at (8,0) {};
       
    \node[] (P) at (2.5,1) {$R_1^Q, R_1, R_2^Q$};
    \node[] (Q) at (2.6,-1) {$R_2$};

    \node[] (P) at (6,1) {$R_1^Q, R_2$};
    \node[] (Q) at (6,-1) {$R_2^Q, R_1$};
    
    \node[] (P) at (9.4,1) {$R_1$};
    \node[] (Q) at (9.4,-1) {$R_1^Q, R_2^Q, R_2$};

    \node[] (P) at (6,2) {$\widetilde{G}$};
    
\end{scope}
\end{scope}

\begin{scope}[xshift=0, yshift=-121]
\begin{scope}[xshift=-10, yshift=-195, scale=1.6]
    \draw[thick, fill=ForestGreen!05] (-1,1.8) -- (5,1.8) -- (5,5.5) -- (-1,5.5) -- cycle;
\end{scope}
\begin{scope}[xshift=-75]
    \draw[ultra thick,color=black!60] (2.5,0.6) .. controls (2.6, 0.6) and (3.6, 0.6) .. (4,0);
    \draw[ultra thick,color=black!60] (4,0)--(5,0);
    \draw[ultra thick,color=black!60, dashed] (2,0.6)--(2.5,0.6);
    \draw[ultra thick,color=black!60] (5,0) .. controls (5.4, 0.6) and (6.6, 0.6) .. (7,0);
    \draw[ultra thick,color=black!60] (8,0) .. controls (8.4, 0.6) and (9.4, 0.6) .. (9.5,0.6);
    \draw[ultra thick,color=black!60, dashed] (9.5,0.6)--(10,0.6);

    \draw[ultra thick,color=black!20] (2.5,-0.6) .. controls (2.6, -0.6) and (3.6, -0.5) .. (4,0);
    \draw[ultra thick,color=black!20, dashed] (2,-0.6)--(2.5,-0.6); 
    \draw[ultra thick,color=black!20] (5,0) .. controls (5.4, -0.6) and (6.6, -0.6) .. (7,0);
    \draw[ultra thick,color=black!60] (7,0)--(8,0);
    \draw[ultra thick,color=black!20] (8,0) .. controls (8.4, -0.6) and (9.4, -0.6) .. (9.5,-0.6);
    \draw[ultra thick,color=black!20, dashed] (9.5,-0.6)--(10,-0.6);

    \node[circle, draw, inner sep=1.6pt, fill] (s) at (4,0) {};
    \node[circle, draw, inner sep=1.6pt, fill] (s) at (5,0) {};

    \node[circle, draw, inner sep=1.6pt, fill] (s) at (7,0) {};
    \node[circle, draw, inner sep=1.6pt, fill] (c) at (8,0) {};
       
    \node[] (P) at (1.5,0.6) {$\mathrm{Top}_1$};
   
    \node[] (P) at (6,1.5) {$\widetilde{S}$};
    \node[] (P) at (1.3,1.5) {\textbf{(1)}};     
\end{scope}
\begin{scope}[xshift=-75, yshift=-70]
    \draw[ultra thick,color=black!20] (2.5,0.6) .. controls (2.6, 0.6) and (3.6, 0.6) .. (4,0);
    \draw[ultra thick,color=black!60] (4,0)--(5,0);
    \draw[ultra thick,color=black!20, dashed] (2,0.6)--(2.5,0.6);
    \draw[ultra thick,color=black!20] (5,0) .. controls (5.4, 0.6) and (6.6, 0.6) .. (7,0);
    \draw[ultra thick,color=black!20] (8,0) .. controls (8.4, 0.6) and (9.4, 0.6) .. (9.5,0.6);
    \draw[ultra thick,color=black!20, dashed] (9.5,0.6)--(10,0.6);

    \draw[ultra thick,color=black!60] (2.5,-0.6) .. controls (2.6, -0.6) and (3.6, -0.5) .. (4,0);
    \draw[ultra thick,color=black!60, dashed] (2,-0.6)--(2.5,-0.6); 
    \draw[ultra thick,color=black!60] (5,0) .. controls (5.4, -0.6) and (6.6, -0.6) .. (7,0);
    \draw[ultra thick,color=black!60] (7,0)--(8,0);
    \draw[ultra thick,color=black!60] (8,0) .. controls (8.4, -0.6) and (9.4, -0.6) .. (9.5,-0.6);
    \draw[ultra thick,color=black!60, dashed] (9.5,-0.6)--(10,-0.6);

    \node[circle, draw, inner sep=1.6pt, fill] (s) at (4,0) {};
    \node[circle, draw, inner sep=1.6pt, fill] (s) at (5,0) {};

    \node[circle, draw, inner sep=1.6pt, fill] (s) at (7,0) {};
    \node[circle, draw, inner sep=1.6pt, fill] (c) at (8,0) {};
       
    \node[] (P) at (1.5,-0.6) {$\mathrm{Bot}_2$};
\end{scope}
\end{scope}


\begin{scope}[xshift=290, yshift=8.5]
\begin{scope}[xshift=-15, yshift=-332, scale=1.6]
    \draw[thick, fill=ForestGreen!05] (-1,-1.17) -- (5,-1.17) -- (5,8.5) -- (-1,8.5) -- cycle;
\end{scope}
\begin{scope}[xshift=-80]
    \draw[ultra thick,color=black!60] (2.5,0.6) .. controls (2.6, 0.6) and (3.6, 0.6) .. (4,0);
    \draw[ultra thick,color=black!60] (4,0)--(5,0);
    \draw[ultra thick,color=black!60, dashed] (2,0.6)--(2.5,0.6);
    \draw[ultra thick,color=black!60] (5,0) .. controls (5.4, 0.6) and (6.6, 0.6) .. (7,0);
    \draw[ultra thick,color=black!20] (8,0) .. controls (8.4, 0.6) and (9.4, 0.6) .. (9.5,0.6);
    \draw[ultra thick,color=black!20, dashed] (9.5,0.6)--(10,0.6);

    \draw[ultra thick,color=black!20] (2.5,-0.6) .. controls (2.6, -0.6) and (3.6, -0.5) .. (4,0);
    \draw[ultra thick,color=black!20, dashed] (2,-0.6)--(2.5,-0.6); 
    \draw[ultra thick,color=black!20] (5,0) .. controls (5.4, -0.6) and (6.6, -0.6) .. (7,0);
    \draw[ultra thick,color=black!60] (7,0)--(8,0);
    \draw[ultra thick,color=black!60] (8,0) .. controls (8.4, -0.6) and (9.4, -0.6) .. (9.5,-0.6);
    \draw[ultra thick,color=black!60, dashed] (9.5,-0.6)--(10,-0.6);

    \node[circle, draw, inner sep=1.6pt, fill] (s) at (4,0) {};
    \node[circle, draw, inner sep=1.6pt, fill] (s) at (5,0) {};

    \node[circle, draw, inner sep=1.6pt, fill] (s) at (7,0) {};
    \node[circle, draw, inner sep=1.6pt, fill] (c) at (8,0) {};
       
    \node[] (P) at (1.5,0.6) {$R_1^Q$};

    \node[] (P) at (6,1.5) {$\widetilde{Z}$};
    \node[] (P) at (1.3,1.5) {\textbf{(2)}};    
\end{scope}
\begin{scope}[xshift=-80, yshift=-70]
    \draw[ultra thick,color=black!60] (2.5,0.6) .. controls (2.6, 0.6) and (3.6, 0.6) .. (4,0);
    \draw[ultra thick,color=black!60] (4,0)--(5,0);
    \draw[ultra thick,color=black!60, dashed] (2,0.6)--(2.5,0.6);
    \draw[ultra thick,color=black!20] (5,0) .. controls (5.4, 0.6) and (6.6, 0.6) .. (7,0);
    \draw[ultra thick,color=black!60] (8,0) .. controls (8.4, 0.6) and (9.4, 0.6) .. (9.5,0.6);
    \draw[ultra thick,color=black!60, dashed] (9.5,0.6)--(10,0.6);

    \draw[ultra thick,color=black!20] (2.5,-0.6) .. controls (2.6, -0.6) and (3.6, -0.5) .. (4,0);
    \draw[ultra thick,color=black!20, dashed] (2,-0.6)--(2.5,-0.6); 
    \draw[ultra thick,color=black!60] (5,0) .. controls (5.4, -0.6) and (6.6, -0.6) .. (7,0);
    \draw[ultra thick,color=black!60] (7,0)--(8,0);
    \draw[ultra thick,color=black!20] (8,0) .. controls (8.4, -0.6) and (9.4, -0.6) .. (9.5,-0.6);
    \draw[ultra thick,color=black!20, dashed] (9.5,-0.6)--(10,-0.6);

    \node[circle, draw, inner sep=1.6pt, fill] (s) at (4,0) {};
    \node[circle, draw, inner sep=1.6pt, fill] (s) at (5,0) {};

    \node[circle, draw, inner sep=1.6pt, fill] (s) at (7,0) {};
    \node[circle, draw, inner sep=1.6pt, fill] (c) at (8,0) {};
       
    \node[] (P) at (1.5,0.6) {$R_1$};
   
\end{scope}
\begin{scope}[xshift=-80, yshift=-140]
    \draw[ultra thick,color=black!60] (2.5,0.6) .. controls (2.6, 0.6) and (3.6, 0.6) .. (4,0);
    \draw[ultra thick,color=black!60] (4,0)--(5,0);
    \draw[ultra thick,color=black!60, dashed] (2,0.6)--(2.5,0.6);
    \draw[ultra thick,color=black!20] (5,0) .. controls (5.4, 0.6) and (6.6, 0.6) .. (7,0);
    \draw[ultra thick,color=black!20] (8,0) .. controls (8.4, 0.6) and (9.4, 0.6) .. (9.5,0.6);
    \draw[ultra thick,color=black!20, dashed] (9.5,0.6)--(10,0.6);

    \draw[ultra thick,color=black!20] (2.5,-0.6) .. controls (2.6, -0.6) and (3.6, -0.5) .. (4,0);
    \draw[ultra thick,color=black!20, dashed] (2,-0.6)--(2.5,-0.6); 
    \draw[ultra thick,color=black!60] (5,0) .. controls (5.4, -0.6) and (6.6, -0.6) .. (7,0);
    \draw[ultra thick,color=black!60] (7,0)--(8,0);
    \draw[ultra thick,color=black!60] (8,0) .. controls (8.4, -0.6) and (9.4, -0.6) .. (9.5,-0.6);
    \draw[ultra thick,color=black!60, dashed] (9.5,-0.6)--(10,-0.6);

    \node[circle, draw, inner sep=1.6pt, fill] (s) at (4,0) {};
    \node[circle, draw, inner sep=1.6pt, fill] (s) at (5,0) {};

    \node[circle, draw, inner sep=1.6pt, fill] (s) at (7,0) {};
    \node[circle, draw, inner sep=1.6pt, fill] (c) at (8,0) {};
       
    \node[] (P) at (1.5,0.6) {$R_2^Q$};
   
\end{scope}
\begin{scope}[xshift=-80, yshift=-210]
    \draw[ultra thick,color=black!20] (2.5,0.6) .. controls (2.6, 0.6) and (3.6, 0.6) .. (4,0);
    \draw[ultra thick,color=black!60] (4,0)--(5,0);
    \draw[ultra thick,color=black!20, dashed] (2,0.6)--(2.5,0.6);
    \draw[ultra thick,color=black!60] (5,0) .. controls (5.4, 0.6) and (6.6, 0.6) .. (7,0);
    \draw[ultra thick,color=black!60] (8,0) .. controls (8.4, 0.6) and (9.4, 0.6) .. (9.5,0.6);
    \draw[ultra thick,color=black!60, dashed] (9.5,0.6)--(10,0.6);

    \draw[ultra thick,color=black!60] (2.5,-0.6) .. controls (2.6, -0.6) and (3.6, -0.5) .. (4,0);
    \draw[ultra thick,color=black!60, dashed] (2,-0.6)--(2.5,-0.6); 
    \draw[ultra thick,color=black!20] (5,0) .. controls (5.4, -0.6) and (6.6, -0.6) .. (7,0);
    \draw[ultra thick,color=black!60] (7,0)--(8,0);
    \draw[ultra thick,color=black!20] (8,0) .. controls (8.4, -0.6) and (9.4, -0.6) .. (9.5,-0.6);
    \draw[ultra thick,color=black!20, dashed] (9.5,-0.6)--(10,-0.6);

    \node[circle, draw, inner sep=1.6pt, fill] (s) at (4,0) {};
    \node[circle, draw, inner sep=1.6pt, fill] (s) at (5,0) {};

    \node[circle, draw, inner sep=1.6pt, fill] (s) at (7,0) {};
    \node[circle, draw, inner sep=1.6pt, fill] (c) at (8,0) {};
       
    \node[] (P) at (1.6,-0.3) {$R_2b_1R_1$};

\end{scope}
\begin{scope}[xshift=-80, yshift=-280]
    \draw[ultra thick,color=black!20] (2.5,0.6) .. controls (2.6, 0.6) and (3.6, 0.6) .. (4,0);
    \draw[ultra thick,color=black!60] (4,0)--(5,0);
    \draw[ultra thick,color=black!20, dashed] (2,0.6)--(2.5,0.6);
    \draw[ultra thick,color=black!60] (5,0) .. controls (5.4, 0.6) and (6.6, 0.6) .. (7,0);
    \draw[ultra thick,color=black!20] (8,0) .. controls (8.4, 0.6) and (9.4, 0.6) .. (9.5,0.6);
    \draw[ultra thick,color=black!20, dashed] (9.5,0.6)--(10,0.6);

    \draw[ultra thick,color=black!60] (2.5,-0.6) .. controls (2.6, -0.6) and (3.6, -0.5) .. (4,0);
    \draw[ultra thick,color=black!60, dashed] (2,-0.6)--(2.5,-0.6); 
    \draw[ultra thick,color=black!20] (5,0) .. controls (5.4, -0.6) and (6.6, -0.6) .. (7,0);
    \draw[ultra thick,color=black!60] (7,0)--(8,0);
    \draw[ultra thick,color=black!60] (8,0) .. controls (8.4, -0.6) and (9.4, -0.6) .. (9.5,-0.6);
    \draw[ultra thick,color=black!60, dashed] (9.5,-0.6)--(10,-0.6);

    \node[circle, draw, inner sep=1.6pt, fill] (s) at (4,0) {};
    \node[circle, draw, inner sep=1.6pt, fill] (s) at (5,0) {};

    \node[circle, draw, inner sep=1.6pt, fill] (s) at (7,0) {};
    \node[circle, draw, inner sep=1.6pt, fill] (c) at (8,0) {};
       
    \node[] (P) at (1.5,-0.6) {$R_2$};
\end{scope}
\begin{scope}[xshift=-80, yshift=-350]
    \draw[ultra thick,color=black!20] (2.5,0.6) .. controls (2.6, 0.6) and (3.6, 0.6) .. (4,0);
    \draw[ultra thick,color=black!60] (4,0)--(5,0);
    \draw[ultra thick,color=black!20, dashed] (2,0.6)--(2.5,0.6);
    \draw[ultra thick,color=black!20] (5,0) .. controls (5.4, 0.6) and (6.6, 0.6) .. (7,0);
    \draw[ultra thick,color=black!60] (8,0) .. controls (8.4, 0.6) and (9.4, 0.6) .. (9.5,0.6);
    \draw[ultra thick,color=black!60, dashed] (9.5,0.6)--(10,0.6);

    \draw[ultra thick,color=black!60] (2.5,-0.6) .. controls (2.6, -0.6) and (3.6, -0.5) .. (4,0);
    \draw[ultra thick,color=black!60, dashed] (2,-0.6)--(2.5,-0.6); 
    \draw[ultra thick,color=black!60] (5,0) .. controls (5.4, -0.6) and (6.6, -0.6) .. (7,0);
    \draw[ultra thick,color=black!60] (7,0)--(8,0);
    \draw[ultra thick,color=black!20] (8,0) .. controls (8.4, -0.6) and (9.4, -0.6) .. (9.5,-0.6);
    \draw[ultra thick,color=black!20, dashed] (9.5,-0.6)--(10,-0.6);

    \node[circle, draw, inner sep=1.6pt, fill] (s) at (4,0) {};
    \node[circle, draw, inner sep=1.6pt, fill] (s) at (5,0) {};

    \node[circle, draw, inner sep=1.6pt, fill] (s) at (7,0) {};
    \node[circle, draw, inner sep=1.6pt, fill] (c) at (8,0) {};
       
    \node[] (P) at (1.6,-0.3) {$R_2b_2R_1$};
\end{scope}
\end{scope}

\begin{scope}
\begin{scope}[xshift=-10, yshift=-380, scale=1.6]
    \draw[thick, fill=ForestGreen!05] (-1,0.07) -- (5,0.07) -- (5,3) -- (-1,3) -- cycle;
    \node[] (P) at (-0.6,2.7) {\textbf{(3)}};    
\end{scope}
\begin{scope}[scale=0.3, xshift=300, yshift=-1180]
    \draw[thick, color=black]  (2,0) -- (5,5) -- (2,10) -- (-2,10) -- (-5,5) -- (-2,0) -- (2,0);		
    \draw[thick, color=black]  (2,0) -- (0,5) -- (5,5);
    \draw[thick, color=black]  (-2,0) -- (0,5) -- (-5,5);
    \draw[thick, color=black]  (-2,10) -- (0,5) -- (2,10);    
    
	\node[circle,fill,inner sep=2pt,color=black] at (2,0)  {};
	\node[circle,fill,inner sep=2pt,color=black]  at (5,5) {};
	\node[circle,fill,inner sep=2pt,color=black]  at (2,10) {};
        \node[circle,fill,inner sep=2pt,color=black]  at (-2,10) {}; 
	\node[circle,fill,inner sep=2pt,color=black]  at (-5,5) {};
	\node[circle,fill,inner sep=2pt,color=black]  at (-2,0) {};
 	\node[circle,fill,inner sep=2pt,color=black]  at (0,5) {};  

	\node[circle,fill,inner sep=1.3pt,color=blue] at (0,1.7)  {};
 	\node[circle,fill,inner sep=1.3pt,color=blue] at (-2.5,3.3)  {};
 	\node[circle,fill,inner sep=1.3pt,color=blue] at (2.5,3.3)  {};
 	\node[circle,fill,inner sep=1.3pt,color=blue] at (-2.5,6.6)  {};  
 	\node[circle,fill,inner sep=1.3pt,color=blue] at (2.5,6.6)  {};
 	\node[circle,fill,inner sep=1.3pt,color=blue] at (0,8.3)  {};  

    \draw[thick,dashed,color=blue,-{Stealth[scale=1]}]  (0,1.7) -- (-2.5,3.3);
    \draw[thick,dashed,color=blue,-{Stealth[scale=1]}]  (0,1.7) -- (2.5,3.3);    
    \draw[thick,dashed,color=blue,-{Stealth[scale=1]}]  (2.5,3.3) -- (2.5,6.6);
    \draw[thick,dashed,color=blue,-{Stealth[scale=1]}]  (-2.5,3.3) -- (-2.5,6.6);    
    \draw[thick,dashed,color=blue,-{Stealth[scale=1]}]  (-2.5,6.6) -- (0,8.3);
    \draw[thick,dashed,color=blue,-{Stealth[scale=1]}]  (2.5,6.6) -- (0,8.3);

    \node[] (a) at (2.4,-1.5) {$R_1^Q$};
    \node[] (a) at (-2.4,-1.5) {$R_2^Q$};
    \node[] (a) at (6.4,5) {$R_2$};
    \node[] (a) at (-6.4,5) {$R_1$};  
    \node[] (a) at (1,5.8) {$S$};
    \node[] (a) at (-3,11.4) {$R_2b_2R_1$}; 
    \node[] (a) at (3,11.4) {$R_2b_1R_1$};  
\end{scope}
\end{scope}

\end{tikzpicture}
    }
    \caption{Steps 1-3 in the strategy applied to Case (2.II.ii).}
    \label{fig_photo_caseII_ii}
    \end{figure}

    It remains to consider Case II.(ii).
    Similarly as in the previous case, since $\Top(R_2^Q,R_2)\in Q$ is incoherent with $R_1$, we have that $\Top(R_2^Q,R_2) = R_1^Q$.
    This case is depicted in Figure~\ref{fig.polygons_and_coherence_graphs}~(bottom).
        
    The outcomes of the steps throughout our strategy are illustrated in~\Cref{fig_photo_caseII_ii}.
    Note, in particular, that the graph $\widetilde{G}$ looks a bit different to the graph in~\Cref{fig.join_proof_cases}, because of the forced requirement $\Top(R_2^Q,R_2) = R_1^Q$.
    
    In Step 1 of our strategy, we find the set $\widetilde{S}=\{\Top_1,\Bot_2\}$; note that $\Top_2=\Top_2(R_2^Q,R_2)=R_1^Q$ and $\Bot_1=\Bot(R_1^Q,R_1)=R_2^Q$ do not appear in this set because $R_1^Q,R_2^Q\notin S$.
    In Step~2, we find the set $\widetilde{Z}=\{R_1^Q, R_1, R_2^Q, R_2b_1R_1, R_2, R_2b_2R_1\}$ of routes in $\widetilde{G}$ that are not in $\widetilde S$ and are coherent with $\widetilde{S}$.
    
    In Step~2a, we need to argue that the routes in $\widetilde{Z}$ are coherent with $S$. Since this holds for $R_1^Q, R_2^Q, R_1, R_2$, we only need to check that $R_2b_1R_1$ and $R_2b_2R_1$ are coherent with $S$. 
    Suppose there exists a route $R\in S$ that is incoherent with either $R_2b_1R_1$ or $R_2b_2R_1$. Using similar arguments as in Step~2a for the pentagonal Case I.(i), one can deduce that $R$ is incoherent with $R_1$ or $R_2$, which is a contradiction. 
    
    In Step~3, we find the coherence graph between the routes in $\widetilde{Z}\cup S$, which is a triangulated hexagon, with a point in the middle representing the set $S$. The dashed arrows represent clockwise rotations from the smaller to the larger maximal clique. The hexagon is also illustrated in more detail in~\Cref{fig.hexagon}.  

    In Step~4, we deduce that $C_{\max}(S)= S\cup \{R_2b_2R_1,R_2b_1R_1\}$ just by looking at the direction of the arrows in the triangulated hexagon to determine the largest maximal clique. 
    
    It remains to show that~$C_1 \vee C_2$ exists and is equal to $C_{\max}(S)$. For this, let $M$ be a maximal clique satisfying $C_1\leq M$ and $C_2 \leq M$. We want to show that $C_{\max}(S)\leq M$. 
    As before, let~$R^M$ be a route in $M$, $R$ be a route in $C_{\max}(S)$, and $v \in R^M\cap R$. 
    By Theorem~\ref{thm.ccw}, it suffices to show that $R^M$ is coherent with~$R$ at $v$ or that $R <_v^{\cw} R^M$. Every route $R$ in $S\cup \{R_1,R_2\}$ satisfies this condition, because $C_1\leq M$ and $C_2\leq M$. Since $C_{\max}(S) = S\cup \{R_2b_2R_1,R_2b_1R_1\}$ then we only need to check $R_2b_2R_1$ and $R_2b_1R_1$. 
    Using similar arguments as in the pentagonal Case I.(i), one can deduce that if this would not be the case then either $R^M <_v^\cw R_1$ or $R^M <_v^\cw R_1$, which is a contradiction.  
    This finishes the argument for Step~4.
\end{proof}

\begin{figure}
    \centering
    \begin{tikzpicture}

\begin{scope}[xshift=0, yshift=0, scale=0.8]
    \draw[ultra thick,color=blue] (2.5,0.6) .. controls (2.6, 0.6) and (3.6, 0.6) .. (4,0);
    \draw[ultra thick,color=blue] (4,0.035)--(5,0.035);
    \draw[ultra thick,color=red] (4,-0.035)--(5,-0.035);
    \draw[ultra thick,color=blue, dashed] (2,0.6)--(2.5,0.6);
    \draw[ultra thick,color=blue] (5,0.035) .. controls (5.4, -0.565) and (6.4, -0.565) .. (6.5,-0.565);
    \draw[ultra thick,color=blue, dashed] (6.5,-0.565)--(7.25,-0.565);
    \draw[ultra thick,color=red] (4.95,-0.035) .. controls (5.4, -0.635) and (6.4, -0.635) .. (6.5,-0.635);
    \draw[ultra thick,color=red, dashed] (6.5,-0.635)--(7.25,-0.635);

    \draw[ultra thick,color=red] (2.5,-0.6) .. controls (2.6, -0.6) and (3.6, -0.6) .. (4,0);
    \draw[ultra thick,color=red, dashed] (1,-0.6)--(1.5,-0.6);
    \draw[ultra thick,color=red] (1.5,-0.6)--(2.5,-0.6);    
    \draw[ultra thick,color=black!30] (5,0) .. controls (5.4, 0.6) and (6.4, 0.6) .. (6.5,0.6);
    \draw[ultra thick,color=black!30, dashed] (6.5,0.6)--(7.25,0.6);

    \draw[ultra thick,color=black!30, dashed] (1,-1) --(1.5,-1);    
    \draw[ultra thick,color=black!30] (1.5,-1) .. controls (1.6, -1) and (1.9, -0.8) .. (2,-0.6);
    \draw[ultra thick,color=black!30] (2.5,-0.6) .. controls (2.6, -0.35) and (2.6, -0.35) .. (2.7,-0.2); 
    \draw[ultra thick,color=black!30, dashed] (2.7,-0.2) .. controls (2.85,0) and (3.3,0.1) .. (3.4,0.1);

    \node[circle, draw, inner sep=1.6pt, fill,label=below:\tiny{$b_2$}] (s) at (2,-0.6) {};
    \node[circle, draw, inner sep=1.6pt, fill,label=below:\tiny{$c_2$}] (s) at (2.5,-0.6) {};

    \node[label=below:\tiny{$b_1$}] (s) at (4,0) {};
    
    \node[circle, draw, inner sep=1.6pt, fill] (s) at (4,0) {};    
    \node[circle, draw, inner sep=1.6pt, fill, label=below:\tiny{$c_1$}] (c) at (5,0) {};
       
    \node[] (P) at (1.6,0.6) {\tiny \textcolor{blue}{$R_1^Q$}};
    \node[] (Q) at (0.6,-0.5) {\tiny \textcolor{red}{$R_2^Q$}};
    
\end{scope}

\begin{scope}[xshift=-120, yshift=60, scale=0.8]

    \draw[ultra thick,color=black!30] (2.5,0.6) .. controls (2.6, 0.6) and (3.6, 0.6) .. (4,0);
    \draw[ultra thick,color=blue] (4,0.035)--(5,0.035);
    \draw[ultra thick,color=red] (4,-0.035)--(5,-0.035);
    \draw[ultra thick,color=black!30, dashed] (2,0.6)--(2.5,0.6);
    \draw[ultra thick,color=red] (5,0.035) .. controls (5.4, -0.565) and (6.4, -0.565) .. (6.5,-0.565);
    \draw[ultra thick,color=red, dashed] (6.5,-0.565)--(7.25,-0.565);

    \draw[ultra thick,color=blue] (2.5,-0.565) .. controls (2.6, -0.565) and (3.6, -0.565) .. (4,0.035);
    \draw[ultra thick,color=red] (2.5,-0.635) .. controls (2.6, -0.635) and (3.6, -0.635) .. (4.05,-0.035);
    \draw[ultra thick,color=blue, dashed] (1,-0.565)--(1.5,-0.565);
    \draw[ultra thick,color=blue] (1.5,-0.565)--(2.5,-0.565);
    \draw[ultra thick,color=red, dashed] (1,-0.635)--(1.5,-0.635);
    \draw[ultra thick,color=red] (1.5,-0.635)--(2.5,-0.635);    
    \draw[ultra thick,color=blue] (5,0) .. controls (5.4, 0.6) and (6.4, 0.6) .. (6.5,0.6);
    \draw[ultra thick,color=blue, dashed] (6.5,0.6)--(7.25,0.6);

    \draw[ultra thick,color=black!30, dashed] (1,-1) --(1.5,-1);    
    \draw[ultra thick,color=black!30] (1.5,-1) .. controls (1.6, -1) and (1.9, -0.8) .. (2,-0.6);
    \draw[ultra thick,color=red] (2,-0.635)--(2.5,-0.635);  
    \draw[ultra thick,color=black!30] (2.5,-0.6) .. controls (2.6, -0.35) and (2.6, -0.35) .. (2.7,-0.2); 
    \draw[ultra thick,color=black!30, dashed] (2.7,-0.2) .. controls (2.85,0) and (3.3,0.1) .. (3.4,0.1);

    \node[circle, draw, inner sep=1.6pt, fill,label=below:\tiny{$b_2$}] (s) at (2,-0.6) {};
    \node[circle, draw, inner sep=1.6pt, fill,label=below:\tiny{$c_2$}] (s) at (2.5,-0.6) {};

    \node[label=below:\tiny{$b_1$}] (s) at (4,0) {};
    
    \node[circle, draw, inner sep=1.6pt, fill] (s) at (4,0) {};    
    \node[circle, draw, inner sep=1.6pt, fill, label=below:\tiny{$c_1$}] (c) at (5,0) {};
       
    \node[] (P) at (7.5,0.6) {\tiny \textcolor{blue}{$R_1$}};
    \node[] (Q) at (7.5,-0.6) {\tiny \textcolor{red}{$R_2^Q$}};
    
\end{scope}

\begin{scope}[xshift=-120, yshift=120, scale=0.8]
    \draw[ultra thick,color=black!30] (2.5,0.6) .. controls (2.6, 0.6) and (3.6, 0.6) .. (4,0);
    \draw[ultra thick,color=blue] (4,0.035)--(5,0.035);
    \draw[ultra thick,color=ForestGreen] (4,-0.035)--(5,-0.035);
    \draw[ultra thick,color=black!30, dashed] (2,0.6)--(2.5,0.6);
    \draw[ultra thick,color=black!30] (4.95,0.035) .. controls (5.4, -0.565) and (6.4, -0.565) .. (6.5,-0.565);
    \draw[ultra thick,color=black!30, dashed] (6.5,-0.565)--(7.25,-0.565);

    \draw[ultra thick,color=blue] (2.5,-0.565) .. controls (2.6, -0.565) and (3.6, -0.565) .. (4,0.035);
    \draw[ultra thick,color=ForestGreen] (2.5,-0.635) .. controls (2.6, -0.635) and (3.6, -0.635) .. (4.05,-0.035);
    \draw[ultra thick,color=blue, dashed] (1,-0.565)--(1.5,-0.565);
    \draw[ultra thick,color=blue] (1.5,-0.565)--(2.5,-0.565);
    \draw[ultra thick,color=blue] (4.95,0.035) .. controls (5.4, 0.635) and (6.4, 0.635) .. (6.5,0.635);
    \draw[ultra thick,color=blue, dashed] (6.5,0.635)--(7.25,0.635);
    \draw[ultra thick,color=ForestGreen] (5,-0.035) .. controls (5.4, 0.565) and (6.4, 0.565) .. (6.5,0.565);
    \draw[ultra thick,color=ForestGreen, dashed] (6.5,0.565)--(7.25,0.565);

    \draw[ultra thick,color=ForestGreen, dashed] (1,-1) --(1.5,-1);    
    \draw[ultra thick,color=ForestGreen] (1.5,-1) .. controls (1.6, -1) and (1.9, -0.8) .. (2,-0.6);
    \draw[ultra thick,color=ForestGreen] (2,-0.635)--(2.5,-0.635);  
    \draw[ultra thick,color=black!30] (2.5,-0.6) .. controls (2.6, -0.35) and (2.6, -0.35) .. (2.7,-0.2); 
    \draw[ultra thick,color=black!30, dashed] (2.7,-0.2) .. controls (2.85,0) and (3.3,0.1) .. (3.4,0.1);

    \node[circle, draw, inner sep=1.6pt, fill,label=below:\tiny{$b_2$}] (s) at (2,-0.6) {};
    \node[circle, draw, inner sep=1.6pt, fill,label=below:\tiny{$c_2$}] (s) at (2.5,-0.6) {};

    \node[label=below:\tiny{$b_1$}] (s) at (4,0) {};
    
    \node[circle, draw, inner sep=1.6pt, fill] (s) at (4,0) {};    
    \node[circle, draw, inner sep=1.6pt, fill, label=below:\tiny{$c_1$}] (c) at (5,0) {};
       
    \node[] (P) at (0.6,-0.5) {\tiny \textcolor{blue}{$R_1$}};
    \node[] (Q) at (0.3,-1) {\tiny \textcolor{ForestGreen}{$R_2b_2R_1$}};
\end{scope}

\begin{scope}[xshift=120, yshift=90, scale=0.8]

    \draw[ultra thick,color=blue] (2.5,0.6) .. controls (2.6, 0.6) and (3.6, 0.6) .. (4,0);
    \draw[ultra thick,color=blue] (4,0)--(5,0);
    \draw[ultra thick,color=blue, dashed] (2,0.6)--(2.5,0.6);
    \draw[ultra thick,color=blue] (5,0.035) .. controls (5.4, -0.565) and (6.4, -0.565) .. (6.5,-0.565);
    \draw[ultra thick,color=blue, dashed] (6.5,-0.565)--(7.25,-0.565);

    \draw[ultra thick,color=black!30] (2.5,-0.6) .. controls (2.6, -0.6) and (3.6, -0.6) .. (4,0);
    \draw[ultra thick,color=black!30, dashed] (1,-0.6)--(1.5,-0.6);
    \draw[ultra thick,color=black!30] (1.5,-0.6)--(2.5,-0.6);    
    \draw[ultra thick,color=black!30] (5,0) .. controls (5.4, 0.6) and (6.4, 0.6) .. (6.5,0.6);
    \draw[ultra thick,color=black!30, dashed] (6.5,0.6)--(7.25,0.6);

    \draw[ultra thick,color=red, dashed] (1,-1) --(1.5,-1);    
    \draw[ultra thick,color=red] (1.5,-1) .. controls (1.6, -1) and (1.9, -0.8) .. (2,-0.6);
    \draw[ultra thick,color=red] (2,-0.6)--(2.5,-0.6);  
    \draw[ultra thick,color=red] (2.5,-0.6) .. controls (2.6, -0.35) and (2.6, -0.35) .. (2.7,-0.2); 
    \draw[ultra thick,color=red, dashed] (2.7,-0.2) .. controls (2.85,0) and (3.3,0.1) .. (3.4,0.1);

    \node[circle, draw, inner sep=1.6pt, fill,label=below:\tiny{$b_2$}] (s) at (2,-0.6) {};
    \node[circle, draw, inner sep=1.6pt, fill,label=below:\tiny{$c_2$}] (s) at (2.5,-0.6) {};

    \node[label=below:\tiny{$b_1$}] (s) at (4,0) {};
    
    \node[circle, draw, inner sep=1.6pt, fill] (s) at (4,0) {};    
    \node[circle, draw, inner sep=1.6pt, fill, label=below:\tiny{$c_1$}] (c) at (5,0) {};
       
    \node[] (P) at (1.6,0.6) {\tiny \textcolor{blue}{$R_1^Q$}};
    \node[] (Q) at (0.6,-1) {\tiny \textcolor{red}{$R_2$}};
    
\end{scope}

\begin{scope}[xshift=0, yshift=185, scale=0.8]

    \draw[ultra thick,color=black!30] (2.5,0.6) .. controls (2.6, 0.6) and (3.6, 0.6) .. (4,0);
    \draw[ultra thick,color=ForestGreen] (4,0)--(5,0);
    \draw[ultra thick,color=black!30, dashed] (2,0.6)--(2.5,0.6);
    \draw[ultra thick,color=black!30] (4.95,0.035) .. controls (5.4, -0.565) and (6.4, -0.565) .. (6.5,-0.565);
    \draw[ultra thick,color=black!30, dashed] (6.5,-0.565)--(7.25,-0.565);

    \draw[ultra thick,color=ForestGreen] (2.5,-0.635) .. controls (2.6, -0.6) and (3.6, -0.6) .. (4,0);
    \draw[ultra thick,color=black!30, dashed] (1,-0.6)--(1.5,-0.6);
    \draw[ultra thick,color=black!30] (1.5,-0.6)--(2,-0.6);
    \draw[ultra thick,color=ForestGreen] (5,0) .. controls (5.4, 0.6) and (6.4, 0.6) .. (6.5,0.6);
    \draw[ultra thick,color=ForestGreen, dashed] (6.5,0.6)--(7.25,0.6);

    \draw[ultra thick,color=ForestGreen, dashed] (1,-1.035) --(1.5,-1.035);
    \draw[ultra thick,color=red, dashed] (1,-0.965) --(1.5,-0.965);        
    \draw[ultra thick,color=ForestGreen] (1.5,-1.035) .. controls (1.6, -1.035) and (1.9, -0.835) .. (2,-0.635);
    \draw[ultra thick,color=red] (1.5,-0.965) .. controls (1.6, -0.965) and (1.9, -0.76) .. (2,-0.52);    
    \draw[ultra thick,color=ForestGreen] (2,-0.635)--(2.5,-0.635);  
    \draw[ultra thick,color=red] (2,-0.565)--(2.5,-0.565);  
    \draw[ultra thick,color=red] (2.5,-0.6) .. controls (2.6, -0.35) and (2.6, -0.35) .. (2.7,-0.2); 
    \draw[ultra thick,color=red, dashed] (2.7,-0.2) .. controls (2.85,0) and (3.3,0.1) .. (3.4,0.1);

    \node[circle, draw, inner sep=1.6pt, fill,label=below:\tiny{$b_2$}] (s) at (2,-0.6) {};
    \node[circle, draw, inner sep=1.6pt, fill,label=below:\tiny{$c_2$}] (s) at (2.5,-0.6) {};

    \node[label=below:\tiny{$b_1$}] (s) at (4,0) {};
    
    \node[circle, draw, inner sep=1.6pt, fill] (s) at (4,0) {};    
    \node[circle, draw, inner sep=1.6pt, fill, label=below:\tiny{$c_1$}] (c) at (5,0) {};
       
    \node[] (P) at (2.4,0) {\tiny \textcolor{red}{$R_2$}};
    \node[] (Q) at (8,0.6) {\tiny \textcolor{ForestGreen}{$R_2b_2R_1$}};
\end{scope}

\begin{scope}[xshift=13, yshift=10, scale = 0.6]
    \node[] (a) at (0,0) {};
    \node[] (b) at (-2,2) {};
    \draw[very thick] (a) -- (b);
\end{scope}
\begin{scope}[xshift=-20, yshift=75, scale = 0.5]
    \node[] (a) at (0,0) {};
    \node[] (b) at (0,2) {};
    \draw[very thick] (a) -- (b);
\end{scope}
\begin{scope}[xshift=170, yshift=10, scale = 0.6]
    \node[] (a) at (0,0) {};
    \node[] (b) at (2,3) {};
    \draw[very thick] (a) -- (b);
\end{scope}
\begin{scope}[xshift=170, yshift=110, scale = 0.6]
    \node[] (a) at (2,0) {};
    \node[] (b) at (0,3) {};
    \draw[very thick] (a) -- (b);
\end{scope}
\begin{scope}[xshift=13, yshift=140, scale = 0.6]
    \node[] (a) at (-2,0) {};
    \node[] (b) at (0,2) {};
    \draw[very thick] (a) -- (b);
\end{scope}

\begin{scope}[xshift=-5, yshift=0, scale = 1.1]
    \node[] (a) at (0,0) {$Q$};
\end{scope}
\begin{scope}[xshift=10, yshift=190, scale = 1.1]
    \node[] (a) at (0,0) {$C_1'\vee C_2$};
\end{scope}
\begin{scope}[xshift=-105, yshift=65, scale = 1.1]
    \node[] (a) at (0,0) {$C_1$};
\end{scope}
\begin{scope}[xshift=-105, yshift=125, scale = 1.1]
    \node[] (a) at (0,0) {$C_1'$};
\end{scope}
\begin{scope}[xshift=300, yshift=90, scale = 1.1]
    \node[] (a) at (0,0) {$C_2$};
\end{scope}
\end{tikzpicture}
    \caption{The pentagonal sub-cases in the proof of~\Cref{prop.intervalJoin}(v).}
    \label{fig.pentagon}
\end{figure}

\begin{figure}
    \centering
    \begin{tikzpicture}
\begin{scope}[xshift=-35, yshift=-20, scale=0.8]
    \draw[ultra thick,color=blue, dashed] (2,0.635)--(2.5,0.635);
    \draw[ultra thick,color=red, dashed] (2,0.565)--(2.5,0.565);        
    \draw[ultra thick,color=blue] (2.5,0.635) .. controls (2.6, 0.635) and (3.6, 0.635) .. (4.05,0.035);
    \draw[ultra thick,color=red] (2.5,0.565) .. controls (2.6, 0.565) and (3.6, 0.565) .. (4,-0.035);    
    \draw[ultra thick,color=blue] (4,0.035)--(5,0.035);
    \draw[ultra thick,color=blue] (5,0) .. controls (5.4, 0.6) and (6.6, 0.6) .. (7,0);
    \draw[ultra thick,color=blue] (7,0.035)--(8,0.035);    
    \draw[ultra thick,color=black!30] (8,0) .. controls (8.4, 0.6) and (9.4, 0.6) .. (9.5,0.6);
    \draw[ultra thick,color=black!30, dashed] (9.5,0.6)--(10.2,0.6);

    \draw[ultra thick,color=black!30, dashed] (2,-0.6)--(2.5,-0.6);
    \draw[ultra thick,color=black!30] (2.5,-0.6) .. controls (2.6, -0.6) and (3.6, -0.5) .. (4,0);
    \draw[ultra thick,color=red] (4,-0.035)--(5,-0.035);     
    \draw[ultra thick,color=red] (5,0) .. controls (5.4, -0.6) and (6.6, -0.6) .. (7,0);
    \draw[ultra thick,color=red] (7,-0.035)--(8,-0.035);     
    \draw[ultra thick,color=blue] (8,0.035) .. controls (8.4, -0.565) and (9.4, -0.565) .. (9.5,-0.565);
    \draw[ultra thick,color=red] (7.95,-0.035) .. controls (8.4, -0.635) and (9.4, -0.635) .. (9.5,-0.635);    
    \draw[ultra thick,color=blue, dashed] (9.5,-0.565)--(10.2,-0.565);
    \draw[ultra thick,color=red, dashed] (9.5,-0.635)--(10.2,-0.635);    

    \node[circle, draw, inner sep=1.6pt, fill,label=above:\tiny{$b_2$}] (s) at (4,0) {};
    \node[circle, draw, inner sep=1.6pt, fill,label=above:\tiny{$c_2$}] (s) at (5,0) {};

    \node[circle, draw, inner sep=1.6pt, fill,label=above:\tiny{$b_1$}] (s) at (7,0) {};
    \node[circle, draw, inner sep=1.6pt, fill, label=above:\tiny{$c_1$}] (c) at (8,0) {};
       
    \node[] (P) at (6,0.8) {\tiny \textcolor{blue}{$R_1^Q$}};
    \node[] (P) at (6,-0.8) {\tiny \textcolor{red}{$R_2^Q$}};    
\end{scope}

\begin{scope}[xshift=-150, yshift=50, scale=0.8]
    \draw[ultra thick,color=blue, dashed] (2,0.635)--(2.5,0.635);
    \draw[ultra thick,color=red, dashed] (2,0.565)--(2.5,0.565);        
    \draw[ultra thick,color=blue] (2.5,0.635) .. controls (2.6, 0.635) and (3.6, 0.635) .. (4.05,0.035);
    \draw[ultra thick,color=red] (2.5,0.565) .. controls (2.6, 0.565) and (3.6, 0.565) .. (4,-0.035);    
    \draw[ultra thick,color=blue] (4,0.035)--(5,0.035);
    \draw[ultra thick,color=black!30] (5,0) .. controls (5.4, 0.6) and (6.6, 0.6) .. (7,0);
    \draw[ultra thick,color=blue] (7,0.035)--(8,0.035);    
    \draw[ultra thick,color=blue] (8,0) .. controls (8.4, 0.6) and (9.4, 0.6) .. (9.5,0.6);
    \draw[ultra thick,color=blue, dashed] (9.5,0.6)--(10.2,0.6);

    \draw[ultra thick,color=black!30, dashed] (2,-0.6)--(2.5,-0.6);
    \draw[ultra thick,color=black!30] (2.5,-0.6) .. controls (2.6, -0.6) and (3.6, -0.5) .. (4,0);
    \draw[ultra thick,color=red] (4,-0.035)--(5,-0.035);     
    \draw[ultra thick,color=red] (4.95,-0.035) .. controls (5.4, -0.635) and (6.6, -0.635) .. (7.05,-0.035);
    \draw[ultra thick,color=blue] (5,0.035) .. controls (5.4, -0.565) and (6.6, -0.565) .. (7,0.035);    
    \draw[ultra thick,color=red] (7,-0.035)--(8,-0.035);     
    \draw[ultra thick,color=red] (8,0.035) .. controls (8.4, -0.565) and (9.4, -0.565) .. (9.5,-0.565);
    \draw[ultra thick,color=red, dashed] (9.5,-0.565)--(10.2,-0.565);

    \node[circle, draw, inner sep=1.6pt, fill,label=above:\tiny{$b_2$}] (s) at (4,0) {};
    \node[circle, draw, inner sep=1.6pt, fill,label=above:\tiny{$c_2$}] (s) at (5,0) {};

    \node[circle, draw, inner sep=1.6pt, fill,label=above:\tiny{$b_1$}] (s) at (7,0) {};
    \node[circle, draw, inner sep=1.6pt, fill, label=above:\tiny{$c_1$}] (c) at (8,0) {};
       
    \node[] (P) at (9.6,0.9) {\tiny \textcolor{blue}{$R_1$}};
    \node[] (P) at (9.6,-0.9) {\tiny \textcolor{red}{$R_2^Q$}};     
\end{scope}

\begin{scope}[xshift=-150, yshift=130, scale=0.8]
    \draw[ultra thick,color=blue, dashed] (2,0.565)--(2.5,0.565);        
    \draw[ultra thick,color=blue] (2.5,0.565) .. controls (2.6, 0.565) and (3.6, 0.565) .. (4,-0.035);    
    \draw[ultra thick,color=blue] (4,0.035)--(5,0.035);
    \draw[ultra thick,color=black!30] (5,0) .. controls (5.4, 0.6) and (6.6, 0.6) .. (7,0);
    \draw[ultra thick,color=blue] (7,0.035)--(8,0.035);    
    \draw[ultra thick,color=blue] (7.95,0.035) .. controls (8.4, 0.635) and (9.4, 0.635) .. (9.5,0.635);
    \draw[ultra thick,color=red] (8,-0.035) .. controls (8.4, 0.565) and (9.4, 0.565) .. (9.5,0.565);    
    \draw[ultra thick,color=blue, dashed] (9.5,0.635)--(10.2,0.635);
    \draw[ultra thick,color=red, dashed] (9.5,0.565)--(10.2,0.565);

    \draw[ultra thick,color=red, dashed] (2,-0.6)--(2.5,-0.6);
    \draw[ultra thick,color=red] (2.5,-0.6) .. controls (2.6, -0.6) and (3.6, -0.5) .. (4,0);
    \draw[ultra thick,color=red] (4,-0.035)--(5,-0.035);     
    \draw[ultra thick,color=red] (4.95,-0.035) .. controls (5.4, -0.635) and (6.6, -0.635) .. (7.05,-0.035);
    \draw[ultra thick,color=blue] (5,0.035) .. controls (5.4, -0.565) and (6.6, -0.565) .. (7,0.035);    
    \draw[ultra thick,color=red] (7,-0.035)--(8,-0.035);     
    \draw[ultra thick,color=black!30] (8,0.035) .. controls (8.4, -0.565) and (9.4, -0.565) .. (9.5,-0.565);
    \draw[ultra thick,color=black!30, dashed] (9.5,-0.565)--(10.2,-0.565);

    \node[circle, draw, inner sep=1.6pt, fill,label=above:\tiny{$b_2$}] (s) at (4,0) {};
    \node[circle, draw, inner sep=1.6pt, fill,label=above:\tiny{$c_2$}] (s) at (5,0) {};

    \node[circle, draw, inner sep=1.6pt, fill,label=above:\tiny{$b_1$}] (s) at (7,0) {};
    \node[circle, draw, inner sep=1.6pt, fill, label=above:\tiny{$c_1$}] (c) at (8,0) {};
       
    \node[] (P) at (2.6,0.9) {\tiny \textcolor{blue}{$R_1$}};
    \node[] (P) at (2.6,-0.9) {\tiny \textcolor{red}{$R_2b_2R_1$}};     
\end{scope}

\begin{scope}[xshift=90, yshift=50, scale=0.8]
    \draw[ultra thick,color=blue, dashed] (2,0.565)--(2.5,0.565);        
    \draw[ultra thick,color=blue] (2.5,0.565) .. controls (2.6, 0.565) and (3.6, 0.565) .. (4,-0.035);    
    \draw[ultra thick,color=blue] (4,0.035)--(5,0.035);
    \draw[ultra thick,color=blue] (4.95,0.035) .. controls (5.4, 0.635) and (6.6, 0.635) .. (7.05,0.035);
    \draw[ultra thick,color=red] (5,-0.035) .. controls (5.4, 0.565) and (6.6, 0.565) .. (7,-0.035);
    \draw[ultra thick,color=blue] (7,0.035)--(8,0.035);    
    \draw[ultra thick,color=black!30] (8,0) .. controls (8.4, 0.6) and (9.4, 0.6) .. (9.5,0.6);
    \draw[ultra thick,color=black!30, dashed] (9.5,0.6)--(10.2,0.6);

    \draw[ultra thick,color=red, dashed] (2,-0.6)--(2.5,-0.6);
    \draw[ultra thick,color=red] (2.5,-0.6) .. controls (2.6, -0.6) and (3.6, -0.5) .. (4,0);
    \draw[ultra thick,color=red] (4,-0.035)--(5,-0.035);     
    \draw[ultra thick,color=black!30] (5,0) .. controls (5.4, -0.6) and (6.6, -0.6) .. (7,0);
    \draw[ultra thick,color=red] (7,-0.035)--(8,-0.035);     
    \draw[ultra thick,color=blue] (8,0.035) .. controls (8.4, -0.565) and (9.4, -0.565) .. (9.5,-0.565);
    \draw[ultra thick,color=red] (7.95,-0.035) .. controls (8.4, -0.635) and (9.4, -0.635) .. (9.5,-0.635);    
    \draw[ultra thick,color=blue, dashed] (9.5,-0.565)--(10.2,-0.565);
    \draw[ultra thick,color=red, dashed] (9.5,-0.635)--(10.2,-0.635);    

    \node[circle, draw, inner sep=1.6pt, fill,label=above:\tiny{$b_2$}] (s) at (4,0) {};
    \node[circle, draw, inner sep=1.6pt, fill,label=above:\tiny{$c_2$}] (s) at (5,0) {};

    \node[circle, draw, inner sep=1.6pt, fill,label=above:\tiny{$b_1$}] (s) at (7,0) {};
    \node[circle, draw, inner sep=1.6pt, fill, label=above:\tiny{$c_1$}] (c) at (8,0) {};
       
    \node[] (P) at (2.6,0.9) {\tiny \textcolor{blue}{$R_1^Q$}};
    \node[] (P) at (2.6,-0.9) {\tiny \textcolor{red}{$R_2$}};
\end{scope}

\begin{scope}[xshift=90, yshift=130, scale=0.8]
    \draw[ultra thick,color=black!30, dashed] (2,0.565)--(2.5,0.565);        
    \draw[ultra thick,color=black!30] (2.5,0.565) .. controls (2.6, 0.565) and (3.6, 0.565) .. (4,-0.035);    
    \draw[ultra thick,color=blue] (4,0.035)--(5,0.035);
    \draw[ultra thick,color=blue] (4.95,0.035) .. controls (5.4, 0.635) and (6.6, 0.635) .. (7.05,0.035);
    \draw[ultra thick,color=red] (5,-0.035) .. controls (5.4, 0.565) and (6.6, 0.565) .. (7,-0.035);
    \draw[ultra thick,color=blue] (7,0.035)--(8,0.035);    
    \draw[ultra thick,color=blue] (8,-0.035) .. controls (8.4, 0.565) and (9.4, 0.565) .. (9.5,0.565);    
    \draw[ultra thick,color=blue, dashed] (9.5,0.565)--(10.2,0.565);    

    \draw[ultra thick,color=blue, dashed] (2,-0.565)--(2.5,-0.565);
    \draw[ultra thick,color=red, dashed] (2,-0.635)--(2.5,-0.635);    
    \draw[ultra thick,color=blue] (2.5,-0.565) .. controls (2.6, -0.565) and (3.6, -0.565) .. (4,0.035);
    \draw[ultra thick,color=red] (2.5,-0.635) .. controls (2.6, -0.635) and (3.6, -0.635) .. (4.05,-0.035);    
    \draw[ultra thick,color=red] (4,-0.035)--(5,-0.035);     
    \draw[ultra thick,color=black!30] (5,0) .. controls (5.4, -0.6) and (6.6, -0.6) .. (7,0);
    \draw[ultra thick,color=red] (7,-0.035)--(8,-0.035);     
    \draw[ultra thick,color=red] (8,0.035) .. controls (8.4, -0.565) and (9.4, -0.565) .. (9.5,-0.565);
    \draw[ultra thick,color=red, dashed] (9.5,-0.565)--(10.2,-0.565);

    \node[circle, draw, inner sep=1.6pt, fill,label=above:\tiny{$b_2$}] (s) at (4,0) {};
    \node[circle, draw, inner sep=1.6pt, fill,label=above:\tiny{$c_2$}] (s) at (5,0) {};

    \node[circle, draw, inner sep=1.6pt, fill,label=above:\tiny{$b_1$}] (s) at (7,0) {};
    \node[circle, draw, inner sep=1.6pt, fill, label=above:\tiny{$c_1$}] (c) at (8,0) {};
       
    \node[] (P) at (9.6,0.9) {\tiny \textcolor{blue}{$R_2b_1R_1$}};
    \node[] (P) at (9.6,-0.9) {\tiny \textcolor{red}{$R_2$}};
\end{scope}

\begin{scope}[xshift=-35, yshift=200, scale=0.8]
    \draw[ultra thick,color=black!30, dashed] (2,0.565)--(2.5,0.565);        
    \draw[ultra thick,color=black!30] (2.5,0.565) .. controls (2.6, 0.565) and (3.6, 0.565) .. (4,-0.035);    
    \draw[ultra thick,color=blue] (4,0.035)--(5,0.035);
    \draw[ultra thick,color=blue] (5,-0.035) .. controls (5.4, 0.565) and (6.6, 0.565) .. (7,-0.035);
    \draw[ultra thick,color=blue] (7,0.035)--(8,0.035);    
    \draw[ultra thick,color=blue] (7.95,0.035) .. controls (8.4, 0.635) and (9.4, 0.635) .. (9.5,0.635);
    \draw[ultra thick,color=red] (8,-0.035) .. controls (8.4, 0.565) and (9.4, 0.565) .. (9.5,0.565);    
    \draw[ultra thick,color=blue, dashed] (9.5,0.635)--(10.2,0.635);
    \draw[ultra thick,color=red, dashed] (9.5,0.565)--(10.2,0.565);    

    \draw[ultra thick,color=blue, dashed] (2,-0.565)--(2.5,-0.565);
    \draw[ultra thick,color=red, dashed] (2,-0.635)--(2.5,-0.635);    
    \draw[ultra thick,color=blue] (2.5,-0.565) .. controls (2.6, -0.565) and (3.6, -0.565) .. (4,0.035);
    \draw[ultra thick,color=red] (2.5,-0.635) .. controls (2.6, -0.635) and (3.6, -0.635) .. (4.05,-0.035);    
    \draw[ultra thick,color=red] (4,-0.035)--(5,-0.035);     
    \draw[ultra thick,color=red] (5,0) .. controls (5.4, -0.6) and (6.6, -0.6) .. (7,0);
    \draw[ultra thick,color=red] (7,-0.035)--(8,-0.035);     
    \draw[ultra thick,color=black!30] (8,0.035) .. controls (8.4, -0.565) and (9.4, -0.565) .. (9.5,-0.565);
    \draw[ultra thick,color=black!30, dashed] (9.5,-0.565)--(10.2,-0.565);

    \node[circle, draw, inner sep=1.6pt, fill,label=above:\tiny{$b_2$}] (s) at (4,0) {};
    \node[circle, draw, inner sep=1.6pt, fill,label=above:\tiny{$c_2$}] (s) at (5,0) {};

    \node[circle, draw, inner sep=1.6pt, fill,label=above:\tiny{$b_1$}] (s) at (7,0) {};
    \node[circle, draw, inner sep=1.6pt, fill, label=above:\tiny{$c_1$}] (c) at (8,0) {};
       
    \node[] (P) at (6,0.8) {\tiny \textcolor{blue}{$R_2b_1R_1$}};
    \node[] (P) at (6,-0.8) {\tiny \textcolor{red}{$R_2b_2R_1$}};
\end{scope}

\begin{scope}[xshift=10, yshift=5, scale = 1]
    \node[] (a) at (0,0) {};
    \node[] (b) at (-1,1) {};
    \draw[very thick] (a) -- (b);
\end{scope}
\begin{scope}[xshift=-25, yshift=70, scale = 1]
    \node[] (a) at (0,0) {};
    \node[] (b) at (0,1.5) {};
    \draw[very thick] (a) -- (b);
\end{scope}
\begin{scope}[xshift=200, yshift=5, scale = 1]
    \node[] (a) at (0,0) {};
    \node[] (b) at (1,1) {};
    \draw[very thick] (a) -- (b);
\end{scope}
\begin{scope}[xshift=230, yshift=70, scale = 1]
    \node[] (a) at (0,0) {};
    \node[] (b) at (0,1.5) {};
    \draw[very thick] (a) -- (b);
\end{scope}
\begin{scope}[xshift=200, yshift=150, scale = 1]
    \node[] (a) at (1,0) {};
    \node[] (b) at (0,1) {};
    \draw[very thick] (a) -- (b);
\end{scope}
\begin{scope}[xshift=-25, yshift=150, scale = 1]
    \node[] (a) at (0,0) {};
    \node[] (b) at (1,1) {};
    \draw[very thick] (a) -- (b);
\end{scope}

\begin{scope}[xshift=0, yshift=-20, scale = 1.1]
    \node[] (a) at (0,0) {$Q$};
\end{scope}
\begin{scope}[xshift=-10, yshift=200, scale = 1.1]
    \node[] (a) at (0,0) {$C_1\vee C_2$};
\end{scope}
\begin{scope}[xshift=-115, yshift=50, scale = 1.1]
    \node[] (a) at (0,0) {$C_1$};
\end{scope}
\begin{scope}[xshift=-115, yshift=130, scale = 1.1]
    \node[] (a) at (0,0) {$C_1'$};
\end{scope}
\begin{scope}[xshift=330, yshift=50, scale = 1.1]
    \node[] (a) at (0,0) {$C_2$};
\end{scope}
\begin{scope}[xshift=330, yshift=130, scale = 1.1]
    \node[] (a) at (0,0) {$C_2'$};
\end{scope}

\end{tikzpicture}
    \caption{The hexagonal sub-case in the proof of~\Cref{prop.intervalJoin}(v).}
    \label{fig.hexagon}
\end{figure}

A \defn{polygon} in a lattice is an interval $[x,y]$ that is the union of two finite maximal chains from $x$ to $y$ that are disjoint except at $x$ and $y$. 
A lattice is said to be \defn{polygonal} if the following two conditions hold:
\begin{itemize}
    \item[1.] If $y_1$ and $y_2$ are distinct and cover an element $x$, then $[x,y_1\vee y_2]$ is a polygon.
    \item[2.] If $y_1$ and $y_2$ are distinct and are covered by an element $x$, then $[y_1\wedge y_2,x]$ is a polygon.
\end{itemize}

\begin{theorem}\label{thm.framing_lattice}
The framing poset $\scrL_{G,F}$ is a polygonal lattice. The polygons are squares, pentagons or hexagons. The two chains between the minimal and maximal element of a polygon are of length 2 or 3. 
\end{theorem}

\begin{proof}
    The fact that $\scrL_{G,F}$ is a lattice follows from the BEZ Lemma~\ref{lem.BEZ} and~\Cref{prop.intervalJoin}(v).
    If $y_1$ and $y_2$ are distinct and cover an element $x$, then~\Cref{prop.intervalJoin}(iv) and (v) shows that $[x,y_1\vee y_2]$ is a polygon, whose two chains are of length 2 or 3.
    For the case when $y_1$ and $y_2$ are distinct and are covered by an element $x$, we can use the fact that the dual of $\scrL_{G,F}$ is a framing lattice (Lemma~\ref{lem.dualizing_operations}), and repeat the previous argument.
\end{proof}

\subsection{Meets and joins}

In the previous section, the join (resp. meet) of two maximal cliques $C$ and $C'$ in $\scrL_{G,F}$ was shown to be $C_{\max}(C\cap C')$ (resp. $C_{\min}(C\cap C')$) in the special case where $C$ and $C'$ cover (resp. are covered by) a common element. Figure~\ref{fig_join_CmaxRem} shows that this is not true in general; it shows two maximal cliques $C$ and $C'$ satisfying $C\vee C'\neq C_{\max}(C\cap C')$.  
In contrast, we also have cases where $C\wedge C'\neq C_{\min}(C\cap C')$.

\begin{figure}[h]
    \centering

    \caption{Two maximal cliques satisfying $C\vee C' \neq C_{\max}(C\cap C')$.}
    \label{fig_join_CmaxRem}
\end{figure}

The goal of this subsection is to introduce two modified algorithms to compute joins and meets. When we compute $C_{\max}(C\cap C')$ using Algorithm~\ref{alg.CmaxS} for the example in~\Cref{fig_join_CmaxRem}, we get a route $R$ that does not belong to $C\vee C'$. This problem can be solved by forbidding this route to be added during the algorithm. 
This idea leads to the definition of a \defn{(ccw) removed set} $\Rem_{\ccw}(C,C')$ that appears naturally in the context of cross-Tamari lattices, a class of framing lattices that will be introduced in~\Cref{sec_crossTamari} 
(the points inside the square grids in~\Cref{fig_join_CmaxRem} represent routes in that context).  
For the computation of the meet we define another \defn{(cw) removed set} $\Rem_{\cw}(C,C')$.  

Let $P$ be a path and let $\ccw(P)$ denote the set of routes ccw from $P$, i.e.
$$ \ccw(P) := \{R \mid  P<_v^{cw} R \text{ for some } v \in G\}. $$ 
Similarly, let $\cw(P)$ denote the set of routes clockwise from $P$, i.e. 
$$ \cw(P) := \{R \mid  R<_v^{cw} P \text{ for some } v \in G\}. $$ 
Note that if $P$ is of length less than $2$, then $\ccw(P)$ and $\cw(P)$ are necessarily empty.
An example is shown in~\Cref{fig_ccwP}.

\begin{figure}[h]
    \centering
    \begin{tikzpicture}

\begin{scope}[scale=0.7, xshift=-50, yshift=0]
    \draw[thick, color=black!30] (0,0) -- (1,0);
    \draw[thick, color=black!30] (1,0) -- (2,0);
    \draw[ultra thick, color=blue] (2,0) -- (3,0);
    \draw[ultra thick, color=blue] (3,0) -- (4,0);
    \draw[thick, color=black!30] (4,0) -- (5,0);
    \draw[thick, color=black!30] (5,0) -- (6,0);
    \draw[thick, color=black!30] (6,0) -- (7,0);

    \draw[thick, color=black!30] (0,0) .. controls (0.4, -0.7) and (1.6, -0.7) .. (2,0);
    \draw[thick, color=black!30] (0,0) .. controls (0.4, -1.3) and (2.6, -1.3) .. (3,0);
    \draw[thick, color=black!30] (5,0) .. controls (5.4, -0.7) and (6.6, -0.7) .. (7,0);
    \draw[ultra thick, color=blue] (4,0) .. controls (4.4, -1.3) and (6.6, -1.3) .. (7,0);    
    
    \node[circle, draw, inner sep=1.4pt, fill] (s) at (0,0) {};
    \node[circle, draw, inner sep=1.4pt, fill] (1) at (1,0) {};
    \node[circle, draw, inner sep=1.4pt, fill] (2) at (2,0) {};
    \node[circle, draw, inner sep=1.4pt, fill] (3) at (3,0) {};    
    \node[circle, draw, inner sep=1.4pt, fill] (33) at (4,0) {};
    \node[circle, draw, inner sep=1.4pt, fill] (22) at (5,0) {};
    \node[circle, draw, inner sep=1.4pt, fill] (11) at (6,0) {};
    \node[circle, draw, inner sep=1.4pt, fill] (t) at (7,0) {};        

    \node[color=blue] at (3.5,-0.6) {$P$};

    \node[] (a) at (0,0.4) {\scriptsize $s$};
    \node[] (a) at (1,0.4) {\scriptsize $1$};  
    \node[] (a) at (2,0.4) {\scriptsize $2$};
    \node[] (a) at (3,0.4) {\scriptsize $3$};
    \node[] (a) at (4,0.4) {\scriptsize $\overline{3}$};
    \node[] (a) at (6,0.4) {\scriptsize $\overline{1}$};
    \node[] (a) at (5,0.4) {\scriptsize $\overline{2}$};  
    \node[] (a) at (7,0.4) {\scriptsize $t$};
\end{scope}	

\begin{scope}[scale=0.65, xshift=200, yshift=0]
    \draw[ultra thick, color=black!30, ForestGreen] (3,0) -- (4,0);
    \draw[ultra thick, color=black!30, ForestGreen] (4,0) -- (5,0);

    \draw[ultra thick, color=ForestGreen] (0,0) .. controls (0.4, -1.3) and (2.6, -1.3) .. (3,0);
    \draw[ultra thick, color=ForestGreen] (5,0) .. controls (5.4, -0.7) and (6.6, -0.7) .. (7,0);
    
    \node[circle, draw, inner sep=1.4pt, fill] (s) at (0,0) {};
    \node[circle, draw, inner sep=1.4pt, fill] (1) at (1,0) {};
    \node[circle, draw, inner sep=1.4pt, fill] (2) at (2,0) {};
    \node[circle, draw, inner sep=1.4pt, fill] (3) at (3,0) {};    
    \node[circle, draw, inner sep=1.4pt, fill] (33) at (4,0) {};
    \node[circle, draw, inner sep=1.4pt, fill] (22) at (5,0) {};
    \node[circle, draw, inner sep=1.4pt, fill] (11) at (6,0) {};
    \node[circle, draw, inner sep=1.4pt, fill] (t) at (7,0) {};        


\end{scope}

\begin{scope}[scale=0.65, xshift=200, yshift=-50]
    \draw[ultra thick, color=black!30, ForestGreen] (3,0) -- (4,0);
    \draw[ultra thick, color=black!30, ForestGreen] (4,0) -- (5,0);
    \draw[ultra thick, color=black!30, ForestGreen] (5,0) -- (6,0);
    \draw[ultra thick, color=black!30, ForestGreen] (6,0) -- (7,0);

    \draw[ultra thick, color=ForestGreen] (0,0) .. controls (0.4, -1.3) and (2.6, -1.3) .. (3,0);
    
    \node[circle, draw, inner sep=1.4pt, fill] (s) at (0,0) {};
    \node[circle, draw, inner sep=1.4pt, fill] (1) at (1,0) {};
    \node[circle, draw, inner sep=1.4pt, fill] (2) at (2,0) {};
    \node[circle, draw, inner sep=1.4pt, fill] (3) at (3,0) {};    
    \node[circle, draw, inner sep=1.4pt, fill] (33) at (4,0) {};
    \node[circle, draw, inner sep=1.4pt, fill] (22) at (5,0) {};
    \node[circle, draw, inner sep=1.4pt, fill] (11) at (6,0) {};
    \node[circle, draw, inner sep=1.4pt, fill] (t) at (7,0) {};        


\end{scope}

\begin{scope}[xshift = 360, yshift = -25, scale=0.5]
    \draw[very thin, color=gray!70] (0,0) grid (2,2);

    \node[circle, draw=black, inner sep=1.3pt, fill, color=ForestGreen]  at (2,1)  {};
    \node[circle, draw=black, inner sep=1.3pt, fill, color=ForestGreen]  at (2,2)  {};

    \node[]  at (0,2.6)  {\footnotesize $1$};
    \node[]  at (1,2.6)  {\footnotesize $2$};
    \node[]  at (2,2.6)  {\footnotesize $3$};
    
    \node[]  at (-.6,0)  {\footnotesize $\overline{3}$};
    \node[]  at (-.6,1)  {\footnotesize $\overline{2}$};
    \node[]  at (-.6,2)  {\footnotesize $\overline{1}$};

    \node[]  at (1,-0.8)  {$\ccw(P)$};
\end{scope}

\begin{scope}[xshift = 275, yshift = -15, scale=0.6]
    \node[]  at (0,0)  {\large $\Biggr\}$};
    \node[]  at (1.6,0)  {$\ccw(P)$};
\end{scope}

\end{tikzpicture}
    \caption{Example of $\ccw(P)$.}
    \label{fig_ccwP}
\end{figure}

For two maximal cliques $C$ and $C'$, we define the following sets of routes: 
\begin{align*}
\Rem_{\ccw}(C,C') &:= \bigcup_{\substack{P \text{ is a path} \\ \ccw(P) \cap (C\cup C') =\emptyset }} \ccw(P)\, , \text{ and} \\  
\Rem_{\mathrm{cw}}(C,C') &:= \bigcup_{\substack{P \text{ is a path} \\ \cw(P) \cap (C\cup C') =\emptyset }}  \mathrm{cw}(P).
\end{align*}

Here $\Rem$ stands for ``removed'', as the routes in $\Rem_{\ccw}(C,C')$ and $\Rem_{\cw}(C,C')$ will be removed from consideration in the construction of the join and meet.
For the join, we construct the clique $C_{\max}^{\Rem}(C,C')$ with an algorithm similar to Algorithm~\ref{alg.CmaxS}, but add an extra condition in line $5$ that $PvQ \notin \Rem_{\ccw}(C,C')$.

Also note that in line $1$ we start with $C_{\max}^{\Rem}(C,C') := C\cap C'$ instead of the empty set, because $S=C\cap C'$ is included in $C_{\max}(S)$ at the end of Algorithm~\ref{alg.CmaxS} anyway, so this small change does not alter the outcome of the algorithm.
The algorithm is given below in Algorithm~\ref{alg.join}, and we will show that $C_{\max}^{\Rem}(C, C') = C\vee C'$. 
An example of $\Rem_{\ccw}(C,C')$ and $C_{\max}^{\Rem}(C, C')$ is shown in~\Cref{fig_join_CmaxRem}. 

\begin{algorithm}
  \caption{The construction of $C_{\max}^{\Rem}(C,C')$}\label{alg.join}
  \begin{algorithmic}[1]
    \State $C_{\max}^{\Rem}(C,C') := C\cap C'$
    \For{$v \in V(G)$ (in increasing order)} 
        \For{$Pv \in \scrI(v)$ (in the order $\leq_{\scrI(v)}^{\rev}$)} \Comment{$Pv$ possibly empty}
            \For{$vQ \in \scrO(v)$ (in the order $\leq_{\scrO(v)}$)} \Comment{$vQ$ possibly empty}
                \If{$PvQ$ is coherent with $C_{\max}^{\Rem}(C,C')$ \textbf{and} $PvQ \notin \Rem_{\ccw}(C,C')$}
                    \State $C_{\max}^{\Rem} (C,C') := C_{\max}^{\Rem}(C,C') \cup \{PvQ\}$
                    \State \textbf{break} \Comment{This terminates the innermost loop}
                \EndIf
            \EndFor
        \EndFor
    \EndFor  
  \end{algorithmic}
\end{algorithm}

We construct a maximal clique $C_{\min}^{\Rem}(C,C')$ similarly, but in Algorithm~\ref{alg.join} we replace $\Rem_{\ccw}(C,C')$ with $\Rem_{\mathrm{cw}}(C,C')$.
We also reverse the orders in which $\scrI(v)$ and $\scrO(v)$ are read. 
That is, we replace $\leq_{\scrI(v)}^{\rev}$ with $\leq_{\scrI(v)}$ and replace $\leq_{\scrO(v)}$ with $\leq_{\scrO(v)}^{\rev}$.
See Algorithm~\ref{alg.meet} below. 

\begin{algorithm}
  \caption{The construction of $C_{\min}^{\Rem}(C,C')$}\label{alg.meet}
  \begin{algorithmic}[1]
    \State $C_{\max}^{\Rem}(C,C') := C\cap C'$
    \For{$v \in V(G)$ (in increasing order)} 
        \For{$Pv \in \scrI(v)$ (in the order $\leq_{\scrI(v)}$)} \Comment{$Pv$ possibly empty}
            \For{$vQ \in \scrO(v)$ (in the order $\leq_{\scrO(v)}^{\rev}$)} \Comment{$vQ$ possibly empty}
                \If{$PvQ$ is coherent with $C_{\max}^{\Rem}(C,C')$ \textbf{and} $PvQ \notin \Rem_{\cw}(C,C')$}
                    \State $C_{\max}^{\Rem} (C,C') := C_{\max}^{\Rem}(C,C') \cup \{PvQ\}$
                    \State \textbf{break} \Comment{This terminates the innermost loop}
                \EndIf
            \EndFor
        \EndFor
    \EndFor  
  \end{algorithmic}
\end{algorithm}

We will show that $C_{\min}^{\Rem}(C, C') = C\wedge C'$.
However, we first verify that $C_{\max}^{\Rem}(C, C') = C\vee C'$ and $C_{\min}^{\Rem}(C, C') = C\wedge C'$ are in fact maximal cliques.

\begin{lemma}\label{lem_CmaxRem_maximal}
    The cliques $C_{\max}^{\Rem}(C,C')$ and $C_{\min}^{\Rem}(C,C')$ are maximal. 
\end{lemma}

\begin{proof}
    We prove only the statement for $C_{\max}^{\Rem}(C,C')$ as the proof for $C_{\max}^{\Rem}(C,C')$ is symmetric. 
    Suppose toward a contradiction that $C_{\max}^{\Rem}(C,C')$ is not maximal.
    Then, there is a route $R^*$ that is coherent with all routes in $C_{\max}^{\Rem}(C,C')$, but not in it. 
    Note that when $v=t$ (the sink), Algorithm~\ref{alg.join} runs over all the routes of the graph ($P$ runs over all routes and $Q$ is empty). Since $R^*$ is not added during the algorithm, then $R^*$ must be in $\Rem_{\ccw}(C,C')$, i.e.~$R^* \in \ccw(P)$ for some path $P$ such that $\ccw(P)\cap (C\cup C') =\emptyset$. 
    In other words, $P <_v^{\cw} R^*$ for some $v$ and any route in $C\cup C'$ is coherent with $P$ or clockwise from $P$.
    We can assume that the only vertices of $P$ not in $P\cap R^*$ are $\min(P)$ and $\max(P)$.

    By the choice of $P$, there cannot be a route in $C\cap C'$ that is clockwise from $P$ at $v$ (otherwise it would be incoherent with $R^*$). 
    Since $\ccw(P)\subseteq \Rem_{\ccw}(C,C')$, it follows that~$P$ is coherent with all routes in $C\cap C'$.
    Since in the construction of $C_{\max}^{\Rem}(C,C')$ we can never add a route clockwise from $P$ (otherwise it would be incoherent with $R^*$) and any route in $\ccw(P)$ is in $\Rem_{\ccw}(C,C')$, it follows that $P$ is coherent with all routes in~$C_{\max}^{\Rem}(C,C')$.
    By Lemma~\ref{lem.path_is_extendable},~$P$ can be extended to a route $R_P$ that is coherent with all routes in $C_{\max}^{\Rem}(C,C')$. 
    We can assume that  $R_Pv$ is minimal with respect to $\leq_{\In(v)}$, and among those routes containing $R_Pv$, $vR_P$ is maximal with respect to $\leq_{\Out(v)}$. See~\Cref{fig_lem_CmaxRem_maximal} for a schematic illustration. 

\begin{figure}
    \centering
\begin{tikzpicture}

\begin{scope}[xshift=0, yshift=0, scale=0.8]
    \draw[ultra thick,color=blue,dashed] (2.5,0.6) .. controls (2.6, 0.6) and (3.2, 0.6) .. (3.5,0.4);
    \draw[ultra thick,color=black] (3.5,0.4)--(4,0);
    \draw[ultra thick,color=black] (4,0)--(6,0);
    \draw[ultra thick,color=black] (6,0)--(6.45,-0.38);
    \draw[ultra thick,color=blue, dashed] (2,0.6)--(2.5,0.6);
    \draw[ultra thick,color=blue, dashed] (6.5,-0.4) .. controls (6.8, -0.6) and (7.4, -0.6) .. (7.5,-0.6);
    \draw[ultra thick,color=blue, dashed] (7.5,-0.6)--(8.25,-0.6);

    \draw[ultra thick,color=ForestGreen] (2.5,-0.6) .. controls (2.6, -0.6) and (3.6, -0.6) .. (4,0);
    \draw[ultra thick,color=ForestGreen, dashed] (2,-0.6)--(2.5,-0.6);    
    \draw[ultra thick,color=ForestGreen] (6,0) .. controls (6.4, 0.6) and (7.4, 0.6) .. (7.5,0.6);
    \draw[ultra thick,color=ForestGreen, dashed] (7.5,0.6)--(8.25,0.6);

    \node[circle, draw, inner sep=1.6pt, fill,label=below:\scriptsize{$v$}] (s) at (5,0) {};
    
    \node[] (RP) at (1.6,0.6) {\scriptsize \textcolor{blue}{$R_P$}};
    \node[] (R*) at (1.6,-0.6) {\scriptsize \textcolor{ForestGreen}{$R^*$}};
    \node[] (RP) at (8.6,-0.6) {\scriptsize \textcolor{blue}{$R_P$}};
    \node[] (R*) at (8.6,0.6) {\scriptsize \textcolor{ForestGreen}{$R^*$}};
    \node[] (P) at (4.5,0.4) {\scriptsize \textcolor{black}{$P$}};
   
\end{scope}

\begin{scope}[xshift=230, yshift=0, scale=0.8]
    \draw[ultra thick,color=blue,dashed] (0.5,0.6)--(1.5,0.6);
    \draw[ultra thick,color=blue,dashed] (2.5,0.6) .. controls (2.6, 0.6) and (3.2, 0.6) .. (3.5,0.4);
    \draw[ultra thick,color=black] (3.5,0.4)--(4,0);
    \draw[ultra thick,color=black] (4,0)--(6,0);
    \draw[ultra thick,color=black] (6,0)--(6.45,-0.38);
    \draw[ultra thick,color=blue, dashed] (6.5,-0.4) .. controls (6.8, -0.6) and (7.4, -0.6) .. (7.5,-0.6);
    \draw[ultra thick,color=blue, dashed] (7.5,-0.6)--(8.25,-0.6);

    \draw[ultra thick,color=ForestGreen] (2.5,-0.6) .. controls (2.6, -0.6) and (3.6, -0.6) .. (4,0);
    \draw[ultra thick,color=ForestGreen, dashed] (2,-0.6)--(2.5,-0.6);    
    \draw[ultra thick,color=ForestGreen] (6,0) .. controls (6.4, 0.6) and (7.4, 0.6) .. (7.5,0.6);
    \draw[ultra thick,color=ForestGreen, dashed] (7.5,0.6)--(8.25,0.6);

    \draw[ultra thick,color=red] (1.1,1)--(1.5,0.6);
    \draw[ultra thick,color=red] (1.5,0.6)--(2.5,0.6);
    \draw[ultra thick,color=red] (2.5,0.6)--(2.9,0.2); 

    \node[circle, draw, inner sep=1.6pt, fill,label=below:\tiny{$w_1$}] (s) at (1.5,0.6) {};
    \node[circle, draw, inner sep=1.6pt, fill,label=below:\tiny{$w$}] (s) at (2,0.6) {};
    \node[circle, draw, inner sep=1.6pt, fill,label=below:\tiny{$w_2$}] (s) at (2.5,0.6) {};

    \node[circle, draw, inner sep=1.6pt, fill,label=below:\scriptsize{$v$}] (s) at (5,0) {};
    
    \node[] (RP) at (0.1,0.6) {\scriptsize \textcolor{blue}{$R_P$}};
    \node[] (R*) at (1.6,-0.6) {\scriptsize \textcolor{ForestGreen}{$R^*$}};
    \node[] (RP) at (8.6,-0.6) {\scriptsize \textcolor{blue}{$R_P$}};
    \node[] (R*) at (8.6,0.6) {\scriptsize \textcolor{ForestGreen}{$R^*$}};
    \node[] (P) at (4.5,0.4) {\scriptsize \textcolor{black}{$P$}};
    \node[] (P) at (2,1) {\scriptsize \textcolor{red}{$P'$}};
\end{scope}

\end{tikzpicture}
    \caption{Some routes and paths involved in the proof of~\Cref{lem_CmaxRem_maximal}.}
    \label{fig_lem_CmaxRem_maximal}
\end{figure}
    
    Since $R_P$ was not added to $C_{\max}^{\Rem}(C,C')$, it must be in $\Rem_{\ccw}(C,C')$.
    Thus there is some path $P'$ such that $\ccw(P')\cap (C\cup C') = \emptyset$ and $P' <_w^{\cw} R_P$ for some $w$. 
    As in the case with~$P$, we assume that the only vertices of $P'$ not in $P'\cap R_P$ are $\min(P')$ and $\max(P')$.
    Let $w_1$ and $w_2$ respectively be the minimal and maximal elements in $P'\cap R_P$. 
    We consider two cases separately: $w\leq v$ and $v\leq w$.
    If $w\leq v$, then since $R_Pv$ is minimal with respect to $\leq_{\In(v)}$, there must exist a route $\widetilde{R} \in C_{\max}^{\Rem}(C,C')$ 
    incoherent with $P'wR_P$.
    This route~$\widetilde{R}$ satisfies 
    $P'w_1 <_{\In(w_1)} \widetilde{R}w_1 <_{\In(w_1)} R_Pw_1$,
    otherwise $\widetilde{R}$ would be incoherent with~$R_P$.
    Since~$\widetilde{R}$ is coherent with $R_P$, we also have 
    $w_1\widetilde{R} \leq_{\Out(w_1)} w_1R_P <_{\Out(w_1)} w_1P'$.
    It follows that~$\widetilde{R} \in \ccw(P')$, and thus $\widetilde{R} \in \Rem_{\ccw}(C,C')$, which contradicts the fact that $\widetilde{R} \in C_{\max}^{\Rem}(C,C')$.
    A contradiction is similarly reached in the case when $v\leq w$.
\end{proof}

\begin{lemma}\label{lem.CmaxRem_is_ccw_maximal}
    The cliques $C_{\max}^{\Rem}(C,C')$ and $C_{\min}^{\Rem}(C,C')$ are the unique maximal cliques with the following properties. 
    If a route $R \notin \Rem(C,C')$ is not in $C\cap C'$ but coherent with all routes in $C\cap C'$, then 
    \begin{enumerate}
        \item for any $R' \in C_{\max}^{\Rem}(C,C')$ and $v \in R\cap R'$ either $R$ and $R'$ are coherent at $v$ or $R <_v^{\cw} R'$.
        \item for any $R' \in C_{\min}^{\Rem}(C,C')$ and $v \in R\cap R'$ either $R$ and $R'$ are coherent at $v$ or $R' <_v^{\cw} R$.
    \end{enumerate} 
\end{lemma}

\begin{proof}
    The proof is nearly identical to that of Lemma~\ref{lem.CmaxS_is_ccw_maximal}. 
    In the case when $v<v'$, one needs to guarantee the existence of $\widetilde{R} \notin \Rem(C,C')$, which follows from the fact that $R \notin \Rem(C,C')$. 
    In the case when $v \geq v'$, one needs to guarantee that $R'vR \notin \Rem(C,C')$, which follows directly from the fact that neither $R$ nor $R'$ is in $\Rem(C,C')$. 
    We leave the full details to the reader.
\end{proof}

\begin{proposition}\label{prop.CmaxRem}
    The following hold:
    \begin{enumerate}
        \item $C_{\max}^{\Rem}(C,C')$ is the unique maximal clique that is bigger in the order $\leq_{\mathrm{rot}}^{\ccw}$ than all the maximal cliques containing $C\cap C'$ and no routes from $\Rem(C,C')$. 
        \item $C_{\min}^{\Rem}(C,C')$ is the unique maximal clique that is smaller in the order $\leq_{\mathrm{rot}}^{\ccw}$ than all the maximal cliques containing $C\cap C'$ and no routes from $\Rem(C,C')$.
        \qed
    \end{enumerate}
\end{proposition}

\begin{proof}
This is a direct consequence of~\Cref{lem.CmaxRem_is_ccw_maximal} and the Characterization Theorem~\ref{thm.ccw}. 
\end{proof}

\begin{theorem}
\label{thm.MeetAndJoin}
    The maximal cliques $C_{\min}^{\Rem}(C,C')$ and $C_{\max}^{\Rem}(C,C')$ are respectively the meet $C\wedge C'$ and the join $C\vee C'$ in $\scrL_{G,F}$. 
\end{theorem}

\begin{proof}
    We show only the case of the join, as the meet is argued analogously.
    It follows from Proposition~\ref{prop.CmaxRem} that $C\leq C_{\max}^{\Rem}(C,C')$ and $C' \leq C_{\max}^{\Rem}(C,C')$.
    Suppose toward a contradiction that $C_{\max}^{\Rem}(C,C')$ is not the join of $C$ and $C'$.
    Then, there is a maximal clique $C^* < C_{\max}^{\Rem}(C,C')$ such that $C\leq C^*$ and $C'\leq C^*$.
    Furthermore, there is a route $R^*\in C^*$ and $R \in C_{\max}^{\Rem}(C,C')$ such that $R^* <_v^{\cw} R$ for some $v$.
    In particular, $R \in \ccw(R^*)$. 
    Note that since $C\leq C^*$ we have that $\ccw(R^*) \cap C = \emptyset$.
    Similarly, $\ccw(R^*)\cap C' = \emptyset$, and hence $\ccw(R^*)\cap (C\cup C') = \emptyset$.
    It follows that $R \in \Rem_{\ccw}(C,C')$, and so $R\notin C_{\max}^{\Rem}(C,C')$, which is a contradiction.
\end{proof}

As a consequence of~\Cref{prop.CmaxRem} and~\Cref{thm.MeetAndJoin}, $C\vee C'$ can be obtained from either~$C$ or $C'$ by counter-clockwise rotations of routes that are not in $C\cap C'$ and do not result in a route in $\Rem(C,C')$.

\subsection{Semidistributivity}

Recall that a lattice $L$ is \defn{join-semidistributive} if any elements $x$, $y$, and $z$ in $L$ with $z \vee x = z\vee y$ also satisfy $z\vee (x\wedge y) = z\vee x$.
Similarly, $L$ is \defn{meet-semidistributive} if the dual condition holds, i.e. $z\wedge (x\vee y) = z \wedge x$ whenever $z\wedge x = z \wedge y$. 
A lattice is \defn{semidistributive} if it is both meet- and join-semidistributive. 

\begin{lemma}(BEZ Lemma for semidistributivity \cite[Lemma 9-2.6]{Rea16}) \label{lem.BEZ_semilattice}
Let $L$ be a finite lattice and let $x$, $y$, $z$ be elements of $L$. 
Consider the following criteria:
\begin{itemize}
    \item[(i)] If $x$ and $y$ are covered by a common element and $z \vee x = z \vee y$, then
$z \vee (x \wedge y) = z \vee x = z \vee y$. 
    \item[(ii)] If $x$ and $y$ cover a common element and $z\wedge x = z\wedge y$, then $z\wedge (x \vee y) = z\wedge x = z \vee y$.
\end{itemize}
If $L$ satisfies (i) then it is join-semidistributive, if $L$ satisfies (ii) then it is meet-semidistributive, and if $L$ satisfies both (i) and (ii) then it is semidistributive. 
\qed
\end{lemma}

To show that framing lattices are semidistributive, we make use of the following lemma. 

\begin{proposition}
    \label{prop.Z}
    Let $C_1\prec C_2$ be maximal cliques such that $C_2 = C_1\setminus R_1 \cup R_2$. 
    Let $P$ be the path at which $R_1$ and $R_2$ are incoherent, and let  $\widetilde{P}_1$ (resp. $\widetilde{P}_2$) be the path obtained from~$P$ by adding the edge of $R_1$ (resp. $R_2$) incoming to $P$ and the edge of $R_1$ (resp. $R_2$) outgoing from $P$.
    Then for a maximal clique~$C^*$ the following hold:
    \begin{enumerate}
        \item $C_1\vee C^* = C_2\vee C^*$ if and only if there exists a route $R^* \in C^*$ such that $\widetilde{P}_1 <_v^{\cw} R^*$ for some $v \in P$. 
        \item $C_1\wedge C^* = C_2\wedge C^*$ if and only if there exists a route $R^* \in C^*$ such that $R^* <_v^{\cw} \widetilde{P}_2$ for some $v \in P$.  
    \end{enumerate}
\end{proposition}

\begin{proof}
   We prove only part (1) because (2) follows by symmetry, reversing the framing. 

   $(\Leftarrow)$ Let $R^* \in C^*$ be a route satisfying $\widetilde{P}_1 <_v^{\cw} R^*$ at some $v \in P$.
    Since $C_1 \leq C_2$, we have that $C_1\vee C^* \leq C_2\vee C^*$. So, it suffices to show that $C_2\vee C^* \leq C_1\vee C^*$, or equivalently that $C_2 \leq C_1\vee C^*$.
    By~\Cref{thm.ccw}, this follows if every route $R'\in C_2$ is coherent with or clockwise from every route $R \in C_1 \vee C^*$.
    Since we know that this property holds for every route in $C_1$, we only need to prove it for $R'=R_2$. 

    Assume toward a contradiction that there is a route $R \in C_1 \vee C^*$ such that $R <_w^{\cw} R_2$ for some $w$.
    Since $\Top(R_1,R_2)$ and $\Bot(R_1,R_2)$ are in $C_1$, they are also coherent with or clockwise from $R$. 
    Thus, $R$ must be a route weakly in between $R_1$ and $R_2$ at $v$. 
    We distinguish two possible cases: (1) $R=R_1$ and (2) $R$ is in between $R_1$ and $R_2$ at $v$.
    
    Case (1). If $R=R_1$ then $R <_v^{\cw} R^*$, which contradicts the fact that $C^*\leq C_1\vee C^*$.
    
    Case (2). If $R$ is in between $R_1$ and $R_2$ at $v$ then 
    \begin{itemize}
        \item[(i)] $R_1v <_{\scrI(v)} Rv <_{\scrI(v)} R_2v$ and $vR_2 \leq_{\scrO(v)} vR \leq_{\scrO(v)} vR_1$; or 
        \item[(ii)] $R_1v \leq_{\scrI(v)} Rv \leq_{\scrI(v)} R_2v$ and $vR_2 <_{\scrO(v)} vR <_{\scrO(v)} vR_1$.
    \end{itemize}
    We consider only case (i) here as (ii) is symmetric.
    Let $P_R$ be the maximal path in $R\cap R_2$ that contains $v$.
    Let $\widetilde{P}_R$ be the path $P_R$ with the additional edge of $R$ incoming to $P_R$.
    Since $R\in C_1 \vee C^*$ and $C_1 \leq C_1\vee C^*$, then any route in $C_1$ must be coherent with or clockwise from~$R$, and hence also $\widetilde{P}_Rv$. 
    However, a route of $C_1$ cannot be clockwise from $\widetilde{P}_Rv$ as it would then also be clockwise from $\Bot(R_1,R_2)$. 
    Thus $\widetilde{P}_Rv$ is coherent with $C_1$, but by Lemma~\ref{lem.path_is_extendable} it extends to a route in $R' \in C_1$. 
    Moreover $vR'=vR_1$, otherwise $R'$ would be incoherent with $\Bot(R_1,R_2)$ or $R_1$. 
    Therefore, $R'\in C_1$ is in between $R_1$ and $R_2$, which is a contradiction.

    $(\Rightarrow)$ Let $C_1\vee C^* = C_2 \vee C^*$. 
    Assume toward a contradiction that every route of $C^*$ is coherent with or clockwise from $\widetilde{P}_1$. 
    Since every route of $C_1$ is coherent with $\widetilde{P}_1$, any route counter-clockwise from $\widetilde{P}_1$ is in $\Rem_{\ccw}(C_1,C^*)$ by construction.
    Therefore, every route in $C_1\vee C^* = C_{\max}^{\Rem}(C_1,C^*)$ is coherent with or clockwise from $\widetilde{P}_1$.
    Since $C_1\vee C^* = C_2\vee C^*$, every route in $C_2\vee C^*$ is coherent with or clockwise from $\widetilde{P}_1$.
    However, since $R_2\in C_2$ and~$\widetilde{P}_1$ is clockwise from $R_2$, it follows from 
    \Cref{lem.ccwExists} that there exists a route $R'\in C_2\vee C^*$ that is counter-clockwise from $\widetilde{P}_1$. This is a contradiction.
\end{proof}

\begin{theorem}
\label{thm.semidistributivity}
    The framing lattice $\scrL_{G,F}$ is semidistributive. 
\end{theorem}

\begin{proof}
    We provide only the proof of join-semidistributivity below using Lemma~\ref{lem.BEZ_semilattice}(i), as meet-semidistributivity is obtained symmetrically.
    We assume that $x$ and $y$ are covered by a common element with $z\vee x = z\vee y$ and show that $z\vee (x\wedge y) = z\wedge x$.

    The interval $[x\wedge y, x\vee y]$ is a square, pentagon, or hexagon. 
    We only write down the details of the hexagonal case here, as the other two cases are nearly identical.
    Let $x'$ and $y'$ be the elements in the hexagon covering $x\wedge y$, such that $x'\leq y$ and $y'\leq x$ (in the cases of the square and pentagon, we would have $x'=y$ and/or $y'=x$).
    It follows from our assumptions that 
    $z\vee x = z\vee y = z \vee (x\vee y) = z \vee (x' \vee y').$
    It may be helpful for the reader to refer to Figure~\ref{fig.hexagon} with $x:= C_1'$, $y:=C_2'$, $x'=C_2$, $y'=C_1$, and $x\wedge y :=Q$.
    Let $P_x$ be the path at which the two routes involved in the rotation from $x$ to $x\vee y$ are incoherent. 
    In Figure~\ref{fig.hexagon}, these routes are $R_1$ and $R_2b_1R_1$, and $P_x$ is the path from $b_2$ to $c_2$.
    Extending the head and tail of~$P_x$ by an edge along the route of $x$ involved the rotation, we obtain the extended path~$\widetilde{P}_x$ (in Figure~\ref{fig.hexagon} this extension is along $R_1$).  
    Since $z\vee x = z\vee (x\vee y)$,~\Cref{prop.Z} implies that there is a route $R_x \in z$ such that $\widetilde{P}_x <_v^{\cw} R_x$ for some $v \in P_x$. 
    Now, $x\wedge y$ is covered by $x'$ and they differ by routes that are also incoherent exactly at $P_x$ with $\widetilde{P}_x$ contained in the route of $x\wedge y$ involved in the rotation to $x'$.
    Therefore, by~\Cref{prop.Z}, we have that $z\vee x' = z\vee (x\wedge y)$. 
    Similarly, we obtain $z\vee y' = z \vee (x\wedge y)$. 
    We can now simply compute:
    $$ z\vee x = z \vee (x'\vee y') = (z\vee x')\vee y' = z\vee (x\wedge y) \vee y' = z\vee y' = z \vee (x\wedge y).$$
    This finishes our proof. 
\end{proof}

\subsection{\texorpdfstring{$\calH\calH$}{}-lattice property and congruence uniformity}

Our next goal is to show that the framing lattice is congruence uniform. This is implied by the, perhaps less known, $\calH\calH$-lattice property that we now recall.

\begin{definition}[{\cite[Definition 10]{CdPBM04}}, Cf.~\cite{CP22}]
    \label{def.HHlattice}
    Given a lattice $\scrL$, let $E(\scrL)$ denote the set of covering relations of $\scrL$. 
    We say that $\scrL$ is an \defn{$\calH\calH$-lattice} if it is finite, semidistributive, polygonal, and there exist a labeling function
    $$ \ell: E(\scrL) \to \calL$$
    where $\calL$ is a set of labels, and a ranking function
    $$r: \calL \to \mathbb{N}$$
    satisfying the following condition on every polygon $[x,y]$ of $\scrL$.

    Let $x_1$ and $x_2$ denote the two elements covering $x$, and let $y_1$ and $y_2$ denote the two elements covered by $y$, such that $x_1$ and $y_1$ (resp. $x_2$ and $y_2$) belong to the same maximal chain. 
    The labeling $\ell$ and rank function $r$
    must satisfy:
    \begin{itemize}
        \item[(1)] $\ell(x,x_1) = \ell(y_2,y)$ and $\ell(x,x_2) = \ell(y_1,y)$;
        \item[(2)] if $t_1,\ldots, t_k$ is a maximal chain in a polygon, then
        $$ r(t_1), r(t_k) < r(t_2), r(t_{k-1}) < \cdots < r(t_{\frac{k+1}{2}}) \,\, \text{ if $k$ is odd; and} $$
        $$ r(t_1), r(t_k) < r(t_2), r(t_{k-1}) < \cdots < r(t_{\frac{k}{2}}),r(t_{\frac{k}{2}+1}) \,\, \text{ if $k$ is even.} $$
    \end{itemize}
\end{definition}

A lattice is \defn{congruence uniform} if it can be obtained from the one element lattice
by a sequence of \textit{doublings} of intervals \cite{Day77}. 
We refer to \cite{Rea03} for more details on congruence uniform lattices. 

\begin{theorem}\cite[Corollary 1]{CdPBM04}
    Every $\calH\calH$-lattice is congruence uniform. 
\end{theorem}

\begin{theorem}
    \label{thm.HHlattice}
    The framing lattice is an $\calH\calH$-lattice, and hence congruence uniform. 
\end{theorem}

\begin{proof}
    We have shown that framing lattices are finite, polygonal, and semidsitributive, so it remains to find a labeling $\ell$ and ranking function $r$ satisfying the necessary conditions in Definition~\ref{def.HHlattice} for each polygon.

    If $(C, C')$ is an edge in $\scrL_{G,F}$ with $C\prec C'$, then $C' = C\setminus R \cup R'$ for some routes $R$ and $R'$ which are incoherent at a path $P$.
    Define $\ell$ by labeling each edge $(C,C')$ by the associated path $P$, and define $r(P) = \max(P)-\min(P)$. 

    By inspection of the polygons in Figures~\ref{fig.square}, \ref{fig.pentagon}, and \ref{fig.hexagon}, condition (1) of Definition~\ref{def.HHlattice} holds for $\ell$. 
    As for the ranking $r$, observe that the maximal chains in the polygons can be of length $2$ or $3$. 
    There is nothing to check when the length is $2$. 
    When the length is $3$, if the chain is labeled by ($t_1$, $t_2$, $t_3$), then we need to have $r(t_1) < r(t_2) > r(t_3)$. 
    In all cases, we have $r(t_1) = \max(P_1)-\min(P_1)$ and $r(t_3) = \max(P_2) - \min(P_2)$ or the other way around.
    Furthermore, $r(t_2) = \max(P_1)-\min(P_2)$. 
    As 
    $\min(P_2) \leq \max(P_2) < \min(P_1) \leq \max(P_1)$,
    the desired inequalities $r(t_1) < r(t_2) > r(t_3)$ follow.
\end{proof}

\section{The framing core label order and noncrossing partitions}\label{sec_core_label_order}

The key ingredient in our proof of semistributivity and the $\calH\calH$-lattice property (and consequently, congruence uniformity) was the labeling of the edges of the lattice by the path~$P$ at which the two routes involved in the rotation are incoherent, as well as the labeling by the extended paths $\widetilde{P}_1$ and $\widetilde{P}_2$ from~\Cref{prop.Z}.

In this section, we define these path and extended path labelings explicitly, and show that they can be used to describe the join- and the meet irreducibles of the framing lattice. We also describe the framing core label order which generalizes the poset of noncrossing partitions in this context.  

\subsection{The path- and extended path labeling}
Let $C_1\prec C_2$ be maximal cliques such that $C_2 = C_1\setminus R_1 \cup R_2$. 
Let $P$ be the path at which $R_1$ and $R_2$ are incoherent, and let  $\widetilde{P}_1$ (resp. $\widetilde{P}_2$) be the path obtained from~$P$ by adding the edge of $R_1$ (resp $R_2$) incoming to~$P$ and the edge of $R_1$ (resp $R_2$) outgoing from $P$. An illustration of these paths is shown in~\Cref{fig_extended_paths}.

\begin{figure}[h]
    \centering
    \begin{tikzpicture}

\begin{scope}[xshift=0, yshift=0, scale=0.8]

    \draw[ultra thick,color=red] (4,0)--(5.5,0);

    \draw[ultra thick,color=black!30, dashed] (2.5,0.6) .. controls (3.1, 0.6) and (3.6, 0.6) .. (4,0);

    \draw[ultra thick,color=black!30, dashed] (2.5,-0.6) .. controls (3.1, -0.6) and (3.6, -0.6) .. (4,0);
    
    \draw[ultra thick,color=black!30, dashed] (5.5,0) .. controls (5.9, 0.6) and (6.4, 0.6) .. (7,0.6);

    \draw[ultra thick,color=black!30, dashed] (5.5,0) .. controls (5.9, -0.6) and (6.4, -0.6) .. (7,-0.6);

    \node[circle, draw, inner sep=1.3pt, fill] (s) at (4,0) {};
    \node[circle, draw, inner sep=1.3pt, fill] (s) at (5.5,0) {};
    
    \node[] (R1) at (2.2,0.6) {\scriptsize \textcolor{black}{$R_1$}};
    \node[] (R2) at (2.2,-0.6) {\scriptsize \textcolor{black}{$R_2$}};
    \node[] (R1) at (7.3,-0.6) {\scriptsize \textcolor{black}{$R_1$}};
    \node[] (R2) at (7.3,0.6) {\scriptsize \textcolor{black}{$R_2$}};
    \node[] (P) at (4.75,0.4) {\scriptsize \textcolor{red}{$P$}};
\end{scope}

\begin{scope}[xshift=160, yshift=0, scale=0.8]

    \draw[ultra thick,color=red] (3.6,0.4)--(4,0);
    \draw[ultra thick,color=red] (4,0)--(5.5,0);
    \draw[ultra thick,color=red] (5.5,0)--(5.9,-0.4);

    \draw[ultra thick,color=black!30, dashed] (2.5,0.6) .. controls (3.1, 0.6) and (3.4, 0.6) .. (3.6,0.4);

    \draw[ultra thick,color=black!30, dashed] (2.5,-0.6) .. controls (3.1, -0.6) and (3.6, -0.6) .. (4,0);
    
    \draw[ultra thick,color=black!30, dashed] (5.5,0) .. controls (5.9, 0.6) and (6.4, 0.6) .. (7,0.6);

    \draw[ultra thick,color=black!30, dashed] (5.9,-0.4) .. controls (6.1, -0.6) and (6.4, -0.6) .. (7,-0.6);

    \node[circle, draw, inner sep=1.3pt, fill] (s) at (4,0) {};
    \node[circle, draw, inner sep=1.3pt, fill] (s) at (5.5,0) {};
    \node[circle, draw, inner sep=1.3pt, fill] (s) at (3.6,0.4) {};
    \node[circle, draw, inner sep=1.3pt, fill] (s) at (5.9,-0.4) {};
    
    \node[] (R1) at (2.2,0.6) {\scriptsize \textcolor{black}{$R_1$}};
    \node[] (R2) at (2.2,-0.6) {\scriptsize \textcolor{black}{$R_2$}};
    \node[] (R1) at (7.3,-0.6) {\scriptsize \textcolor{black}{$R_1$}};
    \node[] (R2) at (7.3,0.6) {\scriptsize \textcolor{black}{$R_2$}};
    \node[] (P) at (4.75,0.4) {\scriptsize \textcolor{red}{$\widetilde{P}_1$}};
\end{scope}

\begin{scope}[xshift=320, yshift=0, scale=0.8]

    \draw[ultra thick,color=red] (3.6,-0.4)--(4,0);
    \draw[ultra thick,color=red] (4,0)--(5.5,0);
    \draw[ultra thick,color=red] (5.5,0)--(5.9,0.4);

    \draw[ultra thick,color=black!30, dashed] (2.5,-0.6) .. controls (3.1, -0.6) and (3.4, -0.6) .. (3.6,-0.4);

    \draw[ultra thick,color=black!30, dashed] (2.5,0.6) .. controls (3.1, 0.6) and (3.6, 0.6) .. (4,0);
    
    \draw[ultra thick,color=black!30, dashed] (5.5,0) .. controls (5.9, -0.6) and (6.4, -0.6) .. (7,-0.6);

    \draw[ultra thick,color=black!30, dashed] (5.9,0.4) .. controls (6.1, 0.6) and (6.4, 0.6) .. (7,0.6);

    \node[circle, draw, inner sep=1.3pt, fill] (s) at (4,0) {};
    \node[circle, draw, inner sep=1.3pt, fill] (s) at (5.5,0) {};
    \node[circle, draw, inner sep=1.3pt, fill] (s) at (3.6,-0.4) {};
    \node[circle, draw, inner sep=1.3pt, fill] (s) at (5.9,0.4) {};
    
    \node[] (R1) at (2.2,0.6) {\scriptsize \textcolor{black}{$R_1$}};
    \node[] (R2) at (2.2,-0.6) {\scriptsize \textcolor{black}{$R_2$}};
    \node[] (R1) at (7.3,-0.6) {\scriptsize \textcolor{black}{$R_1$}};
    \node[] (R2) at (7.3,0.6) {\scriptsize \textcolor{black}{$R_2$}};
    \node[] (P) at (4.75,0.4) {\scriptsize \textcolor{red}{$\widetilde{P}_2$}};
\end{scope}
\end{tikzpicture}
    \caption{The path $P$, and the extended paths $\widetilde{P}_1$ and $\widetilde{P_2}$.}
    \label{fig_extended_paths}
\end{figure}

We define three different edge labelings of the Hasse diagram of the framing lattice: 

\begin{enumerate}
    \item the \defn{path labeling} 
    $\ell(C_1,C_2):=P$,
    \item the \defn{$\cw$-extended path labeling}
    $\widetilde{\ell}_1(C_1,C_2):=\widetilde{P}_1$, and
    \item the \defn{$\ccw$-extended path labeling}
    $\widetilde{\ell}_2(C_1,C_2):=\widetilde{P}_2$.
\end{enumerate}

Since we will be mainly using $\widetilde{\ell}_1$ instead of $\widetilde{\ell}_2$, we often call $\widetilde{\ell}_1$ the \defn{extended path labeling} for simplicity, and make a distinction adding $\cw$ or $\ccw$ when necessary.

As we observed already in the proofs of semidistributivity and of the $\calH\calH$-lattice property, in every polygon of the framing lattice opposite edges containing the minimal or maximal element of the polygon have the same path and extended path labelings. This is stated formally in the following lemma, and is illustrated in~\Cref{fig_opposite_polygon_labels}.

\begin{lemma}\label{lem_opposite_polygon_labels}
Given a polygon $[x,y]$ of the framing lattice $\scrL_{G,F}$,
let $x_1$ and $x_2$ denote the two elements covering $x$, and let $y_1$ and $y_2$ denote the two elements covered by $y$, such that $x_1$ and $y_1$ (resp. $x_2$ and $y_2$) belong to the same maximal chain. Then, 
\begin{enumerate}
    \item $\ell(x,x_1) = \ell(y_2,y)$ and $\ell(x,x_2) = \ell(y_1,y)$,
    \item $\widetilde{\ell}_1(x,x_1) = \widetilde{\ell}_1(y_2,y)$ and $\widetilde{\ell}_1(x,x_2) = \widetilde{\ell}_1(y_1,y)$,
    \item $\widetilde{\ell}_2(x,x_1) = \widetilde{\ell}_2(y_2,y)$ and $\widetilde{\ell}_2(x,x_2) = \widetilde{\ell}_2(y_1,y)$.
\end{enumerate}
That is, opposite edges containing the minimal or maximal element of any polygon have the same path and extended path labelings.
\end{lemma}
\begin{proof}
    The proof follows by inspection of all the polygonal cases described in the proof of~\Cref{prop.intervalJoin}. The corresponding labels are shown in~\Cref{fig_opposite_polygon_labels} for the generic square, pentagon, and hexagonal cases from~\Cref{fig.square,fig.pentagon,fig.hexagon}.  
\end{proof}

\begin{figure}[h]
    \centering

    \caption{Opposite edges containing the minimal or maximal element of a polygon have the same path and extended path labelings.}
    \label{fig_opposite_polygon_labels}
\end{figure}

Explicit examples of the extended path labeling $\widetilde{\ell}_1$ are shown in~\Cref{fig_extended_path_labeling_examples}, corresponding to the examples of framing lattices in~\Cref{fig_weakorder_framinglattice,fig.framing_poset_example,fig.caracol3}. 

\begin{figure}
    \centering
    \scalebox{0.8}{
    %
    }
    \caption{The $\cw$-extended path labeling $\widetilde{\ell}_1$ for the running examples in~\Cref{fig_weakorder_framinglattice,fig.framing_poset_example,fig.caracol3}. The exceptional routes in each maximal clique have been suppressed for clarity.}
    \label{fig_extended_path_labeling_examples}
\end{figure}

The following definitions will help us to characterize the labels that appear in the $\cw$-extended path labeling and the $\ccw$-extended path labeling. 

Let $(G,F)$ be a framed graph and $\widetilde{P}$ be a path containing at least two edges. Let $e_1$ and~$e_2$ be the initial and final edges of $\widetilde{P}$, respectively. Consider the vertices $w_1=\max(e_1)$ and~$w_2=\min(e_2)$ in $\widetilde{P}$.    
We say that $\widetilde{P}$ is a \defn{$\cw$-extended path} (resp. \defn{$\ccw$-extended path}) if the following two properties hold:
\begin{enumerate}
    \item $e_1$ is not maximal (resp. minimal) with respect to the order $\In(w_1)$, and
    \item $e_2$ is not minimal (resp. maximal) with respect to the order $\Out(w_2)$.
\end{enumerate}

\begin{figure}[h]
    \centering
    \begin{tikzpicture}

\begin{scope}[xshift=105, yshift=35, scale=0.8]
    \node[] (k) at (0,0) {\small \text{cw-extended paths}};
\end{scope}

\begin{scope}[xshift=0, yshift=0, scale=0.8]
    \draw[ultra thick,color=red] (3.2,0.6) .. controls (3.3, 0.6) and (3.8, 0.6) ..(4,0);
    
    \draw[ultra thick,color=red] (4,0)--(5.5,0);
    \draw[ultra thick,color=red] (5.5,0) .. controls (5.7, -0.6) and (6.2, -0.6) ..(6.3,-0.6);

    \draw[ultra thick,color=black!30, dashed] (2.5,0.6)-- (3.2,0.6);
    \draw[ultra thick,color=black!30, dashed] (2.5,0)--(4,0); 

    \draw[ultra thick,color=black!30, dashed] (2.5,-0.6) .. controls (3.1, -0.6) and (3.6, -0.6) .. (4,0);
    
    \draw[ultra thick,color=black!30, dashed] (5.5,0) .. controls (5.9, 0.6) and (6.4, 0.6) .. (7,0.6);
    
    \draw[ultra thick,color=black!30, dashed] (6.3,-0.6)--(7,-0.6);   

    \node[circle, draw, inner sep=1.3pt, fill] (s) at (4,0) {};
    \node[circle, draw, inner sep=1.3pt, fill] (s) at (5.5,0) {};
    \node[circle, draw, inner sep=1.3pt, fill] (s) at (3.2,0.6) {};    
    \node[circle, draw, inner sep=1.3pt, fill] (s) at (6.3,-0.6) {};
    
    \node[] (R1) at (4.2,0.2) {\tiny \textcolor{black}{$w_1$}};
    \node[] (R1) at (3.8,0.6) {\tiny \textcolor{black}{$e_1$}};
    \node[] (R2) at (5.25,0.2) {\tiny\textcolor{black}{$w_2$}};
    \node[] (R1) at (6.0,-0.3) {\tiny \textcolor{black}{$e_2$}};
\end{scope}

\begin{scope}[xshift=0, yshift=-50, scale=0.8]
    \draw[ultra thick,color=black!30, dashed] (2.5,0.6) .. controls (3.1, 0.6) and (3.6, 0.6) .. (4,0);

    \draw[ultra thick,color=red] (4,0)--(5.5,0);
    \draw[ultra thick,color=red] (5.5,0) .. controls (5.7, -0.6) and (6.2, -0.6) ..(6.3,-0.6);


    \draw[ultra thick,color=black!30, dashed] (2.5,-0.6) .. controls (3.1, -0.6) and (3.6, -0.6) .. (4,0);
    
    \draw[ultra thick,color=black!30, dashed] (5.5,0) .. controls (5.9, 0.6) and (6.4, 0.6) .. (7,0.6);

    \draw[ultra thick,color=black!30, dashed] (6.3,-0.6)--(7,-0.6);

    \draw[ultra thick,color=red] (3.2,0)--(4,0);
    \draw[ultra thick,color=black!30, dashed] (2.5,0)--(3.2,0);    

    \node[circle, draw, inner sep=1.3pt, fill] (s) at (4,0) {};
    \node[circle, draw, inner sep=1.3pt, fill] (s) at (5.5,0) {};
    \node[circle, draw, inner sep=1.3pt, fill] (s) at (3.2,0) {};    
    \node[circle, draw, inner sep=1.3pt, fill] (s) at (6.3,-0.6) {};
    
    \node[] (R1) at (4.2,0.2) {\tiny \textcolor{black}{$w_1$}};
    \node[] (R1) at (3.6,0.2) {\tiny \textcolor{black}{$e_1$}};
    \node[] (R2) at (5.25,0.2) {\tiny\textcolor{black}{$w_2$}};
    \node[] (R1) at (6.0,-0.3) {\tiny \textcolor{black}{$e_2$}};
\end{scope}

\begin{scope}[xshift=265, yshift=35, scale=0.8]
    \node[] (k) at (0,0) {\small \text{ccw-extended paths}};
\end{scope}

\begin{scope}[xshift=160, yshift=0, scale=0.8]
    \draw[ultra thick,color=red] (3.2,-0.6) .. controls (3.3, -0.6) and (3.8, -0.6) ..(4,0);
    
    \draw[ultra thick,color=red] (4,0)--(5.5,0);
    \draw[ultra thick,color=red] (5.5,0) .. controls (5.7, 0.6) and (6.2, 0.6) ..(6.3,0.6);

    \draw[ultra thick,color=black!30, dashed] (2.5,-0.6)-- (3.2,-0.6);

    \draw[ultra thick,color=black!30, dashed] (2.5,0.6) .. controls (3.1, 0.6) and (3.6, 0.6) .. (4,0);
    
    \draw[ultra thick,color=black!30, dashed] (5.5,0) .. controls (5.9, -0.6) and (6.4, -0.6) .. (7,-0.6);
    
    \draw[ultra thick,color=black!30, dashed] (6.3,0.6)--(7,0.6);

    \draw[ultra thick,color=black!30, dashed] (5.5,0)--(7,0);    

    \node[circle, draw, inner sep=1.3pt, fill] (s) at (4,0) {};
    \node[circle, draw, inner sep=1.3pt, fill] (s) at (5.5,0) {};
    \node[circle, draw, inner sep=1.3pt, fill] (s) at (3.2,-0.6) {};    
    \node[circle, draw, inner sep=1.3pt, fill] (s) at (6.3,0.6) {};
    
    \node[] (R1) at (4.2,0.2) {\tiny \textcolor{black}{$w_1$}};
    \node[] (R1) at (3.55,-0.3) {\tiny \textcolor{black}{$e_1$}};
    \node[] (R2) at (5.25,0.2) {\tiny\textcolor{black}{$w_2$}};
    \node[] (R1) at (5.65,0.6) {\tiny \textcolor{black}{$e_2$}};
\end{scope}

\begin{scope}[xshift=160, yshift=-50, scale=0.8]
    \draw[ultra thick,color=red] (3.2,-0.6) .. controls (3.3, -0.6) and (3.8, -0.6) ..(4,0);
    
    \draw[ultra thick,color=red] (4,0)--(5.5,0);
    \draw[ultra thick,color=red] (5.5,0)--(6.3,0);

    \draw[ultra thick,color=black!30, dashed] (2.5,-0.6)-- (3.2,-0.6);

    \draw[ultra thick,color=black!30, dashed] (2.5,0.6) .. controls (3.1, 0.6) and (3.6, 0.6) .. (4,0);
    
    \draw[ultra thick,color=black!30, dashed] (5.5,0) .. controls (5.9, -0.6) and (6.4, -0.6) .. (7,-0.6);
    
    \draw[ultra thick,color=black!30, dashed] (5.5,0) .. controls (5.9, 0.6) and (6.4, 0.6) .. (7,0.6);

    \draw[ultra thick,color=black!30, dashed] (6.3,0)--(7,0);    

    \node[circle, draw, inner sep=1.3pt, fill] (s) at (4,0) {};
    \node[circle, draw, inner sep=1.3pt, fill] (s) at (5.5,0) {};
    \node[circle, draw, inner sep=1.3pt, fill] (s) at (3.2,-0.6) {};    
    \node[circle, draw, inner sep=1.3pt, fill] (s) at (6.3,0) {};
    
    \node[] (R1) at (4.2,0.2) {\tiny \textcolor{black}{$w_1$}};
    \node[] (R1) at (3.55,-0.3) {\tiny \textcolor{black}{$e_1$}};
    \node[] (R2) at (5.25,0.2) {\tiny\textcolor{black}{$w_2$}};
    \node[] (R1) at (5.97,0.2) {\tiny \textcolor{black}{$e_2$}};
\end{scope}

\begin{scope}[xshift=425, yshift=35, scale=0.8]
    \node[] (k) at (0,0) {\small \text{Neither}};
\end{scope}

\begin{scope}[xshift=320, yshift=0, scale=0.8]
    \draw[ultra thick,color=red] (3.2,0.6) .. controls (3.3, 0.6) and (3.8, 0.6) ..(4,0);
    
    \draw[ultra thick,color=red] (4,0)--(5.5,0);
    \draw[ultra thick,color=red] (5.5,0) .. controls (5.7, 0.6) and (6.2, 0.6) ..(6.3,0.6);

    \draw[ultra thick,color=black!30, dashed] (2.5,0.6)-- (3.2,0.6);
    \draw[ultra thick,color=black!30, dashed] (2.5,0)--(4,0); 

    \draw[ultra thick,color=black!30, dashed] (2.5,-0.6) .. controls (3.1, -0.6) and (3.6, -0.6) .. (4,0);
    
    \draw[ultra thick,color=black!30, dashed] (5.5,0) .. controls (5.9, -0.6) and (6.4, -0.6) .. (7,-0.6);
    
    \draw[ultra thick,color=black!30, dashed] (6.3,0.6)--(7,0.6);

    \draw[ultra thick,color=black!30, dashed] (2.5,0)--(4,0);    

    \node[circle, draw, inner sep=1.3pt, fill] (s) at (4,0) {};
    \node[circle, draw, inner sep=1.3pt, fill] (s) at (5.5,0) {};
    \node[circle, draw, inner sep=1.3pt, fill] (s) at (3.2,0.6) {};    
    \node[circle, draw, inner sep=1.3pt, fill] (s) at (6.3,0.6) {};
    
    \node[] (R1) at (4.2,0.2) {\tiny \textcolor{black}{$w_1$}};
    \node[] (R1) at (3.9,0.6) {\tiny \textcolor{black}{$e_1$}};
    \node[] (R2) at (5.25,0.2) {\tiny\textcolor{black}{$w_2$}};
    \node[] (R1) at (5.6,0.6) {\tiny \textcolor{black}{$e_2$}};
\end{scope}

\begin{scope}[xshift=320, yshift=-50, scale=0.8]
    \draw[ultra thick,color=red] (3.2,-0.6) .. controls (3.3, -0.6) and (3.8, -0.6) ..(4,0);
    
    \draw[ultra thick,color=red] (4,0)--(5.5,0);
    \draw[ultra thick,color=red] (5.5,0) .. controls (5.7, -0.6) and (6.2, -0.6) ..(6.3,-0.6);

    \draw[ultra thick,color=black!30, dashed] (2.5,-0.6)-- (3.2,-0.6);

    \draw[ultra thick,color=black!30, dashed] (2.5,0.6) .. controls (3.1, 0.6) and (3.6, 0.6) .. (4,0);
    
    \draw[ultra thick,color=black!30, dashed] (5.5,0) .. controls (5.9, 0.6) and (6.4, 0.6) .. (7,0.6);
    
    \draw[ultra thick,color=black!30, dashed] (6.3,-0.6)--(7,-0.6);

    \draw[ultra thick,color=black!30, dashed] (2.5,0)--(4,0);    

    \node[circle, draw, inner sep=1.3pt, fill] (s) at (4,0) {};
    \node[circle, draw, inner sep=1.3pt, fill] (s) at (5.5,0) {};
    \node[circle, draw, inner sep=1.3pt, fill] (s) at (3.2,-0.6) {};    
    \node[circle, draw, inner sep=1.3pt, fill] (s) at (6.3,-0.6) {};
    
    \node[] (R1) at (4.2,0.2) {\tiny \textcolor{black}{$w_1$}};
    \node[] (R1) at (3.55,-0.3) {\tiny \textcolor{black}{$e_1$}};
    \node[] (R2) at (5.25,0.2) {\tiny\textcolor{black}{$w_2$}};
    \node[] (R1) at (6.0,-0.3) {\tiny \textcolor{black}{$e_2$}};
\end{scope}
\end{tikzpicture}
    \caption{Visual representations of cw-extended paths and ccw-extended paths obtained from a path from $w_1$ to $w_2$.}
    \label{fig_extended_paths_definition}
\end{figure}

In other words, when the framing is induced by the drawing of the graph from top to bottom, we require that there is an edge entering $w_1$ below (resp. above) $e_1$ and an edge exiting~$w_2$ above (resp. below) $e_2$. We use the term clockwise $\cw$ (resp. counterclockwise $\ccw$) for this reason. Examples are illustrated in~\Cref{fig_extended_paths_definition}.
The following lemma is straightforward by definition. 

\begin{lemma}
    For two maximal cliques $C_1 \prec C_2$ we have
    \begin{enumerate}
        \item $\widetilde{\ell}_1(C_1,C_2)$ is an $\cw$-extended path.
        \item $\widetilde{\ell}_2(C_1,C_2)$ is an $\ccw$-extended path.
    \end{enumerate}
\end{lemma}

In fact, we will see in~\Cref{cor_extendedpath_labels} that every $\cw$-extended path (resp. $\ccw$-path) is the label of an edge in the $\cw$-extended path labeling $\ell_1$ (resp. $\ccw$-extended path labeling $\ell_2$). 
Before doing that we provide a simple bijection between $\cw$-extended paths and $\ccw$-extended paths, which later on will translate into a bijection between meet irreducible and join irreducible elements of the framing lattice. 

Let $\widetilde{P}$ be an $\cw$-extended path (resp. $\ccw$-extended path) and $e_1,e_2$ and $w_1,w_2$ as above.
We define $\mapccw(\widetilde{P})$ (resp. $\mapcw(\widetilde{P})$) as the path obtained from $\widetilde{P}$ by replacing the incoming edge $e_1$ at $w_1$ by its successor (resp. predecessor) in the order $\In(w_1)$, and the outgoing edge~$e_2$ at $w_2$ by its predecessor (resp. successor) in the order $\Out(w_2)$; they exist by definition of extended $cw$-paths (resp. $\ccw$-extended paths).
The following lemma is also straight forward from the definitions.

\begin{lemma}
The following hold:
    \begin{enumerate}
        \item The map $\mapccw$ is a bijection from $\cw$-extended paths to $\ccw$-extended paths.
        \item The map $\mapcw$ is a bijection from $\ccw$-extended paths to $\cw$-extended paths.
    \end{enumerate}
Moreover, $\mapccw$ and $\mapcw$ are inverses of each other.
\end{lemma}

We call $\mapcw$ the \defn{$\cw$-map} and 
$\mapccw$ the \defn{$\ccw$-map}.

\subsection{Join irreducibles and meet irreducibles}
Now that we have introduced the preliminary concepts, we are ready to give a complete and very simple understanding of the join irreducible and meet irreducible elements of the framing lattice. 

\begin{theorem}\label{thm_join_meet_irreducibles}
    The framing lattice $\scrL_{G,F}$ has the following properties:
    \begin{enumerate}
        \item join irreducible elements of $\scrL_{G,F}$ are in bijection with $\ccw$-extended paths of $(G,F)$. 
        \item meet irreducible elements of $\scrL_{G,F}$ are in bijection with $\cw$-extended paths of $(G,F)$. 
    \end{enumerate}
\end{theorem}

\subsubsection{From $\ccw$-extended paths to join irreducible elements}
We will prove~\Cref{thm_join_meet_irreducibles} by providing explicit bijections. 
For a path $P$ we consider the following two routes.
\begin{align*}
    R_{P}^{\cw}: & 
    \text{ clockwise-most route containing } P \\
    R_{P}^{\ccw}: & 
    \text{ counterclockwise-most route containing } P
\end{align*}

\begin{lemma}[Join irreducibles]
\label{lem_join_irreducibles}
    If $\widetilde{P}_2$ is an $\ccw$-extended path of $(G,F)$, then: 
    \begin{enumerate}
        \item $R_{\widetilde{P}_2}^{\cw} \in C_{\min}( \widetilde{P}_2 )$. 
        \item $C_{\min}( \widetilde{P}_2 ) = 
        C_{\min}( R_{\widetilde{P}_2}^{\cw} )$.
        \item $R_{\widetilde{P}_2}^{\cw}$ is the unique route that can be rotated down (i.e. in $\cw$ direction) in $C_{\min}( \widetilde{P}_2 )$.
        \item $C_{\min}( \widetilde{P}_2 )$ is a join irreducible element of the framing lattice $\scrL_{G,F}$. 
    \end{enumerate}
\end{lemma}

\begin{proof}
    (1) Let $R= R_{\widetilde{P}_2}^{\cw}$. 
    We will show that $R$ is coherent with all the routes in $C_{\min}(\widetilde{P}_2)$, which implies that $R\in C_{\min}(\widetilde{P}_2)$.
    Towards a contradiction, assume that there is a route $R' \in C_{\min}(\widetilde{P}_2)$ that is incoherent with $R$ at some vertex $v$.      
    By construction $R$ is coherent with~$\widetilde{P}_2$, and~\Cref{lem.CmaxS_is_ccw_maximal} implies that $R' <_v^{\cw} R$.

    Let $P_v$ be the maximal path containing $v$ at which $R$ and $R'$ are incoherent. If $P_v$ is contained in the interior of $\widetilde{P}_2$ then $R'$ would be incoherent with $\widetilde{P}_2$, which is a contradiction. 
    Therefore, one of the two following cases holds:
    (i) $\min(P_v)\leq \min(\widetilde{P}_2)$ or 
    (ii) $\max(P_v)\geq \max(\widetilde{P}_2)$.  
    In the first case we have that $\widetilde{P}_2\subseteq R'vR$, contradicting that $R$ is the clockwise-most route containing $\widetilde{P}_2$. In the second case we have that $\widetilde{P}_2\subseteq RvR'$, also contradicting that $R$ is the clockwise-most route containing $\widetilde{P}_2$.  

    (2) By~\Cref{cor.interval}, 
    the set of maximal cliques whose routes are coherent with a set of coherent paths $S$ is equal to the interval $I_S:=[C_{\min}(S),C_{\max}(S)]$. Let $S_1=\{\widetilde{P}_2\}$ and~$S_2=\{R_{\widetilde{P}_2}^{\cw}\}$. If all the routes of a maximal clique $C$ are coherent with $S_2$ then they are also coherent with $S_1$, therefore $I_{S_2}\subseteq I_{S_1}$. In particular, $C_{\min}(S_2)\geq C_{\min}(S_1)$.
    
    On the other hand, we proved in (1) that $S_2\subseteq C_{\min}(S_1)$, which implies that the maximal clique $C_{\min}(S_1)\in I_{S_2}$. So, $C_{\min}(S_2)\leq C_{\min}(S_1)$. 
    Combining both inequalities we get $C_{\min}(S_2)= C_{\min}(S_1)$ as desired.

    (3) Let $R= R_{\widetilde{P}_2}^{\cw}$ as before and $C:=C_{\min}(\widetilde{P}_2)$.
    By parts (1) and (2), we know that $R\in C$ and $C=C_{\min}(R)$. 
    Furthermore, for $\widetilde{P}_1:=\mapcw(\widetilde{P}_2)$ and an interior vertex $v$ of $\widetilde{P}_2$
    we have that $\widetilde{P}_1<_v^{\cw} R$. 
    Therefore, by the existence of a $\cw$ rotation~\Cref{lem.rotationExists_cw_version} we know that there exist a route $R_1$ in $C$ that can be rotated down in $\cw$ direction, producing a maximal clique $C_2=(C\setminus R_1) \cup R_2$.
    If $R_1\neq R$, then $R\in C_2$ contradicting that $C$ is the smallest maximal clique containing $R$. Therefore $R_1=R$, and $R$ is the unique route of $C$ that can be rotated in $\cw$ direction. 

    (4) Since $C=C_{\min}(\widetilde{P}_2)$ has exactly one route that can be rotated down (in $\cw$ direction), then $C$ is a join irreducible element.
\end{proof}

\begin{lemma}[Injectivity]\label{lem_join_injectivity}
    Let $\widetilde{P}_2$ and $\widetilde{P}_2'$ be $\ccw$-extended paths. 
    If $C_{\min}(\widetilde{P}_2)=C_{\min}(\widetilde{P}_2')$ then $\widetilde{P}_2=\widetilde{P}_2'$.
\end{lemma}

\begin{proof}
    Let $C=C_{\min}(\widetilde{P}_2)$ and $C'=C_{\min}(\widetilde{P}_2')$. 
    By~\Cref{lem_join_irreducibles}, $\widetilde{P}_2^\cw\in C$ (resp. ${\widetilde{P}_2}^{' \cw}\in C'$) is the unique route that can be rotated down in $\cw$ direction in $C$ (resp. $C'$). Our assumption~$C=C'$ then implies that 
    \[
    \widetilde{P}_2^\cw={\widetilde{P}_2}^{' \cw}.
    \]  
    Moreover, we can recover $\widetilde{P}_2$ from $\widetilde{P}_2^\cw$ as follows: 
    Let $w_1$ be the first vertex in $\widetilde{P}_2^\cw$ whose incoming edge $e_1$ is not the first edge in the order $\In(w_1)$, and let $w_2$ be the last vertex whose outgoing edge $e_2$ is not the last edge in the order $\Out(w_2)$. Then $\widetilde{P}_2=e_1\ w_1\widetilde{P}_2^\cw w_2 \ e_2$.
    Furthermore, We can recover the path $\widetilde{P}_2'$ from ${\widetilde{P}_2}^{' \cw}$ in the same way. 
    Since $ \widetilde{P}_2^\cw={\widetilde{P}_2}^{' \cw}$ then~$\widetilde{P}_2=\widetilde{P}_2'$ as desired.
\end{proof}

\subsubsection{From $\cw$-extended paths to meet irreducible elements}
   
The following two lemmas follow from~\Cref{lem_join_irreducibles} and~\Cref{lem_join_injectivity} by symmetry.

\begin{lemma}[Meet irreducibles]
\label{lem_meet_irreducibles}
    If $\widetilde{P}_1$ is an $\cw$-extended path of $(G,F)$, then: 
    \begin{enumerate}
        \item $R_{\widetilde{P}_1}^{\ccw} \in C_{\max}( \widetilde{P}_1 )$. 
        \item $C_{\max}( \widetilde{P}_1 ) = 
        C_{\max}( R_{\widetilde{P}_1}^{\ccw} )$.
        \item $R_{\widetilde{P}_1}^{\ccw}$ is the unique route that can be rotated up (i.e. in $\ccw$ direction) in $C_{\max}( \widetilde{P}_1 )$.
        \item $C_{\max}( \widetilde{P}_1 )$ is a meet irreducible element of the framing lattice $\scrL_{G,F}$. 
    \end{enumerate}
\end{lemma}

\begin{lemma}[Injectivity]\label{lem_meet_injectivity}
    Let $\widetilde{P}_1$ and $\widetilde{P}_1'$ be $\cw$-extended paths. 
    If $C_{\max}(\widetilde{P}_1)=C_{\max}(\widetilde{P}_1')$ then $\widetilde{P}_1=\widetilde{P}_1'$.
\end{lemma}

\subsubsection{From irreducible elements to extended paths}

In the previous two sections we introduced injective maps 
\[
\widetilde{P}_2\rightarrow C_{\min}(\widetilde{P}_2)
\] 
from $\ccw$-extended paths to join irreducible elements, and 
\[
\widetilde{P}_1\rightarrow C_{\max}(\widetilde{P}_1)
\] 
from $\cw$-extended paths to meet irreducible elements.
Towards the proof of~\Cref{thm_join_meet_irreducibles} it only remains to show that these two maps are surjective. We do this by describing their inverse maps.

\begin{lemma}[Surjectivity]\label{lem_meet_join_surjectivity}
    Let $C_1\prec C_2$ be two maximal cliques with $C_2=(C_1\setminus R_1) \cup R_2$ and $R_1<_v^\cw R_2$. Let $\widetilde{P}_1$ and $\widetilde{P}_2$ be the $\cw$-extended path and $\ccw$-extended path involved in the rotation. Then,
    \begin{enumerate}
        \item If $C_2$ is join irreducible then $C_2=C_{\min}(\widetilde{P}_2)$ and $R_2=R_{\widetilde{P}_2}^\cw$.
        \item If $C_1$ is meet irreducible then $C_1=C_{\max}(\widetilde{P}_1)$ and $R_1=R_{\widetilde{P}_1}^\ccw$.
        
    \end{enumerate}
\end{lemma}
\begin{proof}
    We prove only part (1) because (2) follows by symmetry. 
    Assume that $C_2$ is join irreducible. Then $C_1$ is the only maximal clique that is covered by $C_2$.
    Moreover, since all the routes in $C_2$ are coherent with $\widetilde{P}_2$ then $C_{\min}(\widetilde{P}_2)\leq C_2$. 
    If $C_2\neq C_{\min}(\widetilde{P}_2)$ then by~\Cref{lem.rotationExists2} there would be a maximal clique $C_2'$ whose routes are coherent with~$\widetilde{P}_2$ such that $C_2'\precccwrot C_2$. Then $C_2'=C_1$ and $R_1\in C_2'$, which is a contradiction because $R_1$ is incoherent with $\widetilde{P}_2$.
    Thus $C_2= C_{\min}(\widetilde{P}_2)$.

    Moreover,~\Cref{lem_join_irreducibles} (3) implies that $R_{\widetilde{P}_2}^\cw$ is the unique route that can be rotated down (in the $\cw$ direction) in $C_{\min}(\widetilde{P}_2)=C_2$. 
    Therefore, 
    $R_2=R_{\widetilde{P}_2}^\cw$.
\end{proof}

In other words, the map $C_2\rightarrow \widetilde{P}_2$ (resp. $C_1\rightarrow \widetilde{P}_1$) sends join irreducible (resp. meet irreducible) elements of the framing lattice to $\ccw$-extended paths (resp. $\cw$-extended paths), such that $C_2=C_{\min}(\widetilde{P}_2)$ (resp. $C_1=C_{\max}(\widetilde{P}_1)$).

\begin{corollary}\label{cor_extendedpath_labels}
    Let $\widetilde{P}_1$ be an $\cw$-extended path and $\widetilde{P}_2$ be a $\ccw$-extended path. 
    \begin{enumerate}
        \item If $\frakj=C_{\min}(\widetilde{P}_2)$ is the join irreducible element corresponding to $\widetilde{P}_2$ and $\frakj_*$ is the unique element covered by $\frakj$ then $\widetilde{\ell}_2(\frakj_*,\frakj)=\widetilde{P}_2$. 
        \item If $\frakm=C_{\max}(\widetilde{P}_1)$ is the meet irreducible element corresponding to $\widetilde{P}_1$ and $\frakm^*$ is the unique element covering $\frakm$ then $\widetilde{\ell}_1(\frakm,\frakm^*)=\widetilde{P}_1$. 
    \end{enumerate}
\end{corollary}
\begin{proof}
    We prove only (1) because (2) follows by symmetry. 
    Let $\widetilde{P}_2':=\widetilde{\ell}_2(\frakj_*,\frakj)$.
    By~\Cref{lem_meet_join_surjectivity}, taking $C_1=\frakj_*$ and $C_2=\frakj$ we deduce that $\frakj=C_{\min}(\widetilde{P}_2')$. But by the injectivity~\Cref{lem_join_injectivity} $\frakj=C_{\min}(\widetilde{P}_2)=C_{\min}(\widetilde{P}_2')$ implies that $\widetilde{P}_2=\widetilde{P}_2'$ as desired.
\end{proof}

\begin{proof}[Proof of~\Cref{thm_join_meet_irreducibles}]
By~\Cref{lem_join_irreducibles,lem_join_injectivity}, the map 
\[
\widetilde{P}_2\rightarrow C_{\min}(\widetilde{P}_2)
\] 
is an injective map 
from $\ccw$-extended paths to join irreducible elements. This map is surjective by~\Cref{lem_meet_join_surjectivity} (1).
Its inverse map is given by 
\[
\frakj \rightarrow \widetilde{\ell}_2(\frakj_*,\frakj).
\]
This proves part (1) of the theorem.

For part (2),~\Cref{lem_meet_irreducibles,lem_meet_injectivity} imply that 
\[
\widetilde{P}_1\rightarrow C_{\max}(\widetilde{P}_1)
\] 
is an injective map from $\cw$-extended paths to meet irreducible elements, and surjectivity follows from ~\Cref{lem_meet_join_surjectivity} (2).
Its inverse map is given by 
\[
\frakm \rightarrow \widetilde{\ell}_1(\frakm,\frakm^*).
\]
\end{proof}

\subsubsection{Bijection between join irreducibles and meet irreducibles}

As seen in the proof of~\Cref{thm_join_meet_irreducibles}, join irreducible elements $\frakj$ of the framing lattice can be identified with $\ccw$-extended paths, by taking the $\ccw$-extended path labeling $\widetilde{\ell}_2(\frakj_*,\frakj)$. 
Similarly, meet irreducible elements $\frakm$ correspond to $\cw$-extended paths, by taking the $\cw$-extended path labeling~$\widetilde{\ell}_1(\frakm,\frakm^*)$.  
On the other hand, $\ccw$-extended paths and $\cw$-extended paths are in correspondence via the inverse maps $\mapcw$ and $\mapccw$. So, we can think of these two maps as inverse maps between join irreducible and meet irreducible elements.

The following lemma implies that the join irreducible and meet irreducible elements corresponding to a cover relation are related by these two maps.

\begin{lemma}
\label{lem_extendedpaths_relation}
    Let $C_1\prec C_2$ be two maximal cliques with $C_2=(C_1\setminus R_1) \cup R_2$ and $R_1<_v^\cw R_2$. Let $\widetilde{P}_1$ and $\widetilde{P}_2$ be the $\cw$-extended path and $\ccw$-extended path involved in the rotation. Then,
    \begin{enumerate}
        \item $\widetilde{P}_1=\mapcw(\widetilde{P}_2)$.
        \item $\widetilde{P}_2=\mapccw(\widetilde{P}_1)$. 
    \end{enumerate}
\end{lemma}
\begin{proof}
    We prove only (1) because (2) follows by symmetry. 
    Let $P=\widetilde{P}_1\cap \widetilde{P}_2$ be the path at which $C_1$ and $C_2$ are incoherent, and let $w_1=\min(P)$ and $w_2=\max(P)$.
    Assume towards a contradiction that $\widetilde{P}_1\neq\mapcw(\widetilde{P}_2)$.
    Then, there is either (i) an edge $e$ entering $w_1$ between $\widetilde{P}_1$ and $\widetilde{P}_2$, or (ii) an edge $e$ exiting $w_2$ between $\widetilde{P}_2$ and $\widetilde{P}_1$. In case (i), we know that $e$ belongs to some route $R\in C_1$ by~\Cref{lem.path_is_extendable}. Furthermore, $w_1R=w_1R_1$ otherwise $R$ would be incoherent with either $R_1$ or $\Bot(R_1,R_2)$ in $C_1$.   But then $R\in C_1$ would be a route in between $R_1$ and $R_2$ which is a contradiction. 
    Case (ii) is solved similarly. 
\end{proof}

\begin{theorem}
        The framing lattice satisfies the following properties:
    \begin{enumerate}
        \item The map 
        $
        \frakj \rightarrow 
        \frakm :=C_{\max}(\widetilde{\ell}_1(\frakj_*,\frakj)) 
        $
        is a bijection from join irreducible to meet irreducible elements.
        \item The map 
        $
        \frakm \rightarrow 
        \frakj :=C_{\min}(\widetilde{\ell}_2(\frakm,\frakm^*)) 
        $
        is a bijection from meet irreducible to join irreducible elements.        
        \item These two maps are inverse to each other. 
    \end{enumerate}
\end{theorem}

\begin{proof}
    For (1), consider the composition of maps
        \[
        \frakj \rightarrow 
        \widetilde{P}_2 := \widetilde{\ell}_2(\frakj_*,\frakj) \rightarrow
        \widetilde{P}_1 := \varphi^{\cw}(\widetilde{P_2}) \rightarrow
        \frakm :=C_{\max}(\widetilde{P}_1). 
        \]
    Each map in this composition is a bijection:
    \begin{itemize}
        \item the first from join irreducible elements to $\ccw$-extended paths
        \item the second from $\ccw$-extended paths to $\cw$-extended paths
        \item the third from $\cw$-extended paths to meet irreducible elements 
    \end{itemize}
    Therefore, this composition of maps is a bijection from join irreducible to meet irreducible elements. 
    Moreover,~\Cref{lem_extendedpaths_relation} implies that $\widetilde{P}_1=\widetilde{\ell}_1(\frakj_*,\frakj)$.
    So, we obtain 
    $\frakm =C_{\max}(\widetilde{\ell}_1(\frakj_*,\frakj))$ 
    as desired.

    For parts (2) and (3), we just need to consider the composition of the inverse maps in reverse direction:
        \[
        \frakm \rightarrow 
        \widetilde{P}_1 := \widetilde{\ell}_1(\frakm,\frakm^*) \rightarrow
        \widetilde{P}_2 := \varphi^{\ccw}(\widetilde{P_1}) \rightarrow
        \frakj :=C_{\min}(\widetilde{P}_2). 
        \]
    Again,~\Cref{lem_extendedpaths_relation} implies that $\widetilde{P}_2=\widetilde{\ell}_2(\frakm,\frakm^*)$.
    So, we obtain 
    $\frakj =C_{\min}(\widetilde{\ell}_2(\frakm,\frakm^*))$.    
\end{proof}

As another consequence of~\Cref{lem_extendedpaths_relation} we get the following corollary. 

\begin{corollary} 
\label{cor_extendedpaths_join_meet_irreducibles}
    Let $C_1\prec C_2$ be two maximal cliques with $C_2=(C_1\setminus R_1) \cup R_2$ and $R_1<_v^\cw R_2$. Let $\widetilde{P}_1$ and $\widetilde{P}_2$ be the $\cw$-extended path and $\ccw$-extended path involved in the rotation. Then,
    \begin{enumerate}
        \item If $\frakj=C_{\min}(\widetilde{P}_2)$ then $\widetilde{\ell}_1(\frakj_*,\frakj)=\widetilde{P}_1$ and $\widetilde{\ell}_2(\frakj_*,\frakj)=\widetilde{P}_2$. 
        \item If $\frakm=C_{\max}(\widetilde{P}_1)$ then $\widetilde{\ell}_1(\frakm,\frakm^*)=\widetilde{P}_1$ and $\widetilde{\ell}_2(\frakm,\frakm^*)=\widetilde{P}_2$.
    \end{enumerate}
    In other words, the $\cw$-extended paths and $\ccw$-extended paths involved in the rotations \mbox{$C_1\prec C_2$}, $\frakj_*\prec \frakj$ and $\frakm\prec \frakm^*$ are the same. 
\end{corollary}

\begin{proof}
    We only prove (1) because (2) follows by symmetry. 
    Applying~\Cref{lem_extendedpaths_relation} to the two covering pairs $C_1\prec C_2$ and $\frakj_*\prec \frakj$ we get 
    $\widetilde{P}_1=\mapcw(\widetilde{P}_2)$ and 
    $\widetilde{\ell}_1(\frakj_*,\frakj)=\mapcw(\widetilde{\ell}_2(\frakj_*,\frakj))$, respectively. 
    Furthermore, we already proved that $\widetilde{\ell}_2(\frakj_*,\frakj)=\widetilde{P}_2$ in~\Cref{cor_extendedpath_labels}. 
    Therefore 
    $\widetilde{\ell}_1(\frakj_*,\frakj)=\mapcw(\widetilde{P}_2)=\widetilde{P}_1$.
\end{proof}


\subsection{The core label order of a congruence uniforme lattice}

In~\cite{reading_shard_2011}, Reading introduced the shard intersection order of a finite Coxeter group $W$, a lattice which contains the noncrossing partition lattice $NC(W)$ associated to $W$ as a sublattice. If $W$ is the Coxeter group of type $A_{n-1}$, one recovers the classical noncrossing partition lattice of $[n]$. 

One of Reading's main motivations for considering shards is the study of lattice congruences of the poset of regions of a hyperplane arrangement. Shards are certain closed polyhedral cones obtained by decomposing (or cutting) the hyperplanes of the arrangement. The shard intersection order is the poset of all possible intersections of shards ordered by reverse inclusion.     

In his Chapter~\cite[Section~9-7.4]{Rea16}, Reading proposed a generalization of the shard intersection order for any congruence uniform lattice $\scrL=(L,\leq)$. This is an alternate order structure on $L$ which was coined the name \emph{core label order} by Mühle in~\cite{muehle_core_2019}. If $\scrL$ is the poset of regions of the Coxeter hyperplane arrangement of $W$ one recovers the shard intersection order of $W$. If $\scrL$ is the classical Tamari lattice $\Tam(n)$ then its core label order is the classical lattice of noncrossing partition of $[n]$.    
In this section we recall the definition of the core label order, and propose the core label order of a framing lattice as the natural analog of the poset of noncrossing partitions for framing lattices. 

Let $\scrL=(L,\leq)$ be a finite congruence uniform lattice. 
For a cover relation $u\prec v$ in $L$ let $cg(u,v)$ be the finest lattice congruence of $\scrL$ in which $u$ and $v$ are equivalent. These are precisely the \emph{join irreducible lattice congruences} of $\scrL$. However, different covering pairs $u\prec v$ may lead to the same lattice congruence. In fact, if $\scrL$  is congruence uniform then the map $\frakj \rightarrow cg(\frakj_*,\frakj)$ is a bijection between join irreducible elements and join irreducible lattice congruences of $\scrL$, see~\cite[Section~9-7.4]{Rea16}. In other words, for every covering pair $u\prec v$ there is a unique join irreducible element $\frakj$ such that $cg(u,v)=cg(\frakj_*,\frakj)$. We denote this unique $\frakj$ as~$\frakj_{cg(u,v)}$. We label the edges $u\prec v$ of the Hasse diagram of $\scrL$ by their corresponding join irreducible elements $\frakj_{cg(u,v)}$. 

The following lemma gives a characterization of ~$\frakj_{cg(u,v)}$, cf.~\cite[Lemma~2.6]{graver_oriented_2018}.

\begin{lemma}[Cf.~{\cite[Lemma~2.6]{graver_oriented_2018}}]
\label{lem_join_label}
    Let $\scrL=(L,\leq)$ be a congruence uniform lattice, and $u\prec v$ be a covering pair. Then $\frakj_{cg(u,v)}$ is the unique join irreducible element $\frakj$ satisfying:
    \begin{enumerate}
        \item $\frakj \vee u = v$
        \item $\frakj \wedge u = \frakj_*$
    \end{enumerate}
\end{lemma}

The \defn{core} of an element $x \in L$ is the interval $[x_{\downarrow},x]$ where
\[
x_{\downarrow} = \bigwedge_{y\in L: y\prec x} y
\]
is the meet of all elements covered by $x$.
The \defn{core label set} of $x$ is the set of labels of the edges in the core of $x$:
\[
\psi_{\scrL}(x) := \left\{
\frakj_{cg(u,v)}:\ x^{\downarrow}\leq u\prec v \leq x
\right\}.
\]
When no confusion may arise, we omit the subscript $\scrL$ and just write $\psi(x)$ for simplicity.

The \defn{core label order} $CL(\scrL)=(L,\sqsubseteq )$ is the poset on $L$ with the following relation. For~$x,y\in L$ we say that $x\sqsubseteq y$ if and only if $\psi(x) \subseteq \psi(y)$.

\subsection{The framing core label order}
Our next goal is to give a precise description of the core label order of the framing lattice using $\ccw$-extended paths.  
The first step is to characterize the join irreducible element associated to each covering pair $C_1\prec C_2$.

\begin{proposition}\label{prop_join_labels}
    Let $C_1\prec C_2$ be two maximal cliques with $C_2=(C_1\setminus R_1) \cup R_2$ and $R_1<_v^\cw R_2$. Let $\widetilde{P}_1$ and $\widetilde{P}_2$ be the $\cw$-extended path and $\ccw$-extended path involved in the rotation. 
    Then
    \[
    \frakj_{cg(C_1,C_2)}=C_{\min}(\widetilde{P}_2),     
    \]
    the join irreducible element corresponding to the $\ccw$-extended path $\widetilde{P}_2$.
\end{proposition}
\begin{proof}
    Let $\frakj=C_{\min}(\widetilde{P}_2)$ be the join irreducible element corresponding to the $\ccw$-extended path $\widetilde{P}_2$, and let $\frakj_*$ be the unique element covered by $\frakj$. 
    By~\Cref{lem_join_label}, we just need to show that (1) $\frakj\vee C_1 = C_2$ and (2) $\frakj \wedge C_1 = \frakj_*$. We will see that both properties follow from~\Cref{prop.Z}.
    Let $P=\widetilde{P}_1\cap \widetilde{P}_2$ be the path at which $R_1$ and $R_2$ are incoherent.

    For (1), let $C^*=\frakj$ and consider the route $R^*=R_{\widetilde{P}_2}^{\cw}$. By~\Cref{lem_join_irreducibles}~(1), $R^*\in C^*$. 
    Moreover,  
    $\widetilde{P}_1 <_v^{\cw} R^*$ for any $v \in P$.
    This implies, by~~\Cref{prop.Z}~(1), that  
    \[
    C_1\vee \frakj = C_2\vee \frakj.
    \]
    Now, every route of $C_2$ is coherent with $\widetilde{P}_2$. By the properties of the $C_{\min}$ algorithm in~\Cref{cor_Cmin_properties}, we deduce that 
    $\frakj=C_{\min}(\widetilde{P}_2) \leq C_2$. Therefore
    \[
    C_1\vee \frakj = C_2\vee \frakj = C_2
    \]
    as we wanted to prove.

    For (2), we apply a similar argument for the covering pair $\frakj_*\prec \frakj$, which plays the role of~$C_1\prec C_2$ in~\Cref{prop.Z}. 
    Note that the $\cw$-extended path and $\ccw$-extended path involved in the rotation from $\frakj_*\prec \frakj$
    are also equal to $\widetilde{P}_1$ and $\widetilde{P}_2$, respectively (by~\Cref{cor_extendedpaths_join_meet_irreducibles}).
    Now, we consider $C^*=C_1$ and the route $R^*=R_1\in C^*$. 
    Since $R^* <_v^{\cw} \widetilde{P}_2$ for any $v \in P$,~\Cref{prop.Z} (2) implies that
    \[
    \frakj_*\wedge C_1 = \frakj \wedge C_1.
    \]
    Now we want to apply~\Cref{prop.Z} (2) again for the covering pair $C_1\prec C_2$, $C^*=\frakj_*$ and the route $R^*\in \frakj_*$ obtained by a cw rotation of the route $R_{\widetilde{P}_2}^\cw\in \frakj$. 
    Since $\widetilde{\ell}_1(\frakj_*,\frakj)=\widetilde{P}_1$ (by~\Cref{cor_extendedpaths_join_meet_irreducibles}), then $\widetilde{P}_1$ is a subpath of $R^*$. Thus 
    $R^* <_v^{\cw} \widetilde{P}_2$ for any $v \in P$, and applying~\Cref{prop.Z} (2) we get
    \[
    \frakj_*\wedge C_1 = \frakj_* \wedge C_2.
    \]
    Combining the two previous equations we get 
    \[
    \frakj \wedge C_1 = 
    \frakj_*\wedge C_1 = 
    \frakj_* \wedge C_2.
    \]
    But $\frakj_*\leq \frakj \leq C_2$, and so $\frakj \wedge C_1 = \frakj_* \wedge C_2=\frakj_*$ as desired.
\end{proof}

Instead of labeling the covering pairs (edges) by join irreducible elements (which can be very complicated in general), we can significantly simplify everything by labeling the covering pairs (edges) with the corresponding $\ccw$-extended paths.   

Define the \defn{$\ccw$-core label set} of a maximal clique $x$ in the framing lattice $\scrL=\scrL_{G,F}$ as 
\[
\psi^{\ccw}_{\scrL}(x) := \left\{
\widetilde{\ell}_2(u,v):\ x^{\downarrow}\leq u\prec v \leq x
\right\}.
\]
That is, the set of $\ccw$-extended path labels in the core $[x^{\downarrow},x]$ of $x$. 
When no confusion may arise, we omit the subscript $\scrL$ and simply write $\psi^{\ccw}(x)$ for simplicity.

\begin{figure}
    \centering

    \caption{Examples of the core label order of the framing lattices in~\Cref{fig_extended_path_labeling_examples}}
    \label{fig_corelabel_examples}
\end{figure}

The following result is a powerful and very simple description of the core label order of a framing lattice. 
\begin{theorem}
    \label{thm.core_label_order}
    The core label order of the framing lattice $\scrL=\scrL_{G,F}$ is isomorphic to the poset of $\ccw$-core label sets $\psi^{\ccw}_{\scrL}(x)$ ordered by inclusion. 
\end{theorem}

\begin{proof}
    This follows directly from the description of the join irreducible labels $\frakj_{cg(C_1,C_2)}=C_{\min}(\widetilde{P}_2)$ in~\Cref{prop_join_labels} and the bijective correspondence $C_{\min}(\widetilde{P}_2) \leftrightarrow \widetilde{P}_2=\widetilde{\ell}_2(C_1,C_2)$ between join irreducible elements and $\ccw$-extended paths.
\end{proof}

Figure~\ref{fig_corelabel_examples} shows examples of the core label order for the framing lattices in Figure~\ref{fig_extended_path_labeling_examples}.

\begin{remark}
    We remark that the framing lattice has a richer underlying geometric structure, consisting not only of vertices and edges but also 2-dimensional faces that are squares, pentagons and hexagons, as well as possibly higher dimensional faces. This geometric structure, which we call the framingtope, will be introduced and studied in forthcoming work~\cite{framingtopes}. It is a polytopal complex that is dual to the complex of interior faces of the framed triangulation of a flow polytope. As illustrated on the right of Figure~\ref{fig_corelabel_examples}, it can be topologically decomposed into a disjoint union of half open cells, one for each element of the lattice. If we identify faces with the same set of labels, then the core label order is just the containment poset of equivalence classes of faces (more explicitly their closures) in this cell decomposition.
\end{remark}

\section{Lattice congruences and quotients}\label{sec.quotients}

The purpose of this section is to study certain lattice congruences and quotients of the framing lattice arising from certain operations called $M$-moves. In particular, we will show that the resulting quotients are the framing lattices of the modified framed graphs. 
Two examples for the oruga graph are illustrated in~\Cref{fig_quotients}, and further examples are shown in~\Cref{fig_quotients11223}.

    

\begin{figure}[h]
    \centering
    \begin{tikzpicture}

\begin{scope}[xshift=0, scale=0.4]
    \draw[thick, color=black] (0,0) .. controls (0.4, 0.6) and (1.6, 0.6) .. (2,0);
    \draw[thick, color=black] (0,0) .. controls (0.4, -0.6) and (1.6, -0.6) .. (2,0);    
    \draw[ultra thick, color=ForestGreen] (2,0) .. controls (2.4, 0.6) and (3.6, 0.6) .. (4,0);
    \draw[ultra thick, color=red] (2,0) .. controls (2.4, -0.6) and (3.6, -0.6) .. (4,0);
    \draw[very thick, color=black] (4,0) .. controls (4.4, 0.6) and (5.6, 0.6) .. (6,0);
    \draw[thick, color=black] (4,0) .. controls (4.4, -0.6) and (5.6, -0.6) .. (6,0);        
    \node[circle, draw, inner sep=1.4pt, fill] (1) at (0,0) {};
    \node[circle, draw, inner sep=1.4pt, fill] (2) at (2,0) {};
    \node[circle, draw, inner sep=1.4pt, fill] (3) at (4,0) {};
    \node[circle, draw, inner sep=1.4pt, fill] (4) at (6,0) {}; 

    \node[] (a) at (1,1.5) {$G$};
    \node[] (a) at (3,1) {\textcolor{ForestGreen}{$e$}};
    \node[] (a) at (3,-1) {\textcolor{red}{$e'$}};
\end{scope}
\begin{scope}[xshift=130, scale=0.4]
    \draw[thick, color=black] (0,0) .. controls (0.4, 0.6) and (1.6, 0.6) .. (2,0);
    \draw[thick, color=black] (0,0) .. controls (0.4, -0.6) and (1.6, -0.6) .. (2,0);    
    \draw[thick, color=black] (2,0) .. controls (2.6, 1.0) and (5.4, 1.0) .. (6,0);
    \draw[thick, color=black] (0,0) .. controls (0.6, 1.0) and (3.4, 1.0) .. (4,0);
    
    \draw[ultra thick, color=red] (2,0) .. controls (2.4, -0.6) and (3.6, -0.6) .. (4,0);
    \draw[thick, color=black] (4,0) .. controls (4.4, 0.6) and (5.6, 0.6) .. (6,0);
    \draw[thick, color=black] (4,0) .. controls (4.4, -0.6) and (5.6, -0.6) .. (6,0);        
    \node[circle, draw, inner sep=1.4pt, fill] (1) at (0,0) {};
    \node[circle, draw, inner sep=1.4pt, fill] (2) at (2,0) {};
    \node[circle, draw, inner sep=1.4pt, fill] (3) at (4,0) {};
    \node[circle, draw, inner sep=1.4pt, fill] (4) at (6,0) {}; 

    \node[] (a) at (3,1.5) {\small $M(G,$ \textcolor{ForestGreen}{$e$}$)$};    
\end{scope}
\begin{scope}[xshift=260, yshift=0, scale=0.4]
    \draw[thick, color=black] (0,0) .. controls (0.4, 0.6) and (1.6, 0.6) .. (2,0);
    \draw[thick, color=black] (0,0) .. controls (0.4, -0.6) and (1.6, -0.6) .. (2,0);    
    \draw[thick, color=black] (2,0) .. controls (2.6, -1.0) and (5.4, -1.0) .. (6,0);
    \draw[thick, color=black] (0,0) .. controls (0.6, -1.0) and (3.4, -1.0) .. (4,0);
    
    \draw[ultra thick, color=ForestGreen] (2,0) .. controls (2.4, 0.6) and (3.6, 0.6) .. (4,0);
    \draw[thick, color=black] (4,0) .. controls (4.4, 0.6) and (5.6, 0.6) .. (6,0);
    \draw[thick, color=black] (4,0) .. controls (4.4, -0.6) and (5.6, -0.6) .. (6,0);        
    \node[circle, draw, inner sep=1.4pt, fill] (1) at (0,0) {};
    \node[circle, draw, inner sep=1.4pt, fill] (2) at (2,0) {};
    \node[circle, draw, inner sep=1.4pt, fill] (3) at (4,0) {};
    \node[circle, draw, inner sep=1.4pt, fill] (4) at (6,0) {}; 
    
    \node[] (a) at (3,1.5) {\small $M(G,$ \textcolor{red}{$e'$}$)$}; 
\end{scope}

\begin{scope}[xshift=0,yshift=-100]

\begin{scope}[xshift=20,yshift=0]
\begin{scope}[scale=0.25, xshift=50, yshift=0]
    \draw[thick, color=NavyBlue] (0,0) -- (4,3);
    \draw[ultra thick, color=red] (4,3) -- (4,7);
    \draw[thick, color=NavyBlue] (4,7) -- (0,10);
    \draw[thick, color=NavyBlue] (0,10) -- (-4,7);
    \draw[ultra thick, color=ForestGreen] (-4,7) -- (-4,3);
    \draw[thick, color=NavyBlue] (-4,3)--(0,0);		
	\node[circle, fill,inner sep=1.4pt,color=black] at (0,0)  {};
	\node[circle,fill,inner sep=1.4pt,color=black]  at (4,3) {};
	\node[circle,fill,inner sep=1.4pt,color=black]  at (4,7) {}; 
	\node[circle,fill,inner sep=1.4pt,color=black]  at (0,10) {};
	\node[circle,fill,inner sep=1.4pt,color=black]  at (-4,7) {};
	\node[circle,fill,inner sep=1.4pt,color=black]  at (-4,3) {};
\end{scope}

\begin{scope}[xshift=130,yshift=0]
\begin{scope}[scale=0.25, xshift=50, yshift=0]
    \draw[thick, color=NavyBlue] (0,0) -- (-4,5);
    \draw[thick, color=NavyBlue] (-4,5) -- (0,10);
    \draw[thick, color=NavyBlue] (0,10) -- (4,7);
    \draw[ultra thick, color=red] (4,7) -- (4,3);
    \draw[thick, color=NavyBlue] (4,3)--(0,0);		
	\node[circle,fill,inner sep=1.4pt,color=black] at (0,0)  {};
	\node[circle,fill,inner sep=1.4pt,color=black]  at (-4,5) {};
	\node[circle,fill,inner sep=1.4pt,color=black]  at (0,10) {};
	\node[circle,fill,inner sep=1.4pt,color=black]  at (4,7) {};
	\node[circle,fill,inner sep=1.4pt,color=black]  at (4,3) {};
\end{scope}
\end{scope}

\begin{scope}[xshift=260, yshift=0]
\begin{scope}[scale=0.25, xshift=50, yshift=0]
    xshift=50, yshift=0]
    \draw[thick, color=NavyBlue] (0,0) -- (4,5);
    \draw[thick, color=NavyBlue] (4,5) -- (0,10);
    \draw[thick, color=NavyBlue] (0,10) -- (-4,7);
    \draw[ultra thick, color=ForestGreen] (-4,7) -- (-4,3);
    \draw[thick, color=NavyBlue] (-4,3)--(0,0);		
	\node[circle,fill,inner sep=1.4pt,color=black] at (0,0)  {};
	\node[circle,fill,inner sep=1.4pt,color=black]  at (4,5) {};
	\node[circle,fill,inner sep=1.4pt,color=black]  at (0,10) {};
	\node[circle,fill,inner sep=1.4pt,color=black]  at (-4,7) {};
	\node[circle,fill,inner sep=1.4pt,color=black]  at (-4,3) {};
\end{scope}
\end{scope}
\end{scope}
\end{scope} 

\end{tikzpicture}
    \caption{Some lattice quotients obtained via $M$-moves.}
    \label{fig_quotients}
\end{figure}

\subsection{An equivalence relation via \texorpdfstring{$M$}{}-moves}

We say that an edge $e=(i,j)$ is an \defn{inner edge} of $G$ if~$i$ and~$j$ are inner vertices, i.e. when $e$ is not incident to the source or sink of $G$.
Given an inner edge~$e=(i,j)$ of $G$, define the graph $M(G,e)$ as the graph obtained from $G$ by replacing $e$ with the edges $(s,j)$ and $(i,t)$.
In other words, $M(G,e)$ is obtained from $G$ by cutting the edge $e$ in half and identifying the two sources and identifying the two sinks.
For a framing~$F$ of~$G$, the graph $M(G,e)$ inherits a framing $F_e$ in the natural way, replacing $e$ with $e_s:=(s,j)$ in the order $\leq_{\In(j)}$ and replacing $e$ with $e_t:=(i,t)$ in the order~$\leq_{\Out(i)}$.
In this case, we say that the framed graph $(M(G,e),F_e)$ is obtained from $(G,F)$ by an \defn{$M$-move}. Such $M$-moves are due to Martha Yip and first appear in \cite{Tam23}, where they are used to connect permutree lattices with framed triangulations of flow polytopes. 

\begin{figure}[h]
    \centering
    \begin{tikzpicture}[scale=1]

\begin{scope}[scale=0.6, xshift=0, yshift=0]
    \draw[thick, color=black] (0,0) .. controls (0.4, 0.6) and (1.6, 0.6) .. (2,0);
    \draw[thick, color=black] (0,0) .. controls (0.4, -0.6) and (1.6, -0.6) .. (2,0);    
    \draw[thick, color=black] (2,0) .. controls (2.4, 0.6) and (3.6, 0.6) .. (4,0);
    \draw[thick, color=black] (2,0) .. controls (2.4, -0.6) and (3.6, -0.6) .. (4,0);
    \draw[thick, color=black] (4,0) .. controls (4.4, 0.6) and (5.6, 0.6) .. (6,0);
    \draw[thick, color=black] (4,0) .. controls (4.4, -0.6) and (5.6, -0.6) .. (6,0);        
    \node[circle, draw, inner sep=1.4pt, fill] (1) at (0,0) {};
    \node[circle, draw, inner sep=1.4pt, fill] (2) at (2,0) {};
    \node[circle, draw, inner sep=1.4pt, fill] (3) at (4,0) {};
    \node[circle, draw, inner sep=1.4pt, fill] (4) at (6,0) {};    

    \node[] (a) at (3,0.8) {\scriptsize $e$};    
    \node[] (a) at (2,-0.4) {\scriptsize $i$};
    \node[] (a) at (4,-0.4) {\scriptsize $j$};   
    \node[] (a) at (0,-0.4) {\scriptsize $s$};
    \node[] (a) at (6,-0.4) {\scriptsize $t$};    

    \node[] (a) at (3,-1.5) {$G$};        
\end{scope}	

\begin{scope}[scale=0.6, xshift=240, yshift=0]
    \draw[thick, color=black] (0,0) .. controls (0.4, 0.6) and (1.6, 0.6) .. (2,0);
    \draw[thick, color=black] (0,0) .. controls (0.4, -0.6) and (1.6, -0.6) .. (2,0);    
    \draw[thick, color=black] (2,0) .. controls (2.4, 0.6) and (2.6, 0.6) .. (2.7,0.6);

    \draw[thick, color=black] (3.3,0.6) .. controls (3.4, 0.6) and (3.6, 0.6) .. (4,0);
    
    \draw[thick, color=black] (2,0) .. controls (2.4, -0.6) and (3.6, -0.6) .. (4,0);
    \draw[thick, color=black] (4,0) .. controls (4.4, 0.6) and (5.6, 0.6) .. (6,0);
    \draw[thick, color=black] (4,0) .. controls (4.4, -0.6) and (5.6, -0.6) .. (6,0);        
    \node[circle, draw, inner sep=1.4pt, fill] (1) at (0,0) {};
    \node[circle, draw, inner sep=1.4pt, fill] (2) at (2,0) {};
    \node[circle, draw, inner sep=1.4pt, fill] (3) at (4,0) {};
    \node[circle, draw, inner sep=1.4pt, fill] (4) at (6,0) {};    

    \node[circle, draw, inner sep=1.4pt, fill] (5) at (2.7,0.6) {};
    \node[circle, draw, inner sep=1.4pt, fill] (6) at (3.3,0.6) {};

    \node[] (a) at (2,-0.4) {\scriptsize $i$};
    \node[] (a) at (4,-0.4) {\scriptsize $j$};   
    \node[] (a) at (0,-0.4) {\scriptsize $s$};
    \node[] (a) at (6,-0.4) {\scriptsize $t$};
\end{scope}	

\begin{scope}[scale=0.6, xshift=480, yshift=0]
    \draw[thick, color=black] (0,0) .. controls (0.4, 0.6) and (1.6, 0.6) .. (2,0);
    \draw[thick, color=black] (0,0) .. controls (0.4, -0.6) and (1.6, -0.6) .. (2,0);    
    \draw[thick, color=black] (2,0) .. controls (2.6, 1.0) and (5.4, 1.0) .. (6,0);
    \draw[thick, color=black] (0,0) .. controls (0.6, 1.0) and (3.4, 1.0) .. (4,0);
    
    \draw[thick, color=black] (2,0) .. controls (2.4, -0.6) and (3.6, -0.6) .. (4,0);
    \draw[thick, color=black] (4,0) .. controls (4.4, 0.6) and (5.6, 0.6) .. (6,0);
    \draw[thick, color=black] (4,0) .. controls (4.4, -0.6) and (5.6, -0.6) .. (6,0);        
    \node[circle, draw, inner sep=1.4pt, fill] (1) at (0,0) {};
    \node[circle, draw, inner sep=1.4pt, fill] (2) at (2,0) {};
    \node[circle, draw, inner sep=1.4pt, fill] (3) at (4,0) {};
    \node[circle, draw, inner sep=1.4pt, fill] (4) at (6,0) {};    

    \node[] (a) at (2,-0.4) {\scriptsize $i$};
    \node[] (a) at (4,-0.4) {\scriptsize $j$};   
    \node[] (a) at (0,-0.4) {\scriptsize $s$};
    \node[] (a) at (6,-0.4) {\scriptsize $t$};    

    \node[] (a) at (3,-1.5) {$M(G,e)$};
\end{scope}

\end{tikzpicture}
    \vspace{-1cm}
    \caption{An example of an M-move.}
    \label{fig_M_move_example}
\end{figure}

Given an inner edge $e=(i,j)$ in $(G,F)$, we define the map $\varphi_e$ sending routes in $G$ to sets of routes in $M(G,e)$ as follows.
For a route $R$ in $G$, if $e\notin R$ we define $\varphi_e(R) = R$, and if $e\in R$ we define $\varphi_e(R)$ as the pair of routes $e_s\cup jR$ and $Ri\cup e_t$.
If $A$ is a set of routes in~$G$, we define the \defn{split map} $\Phi_e$ to be the map $A\mapsto \bigcup_{R\in A} \varphi_e(R)$.
Two examples of maximal cliques with the same image are illustrated in~\Cref{fig_phi_map}. 

\begin{figure}[h]
    \centering
    \begin{tikzpicture}[scale=1]


\begin{scope}[scale=0.6, xshift=0, yshift=0]
    \draw[thick, color=black] (0,0) .. controls (0.4, 0.6) and (1.6, 0.6) .. (2,0);
    \draw[thick, color=black] (0,0) .. controls (0.4, -0.6) and (1.6, -0.6) .. (2,0);    
    \draw[ultra thick, color=ForestGreen] (2,0) .. controls (2.4, 0.6) and (3.6, 0.6) .. (4,0);
    \draw[thick, color=black] (2,0) .. controls (2.4, -0.6) and (3.6, -0.6) .. (4,0);
    \draw[thick, color=black] (4,0) .. controls (4.4, 0.6) and (5.6, 0.6) .. (6,0);
    \draw[thick, color=black] (4,0) .. controls (4.4, -0.6) and (5.6, -0.6) .. (6,0);        
    \node[circle, draw, inner sep=1.4pt, fill] (1) at (0,0) {};
    \node[circle, draw, inner sep=1.4pt, fill] (2) at (2,0) {};
    \node[circle, draw, inner sep=1.4pt, fill] (3) at (4,0) {};
    \node[circle, draw, inner sep=1.4pt, fill] (4) at (6,0) {};    

    \node[] (a) at (2.5,0.7) {\scriptsize $e$};    
    \node[] (a) at (3,1.5) {$G$};        
\end{scope}	

\begin{scope}[scale=0.6, xshift=240, yshift=0]
    \draw[thick, color=black] (0,0) .. controls (0.4, 0.6) and (1.6, 0.6) .. (2,0);
    \draw[thick, color=black] (0,0) .. controls (0.4, -0.6) and (1.6, -0.6) .. (2,0);    
    \draw[thick, color=black] (2,0) .. controls (2.6, 1.0) and (5.4, 1.0) .. (6,0);
    \draw[thick, color=black] (0,0) .. controls (0.6, 1.0) and (3.4, 1.0) .. (4,0);
    
    \draw[thick, color=black] (2,0) .. controls (2.4, -0.6) and (3.6, -0.6) .. (4,0);
    \draw[thick, color=black] (4,0) .. controls (4.4, 0.6) and (5.6, 0.6) .. (6,0);
    \draw[thick, color=black] (4,0) .. controls (4.4, -0.6) and (5.6, -0.6) .. (6,0);        
    \node[circle, draw, inner sep=1.4pt, fill] (1) at (0,0) {};
    \node[circle, draw, inner sep=1.4pt, fill] (2) at (2,0) {};
    \node[circle, draw, inner sep=1.4pt, fill] (3) at (4,0) {};
    \node[circle, draw, inner sep=1.4pt, fill] (4) at (6,0) {};    

    \node[] (a) at (3,1.5) {$M(G,e)$};
\end{scope}	

\begin{scope}[scale=0.6, xshift=480, yshift=0]
    \draw[thick, color=black] (0,0) .. controls (0.4, 0.6) and (1.6, 0.6) .. (2,0);
    \draw[thick, color=black] (0,0) .. controls (0.4, -0.6) and (1.6, -0.6) .. (2,0);    
    \draw[ultra thick, color=ForestGreen] (2,0) .. controls (2.4, 0.6) and (3.6, 0.6) .. (4,0);
    \draw[thick, color=black] (2,0) .. controls (2.4, -0.6) and (3.6, -0.6) .. (4,0);
    \draw[thick, color=black] (4,0) .. controls (4.4, 0.6) and (5.6, 0.6) .. (6,0);
    \draw[thick, color=black] (4,0) .. controls (4.4, -0.6) and (5.6, -0.6) .. (6,0);        
    \node[circle, draw, inner sep=1.4pt, fill] (1) at (0,0) {};
    \node[circle, draw, inner sep=1.4pt, fill] (2) at (2,0) {};
    \node[circle, draw, inner sep=1.4pt, fill] (3) at (4,0) {};
    \node[circle, draw, inner sep=1.4pt, fill] (4) at (6,0) {};    

    \node[] (a) at (2.5,0.7) {\scriptsize $e$};    
    \node[] (a) at (3,1.5) {$G$};        
\end{scope}


\begin{scope}[xshift=10,yshift = -50]
\begin{scope}[scale=0.4, xshift=0, yshift=0]
    \draw[thick, color=black] (0,0) .. controls (0.4, -0.6) and (1.6, -0.6) .. (2,0);    
    \draw[thick, color=black] (2,0) .. controls (2.4, -0.6) and (3.6, -0.6) .. (4,0);
    \draw[thick, color=black] (4,0) .. controls (4.4, -0.6) and (5.6, -0.6) .. (6,0);        
    \node[circle, draw, inner sep=1.4pt, fill] (1) at (0,0) {};
    \node[circle, draw, inner sep=1.4pt, fill] (2) at (2,0) {};
    \node[circle, draw, inner sep=1.4pt, fill] (3) at (4,0) {};
    \node[circle, draw, inner sep=1.4pt, fill] (4) at (6,0) {};    

\end{scope}	
\begin{scope}[scale=0.4, xshift=0, yshift=-50]
    \draw[thick, color=black] (0,0) .. controls (0.4, -0.6) and (1.6, -0.6) .. (2,0);    
    \draw[thick, color=black] (2,0) .. controls (2.4, 0.6) and (3.6, 0.6) .. (4,0);
    \draw[thick, color=black] (4,0) .. controls (4.4, -0.6) and (5.6, -0.6) .. (6,0);        
    \node[circle, draw, inner sep=1.4pt, fill] (1) at (0,0) {};
    \node[circle, draw, inner sep=1.4pt, fill] (2) at (2,0) {};
    \node[circle, draw, inner sep=1.4pt, fill] (3) at (4,0) {};
    \node[circle, draw, inner sep=1.4pt, fill] (4) at (6,0) {};    
\end{scope}
\begin{scope}[scale=0.4, xshift=0,yshift=-100]
    \draw[thick, color=black] (0,0) .. controls (0.4, 0.6) and (1.6, 0.6) .. (2,0);
    \draw[thick, color=black] (2,0) .. controls (2.4, 0.6) and (3.6, 0.6) .. (4,0);
    \draw[thick, color=black] (4,0) .. controls (4.4, -0.6) and (5.6, -0.6) .. (6,0);        
    \node[circle, draw, inner sep=1.4pt, fill] (1) at (0,0) {};
    \node[circle, draw, inner sep=1.4pt, fill] (2) at (2,0) {};
    \node[circle, draw, inner sep=1.4pt, fill] (3) at (4,0) {};
    \node[circle, draw, inner sep=1.4pt, fill] (4) at (6,0) {};    
       
\end{scope}
\begin{scope}[scale=0.4, xshift=0,yshift=-150]
    \draw[thick, color=black] (0,0) .. controls (0.4, 0.6) and (1.6, 0.6) .. (2,0);
    \draw[thick, color=black] (2,0) .. controls (2.4, 0.6) and (3.6, 0.6) .. (4,0);
    \draw[thick, color=black] (4,0) .. controls (4.4, 0.6) and (5.6, 0.6) .. (6,0);
    \node[circle, draw, inner sep=1.4pt, fill] (1) at (0,0) {};
    \node[circle, draw, inner sep=1.4pt, fill] (2) at (2,0) {};
    \node[circle, draw, inner sep=1.4pt, fill] (3) at (4,0) {};
    \node[circle, draw, inner sep=1.4pt, fill] (4) at (6,0) {};    
\end{scope}
\begin{scope}[scale=0.5, xshift=70,yshift=-160]
    \node[] (a) at (0,0) {$C$};
    \node[] (a) at (0,-1) {$213$}; 
\end{scope}
\end{scope}


\begin{scope}[xshift=310, yshift = -50]
\begin{scope}[scale=0.4, xshift=0, yshift=0]
    \draw[thick, color=black] (0,0) .. controls (0.4, -0.6) and (1.6, -0.6) .. (2,0);    
    \draw[thick, color=black] (2,0) .. controls (2.4, -0.6) and (3.6, -0.6) .. (4,0);
    \draw[thick, color=black] (4,0) .. controls (4.4, -0.6) and (5.6, -0.6) .. (6,0);        
    \node[circle, draw, inner sep=1.4pt, fill] (1) at (0,0) {};
    \node[circle, draw, inner sep=1.4pt, fill] (2) at (2,0) {};
    \node[circle, draw, inner sep=1.4pt, fill] (3) at (4,0) {};
    \node[circle, draw, inner sep=1.4pt, fill] (4) at (6,0) {};    

\end{scope}	
\begin{scope}[scale=0.4, xshift=0, yshift=-50]
    \draw[thick, color=black] (0,0) .. controls (0.4, -0.6) and (1.6, -0.6) .. (2,0);    
    \draw[thick, color=black] (2,0) .. controls (2.4, 0.6) and (3.6, 0.6) .. (4,0);
    \draw[thick, color=black] (4,0) .. controls (4.4, -0.6) and (5.6, -0.6) .. (6,0);        
    \node[circle, draw, inner sep=1.4pt, fill] (1) at (0,0) {};
    \node[circle, draw, inner sep=1.4pt, fill] (2) at (2,0) {};
    \node[circle, draw, inner sep=1.4pt, fill] (3) at (4,0) {};
    \node[circle, draw, inner sep=1.4pt, fill] (4) at (6,0) {};    
\end{scope}
\begin{scope}[scale=0.4, xshift=0,yshift=-100]
    \draw[thick, color=black] (0,0) .. controls (0.4, -0.6) and (1.6, -0.6) .. (2,0);    
    \draw[thick, color=black] (2,0) .. controls (2.4, 0.6) and (3.6, 0.6) .. (4,0);
    \draw[thick, color=black] (4,0) .. controls (4.4, 0.6) and (5.6, 0.6) .. (6,0);
    \node[circle, draw, inner sep=1.4pt, fill] (1) at (0,0) {};
    \node[circle, draw, inner sep=1.4pt, fill] (2) at (2,0) {};
    \node[circle, draw, inner sep=1.4pt, fill] (3) at (4,0) {};
    \node[circle, draw, inner sep=1.4pt, fill] (4) at (6,0) {};    
       
\end{scope}
\begin{scope}[scale=0.4, xshift=0,yshift=-150]
    \draw[thick, color=black] (0,0) .. controls (0.4, 0.6) and (1.6, 0.6) .. (2,0);
    \draw[thick, color=black] (2,0) .. controls (2.4, 0.6) and (3.6, 0.6) .. (4,0);
    \draw[thick, color=black] (4,0) .. controls (4.4, 0.6) and (5.6, 0.6) .. (6,0);
    \node[circle, draw, inner sep=1.4pt, fill] (1) at (0,0) {};
    \node[circle, draw, inner sep=1.4pt, fill] (2) at (2,0) {};
    \node[circle, draw, inner sep=1.4pt, fill] (3) at (4,0) {};
    \node[circle, draw, inner sep=1.4pt, fill] (4) at (6,0) {};    
\end{scope}
\begin{scope}[scale=0.5, xshift=70,yshift=-160]
    \node[] (a) at (0,0) {$C'$};
    \node[] (a) at (0,-1) {$231$}; 
\end{scope}
\end{scope}


\begin{scope}[xshift=160, yshift = -50]
\begin{scope}[scale=0.4, xshift=0, yshift=0]
    \draw[thick, color=black] (0,0) .. controls (0.4, -0.6) and (1.6, -0.6) .. (2,0);    
    \draw[thick, color=black] (2,0) .. controls (2.4, -0.6) and (3.6, -0.6) .. (4,0);
    \draw[thick, color=black] (4,0) .. controls (4.4, -0.6) and (5.6, -0.6) .. (6,0);        
    \node[circle, draw, inner sep=1.4pt, fill] (1) at (0,0) {};
    \node[circle, draw, inner sep=1.4pt, fill] (2) at (2,0) {};
    \node[circle, draw, inner sep=1.4pt, fill] (3) at (4,0) {};
    \node[circle, draw, inner sep=1.4pt, fill] (4) at (6,0) {};    

\end{scope}	
\begin{scope}[scale=0.4, xshift=0, yshift=-50]
    \draw[thick, color=black] (0,0) .. controls (0.4, -0.6) and (1.6, -0.6) .. (2,0);    

    \draw[thick, color=black] (2,0) .. controls (2.6, 1.0) and (5.4, 1.0) .. (6,0);
    
    \node[circle, draw, inner sep=1.4pt, fill] (1) at (0,0) {};
    \node[circle, draw, inner sep=1.4pt, fill] (2) at (2,0) {};
    \node[circle, draw, inner sep=1.4pt, fill] (3) at (4,0) {};
    \node[circle, draw, inner sep=1.4pt, fill] (4) at (6,0) {};    
\end{scope}
\begin{scope}[scale=0.4, xshift=0,yshift=-100]
    \draw[thick, color=black] (4,0) .. controls (4.4, -0.6) and (5.6, -0.6) .. (6,0);        

    \draw[thick, color=black] (0,0) .. controls (0.6, 1.0) and (3.4, 1.0) .. (4,0);
    
    \node[circle, draw, inner sep=1.4pt, fill] (1) at (0,0) {};
    \node[circle, draw, inner sep=1.4pt, fill] (2) at (2,0) {};
    \node[circle, draw, inner sep=1.4pt, fill] (3) at (4,0) {};
    \node[circle, draw, inner sep=1.4pt, fill] (4) at (6,0) {};    
       
\end{scope}
\begin{scope}[scale=0.4, xshift=0, yshift=-150]
    \draw[thick, color=black] (0,0) .. controls (0.4, 0.6) and (1.6, 0.6) .. (2,0);

    \draw[thick, color=black] (2,0) .. controls (2.6, 1.0) and (5.4, 1.0) .. (6,0);
    
    \node[circle, draw, inner sep=1.4pt, fill] (1) at (0,0) {};
    \node[circle, draw, inner sep=1.4pt, fill] (2) at (2,0) {};
    \node[circle, draw, inner sep=1.4pt, fill] (3) at (4,0) {};
    \node[circle, draw, inner sep=1.4pt, fill] (4) at (6,0) {};    
\end{scope}
\begin{scope}[scale=0.4, xshift=0,yshift=-200]
    \draw[thick, color=black] (4,0) .. controls (4.4, 0.6) and (5.6, 0.6) .. (6,0);

    \draw[thick, color=black] (0,0) .. controls (0.6, 1.0) and (3.4, 1.0) .. (4,0);
    
    \node[circle, draw, inner sep=1.4pt, fill] (1) at (0,0) {};
    \node[circle, draw, inner sep=1.4pt, fill] (2) at (2,0) {};
    \node[circle, draw, inner sep=1.4pt, fill] (3) at (4,0) {};
    \node[circle, draw, inner sep=1.4pt, fill] (4) at (6,0) {};    
\end{scope}

\begin{scope}[scale=0.5, xshift=70,yshift=-200]
    \node[] (a) at (0,0) {$\Phi_e(C)=\Phi_e(C')$};
\end{scope}
\end{scope}

\begin{scope}[xshift=92,yshift=-50]
    \node[] (a) at (0.85,0.5) {$\Phi_e$};
    \draw[color=black, -{Stealth[scale=1]}] (0,0) -- (1.8,0);
    \draw[color=black, -{Stealth[scale=1]}] (0,-0.7) -- (1.8,-0.7);
    \draw[color=black, -{Stealth[scale=1]}] (0,-0.7) -- (1.8,-1.3);
    \draw[color=black, -{Stealth[scale=1]}] (0,-1.4) -- (1.8,-1.4);
    \draw[color=black, -{Stealth[scale=1]}] (0,-1.4) -- (1.8,-2.0);
    \draw[color=black, -{Stealth[scale=1]}] (0,-2.1) -- (1.8,-2.1);
    \draw[color=black, -{Stealth[scale=1]}] (0,-2.1) -- (1.8,-2.7); 
\end{scope}

\begin{scope}[xshift=245,yshift=-50]
    \node[] (a) at (0.85,0.5) {$\Phi_e$};
    \draw[color=black, {Stealth[scale=1]}-] (0,0) -- (1.8,0);
    \draw[color=black, {Stealth[scale=1]}-] (0,-0.7) -- (1.8,-0.7);
    \draw[color=black, {Stealth[scale=1]}-] (0,-1.3) -- (1.8,-0.7);
    \draw[color=black, {Stealth[scale=1]}-] (0,-0.8) -- (1.8,-1.4);
    \draw[color=black, {Stealth[scale=1]}-] (0,-2.6) -- (1.8,-1.4);
    \draw[color=black, {Stealth[scale=1]}-] (0,-2.1) -- (1.8,-2.1);
    \draw[color=black, {Stealth[scale=1]}-] (0,-2.7) -- (1.8,-2.1); 
\end{scope}

\end{tikzpicture}

    \caption{The split map applied to two maximal cliques.}
    \label{fig_phi_map}
\end{figure}

\begin{lemma}
    For an inner edge $e=(i,j)$ of $G$, the split map $\Phi_e$ is a surjection from the set of maximal cliques in $(G,F)$ to the set of maximal cliques in $(M(G,e),F_e)$. 
\end{lemma}

\begin{proof}
    Given a maximal clique $C$ in $(G,F)$, it is clear from the construction that the routes in $\Phi_e(C)$ are pairwise coherent in $(M(G,e),F_e)$.
    Let $R$ be a route in $M(G,e)$ that is coherent with all the routes in $\Phi_e(C)$.
    We show that $R$ is necessarily in $\Phi_e(C)$. 
    If $R$ is in $C$, then $R$ is also in $\Phi_e(C)$ by construction.
    If $R$ is not in $C$, then $R$ must be of the form $e_s\cup jR$ or $Ri\cup e_t$.
    However, in either case there is a path ($(i,j)\cup jR$ or $Ri\cup (i,j)$) in $G$ that is coherent with all routes in $C$.
    By Lemma~\ref{lem.path_is_extendable}, the path can be extended to a route $R'$ that is coherent with all routes in $C$, and hence necessarily also contained in $C$.
    Now by construction $e\in R'$ and $\varphi_e(R')$ contains $R$. So, $R$ is in $\Phi_e(C)$.
\end{proof}

\begin{lemma}\label{lem_fibers}
    Given a maximal clique $D$ in $(M(G,e),F_e)$, the fiber $\Phi_e^{-1}(D)$ is the set of maximal cliques in $(G,F)$ that are coherent with the set $S$ consisting of the following paths: 
    \begin{enumerate}
        \item the routes of $D$ that are also routes in $G$,
        \item the paths of the form $(i,j)\cup jR$ where $R$ is a route in $D$ passing through $e_s$, and
        \item the paths of the form $R'i\cup (i,j)$ where $R'$ is a route in $D$ passing through $e_t$. 
    \end{enumerate}
    In particular, 
    $\Phi_e^{-1}(D)=[C_{\min}(S),C_{\max}(S)]$.       
\end{lemma}

\begin{proof}
    Let $C$ be a maximal clique of $(G,F)$. We need to show that $\Phi_e(C)=D$ if and only if $C$ is coherent with $S$. 

    For the forward direction, if $\Phi_e(C)=D$ then every path $P$ in $S$ is necessarily a subpath of a route in $C$. Thus, every route in $C$ is coherent with $P$.  

    For the backward direction, if $C$ is coherent with $S$ then we need to show that the routes of $D$ corresponding to the paths of types (1), (2) and (3) in $S$ all belong to $\Phi_e(C)$, from which we can deduce that $\Phi_e(C)=D$. 
    For type (1) paths, consider a route $R\in D$ that is also a route in $G$. Since $R$ is coherent with $C$ then $R\in C$ by maximality, and $\Phi_e(R)=R\in D$. 
    For type (2) paths, consider a route $R$ in $D$ passing through $e_s$. Since the path $(i,j)\cup jR$ is coherent with~$C$ then it can be extended to a route $R''\in C$ and $\Phi_e(R'')=R\in D$. 
    Similarly for type (3) paths, consider a route $R'$ in $D$ passing through $e_t$. Since the path $R'i\cup (i,j)$ is coherent with~$C$, it can be extended to a route $R''\in C$, and then $\Phi_e(R'')=R'\in D$. 
    
    This finishes the proof that $C$ is coherent with $S$ if and only if $\Phi_e(C)=D$.
    Finally, $S$ is a set of paths in $(G,F)$ and $\Phi_e^{-1}(D)$ is the set of maximal cliques that are coherent with $S$. By~\Cref{cor.interval}, we deduce that $\Phi_e^{-1}(D)=[C_{\min}(S),C_{\max}(S)]$.  
\end{proof}

Let $\balpha(e)$ be the equivalence relation on $\scrL_{G,F}$ induced by $x\equiv y$ if and only if $x$ and $y$ belong to the same fiber of $\Phi_e$. 

\begin{corollary}\label{cor_classes_intervals}
    Each equivalence class of $\balpha(e)$ is an interval in $\scrL_{G,F}$.
\end{corollary}

\begin{proof}
    Each equivalence class is of the form $\Phi_e^{-1}(D)$ for some maximal clique $D$ in $(M(G,e),F_e)$. As we have just proved $\Phi_e^{-1}(D)=[C_{\min}(S),C_{\max}(S)]$, which is an interval.  
\end{proof}

\subsection{Lattice congruences}
Or next goal is to show that the equivalence relation $\balpha(e)$ is a lattice congruence of the framing lattice~$\scrL_{G,F}$ (\Cref{cor_lattice_congruence}). 
In order to do this we need to understand the projections mapping each element of $\scrL_{G,F}$ to the bottom element and to the top element in its equivalence class. The following lemma will be useful for that purpose. 

\begin{lemma}\label{lem_fiber_coherent_routes}
    Let $C$ be a maximal clique in $(G,F)$, $D=\Phi_e(C)$ and $S$ be the corresponding set of paths as in~\Cref{lem_fibers}. 
    A route in $G$ is coherent with $S$ if and only if it is:
    \begin{enumerate}
        \item a route of $C$ that does not contain the edge $e=(i,j)$, or 
        \item a route of the form $RijR'$ where $R$ and $R'$ are routes in $C$ containing $e$.
    \end{enumerate}
\end{lemma}

\begin{proof}
    We start by proving the forward direction. 
    Let $R''$ be a route in $G$ that is coherent with~$S$. We will show that $R''$ is of type (1) or (2).   

    Type (1): Assume that $R''$ does not contain the edge $e=(i,j)$. If $R''\notin C$ then there would be a route $R \in C$ that is not coherent with $R''$. If the edge $e$ is not in $R$ then $R\in S$, which is a contradiction. If the edge $e$ is in $R$ then $R''$ and $R$ are incoherent outside $e$, and so $R''$ is incoherent with either $(i,j)\cup jR$ or $Ri\cup (i,j)$; since both paths are in $S$, this is also a contradiction. As a consequence $R''\in C$, which means that $R''$ is of type (1). 

    Type (2): Now assume that $R''$ contains the edge $e=(i,j)$. The subpath $R''j$ must be coherent with $C$, otherwise it would be incoherent with $S$. By~\Cref{lem.path_is_extendable}, the path $R''j$ can be extended to a route coherent with $C$, and by maximality this extended route must be equal to a route $R$ of $C$. Therefore, $R''j=Rj$ for some $R\in C$. Similarly, we can show that $iR''=iR'$ for some route $R'\in C$. Thus $R''=RijR'$, which is a route of type (2).      

    For the backward direction we will show that any route of type (1) or (2) is coherent with~$S$. 

    Type (1): The paths in $S$ are either routes in $C$ not containing $e$ or subpaths of routes of~$C$ containing $e$. Since every route of $C$ is coherent with $C$ then it is also coherent with $S$. 

    Type (2): Let $R$ and $R'$ be routes in $C$ containing $e$. Towards a contradiction, assume there is $R''\in S$ that is incoherent with $RijR'$. If $R''$ does not contain $e$ then $R''\in C$. Furthermore, $R''$ must be incoherent with $RijR'$ outside $e$. This implies that $R''$ is incoherent with either $R$ or $R'$, which is a contradiction. If $R''$ contains $e$ then $R''$ is a path that starts or ends with the edge $e$. In particular, the subpath where $R''$ and $RijR'$ are incoherent does not contain the edge $e$, and so $R''$ is incoherent with either $R$ or $R'$. Since $R''$ is a subpath of a route of $C$ we also have a contradiction.   
\end{proof}

We focus our attention on the routes of the form $RijR'$ where $R$ and $R'$ are routes in~$C$ containing $e$. 
We denote by $R_{C,e}^{\cw}=RijR'$ the route obtained when $Rj$ is minimal with respect to~$\leq_{\scrI(j)}$ and $iR'$ is maximal with respect to~$\leq_{\scrO(i)}$; this is the clockwise most route through $e$ that is coherent with~$S$.
Similarly, we denote by $R_{C,e}^{\ccw}=RijR'$ the route obtained when $Rj$ is maximal with respect to $\leq_{\scrI(j)}$ and $iR'$ is minimal with respect to $\leq_{\scrO(i)}$; this is the counterclockwise most route through $e$ that is coherent with $S$.

For the equivalence relation $\balpha(e)$ on $\scrL_{G,F}$, the \defn{down projection $\pi^{\balpha(e)}_{\downarrow}$} (resp. \defn{up projection $\pi^{\balpha(e)}_{\uparrow}$}) is defined as the map sending each element to the bottom (resp. top) element in its equivalence class.

\begin{lemma}[Down projection]\label{lem_downprojection}
    The projection $\pi^{\balpha(e)}_{\downarrow}(C)$ of a maximal clique $C$ in $(G,F)$ can be described equivalently as either:
    \begin{enumerate}
        \item  the unique maximal clique containing the routes of $C$ without $e$ and the route $R_{C,e}^{\cw}$, or
        \item the maximal clique containing the routes of $C$ without $e$, and the routes of the form $R_{C,e}^{\cw}iR$ and $RjR_{C,e}^{\cw}$, where $R$ is a route in $C$ containing $e$.
    \end{enumerate}
\end{lemma}

\begin{proof}
Let $D=\Phi_e(C)$ and $S$ be the corresponding set of paths as in~\Cref{lem_fibers}. 
We know that the fiber containing $C$ is the set
$\Phi_e^{-1}(D)=[C_{\min}(S),C_{\max}(S)]$ 
of maximal cliques that are coherent with $S$, and so $\pi^{\balpha(e)}_{\downarrow}(C)=C_{\min}(S)$. 

Let $C'\in \Phi_e^{-1}(D)$ be a maximal clique in the fiber containing $C$. We will show first that if $R_{C,e}^{\cw} \in C'$ then $C'$ is the maximal clique containing the routes of $C$ without $e$, and the routes of the form $R_{C,e}^{\cw}iR$ and $RjR_{C,e}^{\cw}$, where $R$ is a route in $C$ containing $e$. Our second step in the proof will be to show that this $C'$ is equal to $\pi^{\balpha(e)}_{\downarrow}(C)$. 

For the first step, assume that $R_{C,e}^{\cw} \in C'$.
By~\Cref{lem_fiber_coherent_routes}, the routes of $C'$ are the routes of $C$ without $e=(i,j)$, together with some routes of the form $RijR'$ where $R$ and $R'$ are routes of $C$ containing $e$. But if neither $R$ nor $R'$ is equal to $R_{C,e}^{\cw}$, then $RijR'$ is incoherent with $R_{C,e}^{\cw}$, which is a contradiction. 
Therefore, every route of $C'$ containing $e$ is of the form $R_{C,e}^{\cw}iR$ or $RjR_{C,e}^{\cw}$, where $R$ is a route in $C$ containing $e$. 
Since all of these routes are pairwise coherent, it follows that $C'$ must contain all of them by maximality. 

For the second step, consider a maximal clique $C''\neq C'$ in the fiber $\Phi_e^{-1}(D)$. There must be a route $RijR'$ in $C''\setminus C'$ for some $R,R'\in C$ containing $e$. But then $R\neq R_{C,e}^{\cw}$ and $R'\neq R_{C,e}^{\cw}$, which implies that the route $R_{C,e}^{\cw} \in C'$ is clockwise incoherent from the route $RijR'\in C''$ at $e$. So, $C''$ cannot be smaller than $C'$. Since the fiber has a unique minimum element, this minimum must be equal to $C'$. Thus, $\pi^{\balpha(e)}_{\downarrow}(C)=C'$.  
\end{proof}

By symmetry we get an analogous result for the up projection.

\begin{lemma}[Up projection]\label{lem_upprojection}
    The projection $\pi^{\balpha(e)}_{\uparrow}(C)$ of a maximal clique $C$ in $(G,F)$ can be described equivalently as either:
    \begin{enumerate}
        \item  the unique maximal clique containing the routes of $C$ without $e$ and the route $R_{C,e}^{\ccw}$, or
        \item the maximal clique containing the routes of $C$ without $e$, and the routes of the form $R_{C,e}^{\ccw}iR$ and $RjR_{C,e}^{\ccw}$, where $R$ is a route in $C$ containing $e$.
    \end{enumerate}
\end{lemma}

\begin{corollary}\label{cor_projections_order_preserving}
    The following hold:
    \begin{enumerate}
        \item The map $\pi^{\balpha(e)}_{\downarrow}$ is order preserving.
        \item The map $\pi^{\balpha(e)}_{\uparrow}$ is order preserving.
    \end{enumerate}
\end{corollary}

\begin{proof}
    To prove (1), let $C$ and $C'$ be maximal cliques in $\scrL_{G,F}$ satisfying $C\leq C'$. We need to show that $\pi^{\balpha(e)}_{\downarrow}(C) \leq \pi^{\balpha(e)}_{\downarrow}(C')$. 
    Let $R$ be a route in $\pi^{\balpha(e)}_{\downarrow}(C)$ and $R'$ be a route in~$\pi^{\balpha(e)}_{\downarrow}(C')$.
    By~\Cref{lem_downprojection}, the route $R$ (resp. $R'$) behaves purely as a route of $C$ (resp. $C'$) outside the edge $e=(i,j)$. So, if $R$ and $R'$ are incoherent outside $e$ then $R$ would be clockwise incoherent from $R'$ (because $C\leq C'$), which is what we want to prove.  
    So, it is enough to consider the case where $R$ and~$R'$ are incoherent at a path containing $e$ (in particular, $e\in R$ and $e\in R'$). 
    In this case, it suffices to check that $R_{C,e}^{\cw}j \leq_{\scrI(j)} R_{C',e}^{\cw}j$ and $iR_{C',e}^{\cw} \leq_{\scrO(i)} iR_{C,e}^{\cw}$.
    
    Suppose toward a contradiction that $R_{C',e}^{\cw}j <_{\scrI(j)} R_{C,e}^{\cw}j$.
    Let $P$ be the maximal path in $R_{C,e}^{\cw}j \cap R_{C',e}^{\cw}j$ containing $j$, and let $k$ be the minimal vertex in $P$.
    Consider the path $P'$ formed from $P$ by appending the edge of $R_{C',e}^{\cw}j$ incoming to $k$.
    Now since $R_{C,e}^{\cw}j$ is minimal with respect to $\leq_{\scrI(j)}$, the path $P'$ must be incoherent with some route $R^*$ in $C$ at $k$ (otherwise it can be extended to a route in $C$ that is smaller in the order $\leq_{\scrI(j)}$).
    However, such a route~$R^*$ is counter clockwise from the routes in $C'$ containing $R_{C',e}^{\cw}j$, which contradicts the fact that $C\leq C'$. 
    This shows that $R_{C,e}^{\cw}j \leq_{\scrI(j)} R_{C',e}^{\cw}j$. 
    The proof that $iR_{C',e}^{\cw} \leq_{\scrO(i)} iR_{C,e}^{\cw}$ uses similar arguments.

    The proof of (2), i.e. that $\pi_{\uparrow}^{\balpha(e)}$ is order preserving, is similar to the proof of (1).
\end{proof}

\begin{corollary}\label{cor_lattice_congruence}
    The equivalence relation $\balpha(e)$ is a lattice congruence on $\scrL_{G,F}$.
\end{corollary}

\begin{proof}
    The equivalence relation $\balpha(e)$ is a lattice congruence if and only if each equivalence class in an interval and the maps $\pi^{\balpha(e)}_{\downarrow}$ and $\pi^{\balpha(e)}_{\uparrow}$ are order preserving (see e.g.~\cite[Proposition~9-5.2]{Rea16}). These properties were proven in~\Cref{cor_classes_intervals} and~\Cref{cor_projections_order_preserving}.
\end{proof}

\subsection{Lattice quotients}
Our next goal is to show that the lattice quotient induced by the equivalence relation $\balpha(e)$ on $\scrL_{G,F}$ is isomorphic to the framing lattice $\scrL_{M(G,e),F_e}$ (\Cref{cor_lattice_quotient}).

\begin{proposition}\label{prop_iso_quotient}
    Let $C$ and $C'$ be two minimal elements of the $\balpha(e)$ equivalence classes on $\scrL_{G,F}$, and let $D=\Phi_e(C)$ and $D'=\Phi_e(C')$ be their images on $\scrL_{M(G,e),F_e}$. Then
    $C\leq C'$ if and only $D\leq D'$.
\end{proposition}

\begin{proof}
    For the forward direction we do not need to assume that $C$ and $C'$ are minimal elements in their equivalence classes. 
    If $C$ and $C'$ are maximal cliques of $\scrL_{G,F}$ satisfying $C\leq C'$ then it follows directly that $\Phi_e(C)\leq \Phi_e(C')$. 

    For the backward direction, assume $D\leq D'$ in $\scrL_{M(G,e),F_e}$. We will show that $C\leq C'$. 
    Since $C$ is the minimal element of its class, we have that $C=\pi^{\balpha(e)}_{\downarrow}(C)$. By~\Cref{lem_downprojection}, the route $R_{C,e}^{\cw}\in C$ and every route of $C$ containing $e$ is of the form $R_{C,e}^{\cw}iR$ or $RjR_{C,e}^{\cw}$, where $R$ is a route containing $e$. 
    The similar statement holds for $C'$.

    Let $R_{D,e_s}^{\cw}$ (resp. $R_{D,e_t}^{\cw}$) be the clockwise most route of $D$ passing through the edge $e_s$ (resp. $e_t$), and define these routes similarly for $D'$. Since $D=\Phi_e(C)$ and $D'=\Phi_e(C')$ then  

    \begin{center}        
    \begin{tabular}{cc}
         $R_{C,e}^{\cw}i=R_{D,e_t}^{\cw}i$, \quad &  $R_{C',e}^{\cw}i=R_{D',e_t}^{\cw}i$ \\
         $jR_{C,e}^{\cw}=jR_{D,e_s}^{\cw}$, \quad & $jR_{C',e}^{\cw}=jR_{D',e_s}^{\cw}$
    \end{tabular}
    \end{center}
    
    Since the routes of $C$ (resp. $C'$) behave as routes of $D$ (resp. $D'$) outside of $e$, if a route of $C$ is incoherent with a route of $C'$ outside $e$ then they are incoherent in the right order because $D\leq D'$. Therefore, its is enough to consider the case where they are incoherent at a path containing the edge $e$. 
    In this case, it suffices to check that 
    $R_{C,e}^{\cw}j \leq_{\scrI(j)} R_{C',e}^{\cw}j$ and 
    $iR_{C',e}^{\cw} \leq_{\scrO(i)} iR_{C,e}^{\cw}$ 
    (as in the proof of~\Cref{cor_projections_order_preserving}).
    By the equations above, these two conditions are equivalent to 
    $R_{D,e_t}^{\cw} \leq_{\scrI(t)} R_{D',e_t}^{\cw}$ and 
    $R_{D',e_s}^{\cw} \leq_{\scrO(s)} R_{D,e_s}^{\cw}$.

    Following the same lines as in the proof of~\Cref{cor_projections_order_preserving}, 
    suppose toward a contradiction that $R_{D',e_t}^{\cw} <_{\scrI(t)} R_{D,e_t}^{\cw}$. 
    Let $P$ be the maximal path in $R_{D,e_t}^{\cw} \cap R_{D',e_t}^{\cw}$ containing $t$, and let $k$ be the minimal vertex in $P$.
    Consider the path $P'$ formed from $P$ by appending the edge of $R_{D',e_t}^{\cw}$ incoming to $k$.
    Now since $R_{D,e_t}^{\cw}$ is the minimal route in $D$ containing $e_t$ with respect to $\leq_{\scrI(t)}$, the path $P'$ must be incoherent with some route $R^*$ in $D$ at $k$ (otherwise it can be extended to a route in $D$ that is smaller in the order $\leq_{\scrI(t)}$).
    However, such a route~$R^*$ is counter clockwise from the route $R_{D',e_t}^{\cw}$ in $D'$ at $k$, which contradicts the fact that $D\leq D'$. 
    This shows that $R_{D,e_t}^{\cw} \leq_{\scrI(t)} R_{D',e_t}^{\cw}$. 
    The proof that $R_{D',e_s}^{\cw} \leq_{\scrO(s)} R_{D,e_s}^{\cw}$ uses similar arguments.
 \end{proof}

\begin{corollary}\label{cor_lattice_quotient}
    For any inner edge $e$ of $G$ we have $\scrL_{G,F}/\balpha(e) \cong \scrL_{M(G,e),F_e}$.
\end{corollary}
\begin{proof}
    The lattice quotient $\scrL_{G,F}/\balpha(e)$ is isomorphic to the restriction of $\scrL_{G,F}$ to the set of minimal elements of the equivalence classes (see e.g.~\cite[Proposition~9-5.5]{Rea16}). This restricted poset is isomorphic to $\scrL_{M(G,e),F_e}$ by~\Cref{prop_iso_quotient}. 
\end{proof}

\begin{remark}
    It would be interesting to study all lattice quotients of the framing lattice, not just the ones arising from $M$-moves. In particular, we do not know if every quotient is itself a framing lattice of some modified graph. We leave this as an open question for the interested reader. 
\end{remark}

\subsection{A distributive quotient}
An $M$-move decreases the number of inner edges in a framed graph $(G,F)$ by one, and it can be repeated until the resulting graph has no more inner edges.
We use $M(G)$ to denote the graph obtained from $G$ by repeating $M$-moves until there are no more inner edges.
Note that the induced framing of the graph $M(G)$, which we denote by $M(F)$, is independent of the order of the performed $M$-moves.   
Moreover, since all edges in $M(G)$ are incident to the source or sink, it follows from Lemma~\ref{lem.iso_operations} (parts (3) and (4)) that the framing lattice $\scrL_{M(G),M(F)}$ is independent of the initial framing $F$. Thus we abbreviate it by simply writing $\scrL_{M(G)}$.

In this section, we present a simple combinatorial description of the lattice quotient $\scrL_{M(G)}$ as a product of certain distributive lattices described in terms of lattice paths. In particular, this will imply that $\scrL_{M(G)}$ is a distributive lattice quotient of the framing lattice.

\subsubsection{A distributive lattice on lattice paths}
For each inner vertex $v$ of $M(G)$, let $G_v$ denote the subgraph of $M(G)$ induced by $\{s,v,t\}$. The framing lattice of $G_v$ is independent of its framing, so we can denote it by $\scrL_{G_v}$. 
Furthermore, we have that 
\begin{align}\label{eq_distributive_product}
\scrL_{M(G)} = \prod_{v\in G\setminus \{s,t\}} \scrL_{G_v}.    
\end{align}
Hence, it is enough to concentrate on the lattices $\scrL_{G_v}$.

Let $a+1$ be the in-degree of $v$ in $G_v$, and $b+1$ be its out-degree. 
The flow polytope $\calF_{G_v}$ is the product $\Delta_a\times \Delta_b$ of two simplices, where $\Delta_a:=\conv\{e_1,\dots,e_{a+1}\}\subseteq \mathbb{R}^{a+1}$. The framed triangulation of $\calF_{G_v}$ is a well known triangulation of the product of two simplices, called the staircase triangulation~\cite{triangulations_book}. As we will see in a more general setting in~\Cref{sec.multipermutations}, the framing lattice $\scrL_{G_v}$ is a lattice on the set of multipermutations of $1^a2^b$, with cover relations given by replacing a consecutive pair $12$ by $21$. Two examples are illustrated in~\Cref{fig_lattice_paths_112_1122}.

\begin{figure}[h]
    \centering
    \begin{tikzpicture}

\begin{scope}[xshift = -55, yshift= -20, scale=0.5]

    \draw[thick, color=black] (2,0) -- (4,0);
    \draw[thick, color=black] (2,0) .. controls (2.4, 0.6) and (3.6, 0.6) .. (4,0);
    \draw[thick, color=black] (2,0) .. controls (2.4, -0.6) and (3.6, -0.6) .. (4,0);
    \draw[thick, color=black] (4,0) .. controls (4.4, 0.6) and (5.6, 0.6) .. (6,0);
    \draw[thick, color=black] (4,0) .. controls (4.4, -0.6) and (5.6, -0.6) .. (6,0);     

	\node[circle, draw, inner sep=1pt, fill](2) at (2,0) {};
	\node[circle, draw, inner sep=1pt, fill](3) at (4,0) {};
	\node[circle, draw, inner sep=1pt, fill](4) at (6,0) {};	
    
    \node[] (a) at (2,-0.5) {\scriptsize $s$};
    \node[] (a) at (4,-0.5) {\scriptsize $v$};    
    \node[] (a) at (6,-0.5) {\scriptsize $t$};
\end{scope}

\begin{scope}[xshift = 140, yshift= -20, scale=0.5]

	\draw[thick, color=black] (0,0) -- (2,0);
    \draw[thick, color=black] (2,0) -- (4,0);
	\draw[thick, color=black] (0,0) .. controls (0.4, 0.6) and (1.6, 0.6) .. (2,0);
    \draw[thick, color=black] (0,0) .. controls (0.4, -0.6) and (1.6, -0.6) .. (2,0);    
    \draw[thick, color=black] (2,0) .. controls (2.4, 0.6) and (3.6, 0.6) .. (4,0);
    \draw[thick, color=black] (2,0) .. controls (2.4, -0.6) and (3.6, -0.6) .. (4,0);

    \node[circle, draw, inner sep=1pt, fill](1) at (0,0) {};	
	\node[circle, draw, inner sep=1pt, fill](2) at (2,0) {};
	\node[circle, draw, inner sep=1pt, fill](3) at (4,0) {};
    
    \node[] (a) at (0,-0.5) {\scriptsize $s$};
    \node[] (a) at (2,-0.5) {\scriptsize $v$};    
    \node[] (a) at (4,-0.5) {\scriptsize $t$};
\end{scope}

\begin{scope}[xshift = 0, yshift=0, scale=1.2]
    \node[circle,fill,inner sep=1.5pt,color=NavyBlue] (g) at (0,1)  {};	  
    \node[circle,fill,inner sep=1.5pt,color=NavyBlue] (h) at (0,2)  {};
    \node[circle,fill,inner sep=1.5pt,color=NavyBlue] (i) at (0,3)  {};	      

    \draw[color=NavyBlue] (g) -- (h) -- (i);

    \node[] (112) at (0.5,1) {\small $112$};
    \node[] (121) at (0.5,2) {\small $121$};
    \node[] (211) at (0.5,3) {\small $211$};

    \node[] (path112) at (-0.9,2) {
    \begin{tikzpicture}[scale=0.4]
        \draw[very thin, color=gray!70] (0,0) grid (2,1);
        \draw[very thick, color=RawSienna!70] (0,0) -- (1,0) -- (1,1) -- (2,1);
        \node[circle,fill,inner sep=1pt,color=black] (i) at (0,0)  {};	      
        \node[circle,fill,inner sep=1pt,color=black] (i) at (1,0)  {};	      
        \node[circle,fill,inner sep=1pt,color=black] (i) at (1,1)  {};	      
        \node[circle,fill,inner sep=1pt,color=black] (i) at (2,1)  {};

        \node[] (k) at (0.5,-0.5) {\tiny $1$};
        \node[] (k) at (1.5,1.5) {\tiny $1$};        
        \node[] (k) at (1.4,0.5) {\tiny $2$};        
    \end{tikzpicture}
    };
    \node[] (path211) at (-1,3.1) {
    \begin{tikzpicture}[scale=0.4]
        \draw[very thin, color=gray!70] (0,0) grid (2,1);
        \draw[very thick, color=RawSienna] (0,0) -- (0,1) -- (1,1) -- (2,1);
        \node[circle,fill,inner sep=1pt,color=black] (i) at (0,0)  {};	      
        \node[circle,fill,inner sep=1pt,color=black] (i) at (0,1)  {};	      
        \node[circle,fill,inner sep=1pt,color=black] (i) at (1,1)  {};	      
        \node[circle,fill,inner sep=1pt,color=black] (i) at (2,1)  {};

        \node[] (k) at (0.5,1.5) {\tiny $1$};
        \node[] (k) at (1.5,1.5) {\tiny $1$};        
        \node[] (k) at (-0.5,0.4) {\tiny $2$};        
    \end{tikzpicture}
    };
    \node[] (path112) at (-0.8,0.9) {
    \begin{tikzpicture}[scale=0.4]
        \draw[very thin, color=gray!70] (0,0) grid (2,1);
        \draw[very thick, color=RawSienna!70] (0,0) -- (1,0) -- (2,0) -- (2,1);
        \node[circle,fill,inner sep=1pt,color=black] (i) at (0,0)  {};	      
        \node[circle,fill,inner sep=1pt,color=black] (i) at (1,0)  {};	      
        \node[circle,fill,inner sep=1pt,color=black] (i) at (2,0)  {};	      
        \node[circle,fill,inner sep=1pt,color=black] (i) at (2,1)  {};

        \node[] (k) at (0.5,-0.5) {\tiny $1$};
        \node[] (k) at (1.5,-0.5) {\tiny $1$};        
        \node[] (k) at (2.4,0.5) {\tiny $2$};        
    \end{tikzpicture}
    };
\end{scope}

\begin{scope}[xshift = 170, yshift=20, scale=1]
	\node[circle,fill,inner sep=1.5pt,color=NavyBlue] (a) at (0,0)  {};	
	\node[circle,fill,inner sep=1.5pt,color=NavyBlue] (b) at (0,1)  {};	    
	\node[circle,fill,inner sep=1.5pt,color=NavyBlue] (c) at (-0.7,2)  {};	
	\node[circle,fill,inner sep=1.5pt,color=NavyBlue] (d) at (0.7,2)  {};		
	\node[circle,fill,inner sep=1.5pt,color=NavyBlue] (e) at (0,3)  {};	    
	\node[circle,fill,inner sep=1.5pt,color=NavyBlue] (f) at (0,4)  {};	    

    \draw[color=NavyBlue] (a) -- (b) -- (c) -- (e) -- (f);
	\draw[color=NavyBlue] (b) -- (d) -- (e);

    \node[] (1122) at (0.6,0) {\small $1122$};
    \node[] (1212) at (0.6,1) {\small $1212$};
    \node[] (1221) at (1.3,2) {\small $1221$};
    \node[] (2112) at (-1.3,2) {\small $2112$}; 
    \node[] (2121) at (0.6,3) {\small $2121$};
    \node[] (2211) at (0.6,4) {\small $2211$};    

    \node[] (path1122) at (1.6,-0.4) {
    \begin{tikzpicture}[scale=0.35]
        \draw[very thin, color=gray!70] (0,0) grid (2,2);
        \draw[very thick, color=RawSienna!70] (0,0) -- (1,0) -- (2,0) -- (2,1) -- (2,2);
        \node[circle,fill,inner sep=1pt,color=black] (i) at (0,0)  {};	      
        \node[circle,fill,inner sep=1pt,color=black] (i) at (1,0)  {};	      
        \node[circle,fill,inner sep=1pt,color=black] (i) at (2,0)  {};	      
        \node[circle,fill,inner sep=1pt,color=black] (i) at (2,1)  {};
        \node[circle,fill,inner sep=1pt,color=black] (i) at (2,2)  {};        

        \node[] (k) at (0.5,-0.5) {\tiny $1$};
        \node[] (k) at (1.5,-0.5) {\tiny $1$};        
        \node[] (k) at (2.4,0.5) {\tiny $2$};        
        \node[] (k) at (2.4,1.5) {\tiny $2$};                
    \end{tikzpicture}
    };
    \node[] (path1212) at (1.6,0.7) {
    \begin{tikzpicture}[scale=0.35]
        \draw[very thin, color=gray!70] (0,0) grid (2,2);
        \draw[very thick, color=RawSienna!70] (0,0) -- (1,0) -- (1,1) -- (2,1) -- (2,2);
        \node[circle,fill,inner sep=1pt,color=black] (i) at (0,0)  {};	      
        \node[circle,fill,inner sep=1pt,color=black] (i) at (1,0)  {};	      
        \node[circle,fill,inner sep=1pt,color=black] (i) at (1,1)  {};	      
        \node[circle,fill,inner sep=1pt,color=black] (i) at (2,1)  {};
        \node[circle,fill,inner sep=1pt,color=black] (i) at (2,2)  {};        

        \node[] (k) at (0.5,-0.5) {\tiny $1$};
        \node[] (k) at (1.5,1.5) {\tiny $1$};        
        \node[] (k) at (1.4,0.5) {\tiny $2$};        
        \node[] (k) at (2.4,1.5) {\tiny $2$};                
    \end{tikzpicture}
    };
    \node[] (path1221) at (2.3,2) {
    \begin{tikzpicture}[scale=0.35]
        \draw[very thin, color=gray!70] (0,0) grid (2,2);
        \draw[very thick, color=RawSienna!70] (0,0) -- (1,0) -- (1,1) -- (1,2) -- (2,2);
        \node[circle,fill,inner sep=1pt,color=black] (i) at (0,0)  {};	      
        \node[circle,fill,inner sep=1pt,color=black] (i) at (1,0)  {};	      
        \node[circle,fill,inner sep=1pt,color=black] (i) at (1,1)  {};	      
        \node[circle,fill,inner sep=1pt,color=black] (i) at (1,2)  {};
        \node[circle,fill,inner sep=1pt,color=black] (i) at (2,2)  {};        

        \node[] (k) at (0.5,-0.5) {\tiny $1$};
        \node[] (k) at (1.5,2.5) {\tiny $1$};        
        \node[] (k) at (1.4,0.5) {\tiny $2$};        
        \node[] (k) at (1.4,1.5) {\tiny $2$};                
    \end{tikzpicture}
    };
    \node[] (path2121) at (1.6,3) {
    \begin{tikzpicture}[scale=0.35]
        \draw[very thin, color=gray!70] (0,0) grid (2,2);
        \draw[very thick, color=RawSienna!70] (0,0) -- (0,1) -- (1,1) -- (1,2) -- (2,2);
        \node[circle,fill,inner sep=1pt,color=black] (i) at (0,0)  {};	      
        \node[circle,fill,inner sep=1pt,color=black] (i) at (0,1)  {};	      
        \node[circle,fill,inner sep=1pt,color=black] (i) at (1,1)  {};	      
        \node[circle,fill,inner sep=1pt,color=black] (i) at (1,2)  {};
        \node[circle,fill,inner sep=1pt,color=black] (i) at (2,2)  {};        

        \node[] (k) at (0.5,1.5) {\tiny $1$};
        \node[] (k) at (1.5,2.5) {\tiny $1$};        
        \node[] (k) at (0.4,0.5) {\tiny $2$};        
        \node[] (k) at (1.4,1.5) {\tiny $2$};                
    \end{tikzpicture}
    };
    \node[] (path2211) at (1.6,4.3) {
    \begin{tikzpicture}[scale=0.35]
        \draw[very thin, color=gray!70] (0,0) grid (2,2);        
        \draw[very thick, color=RawSienna!70] (0,0) -- (0,1) -- (0,2) -- (1,2) -- (2,2);
        \node[circle,fill,inner sep=1pt,color=black] (i) at (0,0)  {};	      
        \node[circle,fill,inner sep=1pt,color=black] (i) at (0,1)  {};	      
        \node[circle,fill,inner sep=1pt,color=black] (i) at (0,2)  {};	      
        \node[circle,fill,inner sep=1pt,color=black] (i) at (1,2)  {};
        \node[circle,fill,inner sep=1pt,color=black] (i) at (2,2)  {};        

        \node[] (k) at (0.5,2.5) {\tiny $1$};
        \node[] (k) at (1.5,2.5) {\tiny $1$};        
        \node[] (k) at (0.4,0.5) {\tiny $2$};        
        \node[] (k) at (0.4,1.5) {\tiny $2$};                
    \end{tikzpicture}
    };
    \node[] (path2112) at (-2.3,2) {
    \begin{tikzpicture}[scale=0.35]
        \draw[very thin, color=gray!70] (0,0) grid (2,2);    
        \draw[very thick, color=RawSienna!70] (0,0) -- (0,1) -- (1,1) -- (2,1) -- (2,2);
        \node[circle,fill,inner sep=1pt,color=black] (i) at (0,0)  {};	      
        \node[circle,fill,inner sep=1pt,color=black] (i) at (0,1)  {};	      
        \node[circle,fill,inner sep=1pt,color=black] (i) at (1,1)  {};	      
        \node[circle,fill,inner sep=1pt,color=black] (i) at (2,1)  {};
        \node[circle,fill,inner sep=1pt,color=black] (i) at (2,2)  {};        

        \node[] (k) at (0.5,1.5) {\tiny $1$};
        \node[] (k) at (1.5,1.5) {\tiny $1$};        
        \node[] (k) at (0.4,0.5) {\tiny $2$};        
        \node[] (k) at (2.4,1.5) {\tiny $2$};                
    \end{tikzpicture}
    };
\end{scope}
\end{tikzpicture}
    \caption{Two examples of the distributive lattice $\scrL_{G_v}$}
    \label{fig_lattice_paths_112_1122}
\end{figure}

Alternatively, the elements of $\scrL_{G_v}$ can be described as lattice paths in the plane from $(0,0)$ to $(a,b)$ using unit East steps and unit North steps. Each mutlipermutation of $1^a2^b$ can be transformed in such a lattice path by replacing each $1$ by an $E$ step and each $2$ by a $N$ step. The covering relation can be then described as adding a box to the path. This resulting poset is known to be a distributive lattice.

As a side remark, note that the number of elements of $\scrL_{G_v}$ is then equal to $\binom{a+b}{a}$, which in turn is the volume of the product of two simplices $\Delta_a\times \Delta_b$.

\begin{theorem}\label{thm_distributive_quotient}
    For a framed graph $(G,F)$, the lattice $\scrL_{M(G)}$ is a distributive lattice quotient of $\scrL_{G,F}$ independent of the choice of $F$. 
\end{theorem}

\begin{proof}
    As discussed above, $\scrL_{M(G)}$ is independent of the initial framing $F$ of $G$. So, it remains to show that it is a distributive lattice. Since the product of distributive lattices is distributive, this follows from Equation~\eqref{eq_distributive_product} and the fact that $\scrL_{G_v}$ is distributive. 
\end{proof}

\subsubsection{A Boolean lattice on lattice quotients via $M$-moves}
One may also be interested in the lattice quotients of $\scrL_{G,F}$ obtained by applying $M$-moves to a subset $A$ of internal edges of $G$. We denote by $(M(G,A),F_A)$ be the resulting framed graph; note that the framing $F_A$ is independent of the order in which we apply the $M$-moves for edges in $A$. We can order such lattice quotients by declaring $\scrL_{M(G,A),F_A}\leq \scrL_{M(G,B),F_B}$ whenever $B \subseteq A$. 
The resulting poset is a Boolean lattice, with the original framing lattice $\scrL_{G,F}$ on the top and the distributive quotient $\scrL_{M(G)}$ on the bottom. 
An example is illustrated in~\Cref{fig_quotients11223}.

\begin{figure}[h]
\centering

    \caption{The Boolean lattice of M-moves.}
    \label{fig_quotients11223}
\end{figure}

\begin{example}
    Let $G$ be the graph shown on the top of ~\Cref{fig_quotients11223}. It is a special case of the multioruga graphs discussed in ~\Cref{sec.multipermutations}. 
    The framing lattice $\scrL_{G,F}$ is the lattice of multipermutations of $11223$, with cover relations given by swapping consecutive elements~$i<j$, see~\Cref{fig.lattice_of_multipermutations}. There are three internal edges $e,e',e''$ (corresponding to label~2) colored yellow, red and purple from bottom to top in the figure, respectively. 
    Making an $M$-move on $e$ corresponds to contracting the edges swapping $1$ and $3$ before the first appearance of $2$; these are the four yellow edges $1\mathbf{13}22-1\mathbf{31}22$, $\mathbf{13}122-\mathbf{31}122$, $\mathbf{13}212-\mathbf{31}212$, $\mathbf{13}221-\mathbf{31}221$. For $e'$ we contract the edges swapping $1$ and $3$ in between the first and the second appearance of~$2$ (edges colored red), while for $e''$ we contract edges swapping $1$ and~$3$ after the second appearance of $2$ (edges colored purple). 
\end{example}

\begin{remark}[The multisylvester congruence]
    In general, doing an $M$-move on an internal edge of a multioruga graph corresponds to contracting some edges of the lattice. If the internal edge corresponds to label $b$ and is the $(\ell+1)$th edge from bottom to top, then the edges that are contracted are precisely those swapping any two labels $a$ and $c$ between the $\ell$th and the $(\ell+1)$th appearance of $b$, satisfying $a<b<c$. 
\end{remark}

\begin{remark}[The Sylvester congruence, Tamari lattices and Cambrian lattices]
    As we have seen in~\Cref{ex_weak_via_framing}, the weak order on permutations of $[n]$ is the framing lattice $\scrL_{G,F}$ for the oruga graph $G=\oru{n}$ with the natural framing $F$. Performing $M$-moves to all the lower internal edges gives rise to the lattice congruence identifying permutations $U acV bW \equiv U caV bW$ where $a < b < c$. This is a known lattice congruence of the weak order called the Sylvester congruence, whose corresponding lattice quotient is isomorphic to the Tamari lattice, see~e.g.~\cite[Section~2.1]{PSZ23}.   
    Cambrian lattices of type $A$ can also be obtained as lattice quotients of the weak order via $M$-moves on the oruga graph. Cambrian lattices are indexed by a sequence $\varepsilon\in\{\pm\}^n$ consisting of plus and minuses. We apply one M-move for each sign, except the first and the last, in the graph $\oru{n}$. If the $b$th sign is $+$ (resp.~$-$) then we apply an M-move on the $b$th upper edge (resp. down edge). This operation corresponds to the lattice congruence identifying permutations $U bV acW \equiv U bV caW$ (resp. $U acV bW \equiv U caV bW$) where $a < b < c$.
    The resulting framing lattice after these M-moves is isomorphic the $\varepsilon$-Cambrian lattice.
\end{remark}

\begin{remark}[Permutree lattices]
    \label{rem.permutree_lattices}
    In~\cite{PP18}, Pilaud and Pons introduced a lattice of permutrees associated to a sequence of symbols $\theta=\{\noneCirc{}, \downCirc{}, \upCirc{}, \upDownCirc{}\}^n$. This lattice interpolates between the Boolean lattice and the weak order of permutations of $[n]$, and contains the Tamari lattice and all type $A$ Cambrian lattices as particular examples. They showed that the permutree lattice can be obtained as a lattice quotient of the weak order by the lattice congruence identifying permutations
    \[
    \begin{array}{ll}
    UacVbW \equiv_{\theta} UcaVbW & \text{if } a < b < c \text{ and } \theta_b = \downCirc{} \text{ or } \upDownCirc{}, \\
    UbVacW \equiv_\theta UbVcaW & \text{if } a < b < c \text{ and } \theta_b = \upCirc{} \text{ or } \upDownCirc{}.
    \end{array}
    \]
    As a consequence, the permetree lattice can be obtained as a framing lattice of a graph obtained by applying M-moves to the graph $\oru{n}$. The first and last symbols are irrelevant, and we apply an M-move for every other symbol. Depending on the $b$th symbol we apply an~M-move to the $b$th upper or down edge according to the following rule:   
    \begin{table*}[h]
        \centering
        \begin{tabular}{cc}
           $\noneCirc$  & no M-move  \\
            $\upCirc$ & M-move to the upper edge \\
            $\downCirc$ & M-move to the down edge \\
            $\upDownCirc$ & M-move to the upper and down edge 
        \end{tabular}
    \end{table*}
    
    The first connection between permutrees and triangulations of flow polytopes was discovered by Tamayo in~\cite{Tam23}. 
\end{remark}


\section{Open questions and conjectures}

Note that there are many more lattice quotients of a framing lattice as there are lattice quotients of it obtained via $M$-moves. For instance, the graph $G=\oru{4}$ has four internal edges and the framing lattice $\scrL_{G,F}$ has $2^4=16$ lattice quotients obtained via $M$-moves. However, $\scrL_{G,F}$ is the weak order on permutations of $[4]$, which has 47 lattice congruences (but only 20 up to horizontal and vertical symmetry)~\cite{PS19}. 
As we have shown, every lattice quotient obtained from $M$-moves is itself a framing lattice. But we do not know if every lattice quotient is a framing lattice as well. 
We leave this as an open question, which is interesting even in the simple case of the oruga graph.

\begin{question}
    Is every lattice quotient of a framing lattice isomorphic to a framing lattice?
\end{question}

Given that framing lattices are congruence uniform, they can be constructed by a sequence of Day doublings \cite{Day77}. However, we do not know a framing lattice can be constructed from the singleton lattice as a sequence of doublings such that each lattice in the sequence is itself a framing lattice. In other words, we have the following question.

\begin{question} \label{question.doubling}
    For a framing lattice $\scrL_{G,F}$, is there a sequence of Day doublings $\scrL_0$, $\scrL_1$, $\ldots$, $\scrL_n = \scrL_{G,F}$ such that $\scrL_0$ is the singleton lattice and $\scrL_{i}$ is a framing lattice for each $0\leq i \leq n$? 
\end{question}

Although we expect the following to be a difficult question to answer, it is natural to wonder about an explicit characterization of framing lattices.

\begin{question} \label{question.simplecriterion}
    Is there a simple criterion to decide whether a given lattice is a framing lattice?
\end{question}


We also propose the following purely enumerative conjecture, which is based on extensive computational evidence.

\begin{conjecture}
    \label{conj.linear_intervals}
    If $F_1$ and $F_2$ are framings of the same graph $G$, then $\scrL_{G,F_1}$ and $\scrL_{G,F_2}$ have the same number of linear intervals of length $k$ for  $k=0,1,2,\ldots$. 
\end{conjecture}

For instance,~\Cref{fig.caracol3} shows two different framing lattices of the same caracol graph with two different framings. In both cases,  the number of linear intervals of length zero, one and two are 5,5, and 2 respectively. There are no further intervals of larger length.  
Conjecture~\ref{conj.linear_intervals} is known to hold in the case of alt $\nu$-Tamari lattices. This was shown by~\cite{CC23}, but their methods do not seem to be easily extended even to cross-Tamari lattices.

Our next question is related to the complexity of framing lattices.

\begin{question}
    What is the complexity of finding a shortest path between two maximal cliques in the rotation graph of framing lattices? Is it NP-hard? The complexity problem in the case of the classical Tamari lattice is a widely open renown problem.
\end{question}

\newpage
\part{The zoo of framing lattices}

In this second part, we take a tour through the zoo of framing lattices. 
Framing lattices serve to unify a wide range of lattices previously studied in the literature, namely, as lattices on maximal cliques of coherent routes in framed graphs.
Our goal here is twofold. 
First, we exhibit various combinatorial lattices from the framing lattice perspective, adding to the list of running examples from Part I. 
In particular, we highlight four families: the Boolean lattices, the lattice of multipermutations, $\varepsilon$-Cambrian lattices, and a new family of lattices we call cross-Tamari lattices (see Figure~\ref{fig.graphs_and_lattices}).
Second, we demonstrate the natural taxonomy of these lattice families provided by varying the two underlying parameters of framing lattices: the flow graph and the framing. 
A beautiful example of this is the family of cross-Tamari lattices (Section~\ref{sec_crossTamari}), which arises naturally from the family of caracol graphs and simultaneously unifies many generalizations of the Tamari and Dyck lattices into one large family. 
See Figure~\ref{fig_framing_lattices_zoo} for examples.

\begin{figure}[htb]
    \centering
    \begin{tikzpicture}

\begin{scope}[xshift=-5, yshift=20]
    \node[draw=black, thick, rectangle, rounded corners, inner xsep=60pt, inner ysep=60pt] at (0,0) {
        {}
    };
    \node[] (permutree) at (0,1.7) {\scriptsize \textbf{Permutree lattices}};
    
    \node[draw=black, thick, rectangle, rounded corners, inner xsep=5pt, inner ysep=5pt] at (0,1) {
        {\scriptsize The weak order}
    };
    \node[draw=black, thick, rectangle, rounded corners, inner xsep=5pt, inner ysep=5pt] at (0,0) {
        {\scriptsize The Tamari lattice}
    };
    \node[draw=black, thick, rectangle, rounded corners, inner xsep=50pt, inner ysep=18pt] at (0,-0.25) {
        {}
    };
    \node[] (cambrian) at (0,-0.6) {\scriptsize $\varepsilon$-Cambrian lattices};
    \node[draw=black, thick, rectangle, rounded corners, inner xsep=5pt, inner ysep=5pt] at (0,-1.5) {
        {\scriptsize The Boolean lattice}
    };    
\end{scope}

\begin{scope}[xshift=180]
    \node[draw=black, thick, rectangle, rounded corners, inner xsep=80pt, inner ysep=80pt] at (0,0) {
        {}
    };   
    \node[] (cross) at (0,2.5) {\scriptsize \textbf{Cross-Tamari lattices}};
    
    \node[draw=black, thick, rectangle, rounded corners, inner xsep=60pt, inner ysep=30pt] at (-0.48,1) {
        {}
    };
    \node[] (eIJ) at (-0.47,1.55) {\scriptsize $(\varepsilon,I,J)$-Cambrian lattices};

    \node[draw=black, thick, rectangle, rounded corners, inner xsep=50pt, inner ysep=18pt] at (-0.12,0.58) {
        {}
    };
    \node[] (eIJ) at (-0.2,0.9) {\scriptsize $\varepsilon$-Cambrian lattices};
    \node[] (eIJ) at (0.08,0.3) {\scriptsize The Tamari lattice};
    
    \node[draw=black, thick, rectangle, rounded corners, inner xsep=45pt, inner ysep=20pt] at (0.2,-0.05) {
        {}
    };
    \node[] (eIJ) at (0.25,-0.41) {\scriptsize $\nu$-Tamari lattices};
    
    \node[draw=black, thick, rectangle, rounded corners, inner xsep=5pt, inner ysep=5pt] at (0.3,-1.7) {
        {\scriptsize The Dyck lattice}
    };
    \node[draw=black, thick, rectangle, rounded corners, inner xsep=45pt, inner ysep=18pt] at (0.3,-1.9) {
        {}
    };
    \node[] (eIJ) at (0.3,-2.3) {\scriptsize $\nu$-Dyck lattices};

    \node[draw=black, thick, rectangle, rounded corners, inner xsep=55pt, inner ysep=47pt] at (0.54,-1) {
        {}
    };   
    \node[] (cross) at (0.5,-1) {\scriptsize alt $\nu$-Tamari lattices};
    
\end{scope}

\begin{scope}[xshift=0, yshift=-160]
    \node[draw=black, thick, rectangle, rounded corners, inner xsep=65pt, inner ysep=36pt] at (0,0.65) {
        {}
    };
    \node[] (grid) at (0,1.4) {\scriptsize \textbf{Grid-Tamari lattices}};

    \node[draw=black, thick, rectangle, rounded corners, inner xsep=5pt, inner ysep=5pt] at (0,0) {
        {\scriptsize The Tamari lattice}
    };
    
    \node[draw=black, thick, rectangle, rounded corners, inner xsep=60pt, inner ysep=20pt] at (0,0.2) {
        {}
    };
    \node[] (grassman) at (0,0.6) {\scriptsize Grassman--Tamari lattices};
\end{scope}

\begin{scope}[xshift=0, yshift=-75]
    \node[draw=black, thick, rectangle, rounded corners, inner xsep=65pt, inner ysep=15pt] at (0,0) {
        {}
    };
    \node[] (grid) at (0,0.2) {\scriptsize \textbf{$\tau$-Tilting posets of}};
    \node[] (grid) at (0,-0.2) {\scriptsize \textbf{certain gentle algebras}};
\end{scope}

\begin{scope}[xshift=180, yshift=-150]
    \node[draw=black, thick, rectangle, rounded corners, inner xsep=80pt, inner ysep=43pt] at (0,0.45) {
        {}
    };
    \node[] (grid) at (0,1.45) {\scriptsize \textbf{Weak order generalizations}};

    \node[] (multi) at (0.1,0) {\scriptsize The weak order};    
    \node[draw=black, thick, rectangle, rounded corners, inner xsep=60pt, inner ysep=20pt] at (-0.3,0.4) {
        {}
    };
    \node[draw=black, thick, rectangle, rounded corners, inner xsep=50pt, inner ysep=17pt] at (0.5,-0.3) {
        {}
    };
    \node[] (multi) at (-0.25,0.63) {\scriptsize Multipermutation lattices};
    \node[] (s-weak) at (0.43,-0.53) {\scriptsize The $s$-weak order};
\end{scope}

\end{tikzpicture}
    \caption{Some popular exhibitions at the zoo of framing lattices.}
    \label{fig_framing_lattices_zoo}
\end{figure}

\begin{figure}[h]
    \centering
    \begin{tikzpicture}

\begin{scope}
\begin{scope}[xshift = 0, yshift= 0, scale=0.8]
    \node[circle, draw, inner sep=1pt, fill] (s) at (0,0) {};	
    \node[circle, draw, inner sep=1pt, fill] (1) at (1,0) {};
    \node[circle, draw, inner sep=1pt, fill] (2) at (2,0) {};
    \node[circle, draw, inner sep=1pt, fill] (3) at (3,0) {};	
    \node[circle, draw, inner sep=1pt, fill] (t) at (4,0) {};	
     
    \draw[color=black] (s) .. controls (0.25, 0.3) and (0.75, 0.3) .. (1);
    \draw[color=black] (s) .. controls (0.25, -0.3) and (0.75, -0.3) .. (1);    
    \draw[color=black] (s) .. controls (0.25, 0.7) and (1.75, 0.7) .. (2);
    \draw[color=black] (s) .. controls (0.25, -0.7) and (1.75, -0.7) .. (2);
    \draw[color=black] (s) .. controls (0.25, 1.1) and (2.75, 1.1) .. (3);	
    \draw[color=black] (s) .. controls (0.25, -1.1) and (2.75, -1.1) .. (3);	    
           
    \draw[color=black] (3) .. controls (3.25, 0.3) and (3.75, 0.3) .. (t);
    \draw[color=black] (3) .. controls (3.25, -0.3) and (3.75, -0.3) .. (t);
    \draw[color=black] (2) .. controls (2.25, 0.7) and (3.75, 0.7) .. (t);
    \draw[color=black] (2) .. controls (2.25, -0.7) and (3.75, -0.7) .. (t);
    \draw[color=black] (1) .. controls (1.25, 1.1) and (3.75, 1.1) .. (t);	
    \draw[color=black] (1) .. controls (1.25, -1.1) and (3.75, -1.1) .. (t);

\end{scope}
\begin{scope}[xshift = 47, yshift= -130, scale=1]
    
    \node[circle,fill,inner sep=1pt,color=NavyBlue] at (0,0)  {};
    \node[circle,fill,inner sep=1pt,color=NavyBlue] at (-1,1)  {};
    \node[circle,fill,inner sep=1pt,color=NavyBlue] at (0,1)  {};
    \node[circle,fill,inner sep=1pt,color=NavyBlue] at (1,1)  {};   
    \node[circle,fill,inner sep=1pt,color=NavyBlue] at (-1,2)  {};
    \node[circle,fill,inner sep=1pt,color=NavyBlue] at (0,2)  {};
    \node[circle,fill,inner sep=1pt,color=NavyBlue] at (1,2)  {};   
    \node[circle,fill,inner sep=1pt,color=NavyBlue] at (0,3)  {};   
 
    \draw[color=NavyBlue] (0,0) -- (-1,1) -- (0,2) -- (1,1) -- (0,0);
    \draw[color=NavyBlue] (0,1) -- (-1,2) -- (0,3) -- (1,2) -- (0,1);
    \draw[color=NavyBlue] (-1,1) -- (-1,2);
    \draw[color=NavyBlue] (1,1) -- (1,2);
    \draw[color=NavyBlue] (0,0) -- (0,1);
    \draw[color=NavyBlue] (0,2) -- (0,3);

\end{scope}
\end{scope}

\begin{scope}
\begin{scope}[xshift = 120, yshift= 0, scale=0.4]
	\node[circle, draw, inner sep=1pt, fill](1) at (0,0) {};	
	\node[circle, draw, inner sep=1pt, fill](2) at (2,0) {};
	\node[circle, draw, inner sep=1pt, fill](3) at (4,0) {};
	\node[circle, draw, inner sep=1pt, fill](4) at (6,0) {};	


    \draw[color=black] (0,0) -- (2,0) -- (4,0);
    \draw[color=black] (0,0) .. controls (0.4, 0.8) and (1.6, 0.8) .. (2,0);
    \draw[color=black] (0,0) .. controls (0.4, -0.8) and (1.6, -0.8) .. (2,0);    
    \draw[color=black] (2,0) .. controls (2.4, 0.8) and (3.6, 0.8) .. (4,0);
    \draw[color=black] (2,0) .. controls (2.4, -0.8) and (3.6, -0.8) .. (4,0);
    \draw[color=black] (4,0) .. controls (4.4, 0.8) and (5.6, 0.8) .. (6,0);
    \draw[color=black] (4,0) .. controls (4.4, -0.8) and (5.6, -0.8) .. (6,0);       
\end{scope}	

\tikzmath{\x1 = 0.7; \x2 =0; \y1 = 0.4; \y2 = 0.4; \z1 = 0; \z2 = 0.7; \xy1 = 0.7; \xy2 = 0.2; }
\begin{scope}[xshift = 140, yshift= -140, scale=1, rotate=45]
    \node[circle,fill,inner sep=1pt,color=NavyBlue] (11223) at (0,0)  {};	
    \node[circle,fill,inner sep=1pt,color=NavyBlue] (11232) at (\x1, \x2)  {};
    \node[circle,fill,inner sep=1pt,color=NavyBlue] (11322) at (2*\x1, 2*\x2)  {};	
    \node[circle,fill,inner sep=1pt,color=NavyBlue] (12123) at (\y1, \y2)  {};	
    \node[circle,fill,inner sep=1pt,color=NavyBlue] (12132) at (\x1 + \y1,\x2 + \y2)  {};		
    \node[circle,fill,inner sep=1pt,color=NavyBlue] (13122) at (2*\x1 + \xy1, 2*\x2 + \xy2)  {};		
    \node[circle,fill,inner sep=1pt,color=NavyBlue] (13212) at (2*\x1 + \xy1 + \y1, 2*\x2 + \xy2 + \y2)  {};	
    \node[circle,fill,inner sep=1pt,color=NavyBlue] (12312) at (\x1 + \y1 + \xy1, \x2 + \y2 + \xy2)  {};		
    \node[circle,fill,inner sep=1pt,color=NavyBlue] (12321) at (\x1 + 2*\y1 + \xy1, \x2 + 2*\y2 + \xy2)  {};	
    \node[circle,fill,inner sep=1pt,color=NavyBlue] (13221) at (2*\x1 + 2*\y1 + \xy1, 2*\x2 + 2*\y2 + \xy2)  {};		
    \node[circle,fill,inner sep=1pt,color=NavyBlue] (12213) at (2*\y1, 2*\y2)  {};	
    \node[circle,fill,inner sep=1pt,color=NavyBlue] (12231) at (2*\y1 + \xy1, 2*\y2 + \xy2)  {};	
    \node[circle,fill,inner sep=1pt,color=NavyBlue] (21123) at (\y1 + \z1, \y2 + \z2)  {};	
    \node[circle,fill,inner sep=1pt,color=NavyBlue] (21132) at (\x1 + \y1 + \z1, \x2 + \y2 + \z2)  {};		
    \node[circle,fill,inner sep=1pt,color=NavyBlue] (21312) at (\x1 + \y1 + \xy1 + \z1, \x2 + \y2 + \xy2 + \z2)  {};			
    \node[circle,fill,inner sep=1pt,color=NavyBlue] (21321) at (\x1 + 2*\y1 + \xy1 + \z1, \x2 + 2*\y2 + \xy2 + \z2)  {};	
    \node[circle,fill,inner sep=1pt,color=NavyBlue] (21231) at (2*\y1 + \xy1 + \z1, 2*\y2 + \xy2 + \z2)  {};
    \node[circle,fill,inner sep=1pt,color=NavyBlue] (21213) at (2*\y1 + \z1, 2*\y2 + \z2)  {};	
    \node[circle,fill,inner sep=1pt,color=NavyBlue] (31122) at (2*\x1 + 2*\xy1, 2*\x2 + 2*\xy2)  {};	
    \node[circle,fill,inner sep=1pt,color=NavyBlue] (31212) at (2*\x1 + 2*\xy1 + \y1, 2*\x2 + 2*\xy2 + \y2)  {};	
	\node[circle,fill,inner sep=1pt,color=NavyBlue] (31221) at (2*\x1 + 2*\xy1 + 2*\y1, 2*\x2 + 2*\xy2 + 2*\y2)  {};	
	\node[circle,fill,inner sep=1pt,color=NavyBlue] (32121) at (2*\x1 + 2*\xy1 + 2*\y1 + \z1, 2*\x2 + 2*\xy2 + 2*\y2 + \z2)  {};	
	\node[circle,fill,inner sep=1pt,color=NavyBlue] (32112) at (2*\x1 + 2*\xy1 + \y1 + \z1, 2*\x2 + 2*\xy2 + \y2 + \z2)  {};
	\node[circle,fill,inner sep=1pt,color=NavyBlue] (23112) at (\x1 + 2*\xy1 + \y1 + \z1, \x2 + 2*\xy2 + \y2 + \z2)  {};
	\node[circle,fill,inner sep=1pt,color=NavyBlue] (23121) at (\x1 + 2*\xy1 + 2*\y1 + \z1, \x2 + 2*\xy2 + 2*\y2 + \z2)  {};	
	\node[circle,fill,inner sep=1pt,color=NavyBlue] (22311) at (2*\xy1 + 2*\y1 + 2*\z1, 2*\xy2 + 2*\y2 + 2*\z2)  {};	
	\node[circle,fill,inner sep=1pt,color=NavyBlue] (22131) at (\xy1 + 2*\y1 + 2*\z1, \xy2 + 2*\y2 + 2*\z2)  {};
	\node[circle,fill,inner sep=1pt,color=NavyBlue] (22113) at (2*\y1 + 2*\z1, 2*\y2 + 2*\z2)  {};		
	\node[circle,fill,inner sep=1pt,color=NavyBlue] (23211) at (\x1 + 2*\xy1 + 2*\y1 + 2*\z1, \x2 + 2*\xy2 + 2*\y2 + 2*\z2)  {};		
	\node[circle,fill,inner sep=1pt,color=NavyBlue] (32211) at (2*\x1 + 2*\xy1 + 2*\y1 + 2*\z1, 2*\x2 + 2*\xy2 + 2*\y2 + 2*\z2)  {};			

	\draw[color=NavyBlue] (11223) -- (11232) -- (12132);
	\draw[color=NavyBlue] (11223) -- (12123) -- (12132) -- (12312) -- (13212);
	\draw[color=NavyBlue] (11232) -- (11322) -- (13122) -- (13212);
	\draw[color=NavyBlue] (12312) -- (12321);
	\draw[color=NavyBlue] (12123) -- (12213) -- (12231) -- (12321);		
	\draw[color=NavyBlue] (21123) -- (21132) -- (21312) -- (21321);
	\draw[color=NavyBlue] (21123) -- (21213) -- (21231) -- (21321);
	\draw[color=NavyBlue] (12123) -- (21123);
	\draw[color=NavyBlue] (12213) -- (21213);
	\draw[color=NavyBlue] (12312) -- (21312);
	\draw[color=NavyBlue] (12321) -- (21321);
	\draw[color=NavyBlue] (12231) -- (21231);
	\draw[color=NavyBlue] (12132) -- (21132);			
	\draw[color=NavyBlue] (21213) -- (22113) -- (22131) -- (22311) -- (23211);
	\draw[color=NavyBlue] (21231) -- (22131);			
	\draw[color=NavyBlue] (12321) -- (13221) -- (31221);				
	\draw[color=NavyBlue] (13122) -- (31122) -- (31212);			
	\draw[color=NavyBlue] (13212) -- (31212) -- (32112) -- (32121);
	\draw[color=NavyBlue] (31212) -- (31221) -- (32121);
	\draw[color=NavyBlue] (21312) -- (23112) -- (32112);				
	\draw[color=NavyBlue] (21321) -- (23121) -- (32121);
	\draw[color=NavyBlue] (13221) -- (13212);	
	\draw[color=NavyBlue] (23112) -- (23121) -- (23211) -- (32211);				
	\draw[color=NavyBlue] (23112) -- (32112) -- (32121) -- (32211);	
\end{scope}
\end{scope}

\begin{scope}[xshift = 220, yshift= 0, scale=0.8]
    \begin{scope}[xshift = 0, yshift= 0, scale=0.55]
	\node[circle, draw, inner sep=1pt, fill](1) at (1,0) {};	
	\node[circle, draw, inner sep=1pt, fill](2) at (2,0) {};
	\node[circle, draw, inner sep=1pt, fill](3) at (3,0) {};
	\node[circle, draw, inner sep=1pt, fill](4) at (4,0) {};	
	\node[circle, draw, inner sep=1pt, fill](5) at (5,0) {};	
	\node[circle, draw, inner sep=1pt, fill](6) at (6,0) {};			
	\node[circle, draw, inner sep=1pt, fill](7) at (7,0) {};	

	\draw[] (1) -- (2) -- (3) -- (4) -- (5) -- (6) -- (7);
  
	\draw[color=black] (1) .. controls (1.25, -0.6) and (2.75, -0.6) .. (3);	
	\draw[color=black] (1) .. controls (1.25, -1.0) and (3.75, -1.0) .. (4);	
        \draw[color=black] (1) .. controls (1.25, 1.4) and (4.75, 1.4) .. (5);	
        \draw[color=black] (1) .. controls (1.25, -1.8) and (5.75, -1.8) .. (6);	
        
        \draw[color=black] (5) .. controls (5.25, -0.6) and (6.75, -0.6) .. (7);	
        \draw[color=black] (4) .. controls (4.25, 1.0) and (6.75, 1.0) .. (7);
        \draw[color=black] (3) .. controls (3.25, -1.4) and (6.75, -1.4) .. (7);
        \draw[color=black] (2) .. controls (2.25, -1.8) and (6.75, -1.8) .. (7);

    \end{scope}
    
    \begin{scope}[xshift = 55, yshift= -160, scale=0.8]

	\node[circle,fill,inner sep=1pt,color=NavyBlue] at (0,0)  {};
	\node[circle,fill,inner sep=1pt,color=NavyBlue] at (1,0.5)  {};
	\node[circle,fill,inner sep=1pt,color=NavyBlue] at (0,1)  {};
	\node[circle,fill,inner sep=1pt,color=NavyBlue] at (1,1.5)  {};

	\node[circle,fill,inner sep=1pt,color=NavyBlue] at (-1.5,1.5)  {};
	\node[circle,fill,inner sep=1pt,color=NavyBlue] at (-1.5,3)  {}; 
 	\node[circle,fill,inner sep=1pt,color=NavyBlue] at (-0.5,2.5)  {};
  
	\node[circle,fill,inner sep=1pt,color=NavyBlue] at (1.5,2.5)  {};
 	\node[circle,fill,inner sep=1pt,color=NavyBlue] at (2,1.5)  {};
	\node[circle,fill,inner sep=1pt,color=NavyBlue] at (2.5,3)  {};

 	\node[circle,fill,inner sep=1pt,color=NavyBlue] at (1,3)  {};
 	\node[circle,fill,inner sep=1pt,color=NavyBlue] at (0.5,2.4)  {};
        \node[circle,fill,inner sep=1pt,color=NavyBlue] at (0.5,4.5)  {};
        \node[circle,fill,inner sep=1pt,color=NavyBlue] at (1.5,4)  {};        
 
	\draw[color=NavyBlue] (0,0) -- (1,0.5) -- (1,1.5) -- (0,1) -- (0,0);
 	\draw[color=NavyBlue] (0,0) -- (-1.5,1.5) -- (-1.5,3) -- (-0.5,2.5) -- (0,1);
   	\draw[color=NavyBlue] (1,0.5) -- (2,1.5) -- (2.5,3) -- (1.5,2.5) -- (1,1.5);
        \draw[color=NavyBlue] (-0.5,2.5) -- (1,3) -- (1.5,2.5);
        \draw[color=NavyBlue] (-1.5,1.5) -- (0.5,2.4) -- (2,1.5);
        \draw[color=NavyBlue] (0.5,2.4) -- (0.5,4.5) -- (-1.5,3);
        \draw[color=NavyBlue] (0.5,4.5) -- (1.5,4) -- (1,3);
        \draw[color=NavyBlue] (1.5,4) -- (2.5,3);
    \end{scope}
\end{scope}

\begin{scope}[xshift = 340, yshift= 0, scale=0.8]
    \begin{scope}[xshift = 0, yshift= 0, scale=0.55]
	\node[circle, draw, inner sep=1pt, fill](1) at (1,0) {};	
	\node[circle, draw, inner sep=1pt, fill](2) at (2,0) {};
	\node[circle, draw, inner sep=1pt, fill](3) at (3,0) {};
	\node[circle, draw, inner sep=1pt, fill](4) at (4,0) {};	
	\node[circle, draw, inner sep=1pt, fill](5) at (5,0) {};	
	\node[circle, draw, inner sep=1pt, fill](6) at (6,0) {};			
	\node[circle, draw, inner sep=1pt, fill](7) at (7,0) {};
	\node[circle, draw, inner sep=1pt, fill](8) at (8,0) {};	 
	\node[circle, draw, inner sep=1pt, fill](9) at (9,0) {};	
	\node[circle, draw, inner sep=1pt, fill](10) at (10,0) {};	
 
	\draw[] (1) -- (2) -- (3) -- (4) -- (5) -- (6) -- (7) -- (8) -- (9) -- (10);
  
	\draw[color=black] (1) .. controls (1.25, -0.7) and (2.75, -0.7) .. (3);	
	\draw[color=black] (1) .. controls (1.25, 1.0) and (4.75, 1.0) .. (5);	
        \draw[color=black] (1) .. controls (1.25, 1.5) and (6.75, 1.5) .. (7);	
         
        \draw[color=black] (4) .. controls (4.25, 1.5) and (9.75, 1.5) .. (10);	
        \draw[color=black] (6) .. controls (6.25, -1.2) and (9.75, -1.2) .. (10);
        \draw[color=black] (8) .. controls (8.25, -0.7) and (9.75, -0.7) .. (10);

    \end{scope}
        
    \begin{scope}[xshift = 90, yshift= -50, scale=0.7]

	\node[circle,fill,inner sep=1pt,color=NavyBlue] at (-0.5,0.5)  {};
	\node[circle,fill,inner sep=1pt,color=NavyBlue] at (0.5,0)  {};
	\node[circle,fill,inner sep=1pt,color=NavyBlue] at (1,-1)  {};
	\node[circle,fill,inner sep=1pt,color=NavyBlue] at (-1.4,-1.4)  {};
	\node[circle,fill,inner sep=1pt,color=NavyBlue] at (0,-2)  {};	
	\node[circle,fill,inner sep=1pt,color=NavyBlue] at (2,-2)  {};
	\node[circle,fill,inner sep=1pt,color=NavyBlue] at (1,-3)  {};
	\node[circle,fill,inner sep=1pt,color=NavyBlue] at (-1,-3)  {};	
	\node[circle,fill,inner sep=1pt,color=NavyBlue] at (-2.4,-2.4)  {};
	\node[circle,fill,inner sep=1pt,color=NavyBlue] at (0,-4)  {};	
	\node[circle,fill,inner sep=1pt,color=NavyBlue] at (-0.6,-6)  {};
	\node[circle,fill,inner sep=1pt,color=NavyBlue] at (-1.2,-4)  {};	
	\node[circle,fill,inner sep=1pt,color=NavyBlue] at (-2.2,-5)  {};
	\node[circle,fill,inner sep=1pt,color=NavyBlue] at (0.4,-5)  {};			
								
	\draw[color=NavyBlue] (-0.5,0.5) -- (-1.4,-1.4) -- (0,-2) -- (1,-1) -- (0.5,0) -- (-0.5,0.5);
	\draw[color=NavyBlue] (1,-1) -- (2,-2) -- (1,-3) -- (0,-2);	
	\draw[color=NavyBlue] (-1.4,-1.4) -- (-2.4,-2.4) -- (-1,-3) -- (0,-2);
	\draw[color=NavyBlue] (1,-3) -- (0,-4) -- (-1,-3);
	\draw[color=NavyBlue] (1,-3) -- (0.4,-5);
	\draw[color=NavyBlue] (-2.4,-2.4) -- (-2.2,-5);
	\draw[color=NavyBlue] (-1.4,-1.4) -- (-1.2,-4);
	\draw[color=NavyBlue] (0,-4) -- (-0.6,-6);
	\draw[color=NavyBlue] (-0.6,-6) -- (0.4,-5) -- (-1.2,-4) -- (-2.2,-5) -- (-0.6,-6);		
    \end{scope}

    \begin{scope}[xshift = 30, yshift= -60, scale=0.35]
        \draw[very thin, color=gray] (1,1) grid (3,2);
        \draw[very thin, color=gray] (0,2) grid (3,3);
        \draw[very thin, color=gray] (2,3) grid (3,4);
    \end{scope}
\end{scope}
\end{tikzpicture}

    \caption{Four framed graphs and the Hasse diagrams of their framing lattices. 
    The first is the Boolean lattice~$\scrB_3$. 
    The second is the lattice of multipermutations of $1^22^23$. The third is the $\varepsilon$-cambrian lattice with $\varepsilon = (-,-,+,-)$. 
    The fourth is a cross-Tamari lattice of the cross-shaped grid shown below the right-most graph.}
    \label{fig.graphs_and_lattices}
\end{figure}

A detailed list of the species currently known to us is given below. 

\begin{itemize}
    \item the Tamari lattice (Example~\ref{ex_caracol3_Tamari}); 
    \item $\nu$-Tamari lattices of Pr\'eville-Ratelle and Viennot \cite{PV17};
    \item alt $\nu$-Tamari lattices of Ceballos and Chenevi\`ere \cite{CC23};
    \item type A Cambrian lattices of Reading \cite{Rea06};
    \item $(\varepsilon, I,\overline{J})$-Cambrian lattices (Remark~\ref{rem.eIJ-cambrian}) of Pilaud \cite{Pil20}, which generalize the type A Cambrian lattices;
    \item the Dyck lattice (Example~\ref{ex_caracol3_Dyck});
    \item principal order ideals in Young's lattice (also known as the $\nu$-Dyck lattice \cite{CC23} or Stanley's distributive lattice). This is a direct consequence of \cite{BGMY23};
    \item cross-Tamari lattices (Theorem~\ref{thm_crossTam}) introduced in  Section~\ref{sec_crossTamari}, which generalize the lattice families above;  
    \item the Boolean lattice (Section~\ref{prop.boolean});
    \item the weak order on permutations (Example~\ref{ex_weak_via_framing});
    \item the $s$-weak order of Ceballos and Pons, which is a direct consequence of \cite{GMPTY23};
    \item the lattice of multipermutations (Theorem~\ref{thm.lattice_of_multipermutations});
    \item permutree lattices of Pilaud and Pons \cite{PP18}, which unify the weak order and the Boolean, Tamari, and Cambrian lattices. This is a direct consequence of the connection to flow polytopes given by Tamayo \cite{Tam23}.
    \item $\tau$-tilting posets for certain gentle algebras, which follows directly from \cite{BBBHPSY22};
    \item the Grassmann--Tamari order of Santos--Stump--Welker~\cite{SSW17}.
    \item the Grid-Tamari lattices of McConville \cite{McC2017}, which generalize the Grassman-Tamari order. This is an immediate consequence of Garver--McConville~\cite{GM17}.
\end{itemize}

We give further details on these lattices in the sections to come, and expect new interesting species to emerge.

\section{The Boolean lattice}
\label{sec.boolean}

\begin{figure}[htb]
    \centering
    \begin{tikzpicture}

\begin{scope}[xshift = 0, yshift= 0, scale=1]
    \node[circle, draw, inner sep=1.5pt, fill] (s) at (0,0) {};	
    \node[circle, draw, inner sep=1.5pt, fill] (1) at (1,0) {};
    \node[circle, draw, inner sep=1.5pt, fill] (2) at (2,0) {};
    \node[circle, draw, inner sep=1.5pt, fill] (3) at (3,0) {};	
    \node[circle, draw, inner sep=1.5pt, fill] (t) at (4,0) {};

    \draw[color=black] (s) .. controls (0.25, 0.3) and (0.75, 0.3) .. (1);
    \draw[color=black] (s) .. controls (0.25, -0.3) and (0.75, -0.3) .. (1);    
    \draw[color=black] (s) .. controls (0.25, 0.7) and (1.75, 0.7) .. (2);
    \draw[color=black] (s) .. controls (0.25, -0.7) and (1.75, -0.7) .. (2);
    \draw[color=black] (s) .. controls (0.25, 1.1) and (2.75, 1.1) .. (3);	
    \draw[color=black] (s) .. controls (0.25, -1.1) and (2.75, -1.1) .. (3);	    
           
    \draw[color=black] (3) .. controls (3.25, 0.3) and (3.75, 0.3) .. (t);
    \draw[color=black] (3) .. controls (3.25, -0.3) and (3.75, -0.3) .. (t);
    \draw[color=black] (2) .. controls (2.25, 0.7) and (3.75, 0.7) .. (t);
    \draw[color=black] (2) .. controls (2.25, -0.7) and (3.75, -0.7) .. (t);
    \draw[color=black] (1) .. controls (1.25, 1.1) and (3.75, 1.1) .. (t);	
    \draw[color=black] (1) .. controls (1.25, -1.1) and (3.75, -1.1) .. (t);

    \node[] (a) at (0,-0.4) {\footnotesize $s$};
    \node[] (a) at (1,-0.4) {\footnotesize $1$};
    \node[] (a) at (2,-0.4) {\footnotesize $2$};
    \node[] (a) at (3,-0.4) {\footnotesize $3$};
    \node[] (a) at (4,-0.4) {\footnotesize $t$};
\end{scope}

\begin{scope}[xshift = 270, yshift= -50, scale=1.3]
    
    \node[circle,fill,inner sep=1.5pt,color=NavyBlue] at (0,0)  {};
    \node[circle,fill,inner sep=1.5pt,color=NavyBlue] at (-1,1)  {};
    \node[circle,fill,inner sep=1.5pt,color=NavyBlue] at (0,1)  {};
    \node[circle,fill,inner sep=1.5pt,color=NavyBlue] at (1,1)  {};   
    \node[circle,fill,inner sep=1.5pt,color=NavyBlue] at (-1,2)  {};
    \node[circle,fill,inner sep=1.5pt,color=NavyBlue] at (0,2)  {};
    \node[circle,fill,inner sep=1.5pt,color=NavyBlue] at (1,2)  {};   
    \node[circle,fill,inner sep=1.5pt,color=NavyBlue] at (0,3)  {};   
 
    \draw[color=NavyBlue] (0,0) -- (-1,1) -- (0,2) -- (1,1) -- (0,0);
    \draw[color=NavyBlue] (0,1) -- (-1,2) -- (0,3) -- (1,2) -- (0,1);
    \draw[color=NavyBlue] (-1,1) -- (-1,2);
    \draw[color=NavyBlue] (1,1) -- (1,2);
    \draw[color=NavyBlue] (0,0) -- (0,1);
    \draw[color=NavyBlue] (0,2) -- (0,3);

    \node[] (a) at (0.3,0) {\scriptsize $\emptyset$};
    \node[] (a) at (-0.7,1) {\scriptsize $\{1\}$};
    \node[] (a) at (0.3,1) {\scriptsize $\{2\}$};
    \node[] (a) at (1.3,1) {\scriptsize $\{3\}$};
    \node[] (a) at (-0.6,2) {\scriptsize $\{1,2\}$};
    \node[] (a) at (0.4,2) {\scriptsize $\{1,3\}$};
    \node[] (a) at (1.4,2) {\scriptsize $\{2,3\}$};
    \node[] (a) at (0.5,3) {\scriptsize $\{1,2,3\}$};
\end{scope}
\end{tikzpicture}
    \caption{The graph $G_3$ and the Boolean lattice $\scrB_3$.}
    \label{fig.boolean_lattice}
\end{figure}
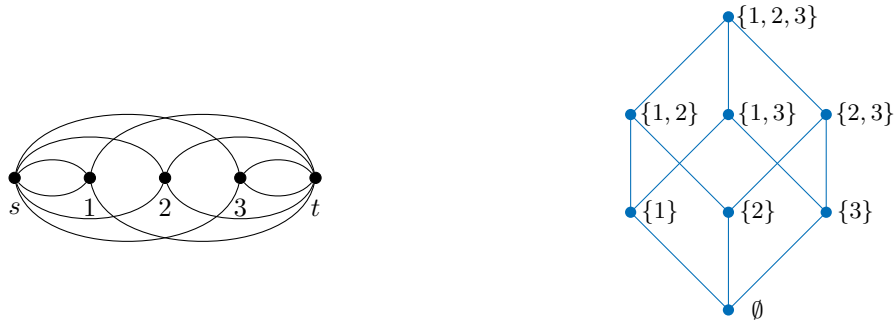

The Boolean lattice $\scrB_n$ is the lattice on the subsets of $[n]$ ordered by inclusion. 
In this section, we describe how to obtain $\scrB_n$ as a framing lattice.
Let $G_n$ be the flow graph with vertex set $\{s,t\}\cup [n]$ and edge set constructed as follows.
For each vertex $i\in [n]$ we add a pair of edges $(s,i)$ and $(s,i)'$ and a pair of edges $(i,t)$ and $(i,t)'$.
In light of Lemma~\ref{lem.iso_operations}, all framings of $G_n$ give isomorphic framing lattices.
For convenience, we choose $F$ to be a framing induced by ordering $(s,i) <_{\scrI(i)} (s,i)'$ and $(i,t) <_{\scrO(i)} (i,t)'$ at each $i \in [n]$.  
See Figure~\ref{fig.boolean_lattice} for an example.

\begin{proposition}
    \label{prop.boolean}
    The framing lattice $\scrL_{G_n,F}$ is the Boolean lattice $\scrB_n$.     
\end{proposition}

\begin{proof}
    A maximal clique of $(G_n,F)$ contains either the route $R_i:=\{(s,i),(i,t)'\}$ or the route $R_i':=\{(s,i)',(i,t)\}$ for each $i\in [n]$.
    For a set $S\subseteq [n]$, define the maximal clique $C_S$ to be the unique maximal clique with routes $R_i'$ with $i\in S$.
    The map $S\mapsto C_S$ is an order preserving bijection between $\scrB_n$ and $\scrL_{G_n,F}$.
\end{proof}

\section{Cambrian lattices of type \texorpdfstring{$A$}{}}
\label{sec.cambrian}
Reading's type $A$ $\varepsilon$-Cambrian lattices~\cite{Rea06} are lattices on triangulations of a polygon. 
The parameter $\varepsilon$ is a map 
$\varepsilon: [n] \to \{\pm\}$
that assigns a positive or negative sign to each element of~$[n]$.  
We define the polygon $P_{\varepsilon}(n)$ as a convex $(n+2)$-gon with vertices $0,1,\dots,n+1$ ordered from left to right, such that $0$ and $n+1$ are on a horizontal line and $i$ is above this line if $\varepsilon(i)=+$, or below if $\varepsilon(i)=-$.   
The \mbox{$\varepsilon$-Cambrian} lattice is the poset on triangulations of~$P_{\varepsilon}(n)$ whose cover relations are increasing slope diagonal flips. The classical Tamari lattice is recovered when $\varepsilon(i)=-$ for all $i$.

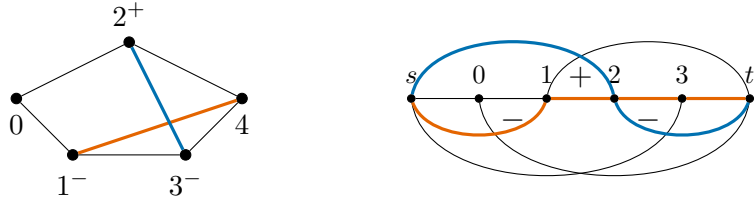
\begin{figure}[htb]
    \centering
    \begin{tikzpicture}[scale=0.75]

    \begin{scope}
    \node[circle, draw, inner sep=1.4pt, fill, label=below:{\small $0$}] (0) at (0,0)  {};
    \node[circle, draw, inner sep=1.4pt, fill, label=below:{\small $1^-$}] (1) at (1,-1)  {};
    \node[circle, draw, inner sep=1.4pt, fill, label=above:{\small $2^+$}] (2) at (2,1)  {};
    \node[circle, draw, inner sep=1.4pt, fill, label=below:{\small $3^-$}] (3) at (3,-1)  {};
    \node[circle, draw, inner sep=1.4pt, fill, label=below:{\small $4$}] (4) at (4,0)  {};

    \draw[] (0) -- (1) -- (3) -- (4);	
    \draw[] (0) -- (2) -- (4);	
 
    \draw[very thick, red] (1) -- (4);
    \draw[very thick, blue] (2) -- (3);
        
    \end{scope}

\begin{scope}[shift={(7,0)}, scale=1.2,
    fnode/.style={circle, draw, inner sep=1pt, fill}]
    
    \node[fnode, label=above:{\footnotesize $s$}] (s) at (0,0) {};
    \node[fnode, label=above:{\footnotesize $0$}] (0) at (1,0) {};
    \node[fnode, label=above:{\footnotesize $1$}] (1) at (2,0) {};
    \node[fnode, label=above:{\footnotesize $2$}] (2) at (3,0) {};
    \node[fnode, label=above:{\footnotesize $3$}] (3) at (4,0) {};
    \node[fnode, label=above:{\footnotesize $t$}] (t) at (5,0) {};
    
    \draw[color=black] (s) -- (t);
    \draw (0) -- (1) node[midway, below]{$-$};
    \draw (1) -- (2) node[midway, above]{$+$};
    \draw (2) -- (3) node[midway, below]{$-$};

    \draw[color=black] (s) .. controls (0.25, -0.7) and (1.75, -0.7) .. (1);
    \draw[color=black] (s) .. controls (0.25, 1.1) and (2.75, 1.1) .. (2);    
    \draw[color=black] (s) .. controls (0.25, -1.5) and (3.75, -1.5) .. (3);

    \draw[color=black] (0) .. controls (1.25, -1.5) and (4.75, -1.5) .. (t);
    \draw[color=black] (1) .. controls (2.25, 1.1) and (4.75, 1.1) .. (t);
    \draw[color=black] (2) .. controls (3.25, -0.7) and (4.75, -0.7) .. (t);

    \draw[very thick, color=blue] (s) .. controls (0.25, 1.1) and (2.75, 1.1) .. (2);    
    \draw[very thick, color=blue] (2) .. controls (3.25, -0.7) and (4.75, -0.7) .. (t);

    \draw[very thick, color=red] (s) .. controls (0.25, -0.7) and (1.75, -0.7) .. (1);
    \draw[very thick, color=red] (1) -- (2) -- (3) -- (t);

\end{scope}

\end{tikzpicture}

    \caption{The polygon $P_\varepsilon(3)$ and the Cambrian caracol graph $\cambrianCaracolGraph{\varepsilon}$ for the parameter $\varepsilon = (-,+,-)$.}
    \label{fig_fpsac_cambrian}
\end{figure}

Let the \defn{Cambrian caracol graph} $\cambrianCaracolGraph{\varepsilon}$ be the graph with vertex set $\{s,0,1,\dots,n,t\}$ and the following three kinds of edges: 
\begin{itemize}
    \item horizontal edges $(s,0),(0,1),(1,2),\dots,(n-1,n),(n,t)$,
    \item positive edges $(s,a)^+,(a-1,t)^+$ when $\varepsilon(a)=+$ (above the horizontal line), and
    \item negative edges $(s,a)^-,(a-1,t)^-$ when $\varepsilon(a)=-$ (below the horizontal line).
\end{itemize}
The graph $\cambrianCaracolGraph{\varepsilon}$ is independent of $\varepsilon$, and coincides with the caracol graph $\car(n+3)$ from~\Cref{ex_caracol3_Tamari}. 
The framing $F_{\varepsilon}$ is the one induced by the drawing, which depends on $\varepsilon$.

\begin{theorem}\label{thm_Cambrian}
    The framing lattice~$\scrL_{\cambrianCaracolGraph{\varepsilon},F_{\varepsilon}}$ is the $\varepsilon$-Cambrian lattice.
\end{theorem}
\begin{proof}
    The routes of $\cambrianCaracolGraph{\varepsilon}$ are in bijection with the diagonals of the polygon~$P_{\varepsilon}(n)$. More precisely, a route is completely determined by its first edge $(s,i)$ entering at $i$ and its last edge $(j-1,t)$ exiting at $j-1$. 
    This route corresponds to the diagonal $ij$ of the polygon. Under this bijection, two routes are coherent if and only if the corresponding diagonals do not cross; see~\Cref{fig_fpsac_cambrian}.
    Therefore, maximal cliques of the framing lattice correspond to triangulations of the polygon. One can also observe that a $\ccw$ rotation corresponds to an increasing slope diagonal flip. As a consequence, the framing lattice~$\scrL_{\cambrianCaracolGraph{\varepsilon},F_{\varepsilon}}$ is isomorphic to the $\varepsilon$-Cambrian lattice.
\end{proof}

An example for $\varepsilon=(-,-,+,-)$ is shown in~\Cref{fig.graphs_and_lattices}.

\begin{remark} \label{rem.eIJ-cambrian}
    The Tamari lattice was generalized to $\nu$-Tamari lattices by Pr\'eville-Ratelle and Viennot \cite{PV17}, which coincide with the $(I,J)$-Tamari lattices described in~\cite{CPS19} for certain subsets $I,J\subseteq \{0\}\sqcup [n]$. Using similar ideas, the $\varepsilon$-Cambrian lattice was generalized to $(\varepsilon,I,J)$-Cambrian lattices by Pilaud in~\cite{Pil20} \footnote{Note that while we use $J\subseteq \{0\}\sqcup [n]$, the convention in  \cite{Pil20} is to use $J\subseteq \{1,\ldots, n+1\}$. Our convention matches that of \cite{CPS19}.}.
    These more general lattices can also be obtained as framing lattices. 
    The $(\varepsilon,I,J)$ graph to be used is obtained from $G_\varepsilon$ by removing the edges $(s,i)$ for $i \notin I$ and the edges $(j-1,t)$ for $j \notin J$. The framing is the one induced by the drawing.    
\end{remark}

\section{The lattice of multipermutations}
\label{sec.multipermutations}

The weak order on $S_n$ can be extended to a weak order on multipermutations. To our knowledge, the resulting lattice of multipermutations was first considered by Bennett and Birkhoff in~\cite{BB94}. For a composition $\mathbf{s} = (s_1,s_2,\ldots, s_n)$ of a positive integer $k$, let $M$ be the set of mutlipermutations of the multiset with $s_i$ copies of $i$ for $1\leq i \leq n$. The lattice of $\mathbf{s}$-multipermutations $\scrM_{\mathbf{s}}$ is the poset induced by the cover relation: $\mu_1 \prec \mu_2$ if and only if the multipermutation $\mu_2$ is obtainable from the multipermutation $\mu_1$ by an increasing transposition of two adjacent numbers. 
Note that when choosing $\mathbf{s} = (1,\ldots, 1)$ with length~$n$, the lattice of multipermutations $\scrM_{\mathbf{s}}$ is the classical weak order on $S_n$. 

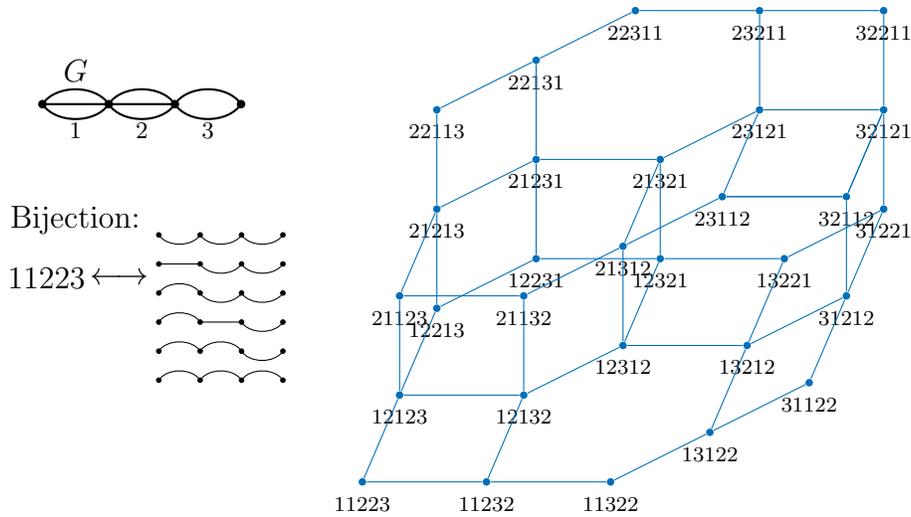
\begin{figure}
\centering
    \begin{tikzpicture}[scale=1.1]

\tikzmath{\x1 = 1; \x2 =0; \y1 = 0.3; \y2 = 0.7; \z1 = 0; \z2 = 0.8; \xy1 = 0.8; \xy2 = 0.4; }

\begin{scope}[xshift = -40, yshift= -20, scale=0.4]
	\node[circle, draw, inner sep=1pt, fill](1) at (0,0) {};	
	\node[circle, draw, inner sep=1pt, fill](2) at (2,0) {};
	\node[circle, draw, inner sep=1pt, fill](3) at (4,0) {};
	\node[circle, draw, inner sep=1pt, fill](4) at (6,0) {};	

	\node (a) at (1,1) {$G$};	

	\draw[thick, color=black] (0,0) -- (2,0) -- (4,0);
	\draw[thick, color=black] (0,0) .. controls (0.4, 0.6) and (1.6, 0.6) .. (2,0);
    \draw[thick, color=black] (0,0) .. controls (0.4, -0.6) and (1.6, -0.6) .. (2,0);    
    \draw[thick, color=black] (2,0) .. controls (2.4, 0.6) and (3.6, 0.6) .. (4,0);
    \draw[thick, color=black] (2,0) .. controls (2.4, -0.6) and (3.6, -0.6) .. (4,0);
    \draw[thick, color=black] (4,0) .. controls (4.4, 0.6) and (5.6, 0.6) .. (6,0);
    \draw[thick, color=black] (4,0) .. controls (4.4, -0.6) and (5.6, -0.6) .. (6,0);       		
    \node[] (a) at (1,-0.8) {\scriptsize $1$};
    \node[] (a) at (3,-0.8) {\scriptsize $2$};    
    \node[] (a) at (5,-0.8) {\scriptsize $3$};
\end{scope}	

\begin{scope}[scale=0.25, xshift=0, yshift=-260]
\begin{scope}[yshift=0]
    \draw[ color=black] (0,0) .. controls (0.4, -0.6) and (1.6, -0.6) .. (2,0);    
    \draw[ color=black] (2,0) .. controls (2.4, -0.6) and (3.6, -0.6) .. (4,0);
    \draw[ color=black] (4,0) .. controls (4.4, -0.6) and (5.6, -0.6) .. (6,0);        
    
    \node[circle, draw, inner sep=0.6pt, fill] (1) at (0,0) {};
    \node[circle, draw, inner sep=0.6pt, fill] (2) at (2,0) {};
    \node[circle, draw, inner sep=0.6pt, fill] (3) at (4,0) {};
    \node[circle, draw, inner sep=0.6pt, fill] (4) at (6,0) {};    
\end{scope}
\begin{scope}[yshift=-40]
    \draw[ color=black] (2,0) .. controls (2.4, -0.6) and (3.6, -0.6) .. (4,0);
    \draw[ color=black] (4,0) .. controls (4.4, -0.6) and (5.6, -0.6) .. (6,0);        

		\draw[color=black] (0,0) -- (2,0);        
    
    \node[circle, draw, inner sep=0.6pt, fill] (1) at (0,0) {};
    \node[circle, draw, inner sep=0.6pt, fill] (2) at (2,0) {};
    \node[circle, draw, inner sep=0.6pt, fill] (3) at (4,0) {};
    \node[circle, draw, inner sep=0.6pt, fill] (4) at (6,0) {};    
\end{scope}
\begin{scope}[yshift=-80]
    \draw[ color=black] (0,0) .. controls (0.4, 0.6) and (1.6, 0.6) .. (2,0);
    \draw[ color=black] (2,0) .. controls (2.4, -0.6) and (3.6, -0.6) .. (4,0);
    \draw[ color=black] (4,0) .. controls (4.4, -0.6) and (5.6, -0.6) .. (6,0);        
    
    \node[circle, draw, inner sep=0.6pt, fill] (1) at (0,0) {};
    \node[circle, draw, inner sep=0.6pt, fill] (2) at (2,0) {};
    \node[circle, draw, inner sep=0.6pt, fill] (3) at (4,0) {};
    \node[circle, draw, inner sep=0.6pt, fill] (4) at (6,0) {};    
\end{scope}
\begin{scope}[yshift=-120]
    \draw[ color=black] (0,0) .. controls (0.4, 0.6) and (1.6, 0.6) .. (2,0);
    \draw[ color=black] (4,0) .. controls (4.4, -0.6) and (5.6, -0.6) .. (6,0);        

	\draw[color=black] (4,0) -- (2,0);        
    
    \node[circle, draw, inner sep=0.6pt, fill] (1) at (0,0) {};
    \node[circle, draw, inner sep=0.6pt, fill] (2) at (2,0) {};
    \node[circle, draw, inner sep=0.6pt, fill] (3) at (4,0) {};
    \node[circle, draw, inner sep=0.6pt, fill] (4) at (6,0) {};    
\end{scope}
\begin{scope}[yshift=-160]
    \draw[ color=black] (0,0) .. controls (0.4, 0.6) and (1.6, 0.6) .. (2,0);
    \draw[ color=black] (2,0) .. controls (2.4, 0.6) and (3.6, 0.6) .. (4,0);
    \draw[ color=black] (4,0) .. controls (4.4, -0.6) and (5.6, -0.6) .. (6,0);        
    
    \node[circle, draw, inner sep=0.6pt, fill] (1) at (0,0) {};
    \node[circle, draw, inner sep=0.6pt, fill] (2) at (2,0) {};
    \node[circle, draw, inner sep=0.6pt, fill] (3) at (4,0) {};
    \node[circle, draw, inner sep=0.6pt, fill] (4) at (6,0) {};    
\end{scope}
\begin{scope}[yshift=-200]
    \draw[ color=black] (0,0) .. controls (0.4, 0.6) and (1.6, 0.6) .. (2,0);
    \draw[ color=black] (2,0) .. controls (2.4, 0.6) and (3.6, 0.6) .. (4,0);
    \draw[ color=black] (4,0) .. controls (4.4, 0.6) and (5.6, 0.6) .. (6,0);
    
    \node[circle, draw, inner sep=0.6pt, fill] (1) at (0,0) {};
    \node[circle, draw, inner sep=0.6pt, fill] (2) at (2,0) {};
    \node[circle, draw, inner sep=0.6pt, fill] (3) at (4,0) {};
    \node[circle, draw, inner sep=0.6pt, fill] (4) at (6,0) {};    
\end{scope}
\begin{scope}[xshift = -140, yshift=-60]
    \node[] (123) at (-0.5,0) {$11223$}; 
    \node[] (a) at (3,0) {$\longleftrightarrow$};           
\end{scope}
\begin{scope}[xshift = -115, yshift=20]
    \node[] (123) at (0,0) {Bijection:}; 
\end{scope}
\end{scope}

\begin{scope}[xshift = 70, yshift= -150, scale=1.5]
	\node[circle,fill,inner sep=1pt,color=NavyBlue, label=below:{\tiny $11223$}] (11223) at (0,0)  {};	
	\node[circle,fill,inner sep=1pt,color=NavyBlue, label=below:{\tiny $11232$}] (11232) at (\x1, \x2)  {};
	\node[circle,fill,inner sep=1pt,color=NavyBlue, label=below:{\tiny $11322$}] (11322) at (2*\x1, 2*\x2)  {};	
	\node[circle,fill,inner sep=1pt,color=NavyBlue, label=below:{\tiny $12123$}] (12123) at (\y1, \y2)  {};	
	\node[circle,fill,inner sep=1pt,color=NavyBlue, label=below:{\tiny $12132$}] (12132) at (\x1 + \y1,\x2 + \y2)  {};		
	\node[circle,fill,inner sep=1pt,color=NavyBlue, label=below:{\tiny $13122$}] (13122) at (2*\x1 + \xy1, 2*\x2 + \xy2)  {};		
	\node[circle,fill,inner sep=1pt,color=NavyBlue, label=below:{\tiny $13212$}] (13212) at (2*\x1 + \xy1 + \y1, 2*\x2 + \xy2 + \y2)  {};		
	\node[circle,fill,inner sep=1pt,color=NavyBlue, label=below:{\tiny $12312$}] (12312) at (\x1 + \y1 + \xy1, \x2 + \y2 + \xy2)  {};			
	\node[circle,fill,inner sep=1pt,color=NavyBlue, label=below:{\tiny $12321$}] (12321) at (\x1 + 2*\y1 + \xy1, \x2 + 2*\y2 + \xy2)  {};	
	\node[circle,fill,inner sep=1pt,color=NavyBlue, label=below:{\tiny $13221$}] (13221) at (2*\x1 + 2*\y1 + \xy1, 2*\x2 + 2*\y2 + \xy2)  {};		
	\node[circle,fill,inner sep=1pt,color=NavyBlue, label=below:{\tiny $12213$}] (12213) at (2*\y1, 2*\y2)  {};	
	\node[circle,fill,inner sep=1pt,color=NavyBlue, label=below:{\tiny $12231$}] (12231) at (2*\y1 + \xy1, 2*\y2 + \xy2)  {};	
	\node[circle,fill,inner sep=1pt,color=NavyBlue, label=below:{\tiny $21123$}] (21123) at (\y1 + \z1, \y2 + \z2)  {};	
	\node[circle,fill,inner sep=1pt,color=NavyBlue, label=below:{\tiny $21132$}] (21132) at (\x1 + \y1 + \z1, \x2 + \y2 + \z2)  {};		
	\node[circle,fill,inner sep=1pt,color=NavyBlue, label=below:{\tiny $21312$}] (21312) at (\x1 + \y1 + \xy1 + \z1, \x2 + \y2 + \xy2 + \z2)  {};			
	\node[circle,fill,inner sep=1pt,color=NavyBlue, label=below:{\tiny $21321$}] (21321) at (\x1 + 2*\y1 + \xy1 + \z1, \x2 + 2*\y2 + \xy2 + \z2)  {};	
	\node[circle,fill,inner sep=1pt,color=NavyBlue, label=below:{\tiny $21231$}] (21231) at (2*\y1 + \xy1 + \z1, 2*\y2 + \xy2 + \z2)  {};
	\node[circle,fill,inner sep=1pt,color=NavyBlue, label=below:{\tiny $21213$}] (21213) at (2*\y1 + \z1, 2*\y2 + \z2)  {};	
	\node[circle,fill,inner sep=1pt,color=NavyBlue, label=below:{\tiny $31122$}] (31122) at (2*\x1 + 2*\xy1, 2*\x2 + 2*\xy2)  {};	
	\node[circle,fill,inner sep=1pt,color=NavyBlue, label=below:{\tiny $31212$}] (31212) at (2*\x1 + 2*\xy1 + \y1, 2*\x2 + 2*\xy2 + \y2)  {};	
	\node[circle,fill,inner sep=1pt,color=NavyBlue, label=below:{\tiny $31221$}] (31221) at (2*\x1 + 2*\xy1 + 2*\y1, 2*\x2 + 2*\xy2 + 2*\y2)  {};	
	\node[circle,fill,inner sep=1pt,color=NavyBlue, label=below:{\tiny $32121$}] (32121) at (2*\x1 + 2*\xy1 + 2*\y1 + \z1, 2*\x2 + 2*\xy2 + 2*\y2 + \z2)  {};	
	\node[circle,fill,inner sep=1pt,color=NavyBlue, label=below:{\tiny $32112$}] (32112) at (2*\x1 + 2*\xy1 + \y1 + \z1, 2*\x2 + 2*\xy2 + \y2 + \z2)  {};
	\node[circle,fill,inner sep=1pt,color=NavyBlue, label=below:{\tiny $23112$}] (23112) at (\x1 + 2*\xy1 + \y1 + \z1, \x2 + 2*\xy2 + \y2 + \z2)  {};
	\node[circle,fill,inner sep=1pt,color=NavyBlue, label=below:{\tiny $23121$}] (23121) at (\x1 + 2*\xy1 + 2*\y1 + \z1, \x2 + 2*\xy2 + 2*\y2 + \z2)  {};	
	\node[circle,fill,inner sep=1pt,color=NavyBlue, label=below:{\tiny $22311$}] (22311) at (2*\xy1 + 2*\y1 + 2*\z1, 2*\xy2 + 2*\y2 + 2*\z2)  {};	
	\node[circle,fill,inner sep=1pt,color=NavyBlue, label=below:{\tiny $22131$}] (22131) at (\xy1 + 2*\y1 + 2*\z1, \xy2 + 2*\y2 + 2*\z2)  {};
	\node[circle,fill,inner sep=1pt,color=NavyBlue, label=below:{\tiny $22113$}] (22113) at (2*\y1 + 2*\z1, 2*\y2 + 2*\z2)  {};		
	\node[circle,fill,inner sep=1pt,color=NavyBlue, label=below:{\tiny $23211$}] (23211) at (\x1 + 2*\xy1 + 2*\y1 + 2*\z1, \x2 + 2*\xy2 + 2*\y2 + 2*\z2)  {};		
	\node[circle,fill,inner sep=1pt,color=NavyBlue, label=below:{\tiny $32211$}] (32211) at (2*\x1 + 2*\xy1 + 2*\y1 + 2*\z1, 2*\x2 + 2*\xy2 + 2*\y2 + 2*\z2)  {};			

	\draw[color=NavyBlue] (11223) -- (11232) -- (12132);
	\draw[color=NavyBlue] (11223) -- (12123) -- (12132) -- (12312) -- (13212);
	\draw[color=NavyBlue] (11232) -- (11322) -- (13122) -- (13212);
	\draw[color=NavyBlue] (12312) -- (12321);
	\draw[color=NavyBlue] (12123) -- (12213) -- (12231) -- (12321);		
	\draw[color=NavyBlue] (21123) -- (21132) -- (21312) -- (21321);
	\draw[color=NavyBlue] (21123) -- (21213) -- (21231) -- (21321);
	\draw[color=NavyBlue] (12123) -- (21123);
	\draw[color=NavyBlue] (12213) -- (21213);
	\draw[color=NavyBlue] (12312) -- (21312);
	\draw[color=NavyBlue] (12321) -- (21321);
	\draw[color=NavyBlue] (12231) -- (21231);
	\draw[color=NavyBlue] (12132) -- (21132);			
	\draw[color=NavyBlue] (21213) -- (22113) -- (22131) -- (22311) -- (23211);
	\draw[color=NavyBlue] (21231) -- (22131);			
	\draw[color=NavyBlue] (12321) -- (13221) -- (31221);				
	\draw[color=NavyBlue] (13122) -- (31122) -- (31212);			
	\draw[color=NavyBlue] (13212) -- (31212) -- (32112) -- (32121);
	\draw[color=NavyBlue] (31212) -- (31221) -- (32121);
	\draw[color=NavyBlue] (21312) -- (23112) -- (32112);				
	\draw[color=NavyBlue] (21321) -- (23121) -- (32121);
	\draw[color=NavyBlue] (13221) -- (13212);	
	\draw[color=NavyBlue] (23112) -- (23121) -- (23211) -- (32211);				
	\draw[color=NavyBlue] (23112) -- (32112) -- (32121) -- (32211);	
\end{scope}
\end{tikzpicture}
    
    \caption{The lattice of multipermutations.}
    \label{fig.lattice_of_multipermutations}
\end{figure}

For $\mathbf{s} = (s_1,s_2,\ldots, s_n)$, let $\oru{\mathbf{s}}$ be the graph on vertex set $[n+1]$ with $s_i+1$ edges between $i$ and $i+1$ for each $i\in[n]$. Its flow polytope is a product of $n$ simplices, $\Delta_{s_1}\times \dots \Delta_{s_n}$.
Let $F$ be the framing induced by a planar drawing (i.e. without edges intersecting) of $\oru{\mathbf{s}}$. For an example see~\Cref{fig.lattice_of_multipermutations}.  

\begin{theorem}
    \label{thm.lattice_of_multipermutations}
    The framing lattice $\scrL_{\oru{\mathbf{s}}, F}$ is the lattice of multipermutations $\scrM_{\mathbf{s}}$.  
\end{theorem}

\begin{proof}
    The argument in the case of the weak order (see Example~\ref{ex.oruga_with_planar_framing}) generalizes naturally to the setting of multipermutations. More precisely, multipermutations are in bijection with maximal cliques as follows. The maximal clique associated to a multipermutation is the set of routes obtained by starting with the bottom-most route and proceeding to bump up the edges corresponding to numbers in the multipermutation (read left to right) until arriving to the top-most route. 
    An example of this generalized bijection is shown in Figure~\ref{fig.lattice_of_multipermutations}.
    Under this bijection, applying an increasing transposition corresponds to a $\ccw$-rotation. Thus, the lattice of multipermutations $\scrM_{\mathbf{s}}$ is isomorphic to the framing lattice $\scrL_{\oru{\mathbf{s}}, F}$.    
\end{proof}

\section{The cross-Tamari lattice}
\label{sec_crossTamari}

In this section we introduce a new family of lattices which we call cross-Tamari lattices. This family includes several well studied lattices in the literature: 
\begin{enumerate}
    \item the $\nu$-Tamari lattices of Pr\'eville--Ratelle and Viennot \cite{PV17} (Cf.~\cite{CPS20}), 
    \item the principal order ideal induced by $\nu$ in Young's lattice~\cite{BGMY23} (also known as the $\nu$-Dyck lattice~\cite{CC23} or Stanley's distributive lattice), 
    \item the alt-$\nu$-Tamari lattices~\cite{CC23}, and
    \item the $\varepsilon$-Cambrian lattices of type $A$~\cite{Rea06}.  
    \item the $(\varepsilon,I,J)$-Cambrian lattices~\cite{Pil20}.
\end{enumerate}

Our starting point is a beautiful connection between the first two examples and triangulations of flow polytopes presented in~\cite{BGMY23}.\footnote{These triangulations appear to be highly special, as they repeatedly arise in the literature under various guises and across a range of different contexts. In particular, they occur as triangulations of integrally equivalent polytopes associated with root polytopes, order polytopes, and products of two simplices. See~\cite[Section~1.4]{CPS19} for more detailed references, as well as~\cite{ceballos_Dyck_2015}.} 
The authors showed that both the Hasse diagram of the $\nu$-Tamari lattice and the principal order ideal induced by $\nu$ in Young's lattice appear as dual graphs of framed triangulations of a single flow polytope $\calF_{\car(\nu)}$, obtained using different framings.
We will introduce the concept of cross-Tamari lattices and show that they are framing lattices of a framed graph induced by a cross-shaped grid (\Cref{thm_crossTam}). 
This immediately implies strong non-trivial consequences that we summarize in~\Cref{cor_crossTam_consequences1,cor_crossTam_consequences2}.   
We will also see that cross-shaped grids that can be obtained from each other by permutations of rows and columns give rise to the same flow polytope but (possibly) different framed triangulations and lattices (\Cref{prop_crossTam_permutations}), recovering the results in~\cite{BGMY23} as a particular case.


\subsection{Cross-shaped grids.}
\begin{figure}[htb]
    \centering
    \begin{tikzpicture}
\begin{scope}[xshift = 0, yshift = 0, scale=0.5]
    \draw[very thin, color=gray!70] (0,0) grid (1,1);
    \draw[very thin, color=gray!70] (0,1) grid (2,2);      
    \draw[very thin, color=gray!70] (0,2) grid (2,3);
    \draw[very thin, color=gray!70] (0,3) grid (4,4);
    \draw[very thin, color=gray!70] (0,4) grid (5,4);

    \node[circle, draw, inner sep=2pt, fill, blue]  at (0,0)  {};
    \node[circle, draw, inner sep=1pt, fill]  at (0,1)  {};
    \node[circle, draw, inner sep=1pt, fill]  at (0,2)  {};
    \node[circle, draw, inner sep=1pt, fill]  at (0,3)  {};    
    \node[circle, draw, inner sep=1pt, fill]  at (0,4)  {};

    \node[circle, draw, inner sep=1pt, fill]  at (1,0)  {};
    \node[circle, draw, inner sep=1pt, fill]  at (1,1)  {};
    \node[circle, draw, inner sep=1pt, fill]  at (1,2)  {};
    \node[circle, draw, inner sep=1pt, fill]  at (1,3)  {};    
    \node[circle, draw, inner sep=1pt, fill]  at (1,4)  {};    

    \node[circle, draw, inner sep=1pt, fill]  at (2,1)  {};
    \node[circle, draw, inner sep=1pt, fill]  at (2,2)  {};    
    \node[circle, draw, inner sep=1pt, fill]  at (2,3)  {};    
    \node[circle, draw, inner sep=1pt, fill]  at (2,4)  {};

    \node[circle, draw, inner sep=1pt, fill]  at (3,3)  {};    
    \node[circle, draw, inner sep=1pt, fill]  at (3,4)  {};

    \node[circle, draw, inner sep=1pt, fill]  at (4,3)  {};    
    \node[circle, draw, inner sep=1pt, fill]  at (4,4)  {};

    \node[circle, draw, inner sep=1pt, fill]  at (5,4)  {};

    \node[]  at (0,4.6)  {\footnotesize {$1$}};
    \node[]  at (1,4.6)  {\footnotesize {$2$}};
    \node[]  at (2,4.6)  {\footnotesize {$3$}};
    \node[]  at (3,4.6)  {\footnotesize {$4$}};    
    \node[]  at (4,4.6)  {\footnotesize {$5$}};
    \node[]  at (5,4.6)  {\footnotesize {$6$}};
    
    \node[]  at (-.6,0)  {\footnotesize {$\overline{5}$}};
    \node[]  at (-.6,1)  {\footnotesize {$\overline{4}$}};
    \node[]  at (-.6,2)  {\footnotesize {$\overline{3}$}};
    \node[]  at (-.6,3)  {\footnotesize {$\overline{2}$}};    
    \node[]  at (-.6,4)  {\footnotesize {$\overline{1}$}};
\end{scope}

\begin{scope}[xshift = 200, yshift = 0, scale=0.5]
    \draw[very thin, color=gray!70] (1,0) grid (3,1);
    \draw[very thin, color=gray!70] (0,1) grid (4,2);
    \draw[very thin, color=gray!70] (4,2) grid (5,2);
    \draw[very thin, color=gray!70] (1,2) grid (3,3);
    \draw[very thin, color=gray!70] (2,3) grid (3,4);

    \node[circle, draw, inner sep=1pt, fill]  at (0,1)  {};
    \node[circle, draw, inner sep=1pt, fill]  at (0,2)  {};

    \node[circle, draw, inner sep=2pt, fill, red]  at (1,0)  {};
    \node[circle, draw, inner sep=1pt, fill]  at (1,1)  {};
    \node[circle, draw, inner sep=1pt, fill]  at (1,2)  {};
    \node[circle, draw, inner sep=1pt, fill]  at (1,3)  {};    

    \node[circle, draw, inner sep=1pt, fill]  at (2,0)  {};
    \node[circle, draw, inner sep=1pt, fill]  at (2,1)  {};
    \node[circle, draw, inner sep=1pt, fill]  at (2,2)  {};    
    \node[circle, draw, inner sep=1pt, fill]  at (2,3)  {};    
    \node[circle, draw, inner sep=1pt, fill]  at (2,4)  {};

    \node[circle, draw, inner sep=1pt, fill]  at (3,0)  {};
    \node[circle, draw, inner sep=1pt, fill]  at (3,1)  {};    
    \node[circle, draw, inner sep=1pt, fill]  at (3,2)  {};
    \node[circle, draw, inner sep=1pt, fill]  at (3,3)  {};    
    \node[circle, draw, inner sep=1pt, fill]  at (3,4)  {};

    \node[circle, draw, inner sep=1pt, fill]  at (4,1)  {};    
    \node[circle, draw, inner sep=1pt, fill]  at (4,2)  {};

    \node[circle, draw, inner sep=1pt, fill]  at (5,2)  {};

    \node[]  at (0,4.6)  {\footnotesize $4$};
    \node[]  at (1,4.6)  {\footnotesize $3$};
    \node[]  at (2,4.6)  {\footnotesize $2$};
    \node[]  at (3,4.6)  {\footnotesize $1$};    
    \node[]  at (4,4.6)  {\footnotesize $5$};
    \node[]  at (5,4.6)  {\footnotesize $6$};
    
    \node[]  at (-.6,0)  {\footnotesize $\overline{4}$};
    \node[]  at (-.6,1)  {\footnotesize $\overline{2}$};
    \node[]  at (-.6,2)  {\footnotesize $\overline{1}$};
    \node[]  at (-.6,3)  {\footnotesize $\overline{3}$};    
    \node[]  at (-.6,4)  {\footnotesize $\overline{5}$};

\end{scope}

\end{tikzpicture}
    \caption{Two cross-shaped grids with proper labelings. Note that they are related by a sequence of row and column commutations that preserve the cross-shaped property. }
    \label{cross-shaped_grid_ex}
\end{figure}
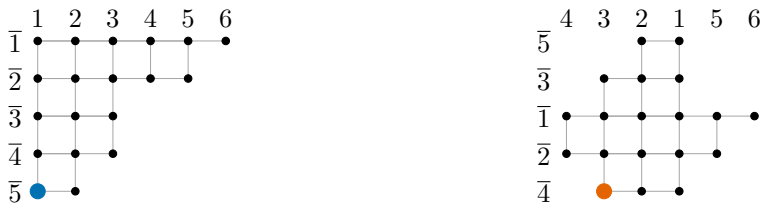

Let $D$ be a set of lattice points in $\bbZ^2$. 
We say that $D$ is \defn{horizontally connected} if for any pair of points $(x,y)$ and $(x',y)$ in $D$ we have $(z,y) \in D$ for all $x< z < x'$. 
Let $\row_D(z)$ denote the set of points in $D$ with $y$-coordinate $z$. 
We say that $D$ is \defn{horizontally nested} if the $x$-coordinates of the points in $\row_D(v)$ are a subset of the $x$-coordinates of the points in $\row_D(w)$ whenever $|\row_D(v)|\leq |\row_D(w)|$.
Similarly, we define \defn{vertically connected} and \defn{vertically nested}. 
A set of lattice points $D \subseteq \bbZ^2$ is a \defn{cross-shaped grid} if it is both horizontally and vertically connected, and horizontally and vertically nested.
Two examples are illustrated in~\Cref{cross-shaped_grid_ex}.

\begin{remark}
    Replacing lattice points by unit boxes, cross-shaped grids coincide with the moon polyominos already used in the literature, see i.e.~\cite{rubey_maximal_2012,serrano_maximal_2012}. We keep the terms ``cross-shaped grid" and ``cross-Tamari" for simplicity.   
\end{remark}

\subsection{The cross-Tamari lattice}

Let $D$ be a cross-shaped grid.
Two distinct points $p,p'\in D$ are \defn{incompatible} if one of them is strictly north-east of the other and every lattice point in the smallest rectangle containing $p$ and $p'$ belongs to $D$.  
Two points are \defn{compatible} if they are not incompatible. 
A \defn{maximal filling} $T$ in $D$ is a maximal set of pairwise compatible points.
If two maximal fillings $T\neq T'$ differ by one single element $T\smallsetminus\{p\}=T'\smallsetminus\{p'\}$ where $p'$ is located strictly north-east of $p$, then we say the $T'$ is obtainable from $T$ by an \defn{increasing rotation}.
The \defn{cross-Tamari order} $\Tam(D)$ is the poset of maximal fillings in $D$ where $T\preceq_D T'$ if $T'$ can be obtained from $T$ by a sequence of increasing rotations. 

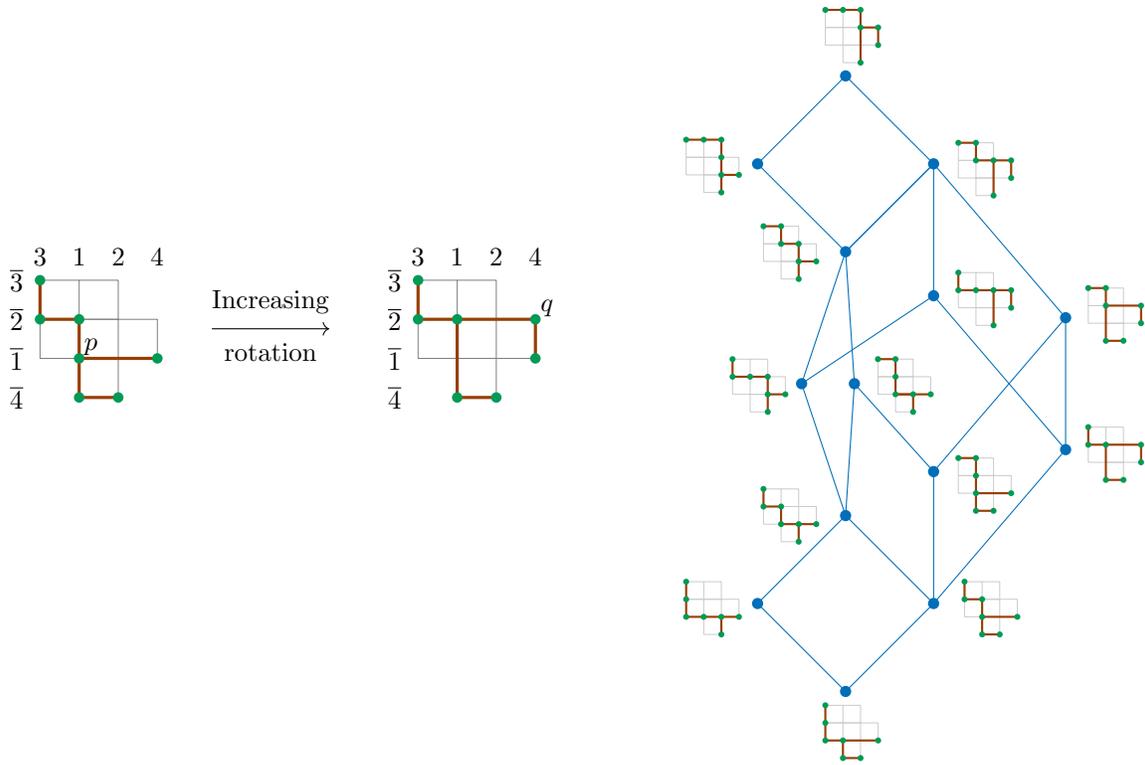
\begin{figure}[htb]
    \centering
    \begin{tikzpicture}
\begin{scope}[xshift = 0, yshift = 0, scale=1.3]

    \begin{scope}[xshift = -40, yshift= -20, scale=0.4]
        \draw[->] (0,0) -- (3,0);
        \node[] (a) at (1.5,0.7) {\footnotesize{Increasing}};
        \node[] (b) at (1.5,-0.6) {\footnotesize{rotation}};
    \end{scope}
    
    \begin{scope}[xshift = -90, yshift= -40, scale=0.4]
        \draw[very thin, color=gray] (1,0) grid (2,1);
        \draw[very thin, color=gray] (0,1) grid (3,2);
        \draw[very thin, color=gray] (0,2) grid (2,3);

        \node[]  at (0,3.6)  {\footnotesize {$3$}};
        \node[]  at (1,3.6)  {\footnotesize {$1$}};
        \node[]  at (2,3.6)  {\footnotesize {$2$}};
        \node[]  at (3,3.6)  {\footnotesize {$4$}};
    
        \node[]  at (-0.6,0)  {\footnotesize {$\overline{4}$}};
        \node[]  at (-0.6,1)  {\footnotesize {$\overline{1}$}};
        \node[]  at (-0.6,2)  {\footnotesize {$\overline{2}$}};
        \node[]  at (-0.6,3)  {\footnotesize {$\overline{3}$}};

        \draw[very thick, color=RawSienna] (1,2) -- (0,2) -- (0,3);
	      \draw[very thick, color=RawSienna] (1,2) -- (1,1) -- (3,1);
	      \draw[very thick, color=RawSienna] (1,1) -- (1,0);
	      \draw[very thick, color=RawSienna] (1,0) -- (2,0);
	             
	      \node[circle,fill,inner sep=1.4pt,color=ForestGreen] at (1,2)  {};
	      \node[circle,fill,inner sep=1.4pt,color=ForestGreen] at (1,1)  {}; 
        \node[circle,fill,inner sep=1.4pt,color=ForestGreen] at (3,1)  {};     
	      \node[circle,fill,inner sep=1.4pt,color=ForestGreen] at (1,0)  {};		     
        \node[circle,fill,inner sep=1.4pt,color=ForestGreen] at (0,2)  {};
	      \node[circle,fill,inner sep=1.4pt,color=ForestGreen] at (0,3)  {};		     
        \node[circle,fill,inner sep=1.4pt,color=ForestGreen] at (2,0)  {};   

        \node[]  at (1.3,1.3)  {\footnotesize {$p$}};
    \end{scope}

    \begin{scope}[xshift = 20, yshift= -40, scale=0.4]
        \draw[very thin, color=gray] (1,0) grid (2,1);
        \draw[very thin, color=gray] (0,1) grid (3,2);
        \draw[very thin, color=gray] (0,2) grid (2,3);

        \node[]  at (0,3.6)  {\footnotesize {$3$}};
        \node[]  at (1,3.6)  {\footnotesize {$1$}};
        \node[]  at (2,3.6)  {\footnotesize {$2$}};
        \node[]  at (3,3.6)  {\footnotesize {$4$}};
    
        \node[]  at (-0.6,0)  {\footnotesize {$\overline{4}$}};
        \node[]  at (-0.6,1)  {\footnotesize {$\overline{1}$}};
        \node[]  at (-0.6,2)  {\footnotesize {$\overline{2}$}};
        \node[]  at (-0.6,3)  {\footnotesize {$\overline{3}$}};

        \draw[very thick, color=RawSienna] (1,2) -- (0,2) -- (0,3);
	      \draw[very thick, color=RawSienna] (1,2) -- (3,2) -- (3,1);
	      \draw[very thick, color=RawSienna] (1,2) -- (1,0);
	      \draw[very thick, color=RawSienna] (1,0) -- (2,0);
	             
	      \node[circle,fill,inner sep=1.4pt,color=ForestGreen] at (1,2)  {};
	      \node[circle,fill,inner sep=1.4pt,color=ForestGreen] at (3,2)  {};		    
        \node[circle,fill,inner sep=1.4pt,color=ForestGreen] at (3,1)  {};     
	      \node[circle,fill,inner sep=1.4pt,color=ForestGreen] at (1,0)  {};		     
        \node[circle,fill,inner sep=1.4pt,color=ForestGreen] at (0,2)  {};
	      \node[circle,fill,inner sep=1.4pt,color=ForestGreen] at (0,3)  {};		     
        \node[circle,fill,inner sep=1.4pt,color=ForestGreen] at (2,0)  {};
    
        \node[]  at (3.3,2.3)  {\footnotesize {$q$}};
    \end{scope}
    
    \begin{scope}[xshift = 170, yshift= -100, scale=0.9]

	\node[circle,fill,inner sep=1.5pt,color=NavyBlue] at (-1,-1)  {};
	  \node[circle,fill,inner sep=1.5pt,color=NavyBlue] at (-2,0)  {};
	\node[circle,fill,inner sep=1.5pt,color=NavyBlue] at (0,0)  {};	
	\node[circle,fill,inner sep=1.5pt,color=NavyBlue] at (-1.5,2.5)  {};
	\node[circle,fill,inner sep=1.5pt,color=NavyBlue] at (-1,1)  {};
	\node[circle,fill,inner sep=1.5pt,color=NavyBlue] at (1.5,1.75)  {};	
	\node[circle,fill,inner sep=1.5pt,color=NavyBlue] at (0,3.5)  {};
	\node[circle,fill,inner sep=1.5pt,color=NavyBlue] at (0,1.5)  {};	
	\node[circle,fill,inner sep=1.5pt,color=NavyBlue] at (-1,4)  {};
	\node[circle,fill,inner sep=1.5pt,color=NavyBlue] at (-0.9,2.5)  {};
	\node[circle,fill,inner sep=1.5pt,color=NavyBlue] at (1.5,3.25)  {};	
	\node[circle,fill,inner sep=1.5pt,color=NavyBlue] at (0,5)  {};
 	\node[circle,fill,inner sep=1.5pt,color=NavyBlue] at (-1,6)  {};
  	\node[circle,fill,inner sep=1.5pt,color=NavyBlue] at (-2,5)  {};
    
	\draw[color=NavyBlue] (-1,-1) -- (-2,0) -- (-1,1) -- (0,0) -- (-1,-1);	
	\draw[color=NavyBlue] (-1,1) -- (-1.5,2.5) -- (0,3.5) -- (1.5,1.75) -- (0,0);

 	\draw[color=NavyBlue] (0,1.5) -- (-0.9,2.5) -- (-1,4) -- (0,5) -- (1.5,3.25) -- (0,1.5);	

  	\draw[color=NavyBlue] (0,0) -- (0,1.5);
  	\draw[color=NavyBlue] (-1,1) -- (-0.9,2.5);
  	\draw[color=NavyBlue] (-1.5,2.5) -- (-1,4);
  	\draw[color=NavyBlue] (0,3.5) -- (0,5);
  	\draw[color=NavyBlue] (1.5,1.75) -- (1.5,3.25);   

   	\draw[color=NavyBlue] (-1,4) -- (0,5) -- (-1,6) -- (-2,5) -- (-1,4);

	\begin{scope}[xshift = -35, yshift= -50, scale=0.20]
        \draw[very thin, color=black!20] (1,0) grid (2,1);
        \draw[very thin, color=black!20] (0,1) grid (3,2);
        \draw[very thin, color=black!20] (0,2) grid (2,3);

	     \draw[thick, color=RawSienna] (0,1) -- (0,3);
	     \draw[thick, color=RawSienna] (0,1) -- (3,1);
	     \draw[thick, color=RawSienna] (1,1) -- (1,0);
	     \draw[thick, color=RawSienna] (1,0) -- (2,0);
	             
	     \node[circle,fill,inner sep=0.8pt,color=ForestGreen] at (0,1)  {};
	     \node[circle,fill,inner sep=0.8pt,color=ForestGreen] at (1,1)  {};		     \node[circle,fill,inner sep=0.8pt,color=ForestGreen] at (3,1)  {};     
	     \node[circle,fill,inner sep=0.8pt,color=ForestGreen] at (1,0)  {};		     \node[circle,fill,inner sep=0.8pt,color=ForestGreen] at (0,2)  {};
	     \node[circle,fill,inner sep=0.8pt,color=ForestGreen] at (0,3)  {};		     \node[circle,fill,inner sep=0.8pt,color=ForestGreen] at (2,0)  {};     
    \end{scope}    
    
	\begin{scope}[xshift = -80, yshift= -10, scale=0.20]
        \draw[very thin, color=black!20] (1,0) grid (2,1);
        \draw[very thin, color=black!20] (0,1) grid (3,2);
        \draw[very thin, color=black!20] (0,2) grid (2,3);

	     \draw[thick, color=RawSienna] (0,1) -- (0,3);
	     \draw[thick, color=RawSienna] (0,1) -- (3,1);
	     \draw[thick, color=RawSienna] (2,1) -- (2,0);

	     \node[circle,fill,inner sep=0.8pt,color=ForestGreen] at (0,1)  {};
	     \node[circle,fill,inner sep=0.8pt,color=ForestGreen] at (2,1)  {};		     \node[circle,fill,inner sep=0.8pt,color=ForestGreen] at (3,1)  {};     
	     \node[circle,fill,inner sep=0.8pt,color=ForestGreen] at (1,1)  {};		     \node[circle,fill,inner sep=0.8pt,color=ForestGreen] at (0,2)  {};
	     \node[circle,fill,inner sep=0.8pt,color=ForestGreen] at (0,3)  {};		     \node[circle,fill,inner sep=0.8pt,color=ForestGreen] at (2,0)  {};     
    \end{scope}        
    
	\begin{scope}[xshift = 10, yshift= -10, scale=0.20]
        \draw[very thin, color=black!20] (1,0) grid (2,1);
        \draw[very thin, color=black!20] (0,1) grid (3,2);
        \draw[very thin, color=black!20] (0,2) grid (2,3);

	     \draw[thick, color=RawSienna] (1,2) -- (0,2) -- (0,3);
	     \draw[thick, color=RawSienna] (1,2) -- (1,1) -- (3,1);
	     \draw[thick, color=RawSienna] (1,1) -- (1,0);
	     \draw[thick, color=RawSienna] (1,0) -- (2,0);
	             
	     \node[circle,fill,inner sep=0.8pt,color=ForestGreen] at (1,2)  {};
	     \node[circle,fill,inner sep=0.8pt,color=ForestGreen] at (1,1)  {};		     \node[circle,fill,inner sep=0.8pt,color=ForestGreen] at (3,1)  {};     
	     \node[circle,fill,inner sep=0.8pt,color=ForestGreen] at (1,0)  {};		     \node[circle,fill,inner sep=0.8pt,color=ForestGreen] at (0,2)  {};
	     \node[circle,fill,inner sep=0.8pt,color=ForestGreen] at (0,3)  {};		     \node[circle,fill,inner sep=0.8pt,color=ForestGreen] at (2,0)  {};     
    \end{scope}            
    
	\begin{scope}[xshift = -55, yshift= 20, scale=0.20]
        \draw[very thin, color=black!20] (1,0) grid (2,1);
        \draw[very thin, color=black!20] (0,1) grid (3,2);
        \draw[very thin, color=black!20] (0,2) grid (2,3);

	     \draw[thick, color=RawSienna] (1,2) -- (0,2) -- (0,3);
	     \draw[thick, color=RawSienna] (1,2) -- (1,1) -- (3,1);
	     \draw[thick, color=RawSienna] (2,1) -- (2,0);

	     \node[circle,fill,inner sep=0.8pt,color=ForestGreen] at (1,2)  {};
	     \node[circle,fill,inner sep=0.8pt,color=ForestGreen] at (1,1)  {};		     \node[circle,fill,inner sep=0.8pt,color=ForestGreen] at (3,1)  {};     
	     \node[circle,fill,inner sep=0.8pt,color=ForestGreen] at (2,1)  {};		     \node[circle,fill,inner sep=0.8pt,color=ForestGreen] at (0,2)  {};
	     \node[circle,fill,inner sep=0.8pt,color=ForestGreen] at (0,3)  {};		     \node[circle,fill,inner sep=0.8pt,color=ForestGreen] at (2,0)  {};     
    \end{scope}

	\begin{scope}[xshift = 50, yshift= 40, scale=0.20]
        \draw[very thin, color=black!20] (1,0) grid (2,1);
        \draw[very thin, color=black!20] (0,1) grid (3,2);
        \draw[very thin, color=black!20] (0,2) grid (2,3);

	     \draw[thick, color=RawSienna] (1,2) -- (0,2) -- (0,3);
	     \draw[thick, color=RawSienna] (1,2) -- (3,2) -- (3,1);
	     \draw[thick, color=RawSienna] (1,2) -- (1,0);
	     \draw[thick, color=RawSienna] (1,0) -- (2,0);
	             
	     \node[circle,fill,inner sep=0.8pt,color=ForestGreen] at (1,2)  {};
	     \node[circle,fill,inner sep=0.8pt,color=ForestGreen] at (3,2)  {};		     \node[circle,fill,inner sep=0.8pt,color=ForestGreen] at (3,1)  {};     
	     \node[circle,fill,inner sep=0.8pt,color=ForestGreen] at (1,0)  {};		     \node[circle,fill,inner sep=0.8pt,color=ForestGreen] at (0,2)  {};
	     \node[circle,fill,inner sep=0.8pt,color=ForestGreen] at (0,3)  {};		     \node[circle,fill,inner sep=0.8pt,color=ForestGreen] at (2,0)  {};     
    \end{scope}                
    
	\begin{scope}[xshift = -65, yshift= 62, scale=0.20]
        \draw[very thin, color=black!20] (1,0) grid (2,1);
        \draw[very thin, color=black!20] (0,1) grid (3,2);
        \draw[very thin, color=black!20] (0,2) grid (2,3);

	     \draw[thick, color=RawSienna] (2,2) -- (0,2) -- (0,3);
	     \draw[thick, color=RawSienna] (2,2) -- (2,1) -- (3,1);
	     \draw[thick, color=RawSienna] (2,1) -- (2,0);
	             
	     \node[circle,fill,inner sep=0.8pt,color=ForestGreen] at (1,2)  {};
	     \node[circle,fill,inner sep=0.8pt,color=ForestGreen] at (2,2)  {};		     \node[circle,fill,inner sep=0.8pt,color=ForestGreen] at (3,1)  {};     
	     \node[circle,fill,inner sep=0.8pt,color=ForestGreen] at (2,1)  {};		     \node[circle,fill,inner sep=0.8pt,color=ForestGreen] at (0,2)  {};
	     \node[circle,fill,inner sep=0.8pt,color=ForestGreen] at (0,3)  {};		     \node[circle,fill,inner sep=0.8pt,color=ForestGreen] at (2,0)  {};     
    \end{scope}                
    
	\begin{scope}[xshift = 8, yshift= 90, scale=0.20]
        \draw[very thin, color=black!20] (1,0) grid (2,1);
        \draw[very thin, color=black!20] (0,1) grid (3,2);
        \draw[very thin, color=black!20] (0,2) grid (2,3);

	     \draw[thick, color=RawSienna] (2,2) -- (0,2) -- (0,3);
	     \draw[thick, color=RawSienna] (2,2) -- (3,2) -- (3,1);
	     \draw[thick, color=RawSienna] (2,2) -- (2,0);
	             
	     \node[circle,fill,inner sep=0.8pt,color=ForestGreen] at (1,2)  {};
	     \node[circle,fill,inner sep=0.8pt,color=ForestGreen] at (2,2)  {};		     \node[circle,fill,inner sep=0.8pt,color=ForestGreen] at (3,1)  {};     
	     \node[circle,fill,inner sep=0.8pt,color=ForestGreen] at (3,2)  {};		     \node[circle,fill,inner sep=0.8pt,color=ForestGreen] at (0,2)  {};
	     \node[circle,fill,inner sep=0.8pt,color=ForestGreen] at (0,3)  {};		     \node[circle,fill,inner sep=0.8pt,color=ForestGreen] at (2,0)  {};     
    \end{scope}                    
    
	\begin{scope}[xshift = 8, yshift= 132, scale=0.20]
        \draw[very thin, color=black!20] (1,0) grid (2,1);
        \draw[very thin, color=black!20] (0,1) grid (3,2);
        \draw[very thin, color=black!20] (0,2) grid (2,3);

	     \draw[thick, color=RawSienna] (0,3) -- (1,3) -- (1,2);
	     \draw[thick, color=RawSienna] (1,2) -- (3,2) -- (3,1);
	     \draw[thick, color=RawSienna] (2,2) -- (2,0);
	             
	     \node[circle,fill,inner sep=0.8pt,color=ForestGreen] at (1,2)  {};
	     \node[circle,fill,inner sep=0.8pt,color=ForestGreen] at (2,2)  {};		    
      \node[circle,fill,inner sep=0.8pt,color=ForestGreen] at (3,1)  {};     
	     \node[circle,fill,inner sep=0.8pt,color=ForestGreen] at (3,2)  {};		     
      \node[circle,fill,inner sep=0.8pt,color=ForestGreen] at (1,3)  {};
	     \node[circle,fill,inner sep=0.8pt,color=ForestGreen] at (0,3)  {};		     
      \node[circle,fill,inner sep=0.8pt,color=ForestGreen] at (2,0)  {};     
    \end{scope}                        
    
	\begin{scope}[xshift = 50, yshift= 85, scale=0.20]
        \draw[very thin, color=black!20] (1,0) grid (2,1);
        \draw[very thin, color=black!20] (0,1) grid (3,2);
        \draw[very thin, color=black!20] (0,2) grid (2,3);

	     \draw[thick, color=RawSienna] (0,3) -- (1,3) -- (1,2);
	     \draw[thick, color=RawSienna] (1,2) -- (3,2) -- (3,1);
	     \draw[thick, color=RawSienna] (1,2) -- (1,0) -- (2,0);
	             
	     \node[circle,fill,inner sep=0.8pt,color=ForestGreen] at (1,2)  {};
	     \node[circle,fill,inner sep=0.8pt,color=ForestGreen] at (1,0)  {};		     \node[circle,fill,inner sep=0.8pt,color=ForestGreen] at (3,1)  {};     
	     \node[circle,fill,inner sep=0.8pt,color=ForestGreen] at (3,2)  {};		     \node[circle,fill,inner sep=0.8pt,color=ForestGreen] at (1,3)  {};
	     \node[circle,fill,inner sep=0.8pt,color=ForestGreen] at (0,3)  {};		     \node[circle,fill,inner sep=0.8pt,color=ForestGreen] at (2,0)  {};     
    \end{scope}                            
    
	\begin{scope}[xshift = 8, yshift= 30, scale=0.20]
        \draw[very thin, color=black!20] (1,0) grid (2,1);
        \draw[very thin, color=black!20] (0,1) grid (3,2);
        \draw[very thin, color=black!20] (0,2) grid (2,3);

	     \draw[thick, color=RawSienna] (0,3) -- (1,3) -- (1,2);
	     \draw[thick, color=RawSienna] (1,1) -- (3,1);
	     \draw[thick, color=RawSienna] (1,2) -- (1,0) -- (2,0);
	             
	     \node[circle,fill,inner sep=0.8pt,color=ForestGreen] at (1,2)  {};
	     \node[circle,fill,inner sep=0.8pt,color=ForestGreen] at (1,0)  {};		     \node[circle,fill,inner sep=0.8pt,color=ForestGreen] at (3,1)  {};     
	     \node[circle,fill,inner sep=0.8pt,color=ForestGreen] at (1,1)  {};		     \node[circle,fill,inner sep=0.8pt,color=ForestGreen] at (1,3)  {};
	     \node[circle,fill,inner sep=0.8pt,color=ForestGreen] at (0,3)  {};		     \node[circle,fill,inner sep=0.8pt,color=ForestGreen] at (2,0)  {};     
    \end{scope}                            
    
	\begin{scope}[xshift = -18, yshift= 62, scale=0.20]
        \draw[very thin, color=black!20] (1,0) grid (2,1);
        \draw[very thin, color=black!20] (0,1) grid (3,2);
        \draw[very thin, color=black!20] (0,2) grid (2,3);

	     \draw[thick, color=RawSienna] (0,3) -- (1,3) -- (1,2);
	     \draw[thick, color=RawSienna] (1,1) -- (3,1);
	     \draw[thick, color=RawSienna] (1,2) -- (1,1) -- (2,1) -- (2,0);
	             
	     \node[circle,fill,inner sep=0.8pt,color=ForestGreen] at (1,2)  {};
	     \node[circle,fill,inner sep=0.8pt,color=ForestGreen] at (2,1)  {};		     \node[circle,fill,inner sep=0.8pt,color=ForestGreen] at (3,1)  {};     
	     \node[circle,fill,inner sep=0.8pt,color=ForestGreen] at (1,1)  {};		     \node[circle,fill,inner sep=0.8pt,color=ForestGreen] at (1,3)  {};
	     \node[circle,fill,inner sep=0.8pt,color=ForestGreen] at (0,3)  {};		     \node[circle,fill,inner sep=0.8pt,color=ForestGreen] at (2,0)  {};     
    \end{scope}                                
    
	\begin{scope}[xshift = -55, yshift= 105, scale=0.20]
        \draw[very thin, color=black!20] (1,0) grid (2,1);
        \draw[very thin, color=black!20] (0,1) grid (3,2);
        \draw[very thin, color=black!20] (0,2) grid (2,3);

	     \draw[thick, color=RawSienna] (0,3) -- (1,3) -- (1,2) -- (2,2);
	     \draw[thick, color=RawSienna] (2,1) -- (3,1);
	     \draw[thick, color=RawSienna] (2,2) -- (2,0);
	             
	     \node[circle,fill,inner sep=0.8pt,color=ForestGreen] at (1,2)  {};
	     \node[circle,fill,inner sep=0.8pt,color=ForestGreen] at (2,1)  {};		     \node[circle,fill,inner sep=0.8pt,color=ForestGreen] at (3,1)  {};     
	     \node[circle,fill,inner sep=0.8pt,color=ForestGreen] at (2,2)  {};		     \node[circle,fill,inner sep=0.8pt,color=ForestGreen] at (1,3)  {};
	     \node[circle,fill,inner sep=0.8pt,color=ForestGreen] at (0,3)  {};		     \node[circle,fill,inner sep=0.8pt,color=ForestGreen] at (2,0)  {};     
    \end{scope}                                    

	\begin{scope}[xshift = -80, yshift= 133, scale=0.20]
        \draw[very thin, color=black!20] (1,0) grid (2,1);
        \draw[very thin, color=black!20] (0,1) grid (3,2);
        \draw[very thin, color=black!20] (0,2) grid (2,3);

	     \draw[thick, color=RawSienna] (0,3) -- (2,3) -- (2,0);
	     \draw[thick, color=RawSienna] (2,1) -- (3,1);
	             
	     \node[circle,fill,inner sep=0.8pt,color=ForestGreen] at (2,3)  {};
	     \node[circle,fill,inner sep=0.8pt,color=ForestGreen] at (2,1)  {};		     \node[circle,fill,inner sep=0.8pt,color=ForestGreen] at (3,1)  {};     
	     \node[circle,fill,inner sep=0.8pt,color=ForestGreen] at (2,2)  {};		     \node[circle,fill,inner sep=0.8pt,color=ForestGreen] at (1,3)  {};
	     \node[circle,fill,inner sep=0.8pt,color=ForestGreen] at (0,3)  {};		     \node[circle,fill,inner sep=0.8pt,color=ForestGreen] at (2,0)  {};     
    \end{scope}                                        
    
	 \begin{scope}[xshift = -35, yshift=175, scale=0.20]
        \draw[very thin, color=black!20] (1,0) grid (2,1);
        \draw[very thin, color=black!20] (0,1) grid (3,2);
        \draw[very thin, color=black!20] (0,2) grid (2,3);

	     \draw[thick, color=RawSienna] (0,3) -- (2,3) -- (2,2) -- (3,2) -- (3,1);
	     \draw[thick, color=RawSienna] (2,2) -- (2,0);
	             
	     \node[circle,fill,inner sep=0.8pt,color=ForestGreen] at (2,3)  {};
	     \node[circle,fill,inner sep=0.8pt,color=ForestGreen] at (3,2)  {};		     \node[circle,fill,inner sep=0.8pt,color=ForestGreen] at (3,1)  {};     
	     \node[circle,fill,inner sep=0.8pt,color=ForestGreen] at (2,2)  {};		     \node[circle,fill,inner sep=0.8pt,color=ForestGreen] at (1,3)  {};
	     \node[circle,fill,inner sep=0.8pt,color=ForestGreen] at (0,3)  {};		     \node[circle,fill,inner sep=0.8pt,color=ForestGreen] at (2,0)  {};     
    \end{scope}                                     
    \end{scope}
    
\end{scope}
\end{tikzpicture}
    \caption{An example of an increasing rotation in a cross-shaped grid (left), and a cross-Tamari lattice (right).}
    \label{cross-tamari_ex}
\end{figure}

\begin{figure}[htb]
    \centering
    \begin{tikzpicture}[scale=0.9]

\begin{scope}[scale=0.4, xshift=0, yshift=40]

    \draw[very thin, color=black] (0,0) grid (0,1);
    \draw[very thin, color=black] (0,1) grid (1,3);
    \draw[very thin, color=black] (1,2) grid (2,3);
    \draw[very thin, color=black] (2,3) grid (3,3);
    
\end{scope}	

\begin{scope}[scale=0.3, xshift=90, yshift=-300]
	\draw[thick, color=NavyBlue]    (0,0) -- (3,5) -- (0,10) -- (-3,7) -- (-3,3) -- (0,0);		
	\node[circle,fill,inner sep=1.5pt,color=NavyBlue] at (0,0)  {};
	\node[circle,fill,inner sep=1.5pt,color=NavyBlue]  at (3,5) {}; 
	\node[circle,fill,inner sep=1.5pt,color=NavyBlue]  at (0,10) {};
	\node[circle,fill,inner sep=1.5pt,color=NavyBlue]  at (-3,7) {};
	\node[circle,fill,inner sep=1.5pt,color=NavyBlue]  at (-3,3) {}; 
\end{scope}

\begin{scope}[xshift=100]
\begin{scope}[scale=0.4, xshift=0, yshift=50]
    \draw[very thin, color=black] (0,0) grid (1,2);
    \draw[very thin, color=black] (1,1) grid (3,2);
\end{scope}
\begin{scope}[scale=0.25, xshift=90, yshift=-400]
	\draw[thick, color=NavyBlue]    (0,0) -- (3,6) -- (0,10) -- (-3,7) -- (-3,3) -- (0,0);		
    \draw[thick, color=NavyBlue]    (3,6) -- (6,10) -- (3,13) -- (-3,7);
    
	\node[circle,fill,inner sep=1.5pt,color=NavyBlue]  at (0,0)  {};
	\node[circle,fill,inner sep=1.5pt,color=NavyBlue]  at (3,6) {}; 
	\node[circle,fill,inner sep=1.5pt,color=NavyBlue]  at (0,10) {};
	\node[circle,fill,inner sep=1.5pt,color=NavyBlue]  at (-3,7) {};
	\node[circle,fill,inner sep=1.5pt,color=NavyBlue]  at (-3,3) {};
	\node[circle,fill,inner sep=1.5pt,color=NavyBlue]  at (6,10) {};    
	\node[circle,fill,inner sep=1.5pt,color=NavyBlue]  at (3,13) {};    
\end{scope}
\end{scope}

\begin{scope}[xshift=210]
\begin{scope}[scale=0.4, xshift=10, yshift=50]
    \draw[very thin, color=black] (1,0) grid (2,1);
    \draw[very thin, color=black] (0,1) grid (3,2);
\end{scope}
\begin{scope}[scale=0.26, xshift=50, yshift=-360]
	\draw[thick, color=NavyBlue]    (0,0) -- (3,3) -- (0,6) -- (-3,3) -- (0,0);		
    \draw[thick, color=NavyBlue]    (3,3) -- (6,8) -- (3,12) -- (0,9) -- (0,6);
    
	\node[circle,fill,inner sep=1.5pt,color=NavyBlue]  at (0,0)  {};
	\node[circle,fill,inner sep=1.5pt,color=NavyBlue]  at (3,3) {}; 
	\node[circle,fill,inner sep=1.5pt,color=NavyBlue]  at (-3,3) {};
	\node[circle,fill,inner sep=1.5pt,color=NavyBlue]  at (0,6) {};
	\node[circle,fill,inner sep=1.5pt,color=NavyBlue]  at (0,9) {};
	\node[circle,fill,inner sep=1.5pt,color=NavyBlue]  at (6,8) {};    
	\node[circle,fill,inner sep=1.5pt,color=NavyBlue]  at (3,12) {}; 
\end{scope}
\end{scope}

\begin{scope}[xshift=320]
\begin{scope}[scale=0.4, xshift=20, yshift=50]
    \draw[very thin, color=black] (1,0) grid (2,3);
    \draw[very thin, color=black] (0,1) grid (3,2);
\end{scope}	
\begin{scope}[scale=0.25, xshift=100, yshift=-420]
	\draw[thick, color=NavyBlue]    (0,0) -- (3,3) -- (0,6) -- (-3,3) -- (0,0);		
    \draw[thick, color=NavyBlue]    (3,3) -- (6,8) -- (3,12) -- (0,9) -- (0,6);
    \draw[thick, color=NavyBlue]    (-3,3) -- (-6,8) -- (-3,12) -- (0,9);
    \draw[thick, color=NavyBlue]    (-3,12) -- (0,15) -- (3,12);    
    
	\node[circle,fill,inner sep=1.5pt,color=NavyBlue]  at (0,0)  {};
	\node[circle,fill,inner sep=1.5pt,color=NavyBlue]  at (3,3) {}; 
	\node[circle,fill,inner sep=1.5pt,color=NavyBlue]  at (-3,3) {};
	\node[circle,fill,inner sep=1.5pt,color=NavyBlue]  at (0,6) {};
	\node[circle,fill,inner sep=1.5pt,color=NavyBlue]  at (0,9) {};
	\node[circle,fill,inner sep=1.5pt,color=NavyBlue]  at (6,8) {};    
	\node[circle,fill,inner sep=1.5pt,color=NavyBlue]  at (3,12) {};
	\node[circle,fill,inner sep=1.5pt,color=NavyBlue]  at (-3,12) {};    
    \node[circle,fill,inner sep=1.5pt,color=NavyBlue]  at (0,15) {};
    \node[circle,fill,inner sep=1.5pt,color=NavyBlue]  at (-6,8) {};    
    
\end{scope}
\end{scope}
\end{tikzpicture}
    \caption{Some examples of cross-Tamari lattices: (1) a Tamari lattice, (2) a $\nu$-Tamari lattice, (3) an alt $\nu$-Tamari lattice, (4) the cross-Tamari lattice of the ``minimal cross shaped grid''.}
    \label{fig.several_cross-tamari_examples}
\end{figure}
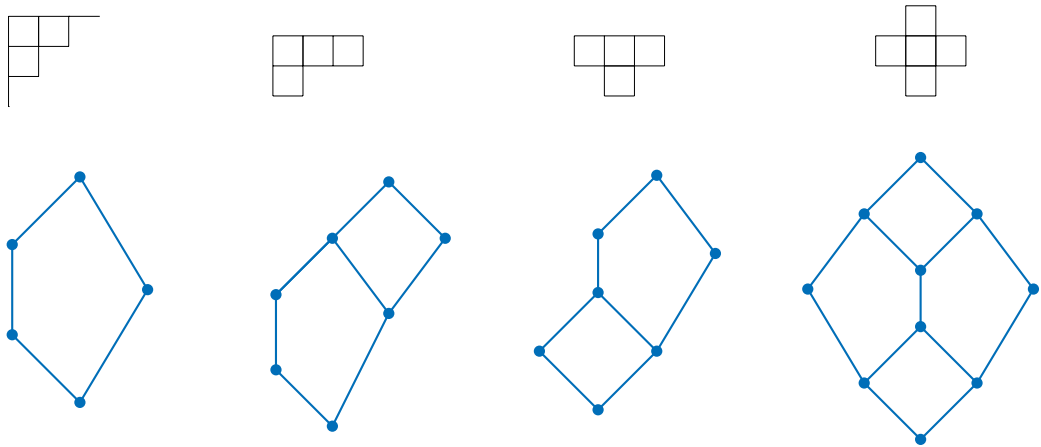

The case where $D$ is the set of lattice points weakly above a staircase shape recovers the classical Tamari lattice. 
If $D$ is the set of lattice points $\mathcal{L}_\nu$ weakly above a given lattice path~$\nu$ then we recover of $\nu$-Tamari lattice of Pr\'eville-Ratelle and Viennot~\cite{PV17} (using the approach of~\cite{CPS20}). Commuting the columns of $\mathcal{L}_\nu$ while preserving the cross-shaped property, gives rise to the alt $\nu$-Tamari lattices studied in~\cite{CC23}. See~\Cref{fig.several_cross-tamari_examples} for some examples.

\subsection{The cross-Tamari lattice as a framing lattice}
Our next objective is to find a framed graph whose framing lattice is isomorphic to the cross-Tamari order $\Tam(D)$ associated to a cross-shaped grid $D$.
If $D$ has $a$ columns and $b$ rows, it is convenient to assign positions to the points in~$D$ according to a relabeling of the columns with the numbers $1,\dots, a$ and the rows with $\overline{1},\dots, \overline{b}$, in some order. 
We identify a point $p\in D$ with its \defn{position} $p=(v,\overline{w})$ where $v$ is the label of column and $\overline{w}$ is the label of the row of the point. 
We denote by $\ell(v)$ (resp.~$\ell(\overline{w})$) the number of elements of $D$ in column $v$ (resp. row~$\overline{w}$). 
A \defn{proper labeling} of the rows and columns of $D$ is a labeling satisfying the following conditions: 
\begin{itemize}
    \item the column labels form a unimodular sequence\footnote{increases and then decreases} and $\ell(v)>\ell(v')$ implies $v<v' \,$;
    \item the row labels form a unimodular sequence and $\ell(\overline{w})>\ell(\overline{w}')$ implies $w<w'$.
\end{itemize}

Intuitively, this means that we label the columns (resp. rows) with the integers 1,2,... in order, from longest to shortest, starting with one of the longest columns (resp. rows) and then adding one column (resp. row) at a time directly to the left or right (resp. above or below) of those previously added.  
Such a labeling is not unique if $D$ has rows or columns of the same length, but any proper labeling will suffice for our purposes. 
See~\Cref{cross-shaped_grid_ex} for two examples of cross-shaped grids with proper labelings of their rows and columns. 
In this figure, the bottom-left corners (colored blue and red) of the two example have positions $(1,\overline{5})$ and $(3,\overline{4})$, respectively.


Let $D$ be a cross-shaped grid and $L$ be a proper labeling of its columns and rows with the numbers $[a]$ and $[\overline b]$. 
We define the \defn{$(D,L)$-caracol graph} $\crossGridGraph{D,L}$ as the graph on the vertex set $\{s,t\}\sqcup [a] \sqcup [\overline b]$, whose edges are given as follows.

First we define a linear order $\prec$ on the vertices, whose minimal element is $s$, maximal element is $t$, and the following three relations hold:
\begin{itemize}
    \item $i_1\prec i_2$ when $i_1<i_2$
    \item $\overline{j}_2\prec \overline{j}_1$ when $j_1<j_2$
    \item $x\prec \overline{y}$ when $(x,\overline{y})\in D$
\end{itemize}

The fact that $\prec$ is a linear order follows from the conditions on $D$ and $L$. We place the vertices $\{s,t\}\sqcup [a] \sqcup [\overline b]$ in a horizontal line following the linear order $\prec$ and draw an edge between each pair of consecutive elements. This looks like $s - 1 - \dots - \overline{1} - t$.
We add additional edges $(s,i)$ and $(\overline{j},t)$ as follows:

\begin{itemize}
    \item For $\overline{1}\neq \overline{j}\in [\overline{b}]$,  we draw an edge $(\overline{j},t)$ below the horizontal line if row label $\overline{j}$ is below row label $\overline{1}$, and above if it is above. 
    \item For $1\neq i \in [a]$, we draw an edge $(s,i)$ below the horizontal line if column label $i$ is on the right of column label $1$, and above if it is on the left. 
\end{itemize}

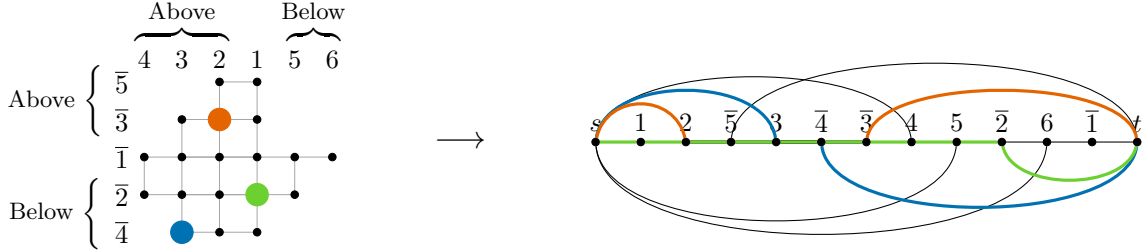
\begin{figure}[htb]
    \centering
    \begin{tikzpicture}

\begin{scope}[xshift = 0, yshift = 0, scale=0.5]
    \draw[very thin, color=gray!70] (1,0) grid (3,1);
    \draw[very thin, color=gray!70] (0,1) grid (4,2);
    \draw[very thin, color=gray!70] (4,2) grid (5,2);
    \draw[very thin, color=gray!70] (1,2) grid (3,3);
    \draw[very thin, color=gray!70] (2,3) grid (3,4);

    \node[circle, draw, inner sep=1pt, fill]  at (0,1)  {};
    \node[circle, draw, inner sep=1pt, fill]  at (0,2)  {};

    \node[circle, draw, inner sep=1pt, fill]  at (1,0)  {};
    \node[circle, draw, inner sep=1pt, fill]  at (1,1)  {};
    \node[circle, draw, inner sep=1pt, fill]  at (1,2)  {};
    \node[circle, draw, inner sep=1pt, fill]  at (1,3)  {};    

    \node[circle, draw, inner sep=1pt, fill]  at (2,0)  {};
    \node[circle, draw, inner sep=1pt, fill]  at (2,1)  {};
    \node[circle, draw, inner sep=1pt, fill]  at (2,2)  {};    
    \node[circle, draw, inner sep=1pt, fill]  at (2,3)  {};    
    \node[circle, draw, inner sep=1pt, fill]  at (2,4)  {};

    \node[circle, draw, inner sep=1pt, fill]  at (3,0)  {};
    \node[circle, draw, inner sep=1pt, fill]  at (3,1)  {};    
    \node[circle, draw, inner sep=1pt, fill]  at (3,2)  {};
    \node[circle, draw, inner sep=1pt, fill]  at (3,3)  {};    
    \node[circle, draw, inner sep=1pt, fill]  at (3,4)  {};

    \node[circle, draw, inner sep=1pt, fill]  at (4,1)  {};    
    \node[circle, draw, inner sep=1pt, fill]  at (4,2)  {};

    \node[circle, draw, inner sep=1pt, fill]  at (5,2)  {};

    \node[]  at (-2,3.5)  {\scriptsize \text{Above}  $\begin{cases} \\ \\ \end{cases}$};

    \node[]  at (-2,0.5)  {\scriptsize \text{Below} $\begin{cases} \\ \\ \end{cases}$};

    \node[] at (1,5.1) {\scriptsize $\overbrace{\color{white}{a+b+c}}^{}$};
    \node[] at (1,5.9) {\scriptsize  \text{Above}};

    \node[] at (4.5,5.1) {\scriptsize $\overbrace{\color{white}{a+b}}^{}$};
    \node[] at (4.5,5.9) {\scriptsize  \text{Below}};

    \node[]  at (0,4.6)  {\footnotesize $4$};
    \node[]  at (1,4.6)  {\footnotesize $3$};
    \node[]  at (2,4.6)  {\footnotesize $2$};
    \node[]  at (3,4.6)  {\footnotesize $1$};    
    \node[]  at (4,4.6)  {\footnotesize $5$};
    \node[]  at (5,4.6)  {\footnotesize $6$};
    
    \node[]  at (-.6,0)  {\footnotesize $\overline{4}$};
    \node[]  at (-.6,1)  {\footnotesize $\overline{2}$};
    \node[]  at (-.6,2)  {\footnotesize $\overline{1}$};
    \node[]  at (-.6,3)  {\footnotesize $\overline{3}$};    
    \node[]  at (-.6,4)  {\footnotesize $\overline{5}$};

    \node[circle, draw, inner sep=3pt, fill, color=blue]  at (1,0)  {};

    \node[circle, draw, inner sep=3pt, fill, color=red]  at (2,3)  {};

    \node[circle, draw, inner sep=3pt, fill, color=green]  at (3,1)  {};

\end{scope}

\begin{scope}[shift={(4.2,1.2)}]
    \node[] at (0,0) {$\longrightarrow$}; 
\end{scope}

\begin{scope}[scale=0.6, shift={(10,2)},
    every node/.style={circle, draw, inner sep=1pt, fill}]

    \node[label=above:{\footnotesize $s$}] (s) at (0,0) {};
    \node[label=above:{\footnotesize $1$}] (a6) at (1,0) {};
    \node[label=above:{\footnotesize $2$}] (a5) at (2,0) {};
    \node[label=above:{\footnotesize $\overline{5}$}] (b1) at (3,0) {};
    \node[label=above:{\footnotesize $3$}] (a4) at (4,0) {};
    \node[label=above:{\footnotesize $\overline{4}$}] (b2) at (5,0) {};
    \node[label=above:{\footnotesize $\overline{3}$}] (b3) at (6,0) {};
    \node[label=above:{\footnotesize $4$}] (a3) at (7,0) {};
    \node[label=above:{\footnotesize $5$}] (a2) at (8,0) {};
    \node[label=above:{\footnotesize $\overline{2}$}] (b4) at (9,0) {};
    \node[label=above:{\footnotesize $6$}] at (10,0) (a1) {};
    \node[label=above:{\footnotesize $\overline{1}$}] (b5) at (11,0) {};
    \node[label=above:{\footnotesize $t$}] (t) at (12,0) {};

    \draw[color=black] (s) -- (t);

    \draw[color=black] (s) .. controls (0.25, -2.7) and (9.75, -2.7) .. (a1);
    \draw[color=black] (s) .. controls (0.25, -2.3) and (7.75, -2.3) .. (a2);    
    \draw[color=black] (s) .. controls (0.25, 1.9) and (6.75, 1.9) .. (a3);
    \draw[color=black] (s) .. controls (0.25, 1.5) and (3.75, 1.5) .. (a4);
    \draw[color=black] (s) .. controls (0.25, 1.1) and (1.75, 1.1) .. (a5);	
               
    \draw[color=black] (b1) .. controls (3.25, 2.3) and (11.75, 2.3) .. (t);
    \draw[color=black] (b2) .. controls (5.25, -1.9) and (11.75, -1.9) .. (t);
    \draw[color=black] (b3) .. controls (6.25, 1.5) and (11.75, 1.5) .. (t);
    \draw[color=black] (b4) .. controls (9.25, -1.1) and (11.75, -1.1) .. (t);

    \draw[color=blue,very thick] (s) .. controls (0.25, 1.5) and (3.75, 1.5) .. (a4);
    \draw[color=blue,very thick] (b2) .. controls (5.25, -1.9) and (11.75, -1.9) .. (t);
    \draw[color=black,very thick] (a4) -- (b2);

    \draw[color=red, very thick] (s) .. controls (0.25, 1.1) and (1.75, 1.1) .. (a5);	
    \draw[color=red, very thick] (b3) .. controls (6.25, 1.5) and (11.75, 1.5) .. (t);
    \draw[color=black,very thick] (a5) -- (b1) -- (a4) -- (b2) -- (b3);

    \draw[color=green, very thick] (b4) .. controls (9.25, -1.1) and (11.75, -1.1) .. (t);
    \draw[color=green,very thick] (s) -- (a6) -- (a5) -- (b1) -- (a4) -- (b2) -- (b3) -- (a3) -- (a2) -- (b4);

\end{scope}
    
\end{tikzpicture}
    \caption{A cross-shaped grid $D$ with a proper labeling $L$ of its rows and columns (left). The $(D,L)$-caracol graph $\crossGridGraph{D,L}$ with the routes corresponding to the marked points in $D$ highlighted (right).}
    \label{cross-Tamari_grid_graph}
\end{figure}

The resulting graph is $\crossGridGraph{D,L}$, and the framing $F_{D,L}$ is the framing induced by our drawing; see~\Cref{cross-Tamari_grid_graph} for an example.

\begin{theorem}\label{thm_crossTam}
    For any proper labeling $L$ of the rows and columns of $D$,
    the framing lattice~$\scrL_{\crossGridGraph{D,L}, F_{D,L}}$ is isomorphic to the cross-Tamari order $\Tam(D)$.  
\end{theorem}

\begin{proof}
    An important feature of $\crossGridGraph{D,L}$ is that its routes can be characterized by two edges, namely the edges of the route incident to the source and sink.
    We can then express the unique route that uses edges $(s,i)$ and $(\overline{j},t)$ by $R_{i,\overline{j}}$.
    The map $(i,\overline{j})\rightarrow R_{i,\overline{j}}$ is a bijection between the lattice points in $D$ and the routes in $\crossGridGraph{D,L}$.
    One can check by inspection that, under this bijection,
    two lattice points are compatible in $D$ if and only if their corresponding routes are coherent (see examples of incompatible/compatible pairs of lattice points and their corresponding incoherent/coherent routes in~\Cref{cross-Tamari_grid_graph}). 
    Therefore, maximal cliques of the framing lattice correspond to maximal fillings. 
    One can also observe that increasing tree rotations correspond to $\ccw$-rotations. 
    As a consequence, the framing lattice $\scrL_{\crossGridGraph{D,L}, F_{D,L}}$ is isomorphic to the cross-Tamari order $\Tam(D)$.
\end{proof}

The following is a direct non-trivial consequence of Theorems~\ref{thm.framing_lattice}, \ref{thm.semidistributivity}, and \ref{thm.HHlattice}.

\begin{corollary}\label{cor_crossTam_consequences1}
    The cross-Tamari order $\Tam(D)$ is an $\mathcal{H}\mathcal{H}$-lattice. Hence it is semidistributive, congruence uniform, and polygonal. Furthermore, its polygons consist only of squares, pentagons or hexagons.
\end{corollary}

The following is another immediate but nontrivial consequence that illustrates the power of the framing lattice approach.

\begin{corollary}\label{cor_crossTam_consequences2}
For a cross-shaped grid $D$ with $a$ columns and $b$ rows the following hold:
\begin{enumerate}
    \item All maximal fillings have the same number of elements, equal to $a+b-1$.
    \item The rotation graph of maximal fillings is connected. 
    \item $\Tam(D)$ has a unique maximal and a unique minimal element.
\end{enumerate}
\end{corollary}
\begin{proof}
    For simplicity, let $G=\crossGridGraph{D,L}$.
    maximal fillings determine the maximal simplices of the framed triangulation of the flow polytope $\calF_{G}$. Thus, the number of elements of any maximal filling $T$ is equal to 
    \begin{align*}
        \dim \calF_{G}+1 &= |E(G)|-|V(G)|+2 \\
        &= (2a+2b-1)-(a+b+2)+2 \\
        &= a+b-1       
    \end{align*}
    This proves item (1). Item (2) follows from the fact that the dual graph of a triangulation of a polytope is connected. Item (3) about $\Tam(D)$ having a unique maximal and a unique minimal element follows from the fact of being a lattice.  
\end{proof}

The following result shows different cross-Tamari lattices obtainable through framed triangulations of the same flow polytope.

\begin{proposition}\label{prop_crossTam_permutations}
    Let $D$ and $D'$ be two cross-shaped grids that are obtainable from each other by a permutation of rows and columns. Then,
    \begin{enumerate}
        \item The graph $\crossGridGraph{D,L}$ is independent of a proper labeling $L$ of the rows and columns of $D$ (but the framing depends on $L$). In particular, the flow polytope $\calF_{\crossGridGraph{D}}:=\calF_{\crossGridGraph{D,L}}$ is independent of $L$.   
        \item The two graphs of $D$ and $D'$ are equal: $\crossGridGraph{D,L}=\crossGridGraph{D',L'}$. 
        \item The Hasse diagrams of the cross-Tamari lattices $\Tam(D)$ and $\Tam(D')$ are the dual graphs of two framed triangulations of the same flow polytope $\calF_{\crossGridGraph{D}}=\calF_{\crossGridGraph{D'}}$.
    \end{enumerate}
\end{proposition}

\begin{proof}
    For item (1), it suffices to check that the horizontal line of $\crossGridGraph{D,L}$ is fixed independent of the proper labeling $L$. 
    Or equivalently, that the linear order $\prec$ on the vertices $\{s,t\}\sqcup [a] \sqcup [\overline b]$ is independent of $L$. 
    This follows from the fact that $(i,\overline{j})\in D$ for $L$ if and only if $(i,\overline{j})\in D$ for $L'$, for any two proper labelings $L$ and $L'$.
    
    For item (2), fix two proper labelings $L$ and $L'$ of $D$ and $D'$ respectively. 
    Since $D$ and $D'$ are obtainable from each other by a permutation of rows and columns, one can observe that $(i,\overline{j})\in D$ if and only if $(i,\overline{j})\in D'$. This implies that $\crossGridGraph{D,L}=\crossGridGraph{D',L'}$. 

    For item (3), since $\crossGridGraph{D,L}=\crossGridGraph{D',L'}$ then their flow polytopes are equal. 
\end{proof}

\begin{proposition}
    Any framing of $G_{D,L}$ produces a framing lattice that is isomorphic to a cross-Tamari lattice.
\end{proposition}

\begin{proof}
    Consider $G_{D,L}$ with framing $F_{D,L}$. 
    Changing the framing by flipping an edge $(s,i)$ with $i\neq 1$ drawn above (resp. below) the horizontal line to being drawn below (resp. above) corresponds with moving the column labeled $i$ in $D$ to the right (resp. left) of the column labeled $1$, while preserving the cross-shaped property.
    Similarly, changing the framing by flipping an edge $(\overline{j},t)$ with $\overline{j}\neq \overline{1}$ drawn above (resp. below) the horizontal line to being drawn below (resp. above) corresponds with moving the column labeled $\overline{j}$ in $D$ below (resp. above) the row labeled $\overline{1}$. 
    Any of these operations gives a new cross-shaped grid $D'$ with proper labeling $L'$, with $G_{D,L} = G_{D',L'}$, while $F_{D,L}$ and $F_{D',L'}$ differ by the drawing of a single edge.
    Since every framing of $G_{D,L}$ can be obtained by iterating these moves, we obtain a corresponding cross-shaped grid and proper labeling for each framing. 
\end{proof}

\begin{figure}[h]
    \centering
    \begin{tikzpicture}[scale=0.75]


\begin{scope}
    \node[circle, draw, inner sep=1.4pt, fill, label=below:{\small $0$}] (0) at (0,0)  {};
    \node[circle, draw, inner sep=1.4pt, fill, label=below:{\small $1^-$}] (1) at (1,-1)  {};
    \node[circle, draw, inner sep=1.4pt, fill, label=above:{\small $2^+$}] (2) at (2,1)  {};
    \node[circle, draw, inner sep=1.4pt, fill, label=below:{\small $3^-$}] (3) at (3,-1)  {};
    \node[circle, draw, inner sep=1.4pt, fill, label=below:{\small $4$}] (4) at (4,0)  {};

    \draw[] (0) -- (1) -- (3) -- (4);	
    \draw[] (0) -- (2) -- (4);	
 
    \draw[very thick, red] (1) -- (4);
    \draw[very thick, blue] (2) -- (3);
        
\end{scope}
\begin{scope}[shift={(7,0)}, scale=1.2,
    fnode/.style={circle, draw, inner sep=1pt, fill}]
    
    \node[fnode, label=above:{\footnotesize $s$}] (s) at (0,0) {};
    \node[fnode, label=above:{\footnotesize $0$}] (0) at (1,0) {};
    \node[fnode, label=above:{\footnotesize $1$}] (1) at (2,0) {};
    \node[fnode, label=above:{\footnotesize $2$}] (2) at (3,0) {};
    \node[fnode, label=above:{\footnotesize $3$}] (3) at (4,0) {};
    \node[fnode, label=above:{\footnotesize $t$}] (t) at (5,0) {};
    
    \draw[color=black] (s) -- (t);
    \draw (0) -- (1) node[midway, below]{$-$};
    \draw (1) -- (2) node[midway, above]{$+$};
    \draw (2) -- (3) node[midway, below]{$-$};

    \draw[color=black] (s) .. controls (0.25, -0.7) and (1.75, -0.7) .. (1);
    \draw[color=black] (s) .. controls (0.25, 1.1) and (2.75, 1.1) .. (2);    
    \draw[color=black] (s) .. controls (0.25, -1.5) and (3.75, -1.5) .. (3);

    \draw[color=black] (0) .. controls (1.25, -1.5) and (4.75, -1.5) .. (t);
    \draw[color=black] (1) .. controls (2.25, 1.1) and (4.75, 1.1) .. (t);
    \draw[color=black] (2) .. controls (3.25, -0.7) and (4.75, -0.7) .. (t);

    \draw[very thick, color=blue] (s) .. controls (0.25, 1.1) and (2.75, 1.1) .. (2);    
    \draw[very thick, color=blue] (2) .. controls (3.25, -0.7) and (4.75, -0.7) .. (t);

    \draw[very thick, color=red] (s) .. controls (0.25, -0.7) and (1.75, -0.7) .. (1);
    \draw[very thick, color=red] (1) -- (2) -- (3) -- (t);

\end{scope}


    
    
    
    

 
        

\begin{scope}[shift={(1,-5.5)}, scale=0.7]
    \draw[very thin, color=gray] (0,1) grid (1,2);
    \draw[very thin, color=gray] (1,0) grid (1,1);
    \draw[very thin, color=gray] (1,1) grid (2,3);
    \draw[very thin, color=gray] (2,2) grid (3,2);    

    \node[]  at (0,3.6)  {\footnotesize {$2$}};
    \node[]  at (1,3.6)  {\footnotesize {$0$}};
    \node[]  at (2,3.6)  {\footnotesize {$1$}};
    \node[]  at (3,3.6)  {\footnotesize {$3$}};
    
    \node[]  at (-0.6,0)  {\footnotesize {$\overline{0}$}};
    \node[]  at (-0.6,1)  {\footnotesize {$\overline{2}$}};
    \node[]  at (-0.6,2)  {\footnotesize {$\overline{3}$}};
    \node[]  at (-0.6,3)  {\footnotesize {$\overline{1}$}};

    \node[circle,fill,inner sep=1.4pt,color=black] at (0,1)  {}; 
    \node[circle,fill,inner sep=1.4pt,color=black] at (0,2)  {};
    \node[circle,fill,inner sep=1.4pt,color=black] at (1,0)  {};     
    \node[circle,fill,inner sep=1.4pt,color=black] at (1,1)  {};		     
    \node[circle,fill,inner sep=1.4pt,color=black] at (1,2)  {};
    \node[circle,fill,inner sep=1.4pt,color=black] at (1,3)  {};		     
    \node[circle,fill,inner sep=1.4pt,color=black] at (2,1)  {};    
    \node[circle,fill,inner sep=1.4pt,color=black] at (2,2)  {};    
    \node[circle,fill,inner sep=1.4pt,color=black] at (2,3)  {}; 
    \node[circle,fill,inner sep=1.4pt,color=black] at (3,2)  {};  

    \node[circle,fill,inner sep=2pt,color=red] at (2,2)  {};      
    \node[circle,fill,inner sep=2pt,color=blue] at (0,1)  {};          
\end{scope}

\begin{scope}[shift={(7,-4)}, scale=1.2,
    fnode/.style={circle, draw, inner sep=1pt, fill}]
    
    \node[fnode, label=above:{\footnotesize $s$}] (s) at (0,0) {};
    
    \node[fnode, label=above:{\footnotesize $0$}] (0) at (0.8,0) {};
    \node[fnode, label=above:{\footnotesize $\overline{0}$}] (0') at (1.2,0) {};
    
    \node[fnode, label=above:{\footnotesize $1$}] (1) at (1.8,0) {};
    \node[fnode, label=above:{\footnotesize $\overline{1}$}] (1') at (2.2,0) {};    
    \node[fnode, label=above:{\footnotesize $2$}] (2) at (2.8,0) {};
    \node[fnode, label=above:{\footnotesize $\overline{2}$}] (2') at (3.2,0) {};    
    \node[fnode, label=above:{\footnotesize $3$}] (3) at (3.8,0) {};
    \node[fnode, label=above:{\footnotesize $\overline{3}$}] (3') at (4.2,0) {};
    \node[fnode, label=above:{\footnotesize $t$}] (t) at (5,0) {};
    
    \draw[color=black] (s) -- (t);
    \draw (0') -- (1);
    \draw (1') -- (2);
    \draw (2') -- (3);

    \node[] (a) at (1.5, -0.17) {$-$};
    \node[] (a) at (2.5, 0.18) {$+$};    
    \node[] (a) at (3.5, -0.17) {$-$};

    \draw[color=black] (s) .. controls (0.25, -0.7) and (1.55, -0.7) .. (1);
    \draw[color=black] (s) .. controls (0.25, 1.1) and (2.55, 1.1) .. (2);    
    \draw[color=black] (s) .. controls (0.25, -1.5) and (3.55, -1.5) .. (3);

    \draw[color=black] (0') .. controls (1.45, -1.5) and (4.75, -1.5) .. (t);
    \draw[color=black] (1') .. controls (2.45, 1.1) and (4.75, 1.1) .. (t);
    \draw[color=black] (2') .. controls (3.45, -0.7) and (4.75, -0.7) .. (t);

    \draw[very thick, color=blue] (s) .. controls (0.25, 1.1) and (2.55, 1.1) .. (2);
    \draw[very thick, color=blue] (2') .. controls (3.45, -0.7) and (4.75, -0.7) .. (t);

    \draw[very thick, color=red] (s) .. controls (0.25, -0.7) and (1.55, -0.7) .. (1);
    \draw[very thick, color=red] (1) -- (2) -- (3) -- (t);
    \draw[very thick, color=blue] (2.8,-0.04) -- (3.2,-0.04);

    \node[fnode, label=above:{\footnotesize $1$}] (1) at (1.8,0) {};
    \node[fnode, label=above:{\footnotesize $\overline{1}$}] (1') at (2.2,0) {};    
    \node[fnode, label=above:{\footnotesize $2$}] (2) at (2.8,0) {};
    \node[fnode, label=above:{\footnotesize $\overline{2}$}] (2') at (3.2,0) {};    
    \node[fnode, label=above:{\footnotesize $3$}] (3) at (3.8,0) {};
    \node[fnode, label=above:{\footnotesize $\overline{3}$}] (3') at (4.2,0) {};
    \node[fnode, label=above:{\footnotesize $t$}] (t) at (5,0) {};

\end{scope}

\end{tikzpicture}
    
    \caption{$\varepsilon$-Cambrian lattices are cross-Tamari lattices (\Cref{rem_Cambrian_cross}): example of a Cambrian caracol graph $G_\varepsilon$ and its corresponding $(D,L)$-caracol graph $G_{D,L}$ (with the order of labels in $[\overline{b}]$ reversed).}
    \label{fig_Cambrian_cross}
\end{figure}
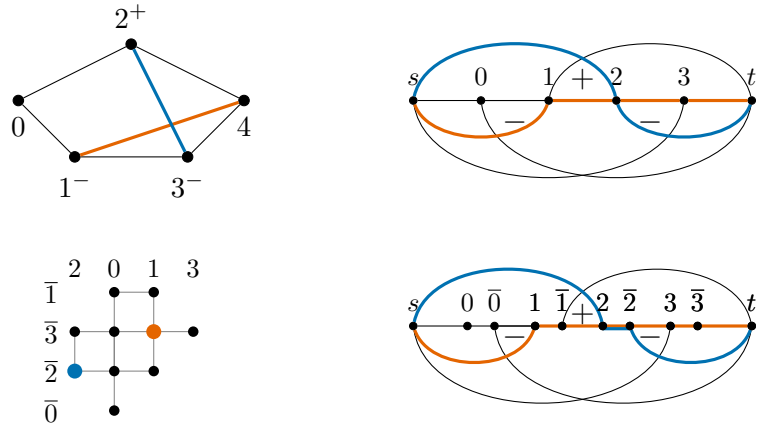

\begin{remark}\label{rem_Cambrian_cross}
    The $\varepsilon$-Cambrian lattices (and the $(\varepsilon,I,J)$-Cambrian lattices of Remark~\ref{rem.eIJ-cambrian}) are special cases of the cross-Tamari lattices.
    This follows from the observation that for any Cambrian caracol graph $G_{\varepsilon}$, one can easily construct a $(D,L)$-caracol graph $G_{D,L}$ that differs from $G_{\varepsilon}$ only by contracting idle edges (thus having isomorphic framing lattices). 
    More precisely, an internal node $i$ of $G_\varepsilon$ is duplicated to a horizontal edge $i-\overline{i}$, and the pair of edges $(s,i),(i,t)$ are transformed to $(s,i),(\overline{i},t)$; Note that in contrast to the cross-Tamari case, the bar labels now appear in increasing order from left to right, so when we translate to the cross-shaped grid the row labels are reversed (labels of rows now \emph{decrease} from longer to shorter, while the labels of columns increase from longer to shorter). A lattice point $(i,\overline{j})$ in the grid corresponds to the diagonal $(i,j+1)$ in the polygon $P_\varepsilon$.  See~\Cref{fig_Cambrian_cross} for an example.
    
    In order to get the $(\varepsilon,I,J)$-Cambrian lattices from~\cite{Pil20}, one just needs to restrict the cross-shaped grid corresponding to $G_\varepsilon$ to a subset of columns and a subset of rows. The restricted cross-shaped grid gives a cross-Tamari lattice isomorphic to a $(\varepsilon,I,J)$-Cambrian lattice. 
\end{remark}

\begin{figure}[h]
    \centering


    \caption{Examples of graphs giving rise to the further species.}
    \label{fig_further_species}
\end{figure}

\section{Further Species}

Further species in the framing lattice zoo include the following.
\begin{itemize}
    \item The $s$-weak order of Ceballos and Pons \cite{CP22}. The connection to triangulated flow polytopes was discovered by Gonz\'ales D'Le\'on et al \cite{GMPTY23}. Given an integer composition $s=(s_1,s_2,\ldots, s_n)$, the $s$-oruga graph $s$-$\oru{n}$ giving rise to the $s$-weak order as a framing lattice is constructed from $\oru{n}$ by adding $s_i$ copies of the edge $(1,n+2-i)$ for each $i\in [n]$. The framing is inherited from $\oru{n}$, with any of the added edges incoming to a vertex $j$ placed between the two edges of $\oru{n}$ incoming to $j$. The order of the added edges incoming to a vertex does not affect the lattice, neither does the order of the edges at the source. See Figure~\ref{fig_further_species} for some examples. 
    \vspace{0.2cm}
    \item The permutree lattices of Pilaud and Pons~\cite{PP18} unify the weak order and the Boolean, Tamari, and Cambrian lattices. Obtaining them as framing lattices is a direct consequence of the connection to flow polytopes given by Tamayo~\cite{Tam23}.
    As mentioned in Remark~\ref{rem.permutree_lattices}, the permutree lattices are the framing lattices obtained from the graph $\oru{n}$ by applying $M$-moves.
    Figure~\ref{fig_further_species} shows some examples of the framed graphs.
    \vspace{0.2cm}
    \item $\tau$-Tilting posets for certain gentle algebras are also obtainable as framing lattices, which is a direct consequence of \cite{BBBHPSY22}. We refer to \cite{BBBHPSY22} for the full details, but give a brief overview here.
    For directed acyclic graphs whose inner vertices have in- and out-degree equal to $2$ (called \emph{full} in \cite{BBBHPSY22}), there exists a framing whose exceptional routes contain all edges of $G$. A framing satisfying this condition (not necessarily unique) is called an \emph{ample framing}, and it has the property that the edges in each exceptional route can all be assigned the value $1$ or $2$ so that the two incoming edges into each inner vertex have different labels (and hence the two outgoing edges also have different labels).
    To any such full $G$ with ample framing, one can associate a quiver $Q$ by removing the source and sink of $G$ along with the edges incident to them, and then reversing the direction of any remaining edges labeled $2$. 
    The resulting quiver $Q$ encodes a gentle algebra~$\Lambda$, whose $\tau$-tilting poset is isomorphic to the framing lattice of $G$ with the ample framing. Therefore, any $\tau$-tilting poset of a gentle algebra whose associated quiver can be obtained from a full $G$ with an ample framing $F$ in this manner is isomorphic to the framing lattice $\scrL_{G,F}$. Figure~\ref{fig_further_species} shows a few examples of graphs with ample framings and their associated quiver. 

    We note that further progress was made by Berggren and Serhiyenko in a slightly more general context in~\cite[Corollary 4.9]{BS24}. 
    They proved the lattice property of framing lattices in the special case of \emph{rooted framed graphs}. 
    These are framed graphs without idle edges such that all \emph{exceptional segments} (i.e. maximal exceptional paths) are incident to a source or a sink.
    These rooted framed graphs generalize the amply framed graphs of \cite{BBBHPSY22}, but do not generally include the framed graphs with more than two edges (or paths) between two inner vertices, such as the second framed graph in Figure~\ref{fig.graphs_and_lattices}.
    
    \vspace{0.2cm}
    \item The Grid-Tamari lattices of McConville \cite{McC2017} are obtained as an immediate consequence of Garver--McConville~\cite{GM17}. 
    Given a grid $\lambda$ (a finite induced subgraph of the~$\mathbb{Z}^2$-lattice), orient the horizontal edges to the right, giving them the label $1$, and the vertical edges downward, giving them the label $2$.
    Then define $G_{\lambda}$ to be the graph obtained from $\lambda$ by adding a source vertex incident to all vertices of $\lambda$ with indegree~$\leq 1$, and a sink vertex incident to all vertices of $\lambda$ with outdegree $\leq 1$. 
    The framing~$F_{\lambda}$ is taken to be induced by the edge labels. 
    See Figure~\ref{fig_further_species} for an example (rotated $45^\circ$ to match our edge orientation convention). The grid-Tamari order $GT(\lambda)$ is isomorphic to the framing lattice $\scrL_{G_{\lambda},F_{\lambda}}$. 
    \vspace{0.2cm}
    \item The Grassmann--Tamari order of Santos--Stump--Welker~\cite{SSW17} is a special case of grid-Tamari orders and hence obtained as a framing lattice as described above. In particular, they are obtained when the grid $\lambda$ is an $k\times (n-k)$ rectangle, and the classical Tamari lattice is obtained by choosing $k=2$.
\end{itemize}

\printbibliography 
\end{document}